\documentclass[10pt]{article}
\usepackage[nohead,margin=1.0in]{geometry}
\usepackage{amssymb, amsmath, amsthm,amsfonts}
\usepackage{graphicx,epsfig}
\usepackage[tight]{subfigure}
\usepackage{cite}
\usepackage{times}
\usepackage{algpseudocode}
\usepackage{float}
\usepackage{appendix}
\usepackage{color}
\usepackage{epstopdf}
\usepackage{mathrsfs}
\usepackage{algorithm}
\usepackage{booktabs}
\usepackage{verbatim}
\usepackage[colorlinks,linktocpage,linkcolor=blue]{hyperref}

\graphicspath{{./pics/}}
\usepackage{booktabs}
\usepackage{threeparttable}

\numberwithin{equation}{section}
\newtheorem{theorem}{Theorem}[section]
\newtheorem{remark}{Remark}[section]

\newtheorem{lemma}{Lemma}[section]

\newtheorem{definition}[theorem]{Definition}
\newtheorem{condition}[theorem]{Condition}
\newtheorem{assumption}[theorem]{Assumption}
\newtheorem{example}{Example}[section]
\allowdisplaybreaks
\def\n{\mathrm{\bf n}}
\def\P01{P_{\mathcal{A}}}

\def\bsgamma{\boldsymbol{\gamma}}

\title{Conductivity Imaging from Internal Measurements with Mixed Least-Squares Deep Neural Networks\thanks{The work of B. Jin is supported by UK EPSRC grant EP/T000864/1 and EP/V026259/1, and a start-up fund from The Chinese University of Hong Kong. The work of  Z. Zhou is supported by Hong Kong Research Grants Council (15303021) and an internal grant of Hong Kong Polytechnic University (Project ID: P0038888, Work Programme: ZVX3).}}

\author{Bangti Jin\thanks{Department of Mathematics, The Chinese University of Hong Kong, Shatin, New Territories, Hong Kong, P.R. China (\texttt{bangti.jin@gmail.com, btjin@math.cuhk.edu.hk}).} \and Xiyao Li\thanks{Department of Computer Science, University College London, Gower Street, London WC1E 6BT, UK (\texttt{xiyao.li.20@ucl.ac.uk})} \and Qimeng Quan\thanks{School of Mathematics and Statistics, Wuhan University, Wuhan 430072, P. R. China (\texttt{quanqm@whu.edu.cn})} \and Zhi Zhou\thanks{Department of Applied Mathematics,
The Hong Kong Polytechnic University, Kowloon, Hong Kong, P.R. China (\texttt{zhizhou@polyu.edu.hk})}}

\begin{document}

\maketitle

\begin{abstract}
In this work we develop a novel approach using deep neural networks to reconstruct the conductivity distribution in elliptic problems from one measurement of the solution over the whole domain. The approach is based on a mixed reformulation of the governing equation and utilizes the standard least-squares objective, with deep neural networks as ansatz functions  to approximate the conductivity and flux simultaneously. We provide a thorough analysis of the deep  neural network approximations of the conductivity for both continuous and empirical losses, including rigorous error estimates that are explicit in terms of the noise level, various penalty parameters and neural network architectural parameters (depth, width and parameter bound). We also provide multiple numerical experiments in two- and multi-dimensions to illustrate distinct features of the approach, e.g., excellent stability with respect to data noise and capability of solving high-dimensional problems.\vskip5pt
\noindent\textbf{Key words}: conductivity imaging, least-squares approach, deep neural network, generalization error, error estimate
\end{abstract}

\section{Introduction}

The conductivity value of an object varies widely with the composition and type of materials and its accurate imaging can provide valuable structural information about the object. This fact underpins several important imaging modalities, e.g., electrical impedance tomography (EIT), current density impedance imaging, and acousto-electrical tomography; see the works \cite{Bal:2013,WidlakScherzer:2012} for comprehensive overviews on mathematical models and relevant theoretical results. In this work, we aim at identifying the conductivity distribution in elliptic problems from internal data using deep neural networks. Let $\Omega\subset \mathbb{R}^d$ be a simply connected open bounded domain with a smooth boundary $\partial\Omega$. Consider the following Neumann boundary value problem for the function $u$
\begin{equation}\label{equ:Neu problem}
	\left\{
	\begin{aligned}
		-\nabla\cdot(q\nabla u) &= f, \ &\mbox{in}&\ \Omega, \\
		q\partial_\n u&=g, \ &\mbox{on}&\ \partial\Omega,
	\end{aligned}
	\right.
\end{equation}
where $\n$ denotes the unit outward normal direction to the boundary $\partial\Omega$ and $\partial_\n$ denotes taking the normal derivative. The functions $f$ and $g$ in \eqref{equ:Neu problem} are the given source and boundary flux, respectively, and satisfy the standard compatibility condition $\int_\Omega f\ {\rm d} x + \int_{\partial\Omega}g\ {\rm d} S = 0$ in order to ensure solvability, and the solution $u$ is unique under a suitable normalization condition, e.g., $\int_\Omega u{\rm d} x=0$. The conductivity $q$ belongs to the admissible set
$\mathcal{A} = \{q\in H^1(\Omega): c_0\leq q\leq c_1 \mbox{ a.e. in } \Omega\}$,
with the constants $ 0<c_0<c_1<\infty $ being the lower and upper bounds on $q$. We use the notation $u(q)$
to indicate the dependence of the solution $u$ to problem \eqref{equ:Neu problem} on the coefficient $q$.

The concerned inverse problem is to recover the conductivity $q$ from one noisy internal measurement $z^\delta$ of the exact data $u(q^\dag)$ (i.e., the solution of problem \eqref{equ:Neu problem} corresponding to the exact conductivity $q^\dag$) over the whole domain $\Omega$.
It has been extensively studied in both engineering and mathematics communities. For
example, the model \eqref{equ:Neu problem} is often used to describe the behavior
of a confined inhomogeneous aquifer, where $u$ represents the piezometric head, $f$
is the recharge, and $q$ is hydraulic conductivity (or transmissivity
in the two-dimensional case); see the works \cite{FrindPinder:1973,Yeh:1986} for extensive discussions
on parameter identifications in hydrology. H\"{o}lder type stability
estimates of the inverse problem have been established under different settings \cite{Richter:1981,
Alessandrini:1986,bonito2017diffusion}. A reconstruction can be
obtained using the regularized output least-squares approach \cite{Acar:1993,
ChenZou:1999,DeckelnickHinze:2012}, equation error approach \cite{Falk:1983,
Karkkainen:1997,AlJamalGockenbach:2012} and mixed type formulation \cite{KohnLowe:1988} etc.
Error bounds on the numerical approximation obtained by the Galerkin finite element method (FEM)
of the regularized formulation were established in \cite{WangZou:2010,jin2021error},
by properly drawing on the stability argument.

In this work, we develop a new approach for recovering the conductivity $q$. It is based on a least-squares mixed-type reformulation of the
governing equation, with an $H^1(\Omega)$
penalty on the unknown conductivity $q$. The mixed least-squares formulation was first proposed by Kohn and Lowe \cite{KohnLowe:1988} (and hence also known as the Kohn-Lowe approach), where the conductivity $q$ and current flux $\sigma:=q\nabla u(q)$ are both approximated using the Galerkin FEM. We approximate both current flux $\sigma$ and conductivity $q$
separately using deep neural networks (DNNs), adopt a least-squares objective for all the equality constraints, and minimize the resulting empirical loss with respect to the DNN parameters. The use of DNNs in place of FEM allows exploiting inductive bias and expressivity of DNNs, which can be highly beneficial for numerical recovery. By leveraging the approximation theory of DNNs \cite{guhring2021approximation}, nonstandard energy argument \cite{KohnLowe:1988,jin2021error} and statistical
learning theory \cite{AnthonyBartlett:1999,BartlettMendelson:2002}, we derive
novel error bounds on the conductivity approximation in terms of the accuracy of
the observational data $z^\delta$, DNN architecture parameters (depth, width and parameter bound), and
the numbers of sampling points in the domain $\Omega$ and on the boundary $\partial\Omega$ etc. This is carried out for both population and empirical losses (resulting from Monte Carlo quadrature of the integrals). These error bounds provide theoretical underpinnings for the approach. The overall analysis relies crucially on the stability of the continuous inverse problem and a suitable discretization via DNNs such that the stability argument can be adapted for the error analysis of the numerical approximation. Moreover,
the proposed approach is easy to implement and robust with respect to noise, especially for partial interior data, and can handle high-dimensional inverse
problems (e.g., $d=5$). For example, the approach can still yield reasonable approximations in the presence of
10\% data noise. These
features make the method highly attractive, and the numerical results
show its significant potential.
The development of the DNN formulation, error analysis and numerical validation are the main contributions of this work.

In recent years, the use of DNNs for solving direct and inverse problems for PDEs has received a lot of attention; see \cite{EHanJentzen:2022,TanyuMaass:2022} for overviews. Existing neural inverse schemes using DNNs roughly fall into two groups: supervised approaches (see, e.g., \cite{SeoKimHarrach:2019,KhooYing:2019,GuoJiang:2021}) and unsupervised approaches (see, e.g., \cite{BaoYeZang:2020,BarSochen:2021,PakravanMistani:2021,XuDarve:2022,jin2022imaging,Pokkunuru:2023}). Supervised methods exploit the availability of (abundant) paired training data to extract problem-specific features, and are concerned with learning approximate inverse operators. In contrast, unsupervised approaches utilize DNNs as function approximators (ansatz functions) for the unknown coefficient and the state, which enjoy excellent expressivity for high-dimensional functions (in favorable situations), and the associated inductive biases \cite{Rahaman:2019}. The works \cite{BaoYeZang:2020}
and \cite{BarSochen:2021} investigated EIT reconstruction using the weak and strong formulations (also with the $L^\infty$ norm consistency for the latter), respectively; See also \cite{Pokkunuru:2023} for further developments of the approach \cite{BarSochen:2021} using energy based models. The work \cite{jin2022imaging} applied the deep Ritz method to a least-gradient reformulation for current density impedance imaging, and derived a generalization error for the loss function. The approach performs reasonably well for both full and partial interior current density data, and shows remarkable robustness against data noise. Pakravan et al \cite{PakravanMistani:2021} developed a hybrid approach, blending high expressivity of DNNs with the accuracy and reliability of traditional numerical methods for PDEs, and showed the approach for recovering the diffusion coefficient in one- and two-dimensional elliptic PDEs; see also \cite{CenJinQuanZhou:2023} for a hybrid DNN-FEM approach for recovering the diffusion coefficient in elliptic and parabolic problems, where the conductivity $q$ and state $u$ are approximated using DNNs and FEM, respectively. The work \cite{CenJinQuanZhou:2023} aims at combining the strengths of neural network and classical Galerkin FEM, i.e., expressivity of DNNs and solid theoretical foundations of the FEM, and provides a thorough theoretical analysis, but it is limited to low-dimensional problems due to the use of FEM for the state $u$. These works have presented encouraging empirical results for a range of PDE inverse problems, and clearly showed enormous potentials of DNNs in solving PDE inverse problems. Our approach follows the unsupervised paradigm, but unlike existing approaches, it employs a mixed formulation of the governing equation and thus differs greatly from the aforementioned ones. In addition, we have established rigorous error bounds on the approximation of the conductivity. Note that  the theoretical analysis of neural inverse schemes is largely elusive, due to the outstanding challenges, mostly associated with nonconvexity of the loss and a lack of linear structure of the DNN function class.

The rest of the paper is organized as follows. In Section \ref{sec:prelim}, we recall preliminaries on DNNs. In Section \ref{sec:Neumann}, we develop the  approach for a Neumann problem \eqref{equ:Neu problem}
based on a mixed formulation of the governing equation. Further we present an error analysis of the approach for both population and empirical losses, using tools from PDEs \cite{evans2010partial} and statistical learning theory \cite{AnthonyBartlett:1999,ShalevBen:20214}.
In Section \ref{sec:Diri}, we describe the extension of the approach to the Dirichlet case.
In Section \ref{sec:numer}, we present numerical experiments to validate the effectiveness of the approach, including large noise, high dimension and partial interior data. We denote by $W^{k,p}(\Omega)$ and $W^{k,p}_0(\Omega)$ the standard Sobolev spaces of order $k$ for any integer $k\geq0$ and real $p\geq1$, equipped with the norm $\|\cdot\|_{W^{k,p}(\Omega)}$. Further, we denote by $W^{-k,p'}(\Omega)$ the dual space of $W^{k,p}_0(\Omega)$, with the pair $(p,p')$ being the H{\"o}lder conjugate exponents. We also write $H^k(\Omega)$ and $H^k_{0}(\Omega)$ with the norm $\|\cdot\|_{H^k(\Omega)}$ if $p=2$ and write $L^p(\Omega)$ with the norm $\|\cdot\|_{L^p(\Omega)}$ if $k=0$. The spaces
on the boundary $\partial\Omega$ are defined similarly. The notation $(\cdot,\cdot)$ denotes the standard $L^2(\Omega)$ inner product. For a Banach space $B$, and the notation $B^d$ represent the $d$-fold product space. We denote by $c$ a generic constant not necessarily the same at each occurrence but it is always independent of the DNN approximation accuracy $\epsilon$, the noise level $\delta$ and the penalty parameters ($\gamma_\sigma$, $\gamma_b$ and $\gamma_q$).

\section{Preliminaries on DNNs}\label{sec:prelim}
In this section, we describe fully connected feedforward
neural networks. Let $\{d_\ell\}_{\ell=0}^L \subset\mathbb{N}$ be fixed with $d_0=d$, and a parameterization $\Theta=
\{(A^{(\ell)},b^{(\ell)})_{\ell=1}^L\}$ consisting of weight matrices and bias vectors, with
$A^{(\ell)}=[W_{ij}^{(\ell)}]\in \mathbb{R}^{d_\ell\times d_{\ell-1}}$ and $b^{(\ell)}=
[b^{(\ell)}_i]\in\mathbb{R}^{d_{\ell}}$ the weight matrix and bias vector at the $\ell$-th layer.
Then a DNN $v_\theta:= v^{(L)}:\Omega\subset\mathbb{R}^d\to\mathbb{R}^{d_L}$
realized by a parameter $\theta\in\Theta$ is defined recursively by
\begin{equation}\label{eqn:NN-realization}
\mbox{DNN realization:}\quad
\left\{\begin{aligned}
v^{(0)}&=x,\quad x\in\Omega\subset\mathbb{R}^d,\\		
v^{(\ell)}&=\rho(A^{(\ell)}v^{(\ell-1)}+b^{(\ell)}),\quad \ell=1,2,\cdots,L-1,\\
		v^{(L)}&=A^{(L)}v^{(L-1)}+b^{(L)},
	\end{aligned}\right.
\end{equation}
where the activation function $\rho:\mathbb{R}\to\mathbb{R}$ is applied
componentwise to a vector, and below we take $\rho\equiv\tanh$: $x\to
\frac{e^x-e^{-x}}{e^x+e^{-x}}$. The DNN $v_\theta$ has a depth $L$ and width $W:=
\max_{\ell=0,\dots,L}(d_{\ell})$.
Given the parametrization $\Theta$ (i.e., architecture), we denote the associated DNN function class $\mathcal{N}_\Theta$ by $\mathcal{N}_\Theta:=
\{v_\theta:\ \theta\in \Theta\}$. The
following approximation result holds \cite[Proposition 4.8]{guhring2021approximation}.
\begin{lemma}\label{lem:tanh-approx}
Let $s\in\mathbb{N}\cup\{0\}$ and $p\in[1,\infty]$ be fixed, and $v\in W^{k,p}
(\Omega)$ with $\mathbb{N}\ni k\geq s+1$. Then for any $\epsilon>0$, there exists at least one
 $\theta\in\Theta$ with depth $O(\log(d+k))$ and number of nonzero weights $O
\big(\epsilon^{-\frac{d}{k-s-\mu (s=2)}}\big)$ and weight parameters bounded by $O(\epsilon^{-2-\frac{2(d/p+d+k+\mu(s=2))+d/p+d}{k-s-\mu(s=2)}})$ in the maximum norm, where $\mu>0$ is arbitrarily small,
such that the DNN realization $v_\theta\in \mathcal{N}_\Theta$ satisfies
$\|v-v_\theta\|_{W^{s,p}(\Omega)} \leq \epsilon.$
\end{lemma}

\begin{remark}
On the domain $\Omega=(0,1)^d$, Guhring and Raslan \cite[pp. 127-128]{guhring2021approximation} proved Lemma \ref{lem:tanh-approx} using two steps. They first approximate a function $v\in W^{k,p}(\Omega)$ by a localized Taylor polynomial $v_{\rm poly}$:
$\|v-v_{\rm poly}\|_{W^{s,p}(\Omega)}\leq c_{\rm poly}N^{-(k-s-\mu(s=2))}$,
where the construction $($see \cite[Definition 4.4]{guhring2021approximation} for details$)$ of $v_{\rm poly}$ relies on an approximate partition of unity $($depending on $N\in\mathbb{N}$$)$ and the constant $c_{\rm poly}=c_{\rm poly}(d, p, k,s)>0$. Next they show that there exists a DNN parameter $\theta$, satisfying the conditions in Lemma \ref{lem:tanh-approx}, such that \cite[Lemma D.5]{guhring2021approximation}:
     \begin{equation*}
        \|v_{\rm poly}-v_\theta\|_{W^{s,p}(\Omega)} \leq c_{\rm NN}\|v\|_{W^{k,p}(\Omega)}\tilde{\epsilon},
     \end{equation*}
    where the constant $c_{\rm NN}=c_{\rm NN}(d,p,k,s)>0$ and  $\tilde{\epsilon}\in(0,\frac12)$.
    Now for small $\epsilon>0$, the desired estimate follows directly from the choice below
    $
        N=(\frac{\epsilon}{2c_{\rm poly}})^{-\frac{1}{k-s-\mu(s=2)}}$ and $\tilde{\epsilon}=\frac{\epsilon}{2c_{\rm NN}\|v\|_{W^{k,p}(\Omega)}}$.
\end{remark}

We denote the set of DNNs of depth $L$, $N_\theta$ nonzero weights, and maximum weight bound $R$ by
\begin{equation*}
    \mathcal{N}(L,N_\theta,R) =: \{v_\theta \mbox{ is a DNN of depth }L: \|\theta\|_{\ell^0}\leq N_\theta, \|\theta\|_{\ell^\infty}\leq R\},
\end{equation*}
where $\|\cdot\|_{\ell^0}$ and $\|\cdot\|_{\ell^\infty}$ denote, respectively, the number of nonzero entries in and the maximum norm of a vector. Furthermore, for any $\epsilon>0$ and $p\geq 1$, we denote the DNN parameter set by $\mathfrak{P}_{p,\epsilon}$ for
\begin{equation*}
\mathcal{N}\Big(c\log(d+2), c\epsilon^{-\frac{d}{1-\mu}}, c \epsilon^{-\frac{4p+3d+3pd}{p(1-\mu)}}\Big).
\end{equation*}

Below, for a vector-valued function, we use the DNNs to approximate its components. This can be achieved by
parallelization, which assembles multiple DNNs into one larger DNN.
Moreover,
the new DNN does not change the depth, and its width equals to the sum of that of subnetworks.
\begin{lemma}\label{lem:NN-paral}
Let $\bar \theta=\{(\bar A^{(\ell)},\bar b^{(\ell)})\}_{\ell=1}^L,\tilde\theta = \{(\tilde A^{(\ell)},\tilde b^{(\ell)})\}_{\ell=1}^L\in\Theta$,  let $\bar v$ and $\tilde v$ be their DNN realizations, and define $\theta=\{( A^{(\ell)}, b^{(\ell)})\}_{\ell=1}^L$ by
	\begin{align*}
		 A^{(1)}&=
		\begin{bmatrix}
		\bar A^{(1)}  \\
		\tilde A^{(1)}
		\end{bmatrix}\in\mathbb{R}^{2d_1\times d_0}, \
	     A^{(\ell)}=
	    \begin{bmatrix}
	   \bar A^{(\ell)} & 0 \\
	    	0            & \tilde A^{(\ell)}
	    \end{bmatrix}\in\mathbb{R}^{2d_\ell\times 2d_{\ell-1}},\quad  \ell=2,\dots,L; \\
		 b^{(\ell)} &=
		\begin{bmatrix}
	\bar b^{(\ell)} \\
		\tilde b^{(\ell)}
		\end{bmatrix}\in \mathbb{R}^{2d_\ell}, \quad \ell=1,\dots,L.
	\end{align*}
Then {$ v_{\theta}=(\bar v_{\bar \theta},\tilde v_{\tilde \theta})^\top $} is the DNN realization of $\theta$,
of depth $L$ and width $2{W}$.
\end{lemma}

\section{Inverse conductivity problem in the Neumann case}
\label{sec:Neumann}

In this section we discuss the inverse conductivity problem for the Neumann problem \eqref{equ:Neu problem}.

\subsection{Mixed formulation and its DNN approximation}
Upon letting $\sigma=q\nabla u$, we recast problem \eqref{equ:Neu problem} into a first-order system
\begin{equation}\label{eqn:mixed-Neum}
	\left\{\begin{aligned}
		\sigma & = q\nabla u,&& \mbox{in}\ \Omega, \\
		-\nabla\cdot\sigma &= f, &&\mbox{in}\ \Omega, \\
		\sigma \cdot \n&=g, &&\mbox{on}\ \partial\Omega.
	\end{aligned}\right.
\end{equation}
To identify $q$, we employ a noisy observation $z^\delta$ (of the exact data
$u(q^\dag)$) in $\Omega$, with
\begin{equation}\label{eqn:delta}
\delta:= \|\nabla(u(q^\dagger)-z^\delta)\|_{L^2(\Omega)}.
\end{equation}
Note that in \eqref{eqn:delta}, we have imposed a mild regularity condition on the noisy
data $z^\delta$, which may be obtained by presmoothing the
raw noisy data beforehand, e.g., via filtering \cite{MurioZhan:1998} or denoising via
DNNs \cite{Ulyanov:2018,jin2022imaging}. This assumption is used in,
e.g., Kohn--Lowe \cite{KohnLowe:1988} and equation error  \cite{AlJamalGockenbach:2012} formulations.

We employ a residual minimization scheme based on the first-order system \eqref{eqn:mixed-Neum},
approximate both $\sigma$ and $q$ using DNNs, and substitute $z^\delta$
for the scalar field $u$. For $q$,
we use a DNN function class (of depth $L_q$ and width $W_q$) with the parametrization $\mathfrak{P}_{p,\epsilon_q}$ (with $p\geq 2$ and tolerance $\epsilon_q$). Similarly, for
$\sigma:\Omega\to\mathbb{R}^d$, we
employ $d$ identical DNN function
classes (of depth $L_{\sigma}$ and width $W_{\sigma}$) with the parametrizations $\mathfrak{P}_{2,\epsilon_\sigma}$ (with tolerance $\epsilon_\sigma$), and stack them into one DNN, cf. Lemma \ref{lem:NN-paral}. Throughout, $\theta$ and $\kappa$ denote the parameters of DNN approximations to $q$ and $\sigma$, respectively.

To enforce the box constraint on $q$,
we employ a cutoff operator $\P01: H^1(\Omega) \rightarrow \mathcal{A}$ defined by
$\P01(v) = \min(\max(c_0,v),c_1)$. Then $P_\mathcal{A}$ is stable  \cite[Corollary 2.1.8]{Ziemer:1989}
\begin{equation}\label{eqn:P01-stab}
\| \nabla \P01(v) \|_{L^p(\Omega)} \le \|\nabla v \|_{L^p(\Omega)},\quad \forall v \in W^{1,p}(\Omega),p\in[1,\infty],
\end{equation}
and moreover, for all $v\in \mathcal{A}$, there holds
\begin{equation}\label{eqn:P01-approx}
\|  \P01(w) - v \|_{L^p(\Omega)} \le \| w - v \|_{L^p(\Omega)},\quad \forall w \in L^{p}(\Omega),~p\in[1,\infty].
\end{equation}
Using the least-squares formulation on the equality constraints, we obtain
\begin{equation}\label{eqn:loss-Neum}
\begin{aligned}
    J_{\bsgamma}(\theta,\kappa)=&\|\sigma_\kappa- \P01(q_\theta)\nabla z^\delta\|_{L^2(\Omega)}^2+\gamma_\sigma\|\nabla\cdot\sigma_\kappa+f\|_{L^2(\Omega)}^2\\
      &+\gamma_b\|\sigma_\kappa\cdot \n-g\|^2_{L^2(\partial\Omega)}+\gamma_q
    \| \nabla q_\theta\|_{L^2(\Omega)}^2.
\end{aligned}
\end{equation}
Then the DNN reconstruction problem reads
\begin{equation}\label{eqn:obj-Neum}
	\min_{(\theta,\kappa)\in (\mathfrak{P}_{p,\epsilon_q},\mathfrak{P}_{2,\epsilon_\sigma}^{\otimes d})} J_{\bsgamma}(\theta,\kappa),
\end{equation}
where the superscript $\otimes d$ denotes the $d$-fold direct product, $\gamma_\sigma>0$, $\gamma_b>0$ and $\gamma_q>0$ are penalty parameters, and $\bsgamma=(\gamma_\sigma,\gamma_b,\gamma_q)\in\mathbb{R}_+^3$.
The $H^1(\Omega)$ seminorm penalty on $q_\theta$ is to stabilize the minimization process. It is essential for overcoming the inherent ill-posedness of the inverse problem \cite{EnglHankeNeubauer:1996,ItoJin:2015}.
The well-posedness of problem
\eqref{eqn:obj-Neum} follows by a standard argument in calculus of variation. Indeed, the
compactness of the parametrizations $\mathfrak{P}_{p,\epsilon_q}$ and $\mathfrak{P}_{2,\epsilon_\sigma}$
holds due to the uniform boundedness on the parameter vectors and finite-dimensionality of the spaces. Meanwhile, the smoothness
of the $\tanh$ activation function implies the continuity of the loss $J_{\bsgamma}$
in the DNN parameters $(\theta,\kappa)$. These two properties imply the existence of a minimizer $(\theta^*,\kappa^*)$.

Note that the objective $J_{\bsgamma}$ involves high-dimensional integrals, and
hence quadrature is needed in practice. This may be achieved using
quadrature rules. In this work, we employ the Monte Carlo method. Let $\mathcal{U}(\Omega)$ and $\mathcal{U}(\partial\Omega)$ be the uniform distributions
over the domain $\Omega$ and the boundary $\partial\Omega$, respectively, and
$(q_\theta,\sigma_\kappa)$ the DNN realization of $(\theta,\kappa)$.
Using the expectation $\mathbb{E}_{\mathcal{U}(\cdot)}[\cdot]$ with respect to $\mathcal{U}(\cdot)$ and likewise $\mathbb{E}_{\mathcal{U}(\partial\Omega)}$, we
can rewrite the loss $J_{\bsgamma}$ as
\begin{align*}
J_{\bsgamma}(\theta,\kappa)&=
|\Omega|\mathbb{E}_{X\sim\mathcal{U}(\Omega)}\Big[\|\sigma_{\kappa}(X)- \P01(q_\theta(X))\nabla z^\delta(X)\|_{\ell^2}^2\Big] \\ &\quad+\gamma_\sigma|\Omega|\mathbb{E}_{X\sim\mathcal{U}(\Omega)}\Big[\big(\nabla\cdot\sigma_\kappa(X)+f(X)\big)^2\Big]\\
&\quad+\gamma_b|\partial\Omega|\mathbb{E}_{Y\sim\mathcal{U}(\partial\Omega)}\Big[\big(\sigma_\kappa(Y)\cdot \n-g(Y)\big)^2\Big]+\gamma_q|\Omega|\mathbb{E}_{X\sim\mathcal{U}(\Omega)}\Big[\|\nabla q_\theta(X)\|_{\ell^2}^2 \Big]   \\
	&=: \mathcal{E}_d(\sigma_{\kappa},q_\theta) + \gamma_\sigma \mathcal{E}_\sigma(\sigma_\kappa)+\gamma_b \mathcal{E}_b(\sigma_\kappa)+\gamma_q\mathcal{E}_q(q_\theta),
\end{align*}
where $|\Omega|$ and $|\partial\Omega|$ denote the Lebesgue measure of $\Omega$ and $\partial\Omega$, respectively, and $\|\cdot\|_{\ell^2}$ denotes the Euclidean norm on $\mathbb{R}^d$. Thus, $J_{\bsgamma}(\theta,\kappa)$ is also known as the population loss. Let
$X=\{X_j\}_{j=1}^{n_r}$ and $Y=\{Y_{j}\}_{j=1}^{n_b}$ be independent and identically
distributed (i.i.d.) samples drawn from the uniform distributions $\mathcal{U}
(\Omega)$ and $\mathcal{U}(\partial\Omega)$, respectively, where $n_r$ and $n_b$ are the numbers of sampling points in the domain $\Omega$ and on the boundary $\partial\Omega$, respectively.
Then the empirical loss $\widehat{J}_{\bsgamma}(\theta,\kappa)$ is given by
\begin{equation}\label{eqn:obj-Neum-dis}
	\widehat{J}_{\bsgamma}(\theta,\kappa)=:\mathcal{\widehat{E}}_d(\sigma_{\kappa},q_\theta) + \gamma_\sigma \mathcal{\widehat{E}}_\sigma(\sigma_\kappa)+\gamma_b\mathcal{\widehat{E}}_b(\sigma_\kappa)+\gamma_q\mathcal{\widehat{E}}_q(q_\theta),
\end{equation}
where $\mathcal{\widehat{E}}_d(\sigma_{\kappa},q_\theta)$, $\mathcal{\widehat{E}}_\sigma(\sigma_\kappa)$,
$\mathcal{\widehat{E}}_b(\sigma_{\kappa})$ and $\mathcal{\widehat{E}}_q(q_\theta)$ are Monte Carlo
approximations of the terms $\mathcal{E}_d(\sigma_{\kappa},q_\theta)$, $\mathcal{E}_\sigma(\sigma_\kappa)$,
$\mathcal{E}_b(\sigma_{\kappa})$ and $\mathcal{E}_q(q_\theta)$, obtained by replacing the expectation with a
sample mean, and are defined by
\begin{align}
\mathcal{\widehat{E}}_d(\sigma_{\kappa},q_\theta)&:=n_r^{-1}|\Omega|\sum_{j=1}^{n_r} \|\sigma_{\kappa}(X_{j})
-\P01(q_\theta(X_{j}))\nabla z^\delta(X_{j})\|_{\ell^2}^2, \label{eqn:loss0}\\
\mathcal{\widehat{E}}_\sigma(\sigma_\kappa)&:=n_r^{-1}|\Omega|\sum_{j=1}^{n_r}\big(\nabla\cdot\sigma_\kappa(X_j)+f(X_j)\big)^2,\label{eqn:loss1} \\
\mathcal{\widehat{E}}_b(\sigma_\kappa)&:=n_b^{-1}|\partial\Omega|\sum_{j=1}^{n_b}\big(\sigma_\kappa(Y_j)\cdot \n-g(Y_j)\big)^2,\label{eqn:loss2} \\
\mathcal{\widehat{E}}_q(q_\theta)&:=n_r^{-1}|\Omega|\sum_{j=1}^{n_r}\|\nabla q_\theta(X_j)\|_{\ell^2}^2.\label{eqn:loss3}
\end{align}
Additionally, we define variants of $\mathcal{E}_b$ and $\widehat{\mathcal{E}}_b$ by
\begin{equation}
  \begin{aligned}
    \mathcal{E}_{b'} &=|\partial\Omega|\mathbb{E}_{\mathcal{U}(\partial\Omega)}[\|\sigma(Y)-q^\dag(Y)\nabla z^\delta(Y)\|_{\ell^2}^2],\\
    \widehat{\mathcal{E}}_{b'} &=n_b^{-1}|\partial\Omega|\sum_{j=1}^{n_b}\|\sigma(Y_j)-q^\dag(Y_j)\nabla z^\delta(Y_j)\|_{\ell^2}^2.
  \end{aligned}
  \label{eqn:loss5}
\end{equation}
The losses $\mathcal{E}_{b'}$ and $\widehat{\mathcal{E}}_{b'}$ will be needed  in Section \ref{sec:Diri}.
The estimation of Monte Carlo errors will be discussed in Section \ref{sec:Neum-MC} below.

\subsection{Error analysis of the population loss}\label{sec:Neumann-pop}
Now we derive (weighted) $L^2(\Omega)$ error bounds on the DNN approximation $q_{\theta}^*$ via the population loss $J_{\bsgamma}(\theta,\kappa)$. The \textit{a priori} regularity in Assumption \ref{ass:Neum} ensures that the exact flux $\sigma^\dag:=q^\dag\nabla u(q^\dag)$ has sufficient regularity and can be approximated via DNNs.

\begin{assumption}\label{ass:Neum}
$q^\dag\in W^{2,p}(\Omega)\cap \mathcal{A}$, $f\in H^1(\Omega)\cap L^{\infty}(\Omega)$ and $g\in H^{\frac32}(\partial\Omega)\cap L^{\infty}(\partial\Omega)$, with $p=\max(2,d+\nu)$ for small $\nu>0$.
\end{assumption}

Then the following \textit{a priori} regularity holds for $u^\dag:=u(q^\dag)$ and $\sigma^\dag$.
\begin{lemma}\label{lem:reg}
Under Assumption \ref{ass:Neum}, the
solution $u^\dag$ to problem \eqref{equ:Neu problem} satisfies $u^\dag\in H^3(\Omega)\cap H_0^1
(\Omega)$ and $\sigma^\dag\in (H^2(\Omega))^d$.
\end{lemma}
\begin{proof}
By Sobolev embedding theorem \cite[Theorem 4.12, p. 85]{AdamsFournier:2003}, $q^\dag\in W^{2,p}(\Omega)\hookrightarrow W^{1,\infty}(\Omega)$ for $p>\max(2,d+\nu)$. Since $f\in L^\infty(\Omega)$ and $g\in L^\infty (\partial\Omega)$, the standard elliptic regularity theory implies $u^\dag\in H^2(\Omega)\cap W^{1,\infty}(\Omega)$. Upon expansion, we have
\begin{equation*}
  \left\{\begin{aligned} -\Delta u^\dag &= \tfrac{f}{q^{\dag}} +\tfrac{\nabla q^\dag\cdot\nabla u^\dag}{q^\dag}, &&\quad \mbox{in }\Omega,\\
  \partial_n u^\dag & = \tfrac{g}{q^\dag}, &&\quad \mbox{on }\partial\Omega,
  \end{aligned}\right.
\end{equation*}
where $\cdot$ denotes Euclidean inner product on $\mathbb{R}^d$.
Since $q^\dag\in W^{1,\infty}(\Omega)\cap \mathcal{A}$ and $f\in L^\infty(\Omega)\cap H^1(\Omega)$, we have for $i=1,\ldots,d$,
$\partial_{x_i}(\frac{f}{q^\dag}) = -\frac{f\partial_{x_i} q^\dag}{(q^\dag)^2} + \frac{\partial_{x_i} f}{q^\dag} \in L^2(\Omega)$,
i.e., $\frac{f}{q^\dag}\in H^1(\Omega)$. Likewise,
since $q^\dag\in W^{1,\infty}(\Omega)\cap W^{2,p}(\Omega)$ and $u^\dag\in W^{1,\infty}(\Omega)\cap H^2(\Omega)$, we have for any $i=1,\ldots,d$,
\begin{equation*}
   \partial_{x_i} \Big(\frac{\nabla q^\dag\cdot\nabla u^\dag}{q^{\dag}}\Big) =
   -\frac{\partial_{x_i}q^\dag (\nabla q^\dag\cdot \nabla u^\dag)}{(q^\dag)^2} + \frac{\nabla (\partial_{x_i} q^\dag) \cdot \nabla u^\dag}{q^\dag} + \frac{\nabla (\partial_{x_i}u^\dag)\cdot \nabla q^\dag}{q^\dag} \in L^2(\Omega),
\end{equation*}
i.e., $\frac{\nabla q^\dag\cdot\nabla u^\dag}{q^\dag}\in H^1(\Omega)$ also. Since $g\in L^{\infty}(\partial\Omega)\cap H^\frac{3}{2}(\partial\Omega)$, by the trace theorem \cite[Theorem 5.36, p. 164]{AdamsFournier:2003}, there exists $\tilde g\in L^{\infty}(\Omega)\cap H^2(\Omega)$ such that $\tilde g|_{\partial\Omega}=g$. Then repeating the above argument gives for any $i,j=1,\ldots,d$
\begin{align*}
    \partial_{x_ix_j}\Big(\frac{\tilde g}{q^\dag}\Big)
    =&\frac{q^\dag\partial_{x_{i}x_j}\tilde g-\partial_{x_i}\tilde g \partial_{x_j}q^\dag}{(q^\dag)^2}-\frac{(\partial_{x_j}\tilde g)(\partial_{x_i}q^\dag) + \tilde g\partial_{x_ix_j}q^\dag}{(q^\dag)^2}\\
    &+2\frac{\tilde g(\partial_{x_i}q^\dag)(\partial_{x_j}q^\dag)}{(q^\dag)^3}\in L^2(\Omega),
\end{align*}
i.e., $\frac{\tilde g}{q^\dag}\in H^2(\Omega)$. This and the trace theorem imply $\frac{g}{q^\dag}\in H^\frac{3}{2}(\partial\Omega)$. The standard elliptic regularity theory yields $u^\dag\in H^3(\Omega)$. Finally,  $q^\dag \in W^{2,p}(\Omega)$ and $u^\dag\in H^3(\Omega)\cap W^{1,\infty}(\Omega)$ imply $\sigma^\dag\equiv q^\dag\nabla u^\dag\in(H^2(\Omega))^d$.
\end{proof}

The next lemma bounds the minimal loss value $J_{\bsgamma}(\theta^*,\kappa^*)$ (and thus its components), which will play a crucial role in the proof of Theorem \ref{thm:pop-loss-Neum}. It provides an \textit{a priori} bound on $\|\nabla q_\theta\|_{L^2(\Omega)}$ or $\|\nabla P_\mathcal{A}(q_\theta)\|_{L^2(\Omega)} $.
\begin{lemma}\label{lem:Neu pri sigma q}
Let Assumption \ref{ass:Neum} hold. Fix small $\epsilon_\sigma$, $\epsilon_q>0$, and let $(\theta^*,\kappa^*)\in( \mathfrak{P}_{p,\epsilon_q}, \mathfrak{P}_{2,\epsilon_\sigma}^{\otimes d})$ be a minimizer of problem
\eqref{eqn:obj-Neum}. Then the following estimate holds
\begin{equation*}
	J_{\bsgamma}(\theta^*,\kappa^*)\leq c\big(\epsilon_q^2+(\gamma_\sigma+\gamma_b+1)\epsilon_\sigma^2+\delta^2+\gamma_q\big).
\end{equation*}
\end{lemma}
\begin{proof}
First, Assumption \ref{ass:Neum} and Lemma \ref{lem:reg} imply $\sigma^\dag\in (H^2(\Omega))^d$. By Lemma \ref{lem:tanh-approx}, there
exists one $(\theta_\epsilon,\kappa_\epsilon) \in (\mathfrak{P}_{p,\epsilon_q},
\mathfrak{P}_{2,\epsilon_\sigma}^{\otimes d})$ such that its DNN realization $(q_{\theta_\epsilon},\sigma_{\kappa_\epsilon})$ satisfies
\begin{equation}\label{eqn:bdd-NN}
   \|q^\dag-q_{\theta_\epsilon}\|_{W^{1,p}(\Omega)}\leq\epsilon_q\quad \mbox{and}\quad  \|\sigma^\dag-\sigma_{\kappa_\epsilon}\|_{H^1(\Omega)}\leq \epsilon_\sigma.
\end{equation}
Then the triangle inequality leads to $\| \nabla q_{\theta_\epsilon} \|_{L^2(\Omega)} \le c $.
The Sobolev embedding $W^{1,p}(\Omega)\hookrightarrow L^\infty(\Omega)$ \cite[Theorem 4.12, p. 85]{AdamsFournier:2003} implies
\begin{equation}\label{eqn:bdd-1}
   \|q^\dag-q_{\theta_\epsilon}\|_{L^\infty(\Omega)}\leq c\epsilon_q.
   \end{equation}
By the minimizing property of $(\theta^*,\kappa^*)$ to $J_{\bsgamma}(\theta,\kappa)$ and triangle inequality, we obtain
\begin{align*}
&J_{\bsgamma}(\theta^*,\kappa^*)\leq J_{\bsgamma}(\theta_\epsilon,\kappa_\epsilon) \\
 =& \|\sigma_{\kappa_\epsilon}- \P01(q_{\theta_\epsilon})\nabla z^\delta\|_{L^2(\Omega)}^2+\gamma_\sigma\|\nabla\cdot\sigma_{\kappa_\epsilon}+f\|_{L^2(\Omega)}^2\\
&+\gamma_b\|\sigma_{\kappa_\epsilon}\cdot \n-g\|^2_{L^2(\partial\Omega)}+\gamma_q\|\nabla q_{\theta_\epsilon}\|_{L^2(\Omega)}^2 \\
\leq& c\big[\|\sigma_{\kappa_\epsilon}-\sigma^\dag\|^2_{L^2(\Omega)}+ \|(q^\dag-\P01(q_{\theta_\epsilon}))\nabla u^\dag\|_{L^2(\Omega)}^2+ \|\P01(q_{\theta_\epsilon})\nabla (u^\dag-z^\delta)\|_{L^2(\Omega)}^2\\
&+ \gamma_\sigma\|\nabla\cdot\sigma_{\kappa_\epsilon}-\nabla\cdot \sigma^\dag\|_{L^2(\Omega)}^2+\gamma_b\|(\sigma_{\kappa_\epsilon}-\sigma^\dag)\cdot \n\|^2_{L^2(\partial\Omega)}+\gamma_q\|\nabla q_{\theta_\epsilon}\|_{L^2(\Omega)}^2\big].
\end{align*}
The continuous embedding $W^{1,p}(\Omega)\hookrightarrow L^\infty(\Omega)$
and the stability estimates \eqref{eqn:P01-stab}--\eqref{eqn:P01-approx} of $\P01$ lead to
\begin{equation*}\label{eqn:bdd-2}
\begin{split}
  \|(q^\dag-\P01(q_{\theta_\epsilon}))\nabla u(q^\dag)\|_{L^2(\Omega)} & \leq \|q^\dag-\P01(q_{\theta_\epsilon})\|_{L^\infty(\Omega)}\|\nabla u^\dag\|_{L^2(\Omega)} \\
   &\leq \|q^\dag-q_{\theta_\epsilon}\|_{L^\infty(\Omega)}\|\nabla u^\dag\|_{L^2(\Omega)},\\
  \|\P01(q_{\theta_\epsilon})\nabla (u^\dag-z^\delta)\|_{L^2(\Omega)} & \leq \|\P01(q_{\theta_\epsilon})\|_{L^{\infty}(\Omega)}\|\nabla(u^\dag-z^\delta)\|_{L^2(\Omega)}\\
   &\le c \| \nabla(u^\dag-z^\delta)\|_{L^2(\Omega)} .
\end{split}
\end{equation*}
The \textit{a priori} estimates in \eqref{eqn:bdd-1} and the trace theorem \cite[Theorem 5.36, p. 164]{AdamsFournier:2003} yield
\begin{align*}
J_{\bsgamma}(\theta^*,\kappa^*)&\leq c\big[\epsilon_\sigma^2 + \|q^\dag-q_{\theta_\epsilon}\|^2_{L^\infty(\Omega)} +  \|\nabla (u^\dag-z^\delta)\|^{2}_{L^2(\Omega)}\\
 &\quad + (\gamma_\sigma+\gamma_b)\|\sigma^\dag-\sigma_{\kappa_\epsilon}\|^2_{H^1(\Omega)}+\gamma_q  \| \nabla q_{\theta_\epsilon} \|^2_{L^2(\Omega)} \big]\\
&\leq c\big(\epsilon_q^2+(\gamma_\sigma+\gamma_b+1)\epsilon_\sigma^2+\delta^2+\gamma_q\big).
\end{align*}
This completes the proof of the lemma.
\end{proof}

Next we give a crucial condition for the error analysis. It holds if $u^\dag\in C^2(\overline{\Omega})$ and $|\nabla u^\dag|
\neq0$ on a smooth domain $\Omega\subset\mathbb{R}^2$; see \cite[Lemma 5]{KohnLowe:1988} for detailed
discussions.
\begin{condition}\label{cond: Neu cond}
For any $\psi\in H^1(\Omega)$, there exists a solution $v_\psi$ of the equation
$\nabla u^\dag\cdot\nabla v_\psi = \psi$ almost everywhere  in $\Omega$ and it satisfies
$\|v_\psi\|_{H^1(\Omega)}\leq c\|\psi\|_{H^1(\Omega)}.$
\end{condition}

Now we can state a first error estimate on the DNN approximation $q^*_\theta$.
\begin{theorem}\label{thm:pop-loss-Neum}
Let Assumption \ref{ass:Neum} and Condition \ref{cond: Neu cond} hold. Fix small
$\epsilon_q,\epsilon_\sigma>0$, and let $(q_\theta^*,\sigma_\kappa^*)$ be the
DNN realization of a minimizer $(\theta^*,\kappa^*)\in (\mathfrak{P}_{p,\epsilon_q},
\mathfrak{P}_{2,\epsilon_\sigma}^{\otimes d})$ of problem \eqref{eqn:obj-Neum}.
Then with $\eta:=\epsilon_q+(\gamma_\sigma^{\frac12}+\gamma_b^{\frac12}+1)\epsilon_\sigma
+\delta+\gamma_q^{\frac12}$, there holds
\begin{equation*}
	\|q^\dag-\P01(q_\theta^*)\|_{L^2(\Omega)}\leq c(1+\gamma_\sigma^{-\frac14}+\gamma_b^{-\frac14} )(1+\gamma_q^{-\frac14} \eta^\frac12) \eta^{\frac12}.
\end{equation*}
\end{theorem}
\begin{proof}
Setting $\psi=q^\dag- \P01(q_\theta^*)\in H^1(\Omega)$ in Condition \ref{cond: Neu cond} and applying
integration by parts yield
\begin{align*}
	 &\quad \|q^\dag- \P01(q_\theta^*)\|^2_{L^2(\Omega)}
	=\big((q^\dag-\P01(q_\theta^*))\nabla u^\dag,\nabla v_\psi\big) = (\sigma^\dag,\nabla v_\psi) - (\P01(q^*_\theta)\nabla u^\dag,\nabla v_\psi)\\
	&=-\big(\nabla\cdot\sigma^\dagger, v_\psi\big)
	- \big(\P01(q_\theta^*)\nabla u^\dag,\nabla v_\psi\big) + (g,v_\psi)_{L^2(\partial\Omega)}  \\
	& = \big(f+\nabla\cdot\sigma_\kappa^*, v_\psi\big) + \big(\sigma_\kappa^*-\P01(q_\theta^*)\nabla z^\delta, \nabla v_\psi \big) + \big( \P01(q_\theta^*)\nabla(z^\delta- u^\dag),\nabla v_\psi\big)\\
&\quad  + (g-\sigma_{\kappa}^*\cdot \n,v_\psi)_{L^2(\partial\Omega)}.
\end{align*}
Then the Cauchy-Schwarz inequality, the trace theorem \cite[Theorem 5.36, p. 164]{AdamsFournier:2003},
Lemma \ref{lem:Neu pri sigma q} and Condition \ref{cond: Neu cond} imply
\begin{align*}
 \|q^\dag-&\P01(q_\theta^*)\|^2_{L^2(\Omega)}
\leq c\big[\|f+\nabla\cdot\sigma_\kappa^*\|_{L^2(\Omega)} + \|\sigma_\kappa^*- \P01(q_\theta^*)\nabla z^\delta\|_{L^2(\Omega)}\\
 & + \|\P01(q_\theta^*)\|_{L^{\infty}(\Omega)}\|\nabla (z^\delta-u^\dag)\|_{L^2(\Omega)}
+\|g-\sigma_{\kappa}^*\cdot \n\|_{L^2(\partial\Omega)}\big]\|v_\psi\|_{H^1(\Omega)} \\
&\leq c \big(1+\gamma_\sigma^{-\frac12}+\gamma_b^{-\frac12} \big) \eta\|q^\dag-\P01(q_\theta^*)\|_{H^1(\Omega)}.
\end{align*}
By the definition of the $(H^1(\Omega))'$-norm, we have
\begin{equation*}
   \|q^\dag- \P01(q_\theta^*)\|_{(H^1(\Omega))'}\leq c \big(1+\gamma_\sigma^{-\frac12}+\gamma_b^{-\frac12} \big) \eta.
\end{equation*}
This, duality pairing and the bound on $\|\nabla q_{\theta}^*\|_{L^2(\Omega)}$ in Lemma \ref{lem:Neu pri sigma q}
and the stability of $\P01$ in \eqref{eqn:P01-stab} imply
\begin{align*}
  \|q^\dag-\P01(q_\theta^*)\|_{L^2(\Omega)}
   \leq&  \|q^\dag-\P01(q_\theta^*)\|^{\frac12}_{(H^{1}(\Omega))'}\|q^\dag-\P01(q_\theta^*)\|^{\frac12}_{H^1(\Omega)} \\
        \leq &c \|q^\dag-\P01(q_\theta^*)\|^{\frac12}_{(H^{1}(\Omega))'} \big(1+\|\nabla \P01(q_{\theta}^*)\|_{L^2(\Omega)}^{\frac12}\big)\\
     \leq &c \|q^\dag-\P01(q_\theta^*)\|^{\frac12}_{(H^{1}(\Omega))'}\big(1+\|\nabla q_{\theta}^*\|_{L^2(\Omega)}^{\frac12}\big)\\
    	  \leq & c(1+\gamma_\sigma^{-\frac14}+\gamma_b^{-\frac14} )(1+\gamma_q^{-\frac14} \eta^\frac12) \eta^{\frac12}.
\end{align*}
Thus we complete the proof of the theorem.
\end{proof}

\begin{remark}\label{rmk:para-Neum}
Theorem \ref{thm:pop-loss-Neum} provides useful guidelines for choosing
the parameters:
$\gamma_\sigma=O(1)$, $\gamma_b=O(1)$, $\gamma_q=O(\delta^2)$, $\epsilon_q=O(\delta)$ and $ \epsilon_\sigma=O(\delta)$.
Then we obtain
$\|q^\dag- \P01(q_\theta^*)\|_{L^2(\Omega)}\leq c\delta^\frac12.$
This bound is consistent with the results using the pure FEM discretization \cite[Theorems 1 and 4]{KohnLowe:1988}. The main tool in both works is the stability estimate of the inverse problem and the approximation property of the ansatz space (DNN versus FEM). However, the meshfree nature of DNNs enables applying the method in high dimensions. This represents one distinct feature of the DNN approach. Of course, like most a priori choice rules, the analysis does not give concrete values for the hyper-parameters, and thus it is of great interest to systematically investigate hyper-parameter choice rules both empirically and theoretically.
\end{remark}
\subsection{Error analysis of the empirical loss}\label{sec:Neum-MC}

Now we analyze the statistical error of the approximation $\widehat{q}_\theta^*$, i.e., the DNN realization
of a minimizer $(\widehat\theta^*,\widehat\kappa^*)$ to the empirical loss $\widehat{J}_{\bsgamma}
(\theta,\kappa)$. The loss $\widehat{J}_{\bsgamma}(\theta,\kappa)$ involves also quadrature errors arising from approximating the integrals via Monte Carlo methods. The analysis requires
the following assumption.
The $L^\infty$ bound is needed in order to apply the standard Rademacher complexity argument (cf. Dudley's formula in Lemma \ref{lem:Dudley}).
\begin{assumption}\label{ass:Neum-gen}
 $f\in L^\infty(\Omega)$, and $z^\delta\in W^{1,\infty}(\Omega)$.
\end{assumption}

The analysis is to bound the error $\sup_{q_\theta\in\mathcal{N}_q,\sigma_\kappa
\in\mathcal{N}_\sigma} \big|J_{\bsgamma}(q_\theta,\sigma_\kappa)-\widehat{J}_{\bsgamma} (q_\theta,
\sigma_\kappa)\big|$ for suitable DNN function classes $\mathcal{N}_q$ and $\mathcal{N}_\sigma$ (corresponding to the sets $\mathfrak{P}_{p,\epsilon_q}$ and $
\mathfrak{P}_{2,\epsilon_\sigma}^{\otimes d}$, respectively), which
are also known as statistical errors in statistical learning theory
\cite{AnthonyBartlett:1999,ShalevBen:20214}. The starting point of the error analysis is the following splitting:
\begin{equation*}
\sup_{q_\theta\in\mathcal{N}_q,\sigma_\kappa\in\mathcal{N}_\sigma} \big|J_{\bsgamma}(q_\theta,\sigma_\kappa)-\widehat{J}_{\bsgamma}
(q_\theta,\sigma_\kappa)\big|\leq \Delta\mathcal{E}_d + \gamma_\sigma \Delta\mathcal{E}_\sigma + \gamma_b\Delta\mathcal{E}_b + \gamma_q \Delta\mathcal{E}_q,
\end{equation*}
with the error components given respectively by
\begin{align*}
&\Delta\mathcal{E}_d: = \sup_{\sigma_{\kappa}\in\mathcal{N}_{\sigma},q_\theta\in\mathcal{N}_q}\big| \mathcal{E}_d(\sigma_{\kappa},q_\theta)-\mathcal{\widehat{E}}_d(\sigma_{\kappa},q_\theta)\big|,
&&\Delta\mathcal{E}_\sigma:= \sup_{\sigma_{\kappa}\in\mathcal{N}_{\sigma}}\big|\mathcal{E}_\sigma(\sigma_\kappa)-\mathcal{\widehat{E}}_\sigma(\sigma_\kappa)\big|, \\ &\Delta\mathcal{E}_b:=\sup_{\sigma_{\kappa}\in\mathcal{N}_{\sigma}}\big|\mathcal{E}_b (\sigma_\kappa)-\mathcal{\widehat{E}}_b(\sigma_\kappa)\big|, &&\Delta\mathcal{E}_q:=\sup_{q_\theta\in\mathcal{N}_q}\big|\mathcal{E}_q(q_\theta)-\mathcal{\widehat{E}}_q(q_\theta)\big|.
\end{align*}
Further we define $\Delta\mathcal{E}_{b'}:=\sup_{\sigma_{\kappa}\in\mathcal{N}_{\sigma}}\big|\mathcal{E}_{b'} (\sigma_\kappa)-\mathcal{\widehat{E}}_{b'}(\sigma_\kappa)\big|.$

Now we state the quadrature error for each term (in high probability). The proof uses PAC-type generalization bounds in Lemma \ref{lem:PAC}, which in turn employs the Rademacher complexity of the function classes $\mathcal{N}_q$ and $\mathcal{N}_\sigma$ via Dudley's formula in Lemma \ref{lem:Dudley}. The overall proof is lengthy and hence it is deferred to the appendix.
\begin{theorem}\label{thm:stat-err}
Let Assumptions \ref{ass:Neum} and \ref{ass:Neum-gen} hold, and $q_\theta\in \mathcal{N}(L_q,N_{\theta},R_q)$ and $\sigma_\kappa \in \mathcal{N}(L_\sigma,N_{\kappa},R_\sigma)$. Fix $\tau\in(0,\frac18)$, and define the following bounds
\begin{align*}
e_d&:=c\frac{R_\sigma^{2}N_\kappa^{2}(N_\kappa+N_\theta)^{\frac12}(\log^\frac12 R_\sigma +\log^\frac12 N_{\kappa} + \log^\frac12 R_q+\log^\frac12 N_\theta+\log^\frac12 n_r)}{\sqrt{n_r}}\\
    &\quad + \tilde c R^2_\sigma N_{\kappa}^2\sqrt{\frac{\log\frac{1}{\tau}}{n_r}},\\
e_\sigma&:=c\frac{R_\sigma^{2L_\sigma}N_{\kappa}^{2L_\sigma-\frac32}
        \big(\log^\frac12R_\sigma+\log^\frac12N_{\kappa}+\log^\frac12n_r\big)}{\sqrt{n_r}}+ \tilde cR_\sigma^{2L_\sigma}N_\kappa^{2L_\sigma-2}\sqrt{\frac{\log\frac{1}{\tau}}{n_r}},\\
e_b&:=c\frac{R_\sigma^2N_{\kappa}^{\frac52}\big(\log^\frac12 R_\sigma+\log^\frac12 N_{\kappa}+\log^\frac12 n_b\big)}{\sqrt{n_b}}+\tilde cR_\sigma^{2}N_\kappa^{2}\sqrt{\frac{\log\frac{1}{\tau}}{n_b}}, \\
e_{b'}&:=c\frac{R_\sigma^2N_{\kappa}^{\frac52}\big(\log^\frac12 R_\sigma+\log^\frac12 N_{\kappa}+\log^\frac12 n_b\big)}{\sqrt{n_b}}+\tilde cR_\sigma^{2}N_\kappa^{2}\sqrt{\frac{\log\frac{1}{\tau}}{n_b}},\\
e_q&:=c\frac{R_q^{2L_q}N_\theta^{2L_q-\frac32}
   \big(\log^\frac12R_q+\log^\frac12N_{\theta}+\log^\frac12n_r\big)}{\sqrt{n_r}}+\tilde cR_q^{2L_q}N_\theta^{2L_q-2}\sqrt{\frac{\log\frac{1}{\tau}}{n_r}},
\end{align*}
where the constants $c$ and $\tilde c$ may depend on  $|\Omega|$, $|\partial\Omega|$, $d$, $\|z^\delta\|_{W^{1,\infty}(\Omega)}$, $\|f\|_{L^\infty(\Omega)}$, and $\|g\|_{L^\infty(\partial\Omega)}$ at most polynomially. Then, with probability
at least $1-\tau$, each of the following statements holds
\begin{equation*}
    \Delta \mathcal{E}_i \leq e_{i}, \quad i\in\{d,\sigma,b,b',q\}.
\end{equation*}
\end{theorem}

Now we state an error estimate on the DNN approximation $\widehat{q}_\theta^*$ via the loss $\widehat{J}_{\bsgamma}$.
\begin{theorem}\label{thm:err-Neum-emp}
Let Assumptions \ref{ass:Neum} and \ref{ass:Neum-gen} and Condition \ref{cond: Neu cond} hold.
Fix small $\epsilon_q$, $\epsilon_\sigma>0$, and let $(\widehat{\theta}^*,\widehat{\kappa}^*)\in
(\mathfrak{P}_{p,\epsilon_q}, \mathfrak{P}_{2,\epsilon_\sigma}^{\otimes d})$ be a
minimizer of the empirical loss \eqref{eqn:obj-Neum-dis}, and $\widehat{q}_\theta^*$ and
$\widehat{\sigma}_\kappa^*$ their NN realizations. Fix $\tau\in (0,\frac18)$, let the bounds $e_d$,
$e_\sigma$, $e_b$ and $e_q$ be defined in Theorem \ref{thm:stat-err}, and further define $\eta$ by
$ \eta:=e_d+\gamma_\sigma e_\sigma+\gamma_be_b+\gamma_qe_q+\epsilon_q^2+(\gamma_\sigma+\gamma_b+1)\epsilon_\sigma^2
+\delta^2+\gamma_q.$
Then with probability at least $1-4\tau$, there holds
\begin{align*}
\|q^\dag-\P01(\widehat{q}_\theta^*)\|_{L^2(\Omega)}\leq& c \big((e_d+\eta)^\frac12+(e_\sigma+\gamma_\sigma^{-1}\eta)^\frac12
+(e_b+\gamma_b^{-1}\eta)^\frac12+(e_q+\gamma_q^{-1}\eta)^\frac12\delta\big)^\frac12\\
&\times\big(1+(e_q+\gamma_q^{-1}\eta)^\frac14\big).
\end{align*}
\end{theorem}
\begin{proof}
Let $(\theta^*,\kappa^*)\in (\mathfrak{P}_{p,\epsilon_q}, \mathfrak{P}_{2,\epsilon_\sigma}^{\otimes d})$ be a minimizer of problem \eqref{eqn:obj-Neum}. Then the
minimizing property of $(\widehat{\theta}^*,\widehat{\kappa}^*)$ to the empirical
loss $\widehat J_{\bsgamma}(\theta,\kappa)$ implies
\begin{align*}
  \widehat J_{\bsgamma}(\widehat\theta^*,\widehat\kappa^*)
   & \leq [\widehat J_{\bsgamma}(\theta^*,\kappa^*) - J_{\bsgamma}(\theta^*,\kappa^*)] + J_{\bsgamma}(\theta^*,\kappa^*)\\
   &\leq |J_{\bsgamma}(\theta^*,\kappa^*) - \widehat J_{\bsgamma}(\theta^*,\kappa^*)| +  J_{\bsgamma}(\theta^*,\kappa^*).
\end{align*}
Consequently, we deduce
\begin{equation*}
   \widehat{J}_{\bsgamma}(\widehat\theta^*,\widehat\kappa^*)  \leq J_{\bsgamma}(\theta^*,\kappa^*) + \sup_{(\theta,\kappa)\in (\mathfrak{P}_{p,\epsilon_q}, \mathfrak{P}_{2,\epsilon_\sigma}^{\otimes d})} |J_{\bsgamma}(\theta,\kappa) - \widehat J_{\bsgamma}(\theta,\kappa)|.
\end{equation*}
The two terms represent the approximation error and statistical error, respectively,
and the former was already bounded in Section \ref{sec:Neumann-pop}. By Lemma \ref{lem:Neu pri sigma q} and Theorem \ref{thm:stat-err},
with a probability at least $1-4\tau$, $\widehat{J}_{\bsgamma}(\widehat{\theta}^*,\widehat{\kappa}^*)\leq c\eta$.
This and the triangle inequality imply that with probability at least $1-4\tau$,
\begin{equation*}
\|\widehat{\sigma}_\kappa^*-\P01(\widehat{q}_\theta^*)\nabla z^\delta\|_{L^2(\Omega)}^2\leq \big[\mathcal{E}_d(\widehat{\sigma}^*_{\kappa},\widehat{q}^*_\theta) -\mathcal{\widehat{E}}_d(\widehat{\sigma}^*_{\kappa},\widehat{q}^*_\theta)\big]+  \mathcal{\widehat{E}}_d(\widehat{\sigma}^*_{\kappa},\widehat{q}^*_\theta)\leq c(e_d+\eta).
\end{equation*}
Similarly, the following estimates hold simultaneously with a probability at least $1-4\tau$,
$\|f+\nabla\cdot\widehat{\sigma}_\kappa^*\|_{L^2(\Omega)}^2\leq c (e_\sigma+\gamma_\sigma^{-1}\eta)$,
$\|g-\widehat{\sigma}_{\kappa}^*\cdot \n\|_{L^2(\partial\Omega)}^2\leq c(e_b+\gamma_b^{-1}\eta)$,
$\|\nabla \widehat{q}_\theta^*\|^2_{L^2(\Omega)}\leq c(e_q+\gamma_q^{-1}\eta)$.
Replacing $q_\theta^*$ by $\widehat{q}_\theta^*$ and repeating the argument of Theorem \ref{thm:pop-loss-Neum} give
\begin{align*}
\|q^\dag-&\P01(\widehat{q}_\theta^*)\|^2_{L^2(\Omega)}\le c\big(\|f+\nabla\cdot\widehat{\sigma}_\kappa^*\|_{L^2(\Omega)} + \|\widehat{\sigma}_\kappa^*-\P01(\widehat{q}_\theta^*)\nabla z^\delta\|_{L^2(\Omega)}\\
&+ \|\P01(\widehat{q}_\theta^*)\|_{L^{\infty}(\Omega)}\|\nabla (z^\delta-u^\dagger)\|_{L^{2}(\Omega)}+\|g-\widehat{\sigma}_{\kappa}^*\cdot \n\|_{L^2(\partial\Omega)}\big)\|v_\psi\|_{H^1(\Omega)}.
\end{align*}
This, the preceding estimates and Condition \ref{cond: Neu cond} imply (with $v_\psi$ from Condition \ref{cond: Neu cond} for $\psi=q^\dag-\P01(\widehat{q}_\theta^*)$)
\begin{equation*}
\|q^\dag-\P01(\widehat{q}_\theta^*)\|_{H^{-1}(\Omega)}\leq c\big((e_d+\eta)^\frac12+(e_\sigma+\gamma_\sigma^{-1}\eta)^\frac12+(e_b+\gamma_b^{-1}\eta)^\frac12 +(e_q+\gamma_q^{-1}\eta)^\frac12\delta\big).
\end{equation*}
The desired estimate now follows from the argument of Theorem \ref{thm:pop-loss-Neum} and the \textit{a priori} bound on the NN approximation $\|\nabla\widehat{q}_\theta^*\|^2_{L^2(\Omega)}$.
\end{proof}

\begin{remark}\label{rmk:para-Neum-emp}
Comparing Theorem \ref{thm:pop-loss-Neum} with Theorem \ref{thm:err-Neum-emp}, especially the expression of $\eta$, implies that the difference between the error of the approximations by the population loss $J_{\bsgamma}$  and empirical loss $\widehat{J}_{\bsgamma}$ lies in the quadrature error arising from the Monte Carlo approximation of the integrals. Similar to Remark \ref{rmk:para-Neum}, Theorem \ref{thm:err-Neum-emp} gives guidelines for choosing the number of sampling points in the domain $\Omega$ and on the boundary $\partial\Omega$.
Under the assumptions of Theorem \ref{thm:stat-err} and the parameter choices in Remark \ref{rmk:para-Neum}, we may choose the numbers $n_r$ and $n_b$ of sample points in the domain $\Omega$ and on the boundary $\partial\Omega$ by
$n_r = O( \max(\frac{R_\sigma^4N_\kappa^4(N_\theta+N_\kappa)}{\delta^{4+s}}, \frac{R_\sigma^{4L_\sigma}N_\kappa^{4L_\sigma-3}}{\delta^{4+s}}, \frac{R_q^{4L_q}N_\kappa^{4L_q-3}}{\delta^{s}}))$ and $n_b=O(\frac{R_\sigma^4N_\kappa^5}{\delta^{4+s}})$,
where the exponent $s>0$ is to absorb the $\log$ factor. Note that by Lemma \ref{lem:tanh-approx}, $R_q=O(\delta^{-\frac{4p+3d+3pd}{p(1-\mu)}})$, $N_\theta=O(\delta^{-\frac{d}{1-\mu}})$, $R_\sigma=O(\delta^{-\frac{8+9d}{2(1-\mu)}})$ and $N_\kappa=O(\delta^{-\frac{d}{1-\mu}})$. Then with probability at least $1-4\tau$, we have
    \begin{equation*}
    \|q^\dag- \P01(\widehat{q}_\theta^*)\|_{L^2(\Omega)}\leq c\delta^\frac12.
\end{equation*}
\end{remark}

\section{Inverse conductivity problem in the Dirichlet case}
\label{sec:Diri}

In this section, we extend the approach to the following Dirichlet problem:
\begin{equation}\label{eqn:Diri}
	\left\{\begin{aligned}
		-\nabla\cdot(q\nabla u) &= f, \ &\mbox{in}&\ \Omega, \\
		u&=0, \ &\mbox{on}&\ \partial\Omega,
	\end{aligned}\right.
\end{equation}
and provide relevant error analysis of the reconstruction scheme.

\subsection{Mixed formulation and its DNN approximation}
We recast problem \eqref{eqn:Diri} into a first-order system:
\begin{equation}\label{eqn:mixed-Diri}
	\left\{\begin{aligned}
		\sigma & = q\nabla u, &&\mbox{in }\ \Omega, \\
		-\nabla\cdot\sigma &= f, &&\mbox{in}\ \Omega, \\
		u&=0,  &&\mbox{on}\ \partial\Omega.
	\end{aligned}\right.
\end{equation}
To recover $q$ , we use a noisy measurement $z^\delta$ (of the exact data $u(q^\dagger)$) in $\Omega$ with a noise level
$\delta:=\|u(q^\dagger)-z^\delta\|_{W^{\frac32,2}(\Omega)}$. Then we discretize $q$ and $\sigma$ by two DNN function classes $\mathfrak{P}_{p,\epsilon_q}$ and $ \mathfrak{P}_{2,\epsilon_\sigma}^{\otimes d}$. Following Section \ref{sec:Neumann}, we may construct a DNN approximation via the population loss
\begin{equation}\label{eqn:obj-Diri}
J_{\bsgamma}(\theta,\kappa)=\|\sigma_\kappa-\P01(q_\theta)\nabla z^\delta\|_{L^2(\Omega)}^2+\gamma_\sigma\|\nabla\cdot\sigma_\kappa+f\|_{L^2(\Omega)}^2+\gamma_q\|\nabla q_\theta\|_{L^2(\Omega)}^2,
\end{equation}
In practice, however, the formulation
\eqref{eqn:obj-Diri} does not lend itself to high-quality reconstructions for highly noisy data $z^\delta$. This is
attributed to a lack of knowledge of the flux $\sigma$ on the boundary $\partial\Omega$ so that the first term does not allow learning $\sigma$ accurately. Hence, we augment the loss $\mathcal{J}_{\bsgamma}$ in \eqref{eqn:obj-Diri} with an additional boundary term
\begin{equation}\label{eqn:obj-Diri1}
\begin{aligned}
J_{\bsgamma}(\theta,\kappa)
=&\|\sigma_\kappa-\P01(q_\theta)\nabla z^\delta\|_{L^2(\Omega)}^2+\gamma_\sigma\|\nabla\cdot\sigma_\kappa+f\|_{L^2(\Omega)}^2\\
&+\gamma_q\|\nabla q_\theta\|_{L^2(\Omega)}^2 +\gamma_b\|\sigma_\kappa-q^\dag\nabla z^\delta\|^2_{L^2(\partial\Omega)}.
\end{aligned}
\end{equation}
 Note that the
assumption on a knowledge $q^\dag|_{\partial\Omega}$ is often made in existing mathematical analysis \cite{Alessandrini:1986,Richter:1981} and numerical studies \cite{BaoYeZang:2020,BarSochen:2021}. Similar to the Neumann case, one can ensure the existence of a minimizer of the loss in \eqref{eqn:obj-Diri1}. Let $(\theta^*,\kappa^*)$ be a minimizer of the loss
\eqref{eqn:obj-Diri1} and  $(q_\theta^*,\sigma_{\kappa}^*)$ be its DNN realization. In practice, we approximate the integrals using Monte Carlo methods. Using the uniform
distributions $\mathcal{U}(\Omega)$ and $\mathcal{U}(\partial\Omega)$, we rewrite the population loss \eqref{eqn:obj-Diri1} as
\begin{align*}
J_{\bsgamma}(\theta,\kappa)&=|\Omega|\mathbb{E}_{X\sim\mathcal{U}(\Omega)} \Big[\|\sigma_{\kappa}(X)-\P01(q_\theta(X))\nabla z^\delta(X)\|_{\ell^2}^2 \Big] \\ &\quad +\gamma_\sigma|\Omega|\mathbb{E}_{X\sim\mathcal{U}(\Omega)}\Big[\big(\nabla\cdot\sigma_\kappa(X)+f(X)\big)^2\Big]\\
&\quad+\gamma_b|\partial\Omega|\mathbb{E}_{Y\sim\mathcal{U}(\partial\Omega)}\Big[\|\sigma_{\kappa}-q^\dagger\nabla z^\delta(Y)\|_{\ell^2}^2\Big]+\gamma_q|\Omega|\mathbb{E}_{X\sim\mathcal{U}(\Omega)}\Big[ \|\nabla q_\theta(X)\|_{\ell^2}^2 \Big] \\
	&=: \mathcal{E}_d(\sigma_{\kappa},q_\theta) + \gamma_\sigma \mathcal{E}_{\sigma}(\sigma_\kappa)+\gamma_b\mathcal{E}_{b'}(\sigma_{\kappa})+\gamma_q\mathcal{E}_q(q_\theta).
\end{align*}
Upon drawing i.i.d. samples $X=\{X_{j}\}_{j=1}^{n_r}\sim \mathcal{U}(\Omega)$ and $Y=\{Y_{j}
\}_{j=1}^{n_b}\sim \mathcal{U}(\partial\Omega)$, we obtain
the empirical loss
\begin{equation}\label{equ: Diri empir loss}
	\widehat{J}_{\bsgamma}(\theta,\kappa):= \mathcal{\widehat{E}}_d(\sigma_{\kappa},q_\theta) + \gamma_\sigma \mathcal{\widehat{E}}_\sigma(\sigma_\kappa)+\gamma_q\mathcal{\widehat{E}}_q(q_\theta)+\gamma_b\mathcal{\widehat{E}}_{b'}(\sigma_{\kappa}),
\end{equation}
where $\mathcal{\widehat{E}}_d(\sigma_{\kappa},q_\theta)$, $\mathcal{\widehat{E}}_\sigma(\sigma_\kappa)$
and $\mathcal{\widehat{E}}_q(q_\theta)$ are given by \eqref{eqn:loss0}--\eqref{eqn:loss2}, and
$\mathcal{\widehat{E}}_{b'}(\sigma_{\kappa,i})$ in \eqref{eqn:loss5}.

\begin{remark}\label{rmk:loss}
Instead of \eqref{eqn:obj-Diri1}, there are alternative formulations of the concerned inverse problem. For example, one may use the population loss \begin{equation}\label{eqn:diriloss1}
\begin{aligned}
J_{\bsgamma}(\theta,\kappa)=&\|\sigma_\kappa-\P01(q_\theta)\nabla z^\delta\|_{L^2(\Omega)}^2+\gamma_\sigma\|\nabla\cdot\sigma_\kappa+f\|_{L^2(\Omega)}^2\\
  &+\gamma_q\|\nabla q_\theta\|_{L^2(\Omega)}^2 + \gamma_b\|q_\theta-q^\dag\|_{L^2(\partial\Omega)}^2,
\end{aligned}
\end{equation}
which enforces the boundary condition $q_\theta=q^\dag$ on the boundary $\partial\Omega$ directly. It can be analyzed analogously.
However, it appears less robust than \eqref{eqn:obj-Diri1} for noisy data.
See Fig. \ref{fig:diri1sb} for numerical illustrations.
\end{remark}

\subsection{Error analysis of the population loss}

First we analyze the error of the DNN realization $q_\theta^*\in\mathcal{N}_q$ of a minimizer $\theta^*$ to the population loss \eqref{eqn:obj-Diri1}.
\begin{assumption}\label{ass:Diri}
$q^\dag\in W^{2,p}(\Omega)\cap \mathcal{A}$, $p=\max(2,d+\nu)$ for some $\nu>0$, and $f\in H^1(\Omega)\cap L^{\infty}(\Omega)$.
\end{assumption}

The following regularity result holds for $u^\dag:=u(q^\dag)$ and $\sigma^\dag:=q^\dag \nabla u^\dag$. The proof follows exactly as Lemma \ref{lem:reg}, and hence it is omitted.
\begin{lemma}\label{lem:reg-Neum}
Let Assumption \ref{ass:Diri} hold. Then the solution $u^\dag$ to problem \eqref{eqn:Diri} satisfies $u^\dag\in H^3(\Omega)\cap W^{2,p}(\Omega)\cap H_0^1(\Omega)$ and $\sigma^\dag\in (H^2(\Omega))^d$.
\end{lemma}

The next result gives an \textit{a priori} estimate on the minimal loss value $J_{\bsgamma}(\theta^*,\kappa^*)$.
\begin{lemma}\label{lem: Diri pri sigma q}
Let Assumption \ref{ass:Diri} hold. Fix small $\epsilon_q$, $\epsilon_\sigma>0$, and let
$(\theta^*,\kappa^*)\in (\mathfrak{P}_{p,\epsilon_q}, \mathfrak{P}_{2,\epsilon_\sigma
}^{\otimes d})$ be a minimizer of the loss $J_{\bsgamma}(\theta,\kappa)$ in \eqref{eqn:obj-Diri1}. Then there holds
\begin{equation*}
	J_{\bsgamma}(\theta^*,\kappa^*)\leq c\big(\epsilon_q^2+(1+\gamma_\sigma+\gamma_b)\epsilon_\sigma^2+(1+\gamma_\sigma)\delta^2+\gamma_q\big).
\end{equation*}
\end{lemma}
\begin{proof}
Assumption \ref{ass:Neum} and Lemma \ref{lem:reg-Neum} imply $\sigma^\dag\in (H^2(\Omega))^d$. Then by Lemma \ref{lem:tanh-approx}, there exists at least one $(\theta_\epsilon,\kappa_\epsilon)\in
(\mathfrak{P}_{p,\epsilon_q}, \mathfrak{P}_{2,\epsilon_\sigma
}^{\otimes d})$ such that its DNN realization $(q_{\theta_\epsilon},\sigma_{\kappa_\epsilon})$ satisfies
\begin{equation*}
\|q^\dag-q_{\theta_\epsilon}\|_{W^{1,p}(\Omega)}\leq\epsilon_q\quad \mbox{and}\quad  \|\sigma^\dag-\sigma_{\kappa_\epsilon}\|_{H^1(\Omega)}\leq \epsilon_\sigma.
\end{equation*}
Then by the minimizing property of $(\theta^*,\kappa^*)$ and the triangle inequality, we have
\begin{align*}
	& J_{\bsgamma}(\theta^*,\kappa^*)\leq J_{\bsgamma}(\theta_\epsilon,\kappa_\epsilon) \\
	 =& \|\sigma_{\kappa_\epsilon}-\P01(q_{\theta_\epsilon})\nabla z^\delta\|_{L^2(\Omega)}^2+\gamma_\sigma\|\nabla\cdot\sigma_{\kappa_\epsilon}+f\|_{L^2(\Omega)}^2\\
  &+\gamma_b\|\sigma_{\kappa_\epsilon}-q^\dag\nabla z^\delta\|^2_{L^2(\partial\Omega)} +\gamma_q\|\nabla q_{\theta_\epsilon}\|_{L^2(\Omega)}^2 \\
\leq& c\big(\|\sigma_{\kappa_\epsilon}-\sigma^\dag\|^2_{L^2(\Omega)}+ \|(q^\dag-\P01(q_{\theta_\epsilon}))\nabla u^\dag\|_{L^2(\Omega)}^2+ \|\P01(q_{\theta_\epsilon})\nabla (u^\dag-z^\delta)\|_{L^2(\Omega)}^2 \\
 &+\gamma_\sigma\|\nabla\cdot(\sigma_{\kappa_\epsilon}- \sigma^\dag)\|_{L^2(\Omega)}^2+\gamma_q\|\nabla q_{\theta_\epsilon}\|_{L^2(\Omega)}^2\\
 &+\gamma_b\|\sigma_{\kappa_\epsilon}-\sigma^\dag\|^2_{H^1(\Omega)} +\gamma_b\|q^\dag\nabla(u^\dag-z^\delta)\|^2_{L^2(\partial\Omega)}\big).
\end{align*}
Now appealing to Assumption \ref{ass:Diri} yields
\begin{align*}
	J_{\bsgamma}(\theta^*,\kappa^*)	&\leq c\big(\epsilon_\sigma^2 +
\|q^\dag-\P01(q_{\theta_\epsilon})\|^2_{L^\infty(\Omega)}\|\nabla u^\dag\|^2_{L^2(\Omega)} \\
 &\quad + \|\P01(q_{\theta_\epsilon})\|^{2}_{L^\infty(\Omega)}\|\nabla(u^\dag-z^\delta)\|^{2}_{L^2(\Omega)}+ \gamma_\sigma\epsilon_\sigma^2+\gamma_q+\gamma_b\epsilon_\sigma^2+\gamma_b\delta^2\big)\\
		&\leq c\big(\epsilon_q^2+(1+\gamma_\sigma+\gamma_b)\epsilon_\sigma^2+(1+\gamma_b)\delta^2+\gamma_q\big).
\end{align*}
This completes the proof of the lemma.
\end{proof}

To derive an $L^2(\Omega)$ error bound, we need the following positivity condition \cite{jin2021error}. It holds with $\beta=2$ if $q^\dag\in\mathcal{A}$, and $f \in  L^2(\Omega)$ with $f\geq c_f>0$  (with $c_f\in\mathbb{R}$) over a Lipschitz domain $\Omega$, and with $\beta=0$ when $q^\dagger\in C^{1,\alpha}(\overline{\Omega})\cap\mathcal{A}$, and $f\in C^{0,\alpha}(\overline{\Omega})$  and $f\geq c_f>0$ on a $C^{2,\alpha}$ domain $\Omega$ for some $\alpha>0$ \cite[Lemmas 3.3 and 3.7]{bonito2017diffusion}.
\begin{condition}\label{Cond:P}
There exist some $\beta>0$ and $c$ such that
$q^\dag |\nabla u^\dag|^2+fu^\dag\geq c\ {\rm dist}(x,\partial\Omega)^\beta$ a.e. in $\Omega$.
\end{condition}

Now we can state an error estimate on the DNN approximation $q_\theta^*$ constructed via the population loss $J_{\bsgamma}(\theta,\gamma)$ in \eqref{eqn:obj-Diri1}. The proof of the theorem differs from that in the Neumann case in Theorem \ref{thm:pop-loss-Neum}.
\begin{theorem}\label{thm: Diri error estimate on conduc}
Let Assumption \ref{ass:Diri} hold. Fix small $\epsilon_q$, $\epsilon_\sigma>0$, and let
$(\theta^*,\kappa^*) \in (\mathfrak{P}_{p,\epsilon_q}, \mathfrak{P}_{2,\epsilon_\sigma
}^{\otimes d})$ be a minimizer of the loss $J_{\bsgamma}(\theta,\kappa)$ in \eqref{eqn:obj-Diri1} and $q_\theta^*$  the DNN realization
of $\theta^*$. Then with
$\eta:=\big(\epsilon_q^2+(1+\gamma_\sigma+\gamma_b)\epsilon_\sigma^2+(1+\gamma_b)\delta^2+\gamma_q\big)^\frac12$,
the following weighted error estimate holds
\begin{equation*}
\int_{\Omega}\Big(\frac{q^\dag-\P01(q_\theta^*)}{q^\dag}\Big)^2\big(q^\dag |\nabla u^\dag|^2+fu^\dag\big)\ \mathrm{d}x\leq  c\big(\gamma_\sigma^{-\frac12}+\gamma_q^{-\frac12}\eta+1\big) \eta.
\end{equation*}
Moreover, if Condition \ref{Cond:P} holds, then
\begin{equation*}
\|q^\dag-\P01(q_\theta^*)\|_{L^2(\Omega)}\leq c[\big(\gamma_\sigma^{-\frac12}+\gamma_q^{-\frac12}\eta+1\big) \eta]^{\frac{1}{2(\beta+1)}}.
\end{equation*}
\end{theorem}
\begin{proof}
For any test function $\varphi\in H_0^1(\Omega)$, by integration by parts, we have
\begin{align*}
&\big((q^\dag-\P01(q_\theta^*))\nabla u^\dag,\nabla\varphi\big) =\big (\sigma^\dag-\sigma_\kappa^*,\nabla\varphi\big)+\big(\P01(q_\theta^*)\nabla(z^\delta-u^\dag),\nabla\varphi\big)\\
   &\quad+\big(\sigma_\kappa^*-\P01(q_\theta^*)\nabla z^\delta,\nabla\varphi\big)  \\
		=& -\big (\nabla\cdot(\sigma^\dag-\sigma_\kappa^*),\varphi\big)+\big(\P01(q_\theta^*)\nabla(z^\delta-u^\dag), \nabla\varphi\big)
       +\big(\sigma_\kappa^*-\P01(q_\theta^*)\nabla z^\delta,\nabla\varphi\big).
\end{align*}
Let $\varphi\equiv\frac{q^\dag-\P01(q_\theta^*)}{q^\dag}u^\dag$. Then by direct computation, we have
\begin{equation*}
	\nabla\varphi = \frac{\nabla(q^\dag-\P01(q_\theta^*))}{q^\dag}u^\dag+\frac{(q^\dag-\P01(q_\theta^*))}{q^\dag}\nabla u^\dag- \frac{(q^\dag-\P01(q_\theta^*))}{(q^\dag)^2}(\nabla q^\dag) u^\dag.
\end{equation*}
Using Assumption \ref{ass:Diri}, the box constraint on $q^\dagger$ and $\P01(q_\theta^*)$, and the stability estimate \eqref{eqn:P01-stab} of the projection operator $P_\mathcal{A}$, we arrive at
\begin{align*}
&\Big\|\frac{\nabla(q^\dag-\P01(q_\theta^*))}{q^\dag}u^\dag\Big\|_{L^2(\Omega)}
\leq c\big(1+\|\nabla \P01(q_\theta^*)\|_{L^2(\Omega)}\big)
\leq c\big(1+\|\nabla q_\theta^*\|_{L^2(\Omega)}\big),\\
&\Big\|\frac{(q^\dag-\P01(q_\theta^*))}{q^\dag}\nabla u^\dag\Big\|_{L^2(\Omega)}+\Big\|\frac{(q^\dag-\P01(q_\theta^*))}{(q^\dag)^2}(\nabla q^\dag) u^\dag\Big\|_{L^2(\Omega)}\leq c .
\end{align*}
This implies $\varphi\in H_0^1(\Omega)$ with the following \textit{a priori} bound
\begin{equation*}
	\|\varphi\|_{L^2(\Omega)} \leq c \quad \text{and}\quad  \|\nabla\varphi\|_{L^2(\Omega)}\leq c(1+\|\nabla q_\theta^*\|_{L^2(\Omega)}).
\end{equation*}
By Lemma \ref{lem: Diri pri sigma q} and the Cauchy-Schwarz inequality, we have
\begin{equation*}
|(\nabla\cdot(\sigma^\dagger-\sigma_\kappa^*),\varphi\big)|\leq  \|\nabla\cdot(\sigma^\dagger-\sigma_\kappa^*)\|_{L^2(\Omega)}\|\varphi\|_{L^2(\Omega)}\leq c\gamma_\sigma^{-\frac12}  \eta.
\end{equation*}
Similarly, we deduce
\begin{equation*}
		|\big(\sigma_\kappa^*-\P01(q_\theta^*)\nabla z^\delta,\nabla\varphi\big)|\leq \|\sigma_{\kappa}^*-P_\mathcal{A}(q_\theta^*)\nabla z^\delta\|_{L^2(\Omega)}\|\nabla\varphi\|_{L^2(\Omega)}\leq c (1+\gamma_q^{-\frac12}\eta)\eta.
\end{equation*}
Meanwhile, by the Cauchy--Schwartz inequality,
\begin{align*}
|\big(\P01(q_\theta^*)\nabla(z^\delta-u^\dagger), \nabla\varphi\big)|&\leq c\|\P01(q_\theta^*)\|_{L^{\infty}(\Omega)}
\|\nabla(z^\delta-u^\dagger)\|_{L^2(\Omega)}\|\nabla\varphi\|_{L^2(\Omega)} \\
 &\leq c(1+\gamma_q^{-\frac12}\eta)\delta.
\end{align*}
Upon repeating the argument in \cite{bonito2017diffusion,jin2021error}, we obtain
\begin{equation*}
	\big((q^\dag-\P01(q_\theta^*))\nabla u^\dag,\nabla\varphi\big) = \frac12\int_{\Omega}\Big(\frac{q^\dag-\P01(q_\theta^*)}{q^\dag}\Big)^2\big(q^\dag|\nabla u^\dagger|^2+fu^\dagger\big)\ \mathrm{d}x.
\end{equation*}
Combining the preceding estimates yields the first assertion.
Next we bound $\|q^\dagger- \P01(q_\theta^*)\|_{L^2(\Omega)}$. Upon fixing $\rho>0$, we can split the domain $\Omega$ into two disjoint sets $\Omega=\Omega_\rho\cup \Omega_\rho^c$, with $\Omega_\rho=\{x\in\Omega:{\rm dist}(x,\partial\Omega)\geq\rho\}$ and $\Omega_\rho^c=\Omega\setminus\Omega_\rho$. Then by the box constraint $q^\dag\in\mathcal{A}$, we have
\begin{align*}
	\|q^\dagger-\P01(q_\theta^*)\|_{L^2(\Omega_\rho)}^2 &= \rho^{-\beta}\int_{\Omega_\rho }(q^\dagger-\P01(q_\theta^*))^2\rho^{\beta}{\rm d}x\\
    &\leq \rho^{-\beta}\int_{\Omega_\rho }(q^\dag-\P01(q_\theta^*))^2\mathrm{dist}(x,\partial\Omega)^{\beta}{\rm d}x\\
	& \leq c\rho^{-\beta}\int_{\Omega_\rho }\Big(\frac{q^\dag-\P01(q_\theta^*)}{q^\dag}\Big)^2\big(q^\dag |\nabla u^\dag|^2+fu^\dag\big){\rm d}x\\
    &\leq c\rho^{-\beta}\big(\gamma_\sigma^{-\frac12}+\gamma_q^{-\frac12}\eta+1\big)\eta.
\end{align*}
Meanwhile, using the box constraint  $q^\dag,\P01(q_\theta^*)\in\mathcal{A}$, we obtain
\begin{equation*}
\|q^\dag-\P01(q_\theta^*)\|_{L^2(\Omega^c_\rho)}^2
\leq  c \int_{\Omega_\rho^c}1\ {\rm d}x \|q^\dag-\P01(q_\theta^*)\|^{2}_{L^\infty(\Omega_\rho^c)}
\leq c\rho.
\end{equation*}
By combining the last two estimates and then optimizing with respect to the scalar $\rho$, we get
the desired bound on $ \|q^\dag-\P01(q_\theta^*)\|_{L^2(\Omega)}$.
\end{proof}

\subsection{Error analysis of the empirical loss}

Now we analyze the impact of quadrature errors on the reconstruction, under the following condition.
\begin{assumption}\label{ass:Diri-gen}
 $f\in L^\infty(\Omega)$ and $\nabla z^\delta\in L^{\infty}(\Omega)\cap L^{\infty}(\partial\Omega)$.
\end{assumption}

We have the following error bound on the DNN realization $\widehat{q}_\theta^*$ of a minimizer
$\widehat\theta^*$ of the loss \eqref{equ: Diri empir loss}.
\begin{theorem}\label{thm:err-Diri-emp}
Let Assumptions \ref{ass:Diri} and \ref{ass:Diri-gen} hold. Fix small $\epsilon_q$, $\epsilon_\sigma>0$, and let the tuple $(\widehat{\theta}^*,\widehat{\kappa}^*)\in (\mathfrak{P}_{p,\epsilon_q}, \mathfrak{P}_{\infty,\epsilon_\sigma}^{\otimes d})$ be a minimizer of the loss in \eqref{equ: Diri empir loss} and $\widehat{q}_\theta^*$ the
DNN realization of $\widehat{\theta}^*$. Fix $\tau\in(0,\frac18)$, let the quantities $e_d$, $e_\sigma$, $e_{b'}$ and
$e_q$ be defined in Theorem \ref{thm:stat-err}, and let
$\eta:=e_d+\gamma_\sigma e_\sigma+\gamma_be_{b'}+\gamma_qe_q
+\epsilon_q^2+ (1+\gamma_\sigma+\gamma_b)\epsilon_\sigma^2+(1+\gamma_b) \delta^2+\gamma_q.$
Then with probability at least $1-4\tau$,
there holds
\begin{align*}
&\int_{\Omega}\Big(\frac{q^\dag-\P01(\widehat{q}_\theta^*)}{q^\dag}\Big)^2\big(q^\dag |\nabla u^\dag|^2+fu^\dag\big)\ \mathrm{d}x\\
\leq& c\big((e_d+\eta)^\frac12+(e_\sigma+\gamma_\sigma^{-1}\eta)^\frac12+\delta(1+e_q+\gamma_q^{-1}\eta)^\frac12\big)\big(1+ e_q+\gamma_q^{-1}\eta\big)^\frac12,
\end{align*}
where the constant $c>0$ depends on  $|\Omega|$, $|\partial\Omega|$, $d$,  $\|z^\delta\|_{W^{1,\infty} (\partial\Omega)}$, $\|f\|_{L^\infty(\Omega)}$ and $\|q^\dagger\|_{L^\infty(\partial\Omega)}$ at most polynomially.
Moreover, if Condition \ref{Cond:P} holds, then with probability at least $1-4\tau$,
$\|q^\dagger-\P01(\widehat{q}_\theta^*)\|_{L^2(\Omega)}\leq c\big(\big((e_d+\eta)^\frac12+(e_\sigma+\gamma_\sigma^{-1}\eta)^\frac12+\delta(1+e_q+\gamma_q^{-1}\eta)^\frac12\big)\big(1+ e_q+\gamma_q^{-1}\eta\big)^\frac12\big)^{\frac{1}{2(\beta+1)}}$.
\end{theorem}
\begin{proof}
The proof is similar to Theorem \ref{thm:err-Neum-emp}, using the estimate on $e_{b'}$
instead of the estimate $e_b$. Indeed, the following splitting of the statistical error holds
\begin{equation*}
\sup_{(\theta,\kappa)\in(\mathfrak{P}_{p,\epsilon_q},\mathfrak{P}_{\infty,\epsilon_\sigma}^{\otimes d})} \big|J_{\bsgamma}(q_\theta,\sigma_\kappa)-\widehat{J}_{\bsgamma}(q_\theta,\sigma_\kappa)\big|\leq
\Delta\mathcal{E}_d+\gamma_\sigma\Delta \mathcal{E}_\sigma+\gamma_{b}\Delta \mathcal{E}_{b'}+\gamma_q\Delta \mathcal{E}_q.
\end{equation*}
Then for any minimizer $(\widehat{\theta}^*,\widehat{\kappa}^*)$ of
the empirical loss $\widehat J_{\bsgamma}(\theta,\kappa)$ \eqref{equ: Diri empir loss}, with probability at least $1-4\tau$,
\begin{equation*}
\widehat{J}_{\bsgamma}(\widehat{\theta}^*,\widehat{\kappa}^*)\leq c\big(e_d+\gamma_\sigma e_\sigma+\gamma_be_{b'}+\gamma_qe_q+\epsilon_q^2+(1+\gamma_\sigma+\gamma_b)\epsilon_\sigma^2+(1+\gamma_b)\delta^2+\gamma_q\big).
\end{equation*}
Then by repeating the argument from Theorem \ref{thm:err-Neum-emp} and
replacing $q_\theta^*$ by $\widehat{q}_\theta^*$, we deduce that
with probability at least $1-4\tau$, the following three estimates hold simultaneously
$\|\widehat{\sigma}_\kappa^*-\P01(\widehat{q}_\theta^*)\nabla z^\delta\|_{L^2(\Omega)}^2\leq c(e_d+\eta)$, $\|f+\nabla\cdot\hat{\sigma}_\kappa^*\|_{L^2(\Omega)}^2\leq c (e_\sigma+\gamma_\sigma^{-1}\eta)$, and
$\|\nabla\widehat{q}_\theta^*\|^2_{L^2(\Omega)}\leq c(e_q+\gamma_q^{-1}\eta)$.
Last, the desired estimates follow from the argument of Theorem
\ref{thm: Diri error estimate on conduc}.
\end{proof}
\begin{remark}
Under the assumptions in Theorem \ref{thm:err-Diri-emp} and the choice of the numbers of sampling points in Remarks \ref{rmk:para-Neum} and \ref{rmk:para-Neum-emp}, there holds with probability at least $1-4\tau$, $
    \|q^\dag- \P01(\widehat{q}_\theta^*)\|_{L^2(\Omega)}\leq c\delta^\frac{1}{2(1+\beta)}$.
This estimate is largely comparable with that by the FEM \cite[Corollary 3.3]{jin2021error}.
\end{remark}

\section{Numerical experiments and discussions}\label{sec:numer}
Now we showcase the performance of the proposed mixed DNN approach. All computations are carried out on TensorFlow
1.15.0 using Intel Core i7-11700K Processor with 16 CPUs. We measure the accuracy of a reconstruction $\hat q$ (with respect to the exact one $q^\dag$) by the relative $L^2(\Omega)$ error $e(\hat q)$ defined by
$e(\hat q)=\|q^\dag-\hat q\|_{L^2(\Omega)}/\|q^\dag\|_{L^2(\Omega)}.$
For an elliptic problem in $\mathbb{R}^d$, we use DNNs with an output dimension $1$ and $d$ to approximate the conductivity
$q$ and flux $\sigma$, respectively. Unless otherwise stated, both DNNs have 4 hidden layers (i.e., depth 5) with 26, 26, 26, and 10
neurons on each layer. Increasing the depth / width enhances the expressivity of the DNN, but also increases the total number of DNN parameters to be learned. Empirically this choice of depth and width balances well the expressivity and efficiency of the DNN, as was similarly adopted  \cite{XuDarve:2022,BarSochen:2021}. Table \ref{table3} below shows the impact of different choices of depth and width on the error $e(\hat{q})$ of the reconstruction $\hat q$. The numbers of sampling points in the domain $\Omega$ and
on the boundary $\partial\Omega$ are denoted by $n_r$ and $n_b$, respectively. The penalty parameters $\gamma_\sigma$, $\gamma_b$ and $\gamma_q$ are for the divergence term, boundary term and $H^1(\Omega)$-penalty
term, respectively, in the losses \eqref{eqn:obj-Neum} and \eqref{eqn:obj-Diri1}. The empirical loss $\widehat{J}_{\bsgamma}$
is  minimized by the popular ADAM \cite{KingmaBa:2015}. We adopt a stagewise decaying learning rate schedule, which is determined by the starting learning rate (lr), decay rate (dr) and the epoch number at which the decay takes place (step). The term epoch refers to the total number of epochs used for the reconstruction. The hyper-parameters of the DNN approach include penalty parameters and optimizer parameters, and are determined in a trial-and-error manner. That is, the algorithm is run with different hyper-parameter settings, and then the set of hyper-parameters attaining the smallest reconstruction error $e(\hat{q})$ is reported. This strategy is often used for tuning hyper-parameters in deep learning \cite{YuZhu:2020}. A systematic study of the hyper-parameter selection is beyond the scope of this work. Tables \ref{tab:algpara} and \ref{tab:algpara2}
summarize the algorithmic parameters for the experiments where the numbers in the brackets indicate the parameters used for noisy data ($\delta=10\%$). For most examples, $q^\dag$ and $u^\dag$ are known explicitly and the values of $\nabla u^\dag$, $f$ and $g$ at the sampling points are computed directly. Then Gaussian noise is added pointwise to $\nabla u^\dag$ at the sampling points drawn uniformly in the domain (and also on the boundary in the Dirichlet case) to form the observation $
    \nabla z^\delta(x)=\nabla u^\dag(x)+\delta\cdot\iota(x)\cdot\max_{x\in\Omega}\|\nabla u^\dag(x)\|_{\ell^\infty}$,
where $\iota(x)\in\mathbb{R}^d$ follows the standard
Gaussian distribution (with zero mean and the covariance being the identity matrix). For Example \ref{exam:discon}, the solution $u^\dag$ and the gradient $\nabla u^\dag$ are approximated using the standard FEM solver FreeFEM++ \cite{Hecht:2012}, and the observation $\nabla z^\delta$ is obtained by adding Gaussian noise pointwise. We present also results by the FEM for several 2D examples, which are computed over a $64\times64$ mesh, with the spaces $P_1$ and $[P_1]^2$ for approximating the conductivity $q$ and the flux $\sigma$, respectively. We employ the \textit{discretize-then-optimize} strategy \cite{KohnLowe:1988}, explicitly forming the finite-dimensional optimization problem and minimizing the resulting problem via MATLAB optimization toolbox. The python codes for reproducing the results will be made available at \url{https://github.com/jlqzprj/mixed-DNN}.

\begin{table}
\begin{center}
\caption{The algorithmic parameters used for the examples. The top block is for Neumann problems, and the bottom block for Dirichlet problems.}\label{tab:algpara}
{
\setlength{\tabcolsep}{2pt}
\begin{tabular}{c|ccccc}
\toprule
  $\backslash$No.& \ref{exam:neu1} & \ref{exam:discon}& \ref{exam:neu2}& \ref{exam:neudim5} &\ref{exam:neupartial3d} \\
\midrule
     $\gamma_\sigma$&10(100)&10(100)&10&1(10)&10 \\
     $\gamma_b$&10(50)&10(50)&10(15)&1(10)&10 \\
     $\gamma_q$&1e-5&1e-5&1e-5&1e-5&1e-5\\
     $n_r$&4e4&4e4&4e4&6e4&6.3e4\\
     $n_b$&4e3&4e3&6e3&3e4&6e3\\
     lr&2e-3&1e-3&3e-3&3e-3&3e-3\\
     dr & 0.7 & 0.75 &0.7 & 0.75&0.8 \\
     step & 2000 & 1500 & 2500 & 3000 &3000\\
     epoch   & 6e4(3e4) &5e4(1.6e4) &6e4(3e4) & 6e4(2e4)&7e4(2e4)\\
\midrule
     &\ref{exam:diri1}&\ref{exam:diridisctn}& \ref{exam:diri2}& \ref{exam:diridim5}&\ref{exam:diripartial3d}\\
     $\gamma_\sigma$ & 10&5(1)&10&10&10\\
     $\gamma_b$ &10&10&10&1(5)&10\\
     $\gamma_q$&1e-5&1e-5&1e-5&1e-5&1e-5\\
     $n_r$&4e4&4e4&4e4&6e4&6.3e4\\
     $n_b$&4e3&4e3&6e3&3e4&6e3\\
     lr&2e-3&3e-3&3e-3&3e-3&3e-3\\
     dr  & 0.7 &0.8& 0.8 & 0.75&0.8\\
     step & 2000 &2500& 3000 & 3000&3000\\
     epoch &6e4(2e4)&8e4(2e4)&6e4(3e4)& 6e4(2e4)&7e4(2e4)\\
\bottomrule
\end{tabular}}
\end{center}
\end{table}

\begin{table}
\begin{center}
\caption{The algorithmic parameters used for the examples.}\label{tab:algpara2}
{
\setlength{\tabcolsep}{.1cm}
\begin{tabular}{c|ccc|ccc}
\toprule
  \multicolumn{1}{c|}{$\backslash$No.}&\multicolumn{3}{c}{\ref{exam:neupartial2d}}  & \multicolumn{3}{c}{\ref{exam:diripartial2d}}\\
  \cmidrule(lr){2-4} \cmidrule(lr){5-7}
  & exact&$5\%$  &$10\%$ & exact& $5\%$  &$10\%$ \\
\midrule
     $\gamma_\sigma$&10&10&10&10&10&10\\
     $\gamma_b$&10&10&10&10&10&10\\
     $\gamma_q$&1e-5&1e-5&1e-5&1e-5&1e-5&1e-5\\
     $n_r$&5e4&5e4&5e4&5e4&5e4&5e4\\
     $n_b$&4e3&4e3&4e3&4e3&4e3&4e3\\
     lr&3e-3&3e-3&3e-3&3e-3&3e-3&3e-3\\
     dr&0.8&0.8&0.8&0.8&0.8&0.8\\
     step&2500&2500&2500&2500&2500&2500\\
     epoch &8e4&3e4&3e4&8e4&3e4&1e4  \\
\bottomrule
\end{tabular}}
\end{center}
\end{table}

\subsection{The Neumann problem}\label{sec:neu}
First we illustrate the approach for the inverse problem in the Neumann case.
The first example is about recovering a smooth conductivity $q^\dag$ with three modes.
\begin{example}
The domain $\Omega = (-1,1)^2$, $q^\dag= 1 +s_1(x_1, x_2) + s_2(x_1, x_2) + s_3(x_1, x_2),$
with $s_1=0.3e^{-20(x_1-0.3)^2-15(x_2-0.3)^2}$,
$s_2= -0.3e^{-10x_1^2-10(x_2+0.5)^2}$ and $s_3 = 0.2e^{-15(x_1+0.4)^2-15(x_2-0.35)^2}$, and $u^\dag=x_1+x_2+\frac{1}{3}(x_1^3+x_2^3)$.
\label{exam:neu1}
\end{example}

\begin{figure}[htbp]
\centering
\setlength{\tabcolsep}{0em}
\begin{tabular}{ccc}
\includegraphics[width=0.32\textwidth]{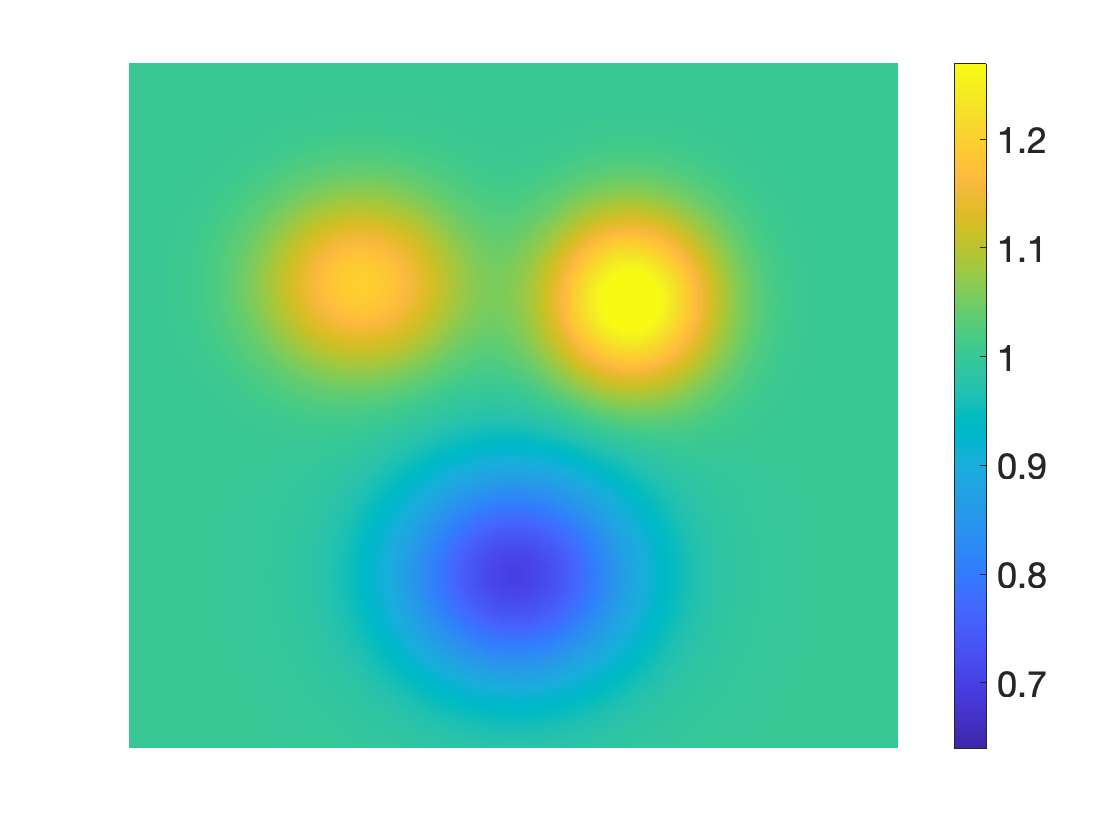} &
\includegraphics[width=0.32\textwidth]{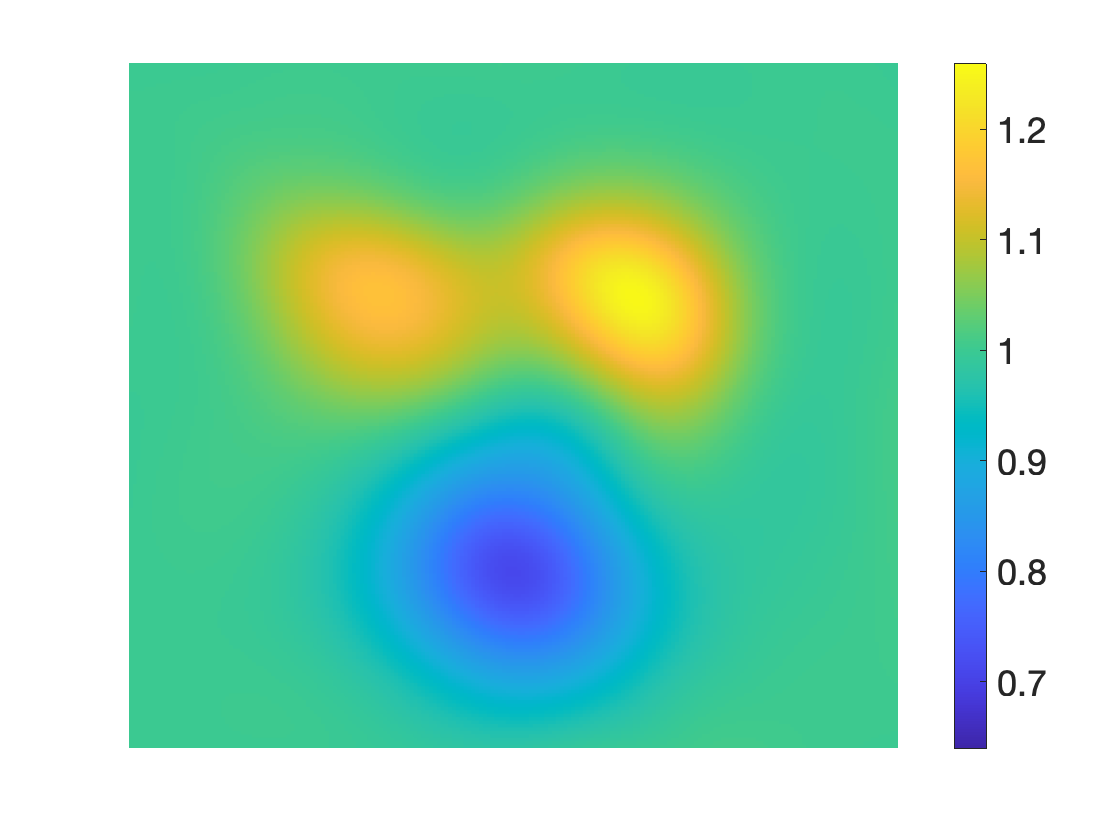} &
\includegraphics[width=0.32\textwidth]{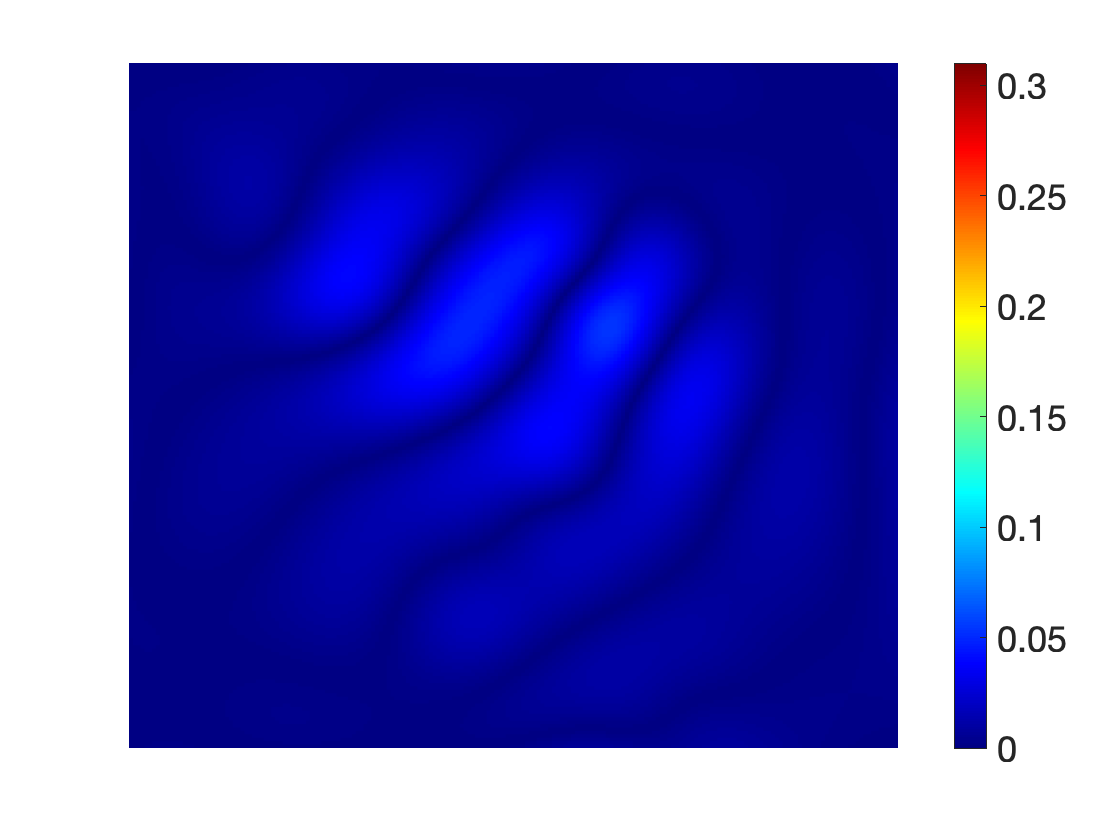}\\
\includegraphics[width=0.32\textwidth]{neu1ex.png} &
\includegraphics[width=0.32\textwidth]{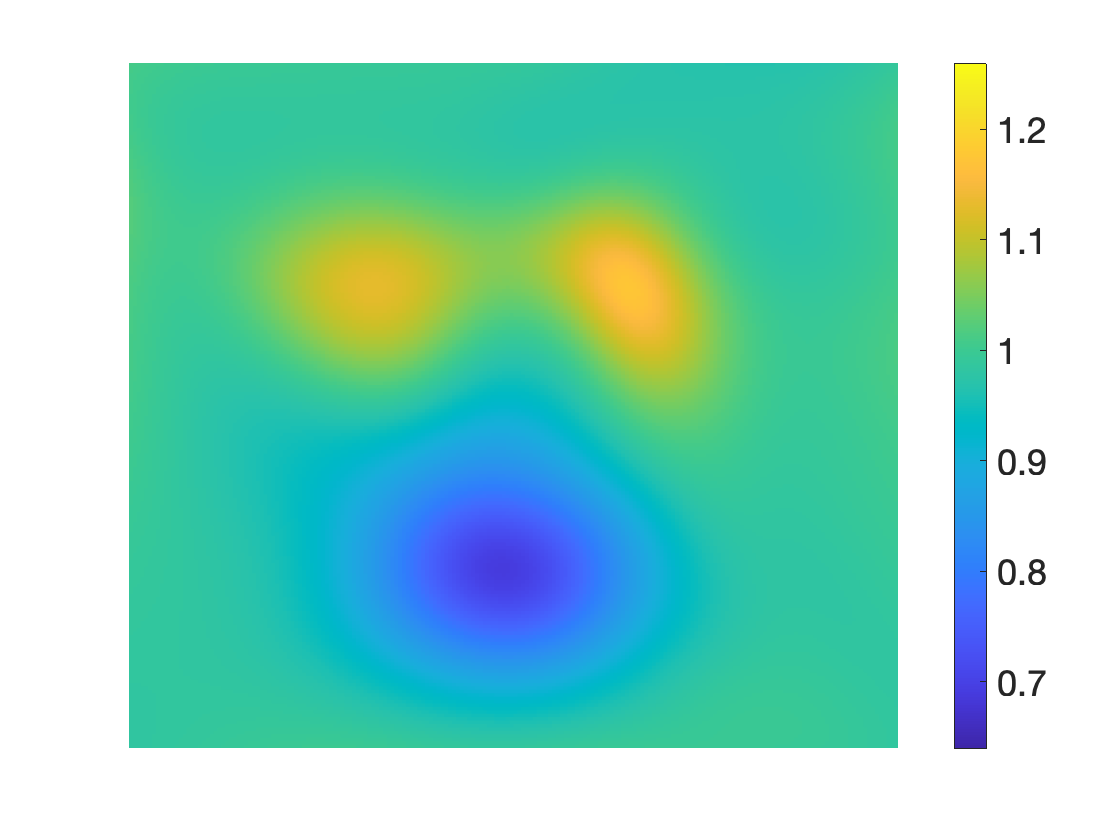} &
\includegraphics[width=0.32\textwidth]{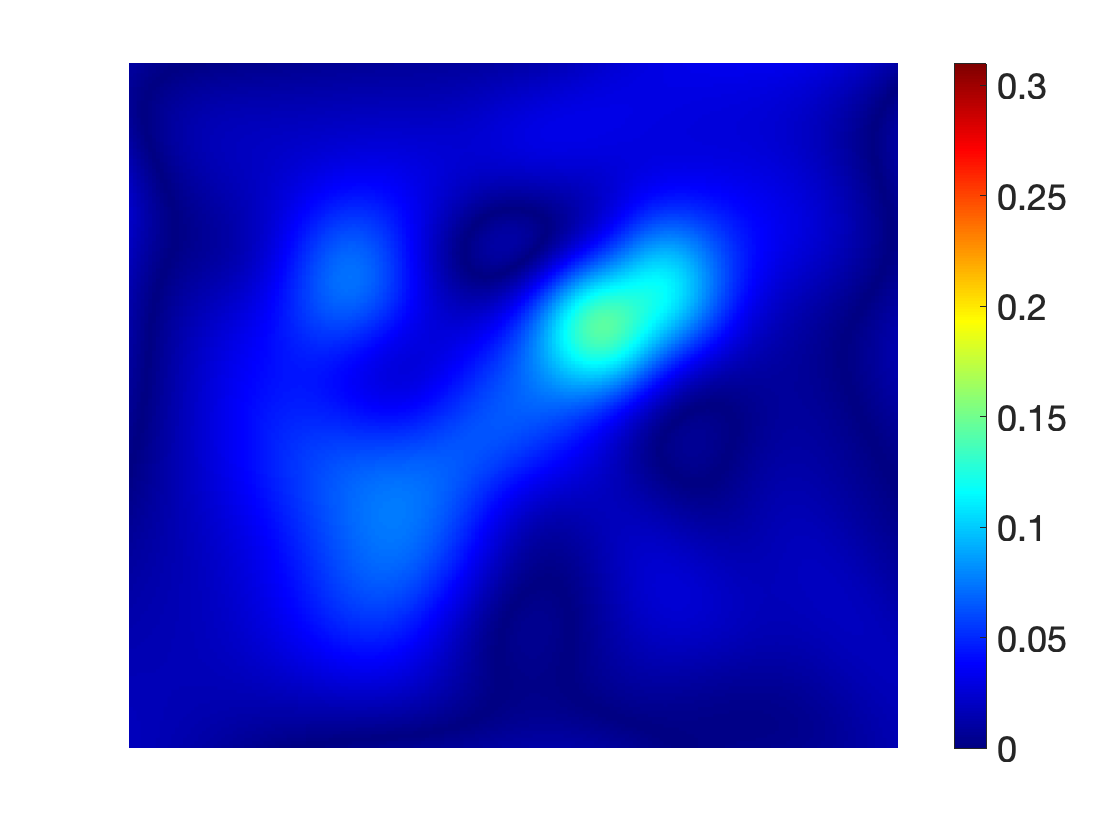}\\
\includegraphics[width=0.32\textwidth]{neu1ex.png} &
\includegraphics[width=0.32\textwidth]{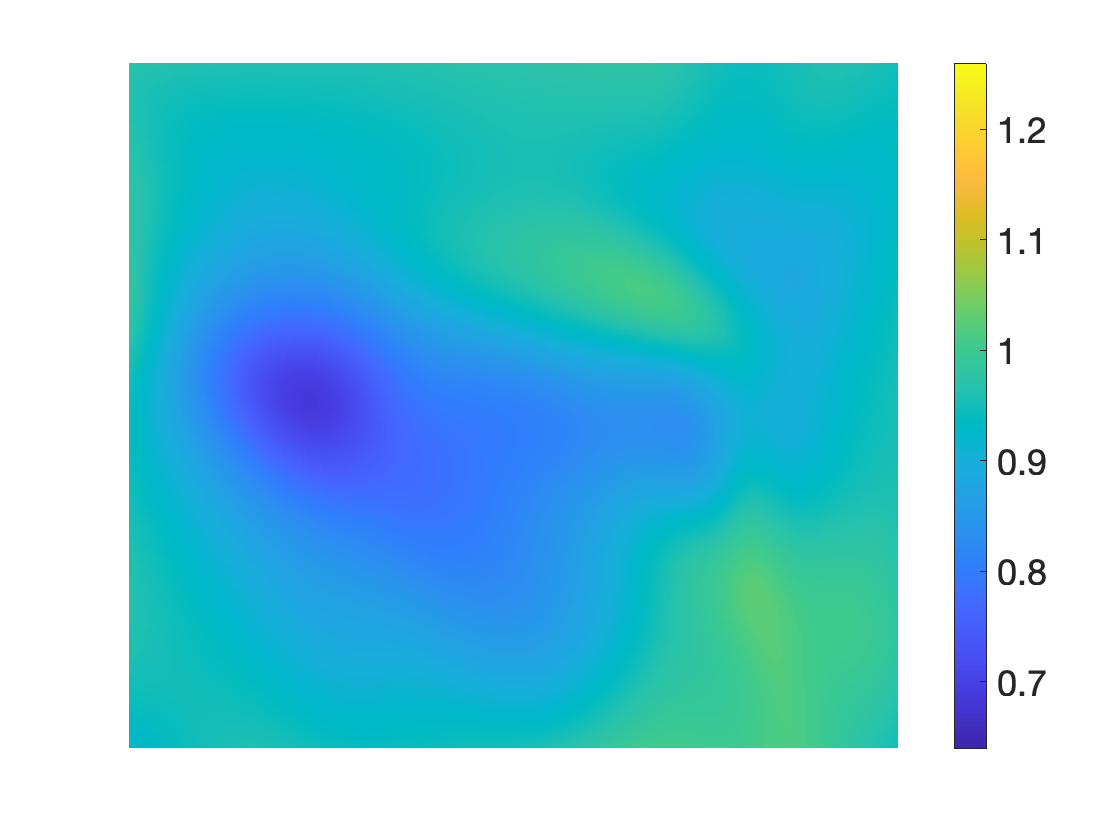} &
\includegraphics[width=0.32\textwidth]{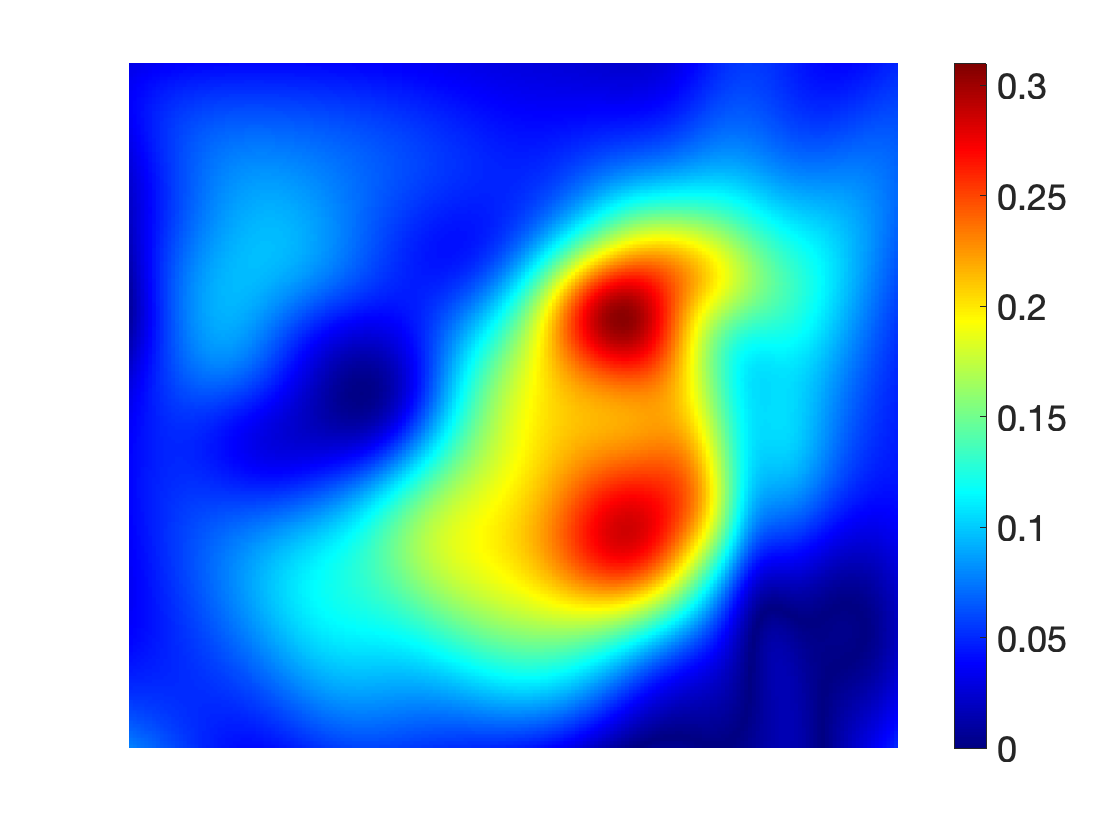}\\
(a) $q^\dag$  & (b) $\hat q$ & (c) $|\hat q-q^\dag|$
\end{tabular}
\caption{The reconstructions for Example \ref{exam:neu1} using DNNs with exact data $($top$)$ and noisy data $(\delta=10\%,\ 20\%$, middle, bottom$)$.}
\label{fig:neu1}
\end{figure}

\begin{figure}[htbp]
\centering
\setlength{\tabcolsep}{0em}
\begin{tabular}{ccc}
\includegraphics[width=0.32\textwidth]{neu1ex.png} &
\includegraphics[width=0.32\textwidth]{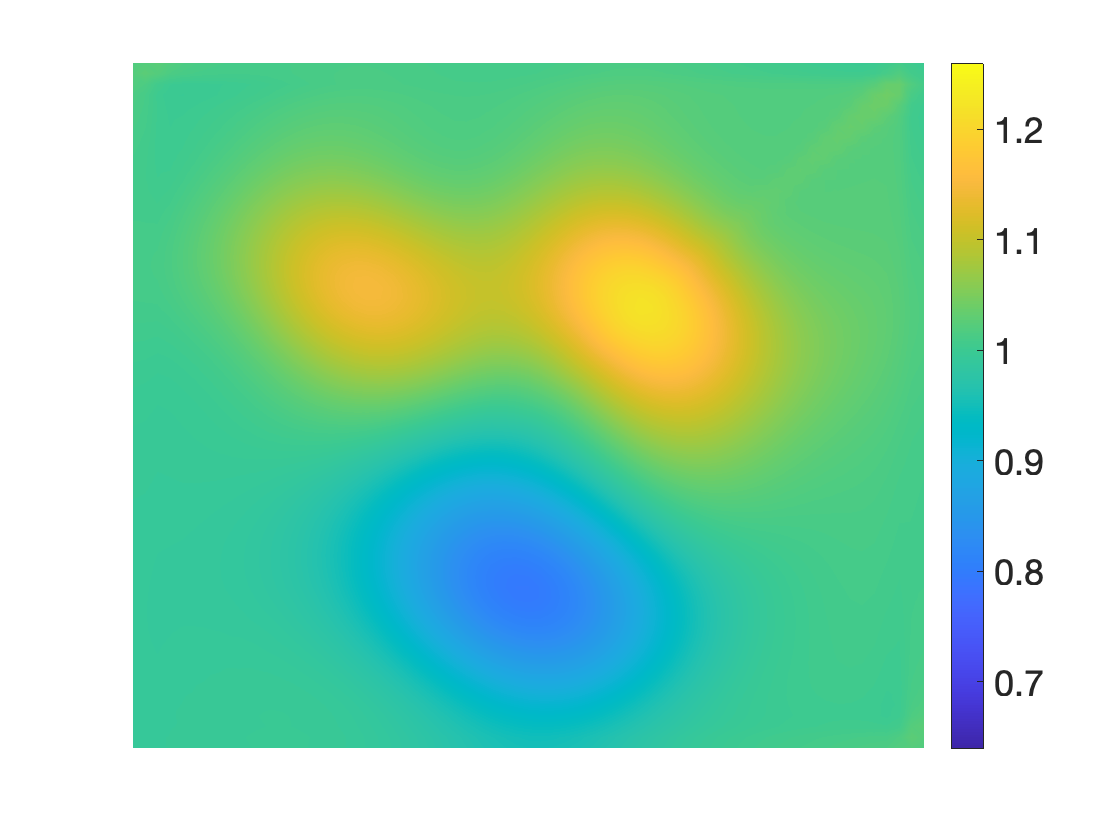} &
\includegraphics[width=0.32\textwidth]{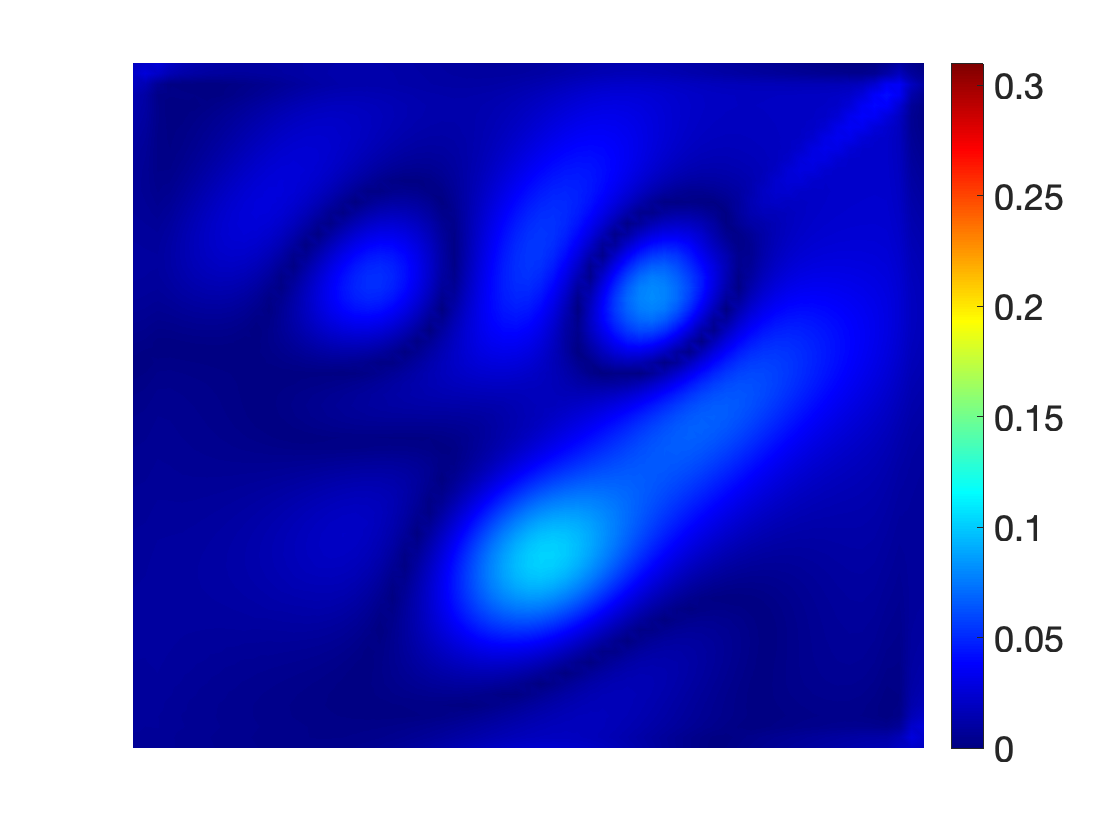}\\
\includegraphics[width=0.32\textwidth]{neu1ex.png} &
\includegraphics[width=0.32\textwidth]{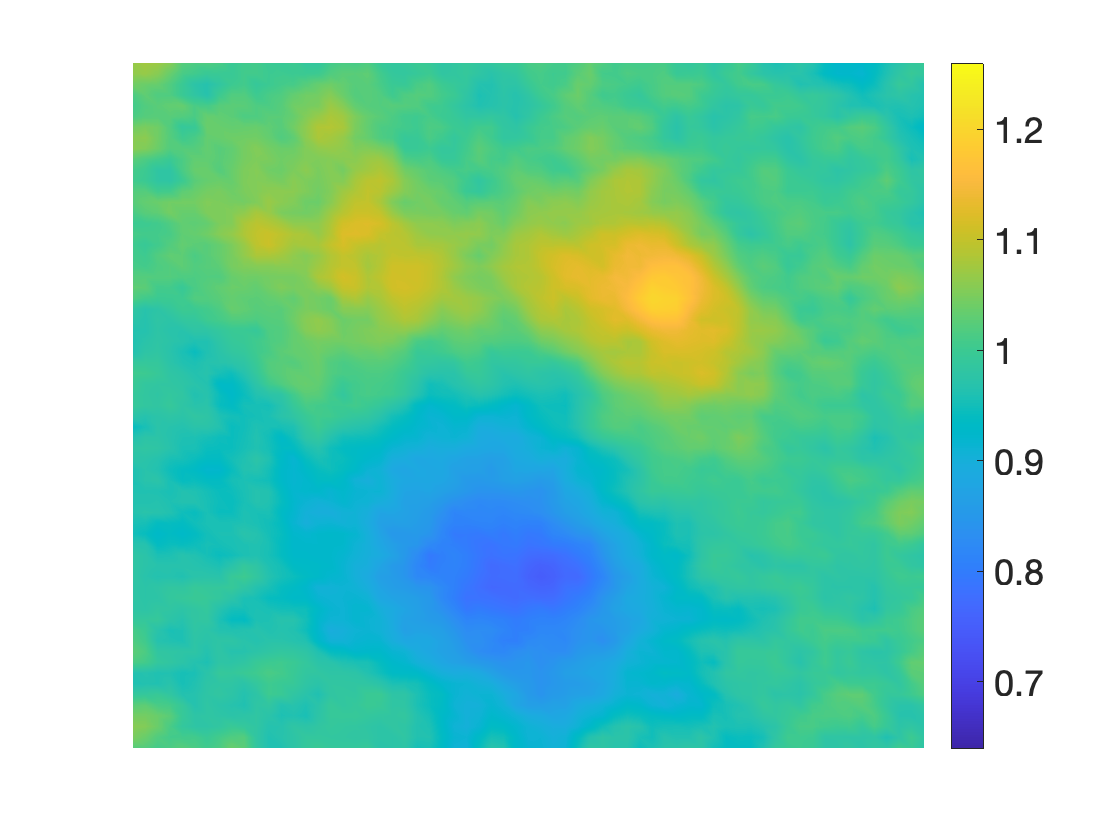} &
\includegraphics[width=0.32\textwidth]{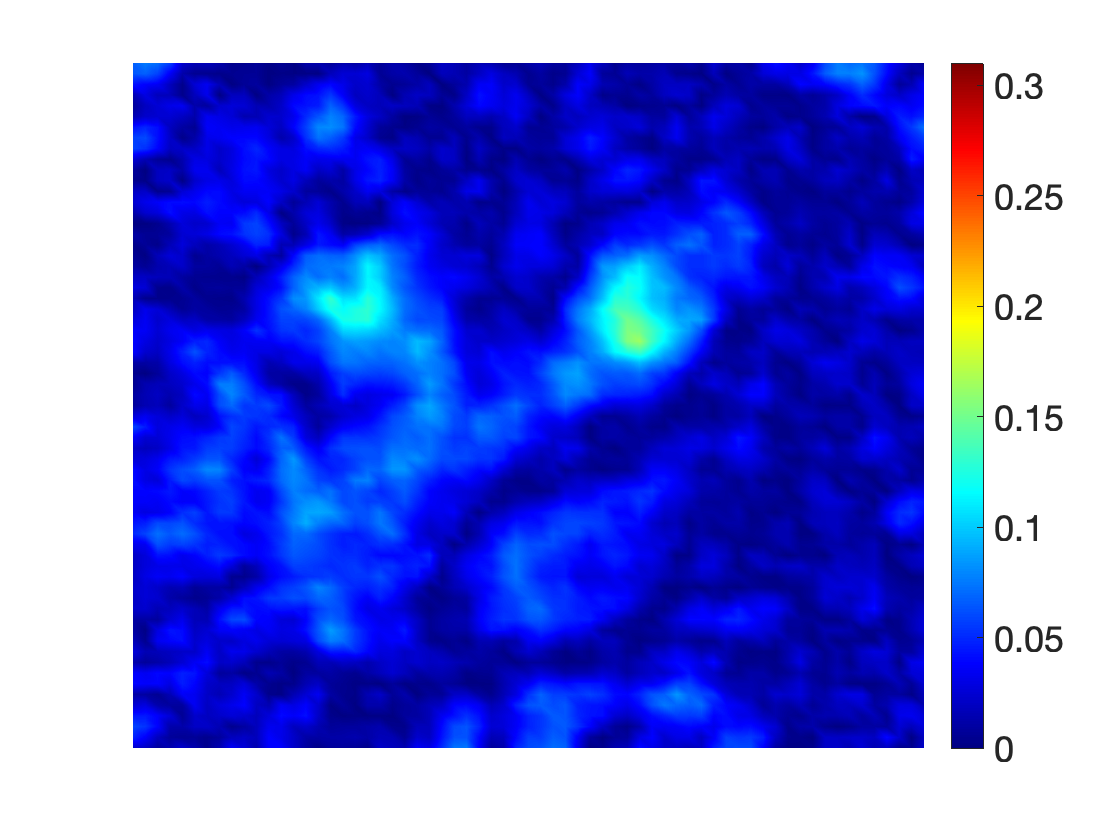}\\
\includegraphics[width=0.32\textwidth]{neu1ex.png} &
\includegraphics[width=0.32\textwidth]{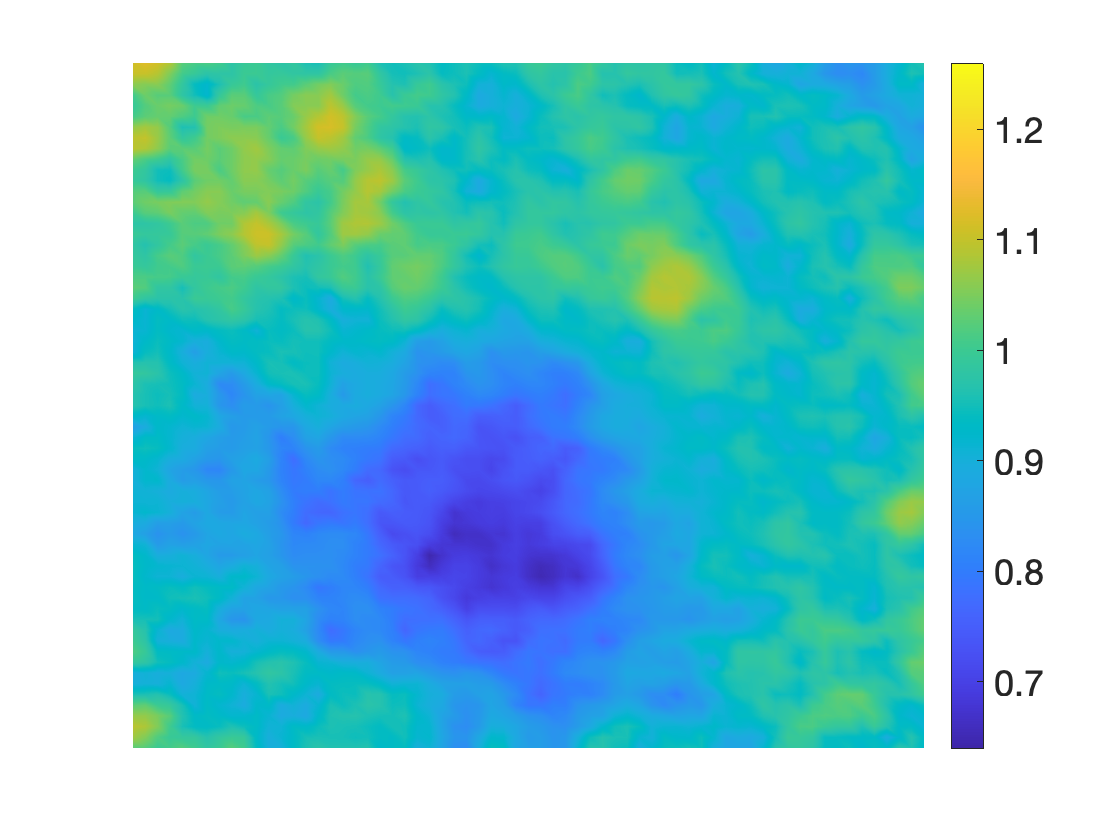} &
\includegraphics[width=0.32\textwidth]{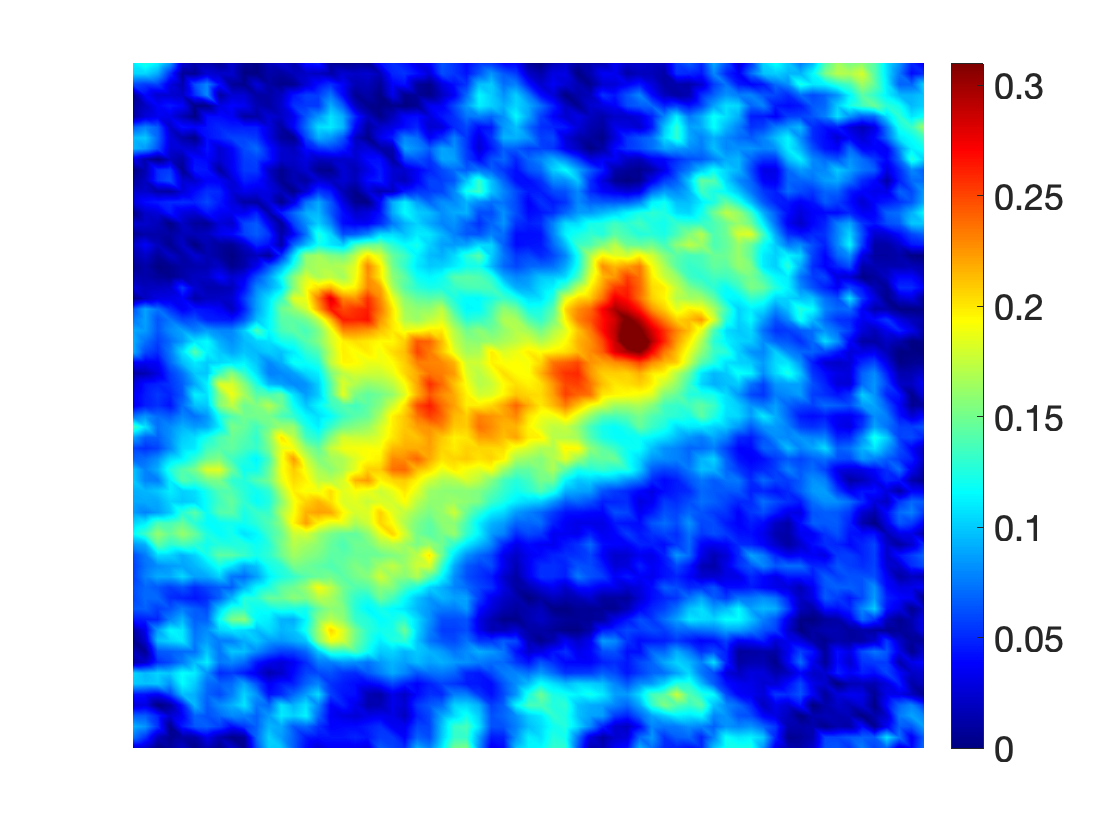}\\
(a) $q^\dag$  & (b) $\hat q$ & (c) $|\hat q-q^\dag|$
\end{tabular}
\caption{The reconstructions for Example \ref{exam:neu1} using FEM  with exact data $($top$)$ and noisy data $(\delta=10\%,\ 20\%$, middle, bottom$)$. }
\label{fig:neu1fem}
\end{figure}

\begin{figure}[htbp]
\centering
\setlength{\tabcolsep}{0em}
\begin{tabular}{ccc}
\includegraphics[width=0.32\textwidth]{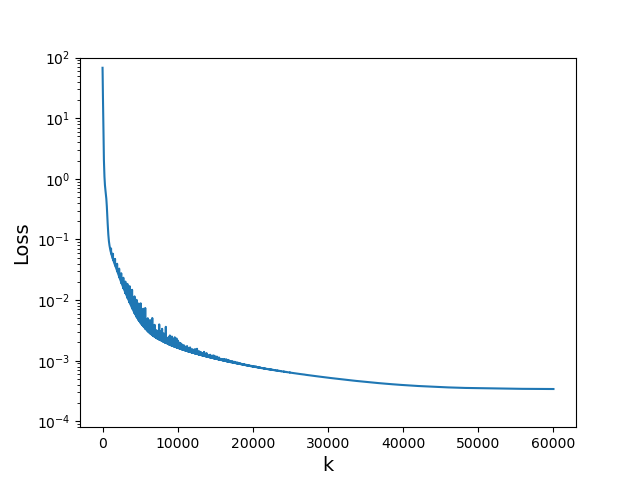} &
\includegraphics[width=0.32\textwidth]{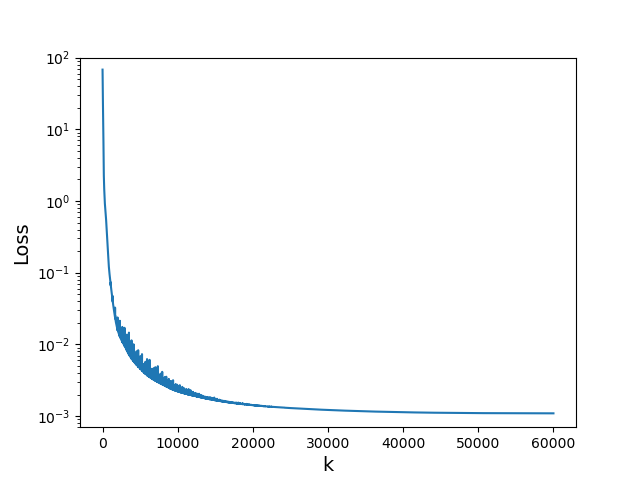} &
\includegraphics[width=0.32\textwidth]{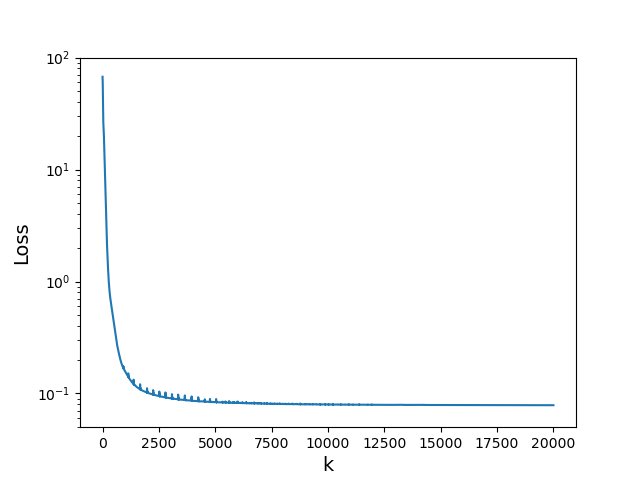}\\
\includegraphics[width=0.32\textwidth]{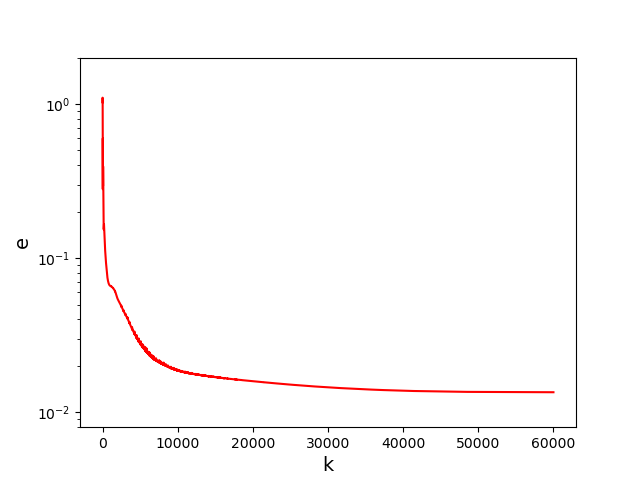} &
\includegraphics[width=0.32\textwidth]{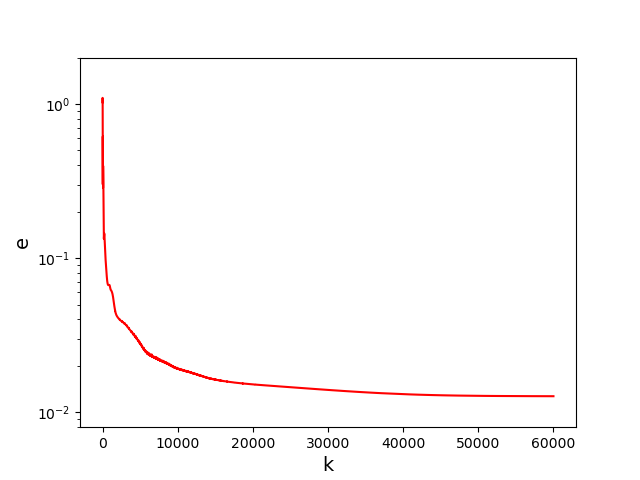} &
\includegraphics[width=0.32\textwidth]{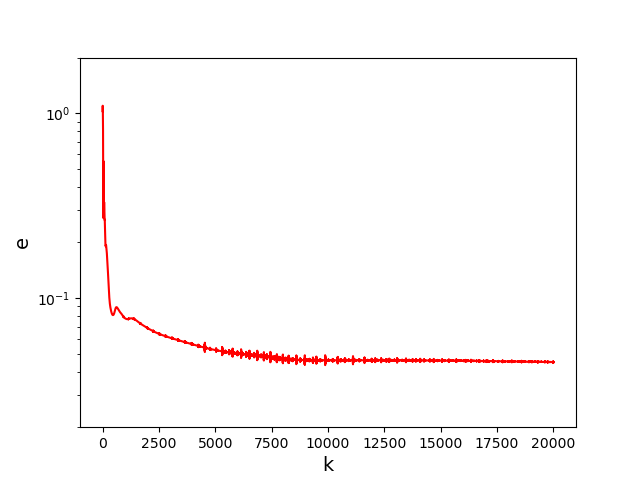}\\
(a) $\delta=0\%$  & (b) $\delta=1\%$ & (c) $\delta=10\%$
\end{tabular}
\caption{The variation of the loss $($top$)$ and the relative $L^2(\Omega)$ error $e$ $($bottom$)$ during the training process for Example \ref{exam:neu1} at different noise levels.}
\label{fig:neu1losse}
\end{figure}

The choice of $u^\dag$ ensures that both $\nabla u^\dag$ and $\Delta u^\dag$ do not vanish over the domain $\Omega$, which is beneficial for both numerical reconstruction and theoretical analysis \cite{Alessandrini:1986}. We employ both DNN and classic FEM to discretize the loss \eqref{eqn:obj-Neum}.
Fig. \ref{fig:neu1} shows the DNN
reconstructions for exact and noisy data, where
the pointwise errors $|\hat q-q^\dag|$ are also shown. The three Gaussian bumps are well resolved, except that the
middle part between the top two bumps is slightly bridged. The DNN approach is robust with respect to data noise, and the reconstruction is accurate apart from the slight under-estimation in peak values in the presence of 10\% noise. However, higher levels of noise (e.g., 20\%) pose challenges to the accurate reconstruction: the bumps can still be recognized but the magnitude is quite off. The DNN reconstructions are largely comparable with that by the classical FEM in Fig. \ref{fig:neu1fem}, but the former tend to be smoother in the reconstruction and error plots. The DNN approach seems fairly stable with respect to the iteration index, cf. Fig. \ref{fig:neu1losse}. Indeed, the presence of the noise does not affect much the convergence behavior of the loss and error, and the final errors are close to each other for all noise levels.
This observation agrees well with that for the pointwise errors in Fig. \ref{fig:neu1}.

There are several algorithmic parameters influencing the overall accuracy of the DNN approximation $\hat q$, e.g., the number of sampling points ($n_r$ and $n_b$), and DNN architectural parameters {(width, depth, and activation function)}. A practical guidance
for choosing these parameters is still missing, and we explore the issue empirically.
In Table \ref{table1}, we study the impact of the regularization parameter $\gamma_q$ on the reconstruction quality. The DNN method is observed to be very robust with respect to the presence
of data noise - the results remain fairly accurate even for up to 10\%  noise, and do
not vary much with $\gamma_q$. That is, the parameter $\gamma_q$ plays a less crucial role as in theoretical analysis. This is attributed to the inductive bias of the DNN representation which prefers smooth functions and thus the formulation has built-in regularizing effect, which is absent from the analysis. Note also that there is an inherent accuracy limitation of the DNN approach, i.e., the reconstruction cannot be made
arbitrarily accurate for exact data $\nabla u^\dag$. This may be attributed to the optimization
error: due to the nonconvexity of the loss landscape, the optimizer may fail to find a global minimizer of the empirical loss $\widehat{J}_{\bsgamma}$ but instead only an approximate local minimizer. This phenomenon has been observed across a broad range of neural solvers based on DNNs
\cite{RaissiPerdikarisKarniadakis:2019,EYu:2018,jin2022imaging}. Tables \ref{table3}-\ref{table4}
show that the $L^2(\Omega)$ error $e(\hat q)$ of the reconstruction $\hat q$ does not vary much
with different DNN architectures and numbers of sampling points. This
agrees with the convergence behavior of ADAM in Fig.
\ref{fig:neu1losse}: both loss $\hat J_{\bsgamma}$ and error $e(\hat q)$ stagnate
at a certain level.

\begin{table}[htp!]
  \centering
  \caption{The variation of the relative $L^2(\Omega)$ error $e(\hat q)$ with respect to various algorithmic parameters. }
\begin{threeparttable}
\subfigure[$e$ v.s. $\gamma_q$ and $\delta$ using DNN method\label{table1}]{\begin{tabular}{c|ccc}
\toprule
  $\gamma_q\backslash\delta$&   0\% &  1\%& 10\%\\
\midrule
     1.00e-2 &  2.03e-2 & 1.95e-2 & 4.36e-2 \\
     1.00e-3 &  1.10e-2 & 8.86e-3 & 4.30e-2\\
     1.00e-4 &  1.01e-2 & 8.52e-3 & 4.30e-2\\
     1.00e-5 &  1.35e-2 & 1.27e-2 & 4.30e-2\\
\bottomrule
\end{tabular}} \\
\subfigure[$e$ v.s. $W$ and  $L$\label{table3}]{
\begin{tabular}{c|ccc}
\toprule
${W}\backslash L$&   5 &  10& 20\\
\midrule
     4&8.58e-2&1.11e-1&2.77e-1 \\
     12&2.27e-2&3.32e-2&9.91e-3\\
     26&1.35e-2&3.88e-3&1.38e-2\\
     40&1.15e-2&1.26e-2&5.17e-3\\
\bottomrule
\end{tabular}}\quad
\subfigure[$e$ v.s. $n_r$ and $n_b$\label{table4}]{
\begin{tabular}{c|cccc}
\toprule
$n_b\backslash n_r$&   5000 &  10000& 20000&40000\\
\midrule
     500&  1.25e-2&1.52e-2&1.23e-2&9.14e-3 \\
     1000& 9.87e-3&2.79e-2&2.91e-2&1.34e-2\\
     4000& 9.71e-3&1.55e-2&1.39e-2&1.35e-2\\

\bottomrule
\end{tabular}}
\end{threeparttable}
\end{table}

\begin{figure}[htbp]
\centering
\setlength{\tabcolsep}{0em}
\begin{tabular}{ccc}
\includegraphics[width=0.32\textwidth]{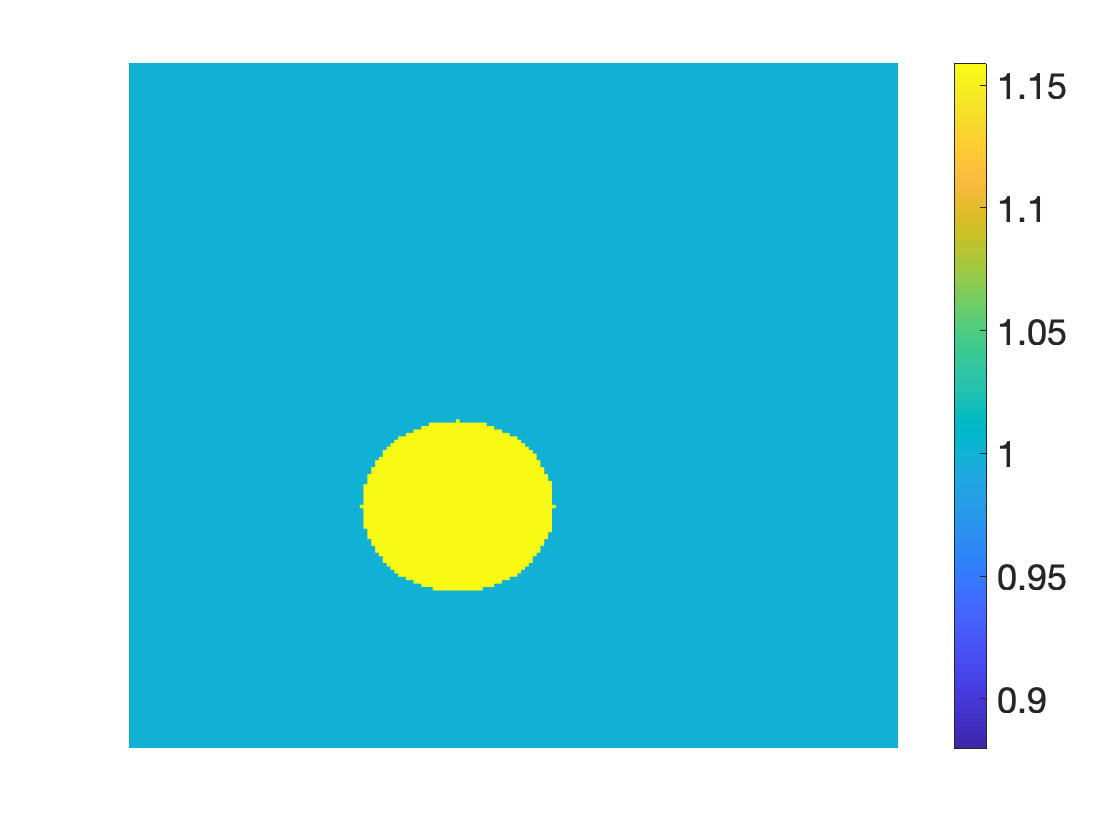} &
\includegraphics[width=0.32\textwidth]{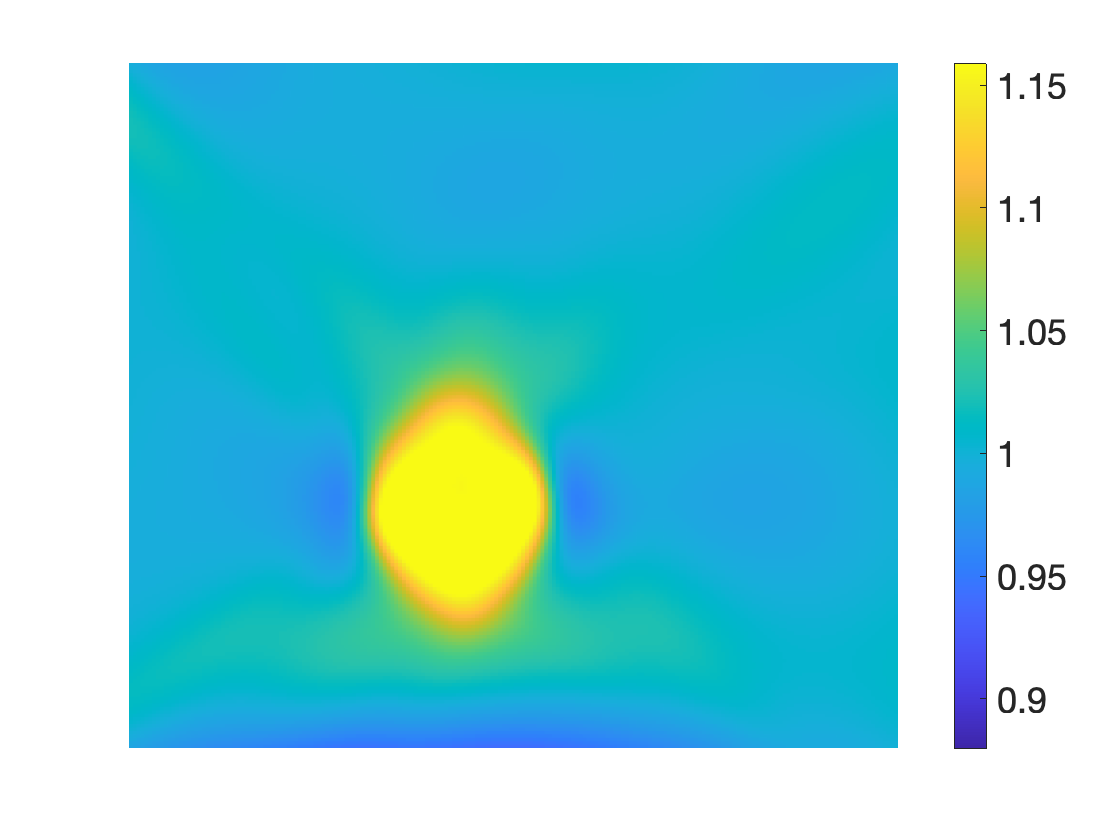} &
\includegraphics[width=0.32\textwidth]{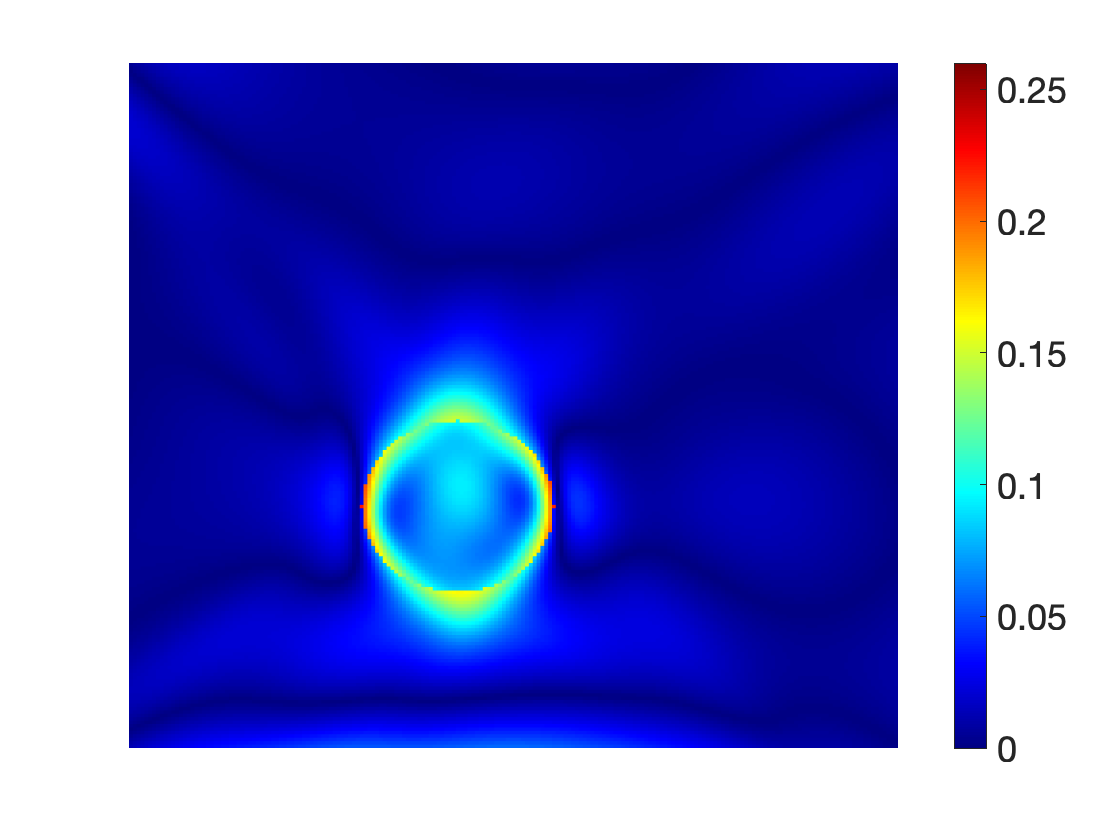}\\
\includegraphics[width=0.32\textwidth]{neudisctnex.png} &
\includegraphics[width=0.32\textwidth]{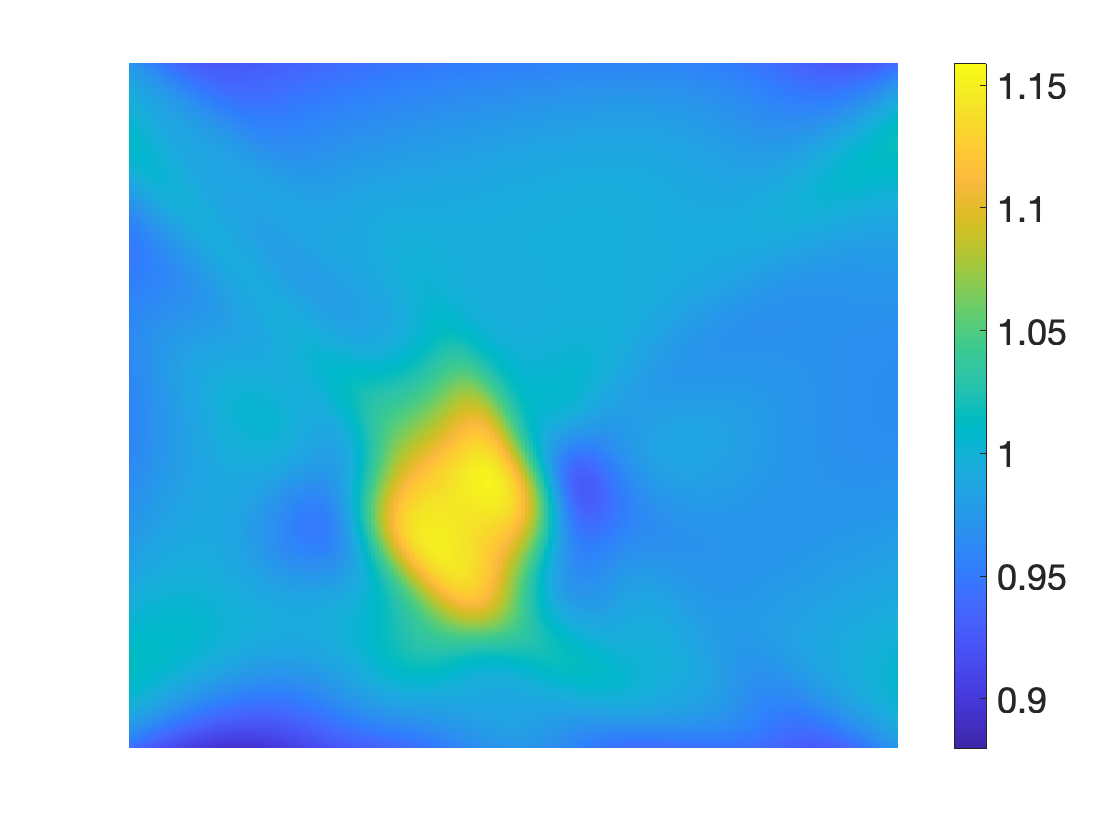} &
\includegraphics[width=0.32\textwidth]{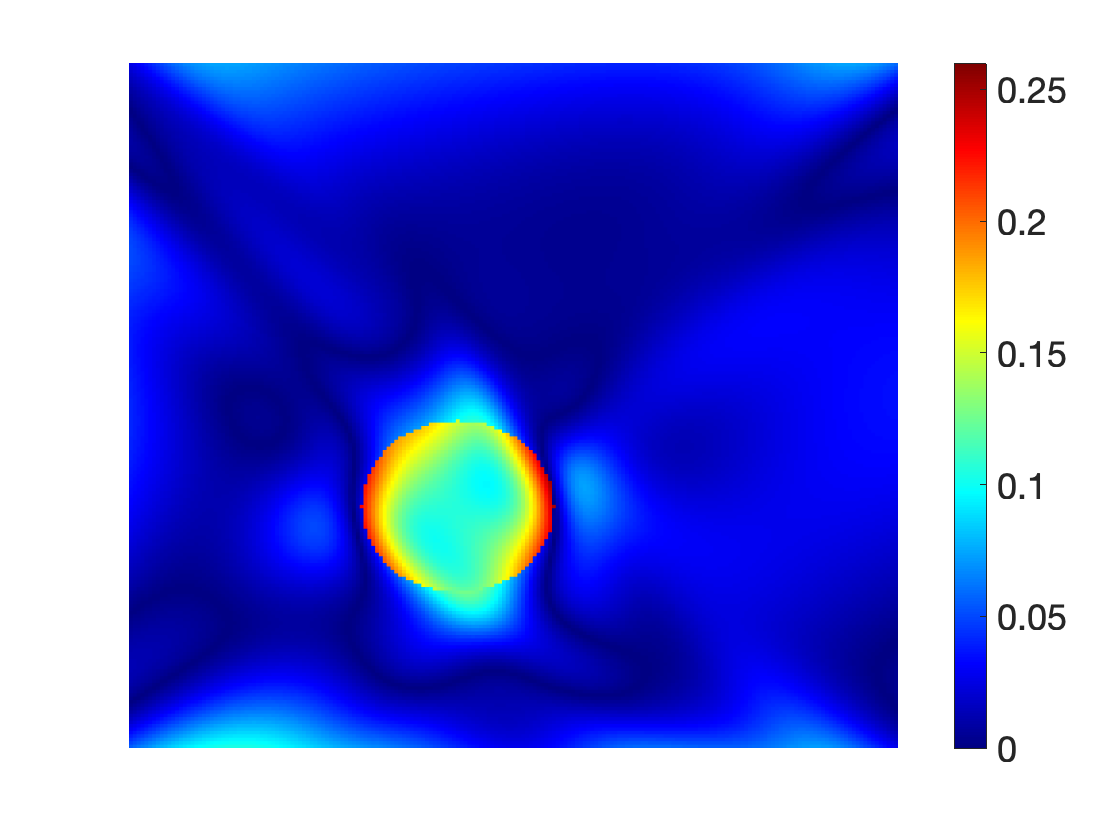}\\
(a) $q^\dag$  & (b) $\hat q$ & (c) $|\hat q-q^\dag|$
\end{tabular}
\caption{The reconstructions for Example \ref{exam:discon} with exact data $($top$)$ and noisy data $($$\delta=10\%$, bottom$)$.}
\label{fig:neudisctn}
\end{figure}

\begin{figure}[htbp]
\centering
\setlength{\tabcolsep}{0em}
\begin{tabular}{ccc}
\includegraphics[width=0.32\textwidth]{neudisctnex.png} &
\includegraphics[width=0.32\textwidth]{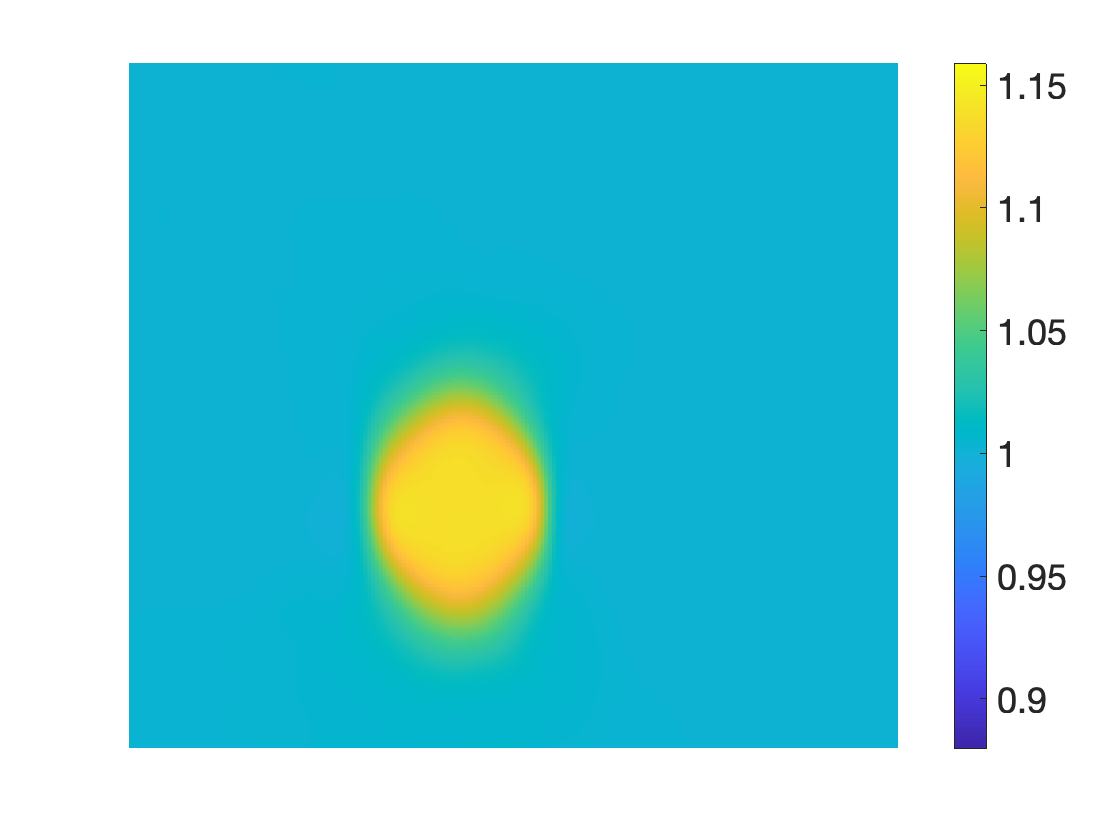} &
\includegraphics[width=0.32\textwidth]{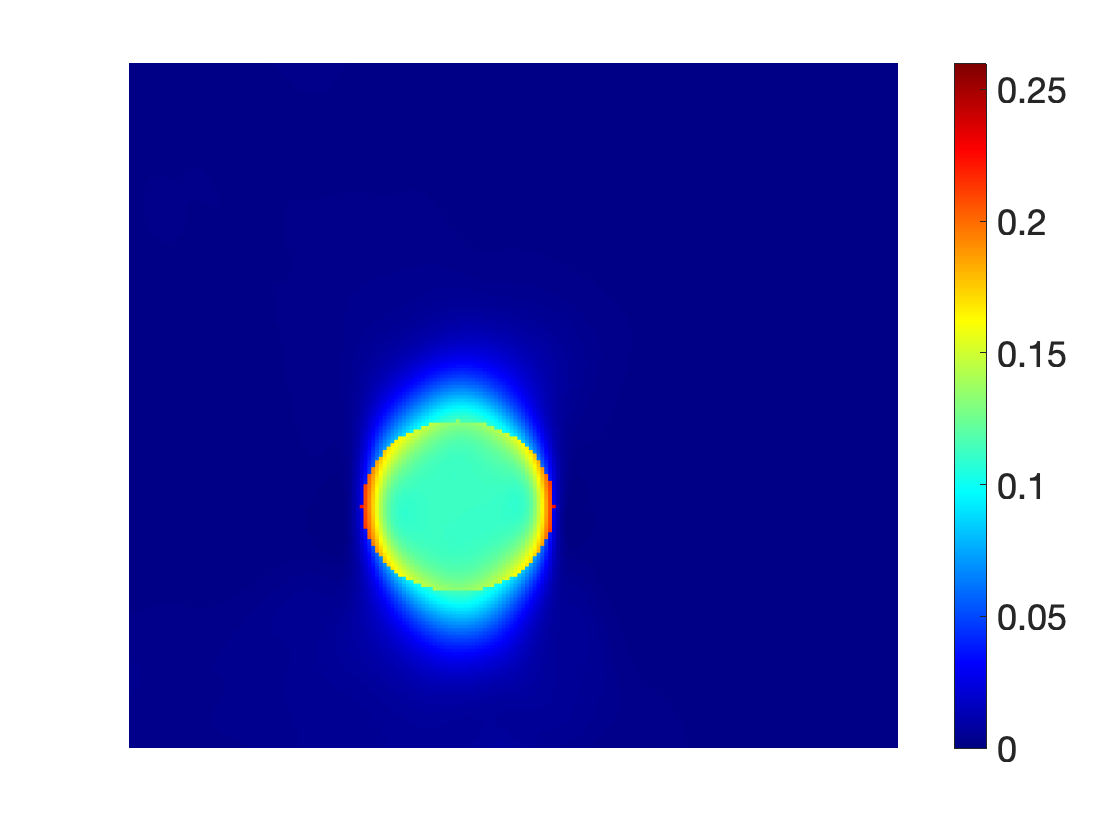}\\
\includegraphics[width=0.32\textwidth]{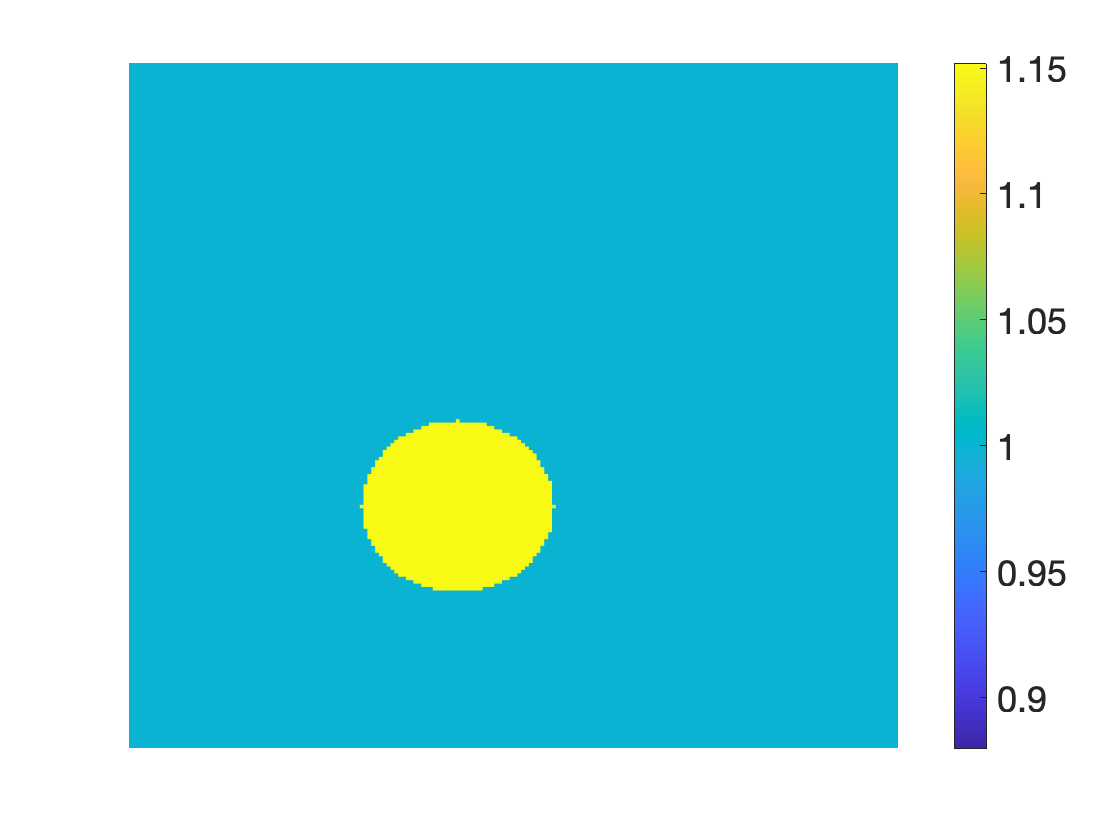} &
\includegraphics[width=0.32\textwidth]{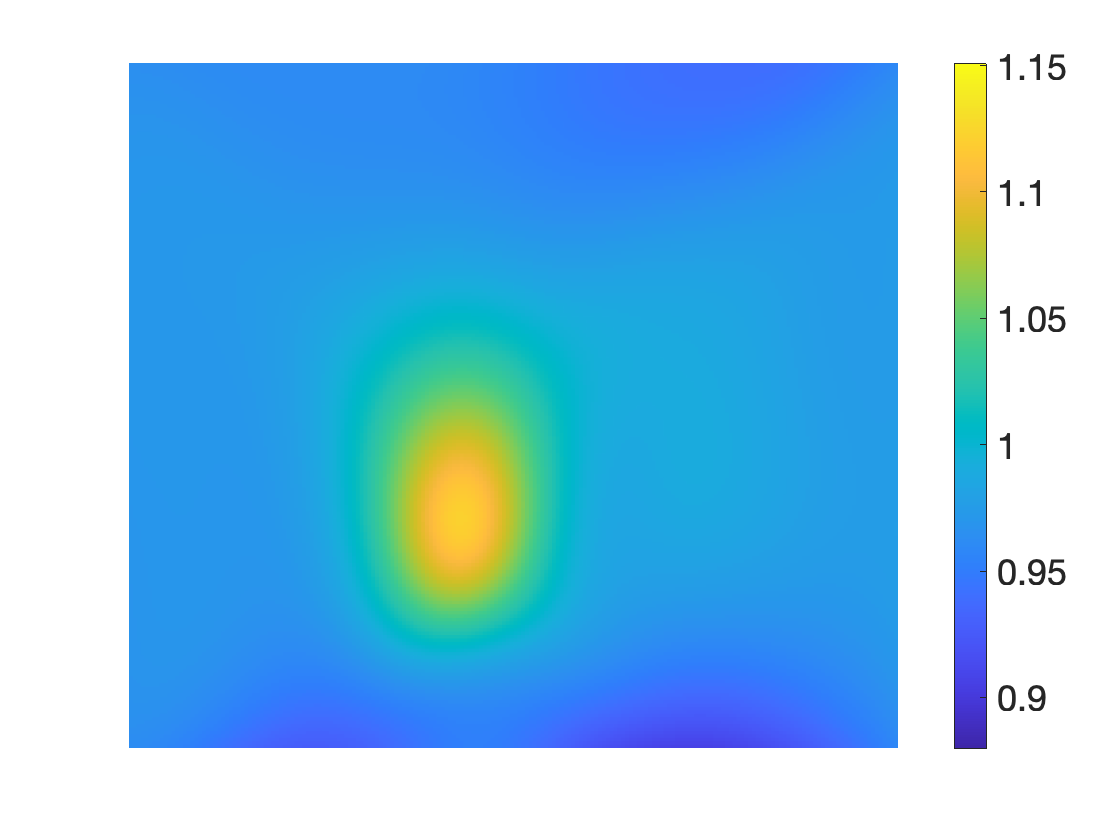} &
\includegraphics[width=0.32\textwidth]{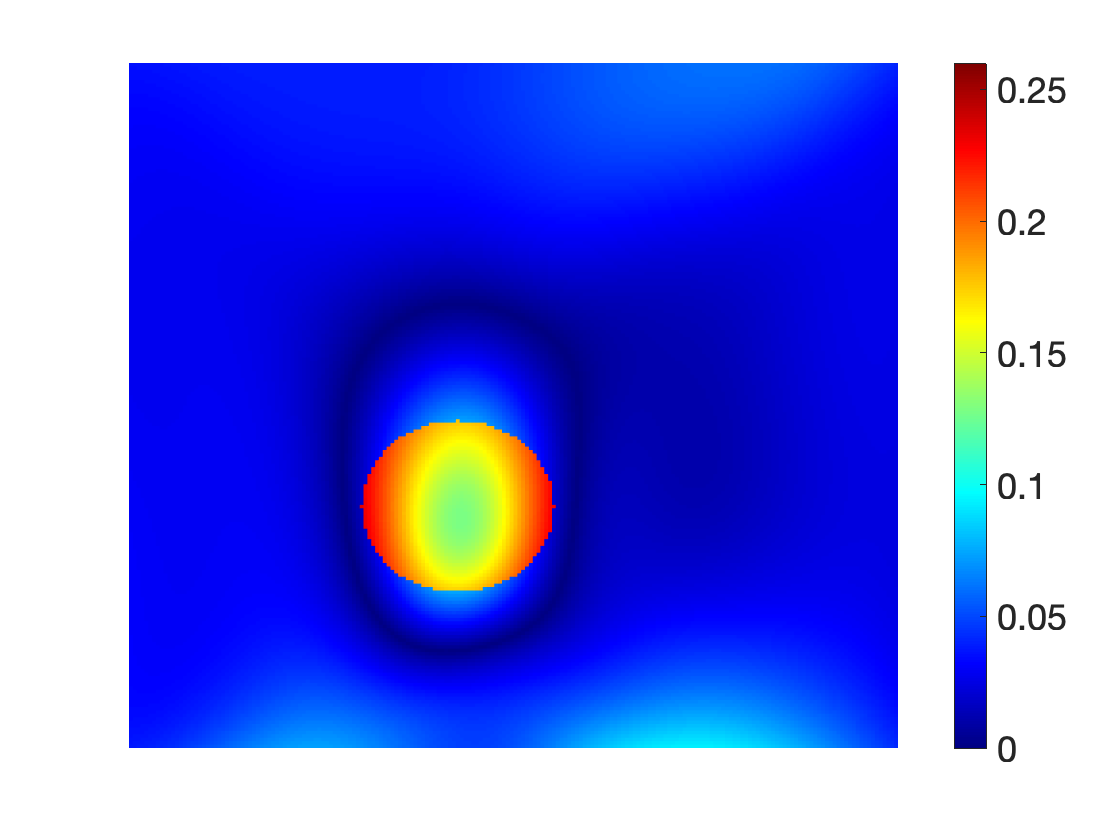}\\
(a) $q^\dag$  & (b) $\hat q$ & (c) $|\hat q-q^\dag|$
\end{tabular}
\caption{The reconstructions for Example \ref{exam:discon} using a loss function including the total variation term, with exact data $($top$)$ and noisy data $($$\delta=10\%$, bottom$)$.}
\label{fig:neudisctntv}
\end{figure}

The second example is about recovering a discontinuous conductivity $q^\dag$.
The notation $\chi_S$ denotes the characteristic function of a set $S$.
\begin{example} \label{exam:discon}
$\Omega = (-1,1)^2$,
$q^\dag = 1+0.25\cdot\chi_{\{(x_1+0.15)^2+(x_2+0.3)^2\leq0.25^2\}}$, $f\equiv0$ and $g=x_1$.
\end{example}

Since $q^\dag$
is piecewise constant, we may also include the  total variation penalty \cite{RudinOsherFatemi:1992},
i.e., $\gamma_{tv}| q|_{\rm TV}$, to the loss $J_{\bsgamma}(\theta,\kappa)$ to
promote piecewise constancy of the reconstruction. Fig. \ref{fig:neudisctn}
shows the error plots without the
total variation term, and Fig.  \ref{fig:neudisctntv} including the total variation
term ($\gamma_{tv}=0.01$ for both exact and noisy data). The reconstruction quality is improved by including the total variation term in the loss so that the reconstructed conductivity better restores the piecewise constancy. However, the reconstruction $\hat q$ still exhibits slight blurring around the discontinuous interface. Moreover, in the presence of $10\%$ noise, the discontinuous bump can be well detected except for an under-estimation in the peak conductivity value.

\begin{figure}[htbp]
\centering
\setlength{\tabcolsep}{0em}
\begin{tabular}{ccc}
\includegraphics[width=0.33\textwidth]{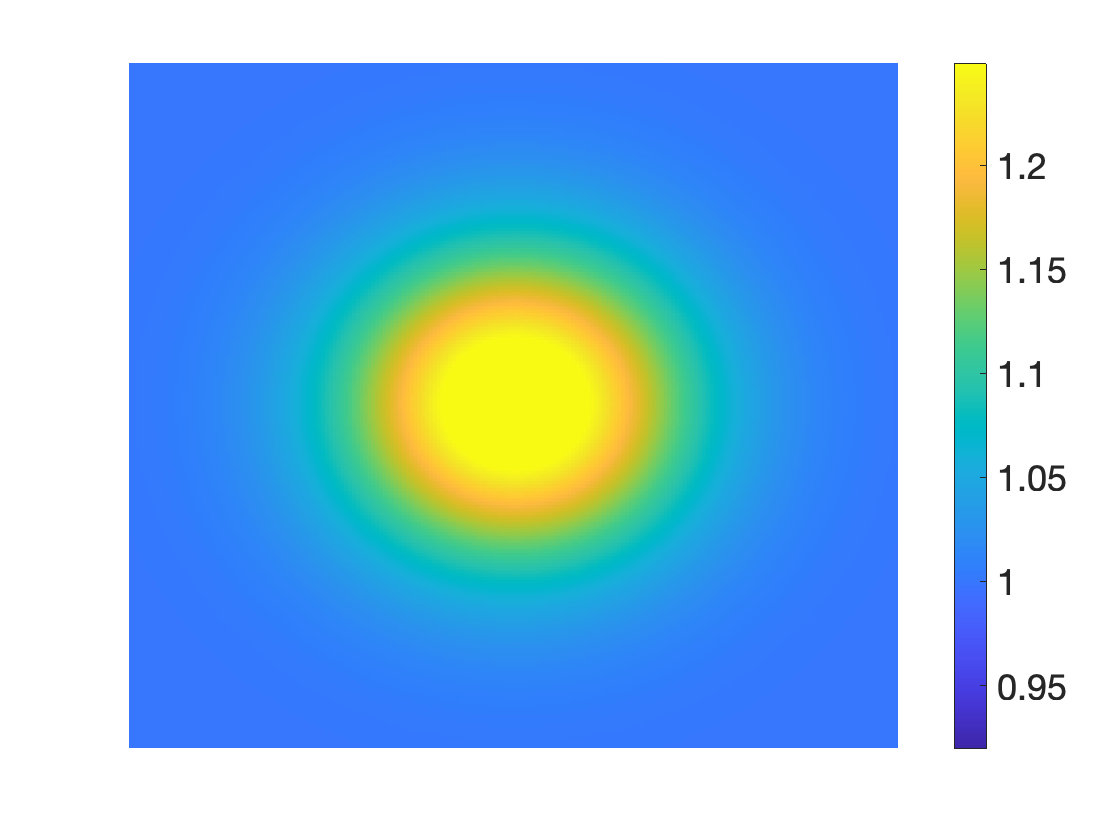} &
\includegraphics[width=0.33\textwidth]{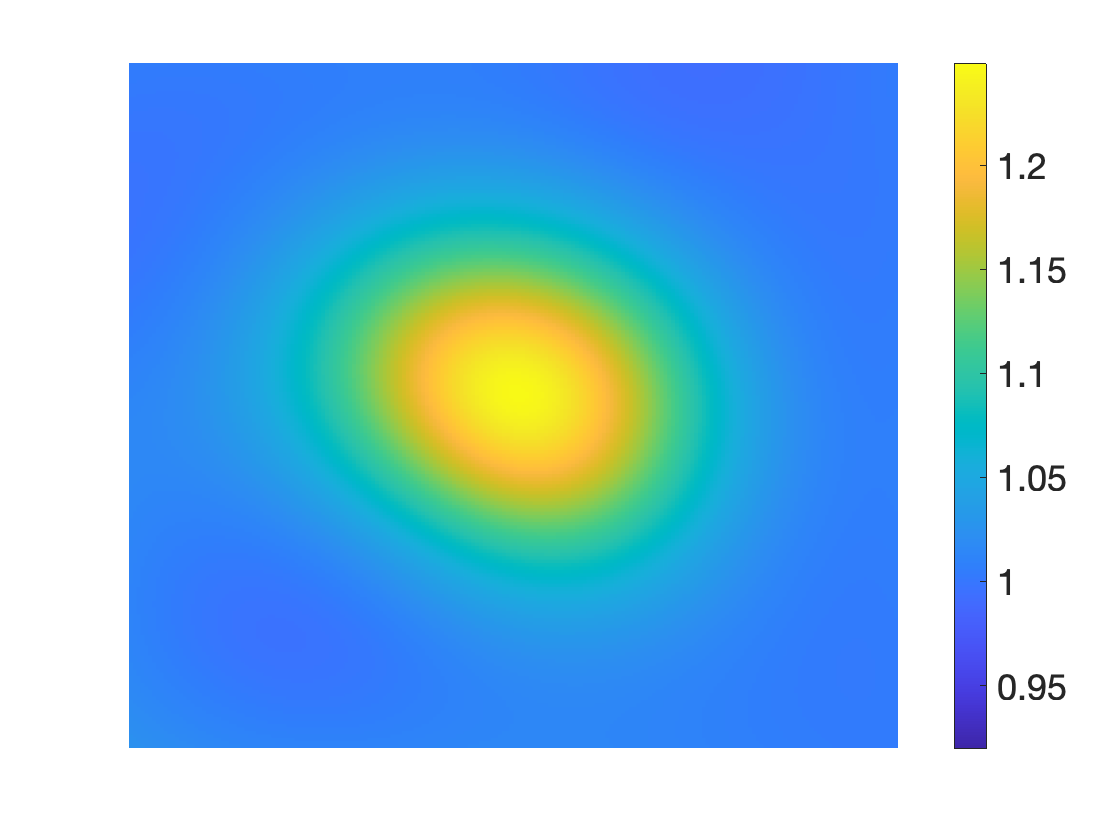} &
\includegraphics[width=0.33\textwidth]{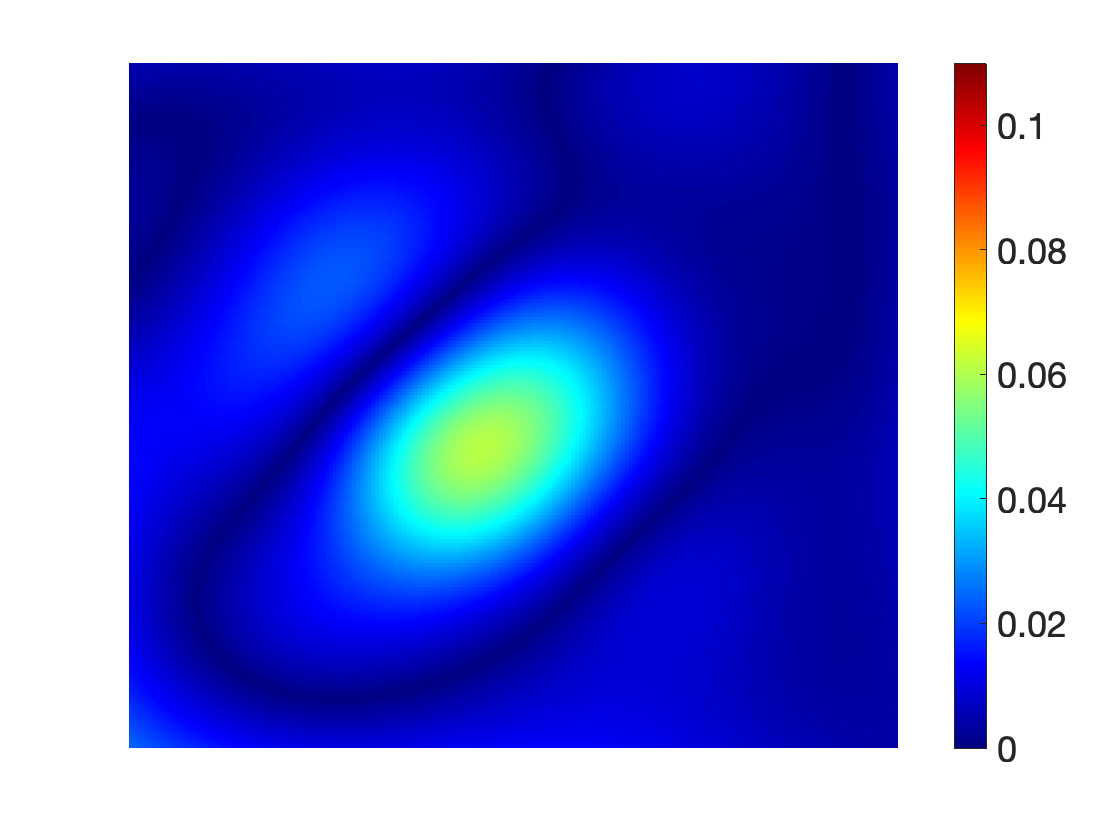}\\
\includegraphics[width=0.33\textwidth]{neu2ex.png} &
\includegraphics[width=0.33\textwidth]{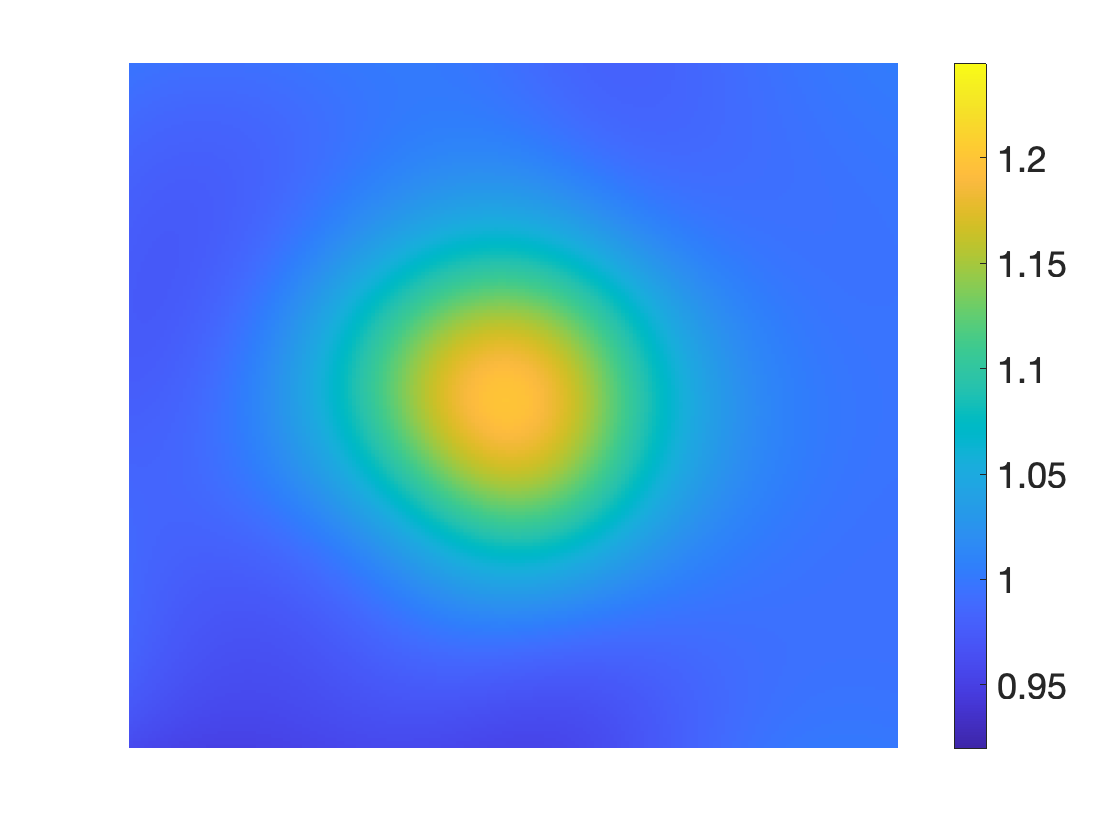} &
\includegraphics[width=0.33\textwidth]{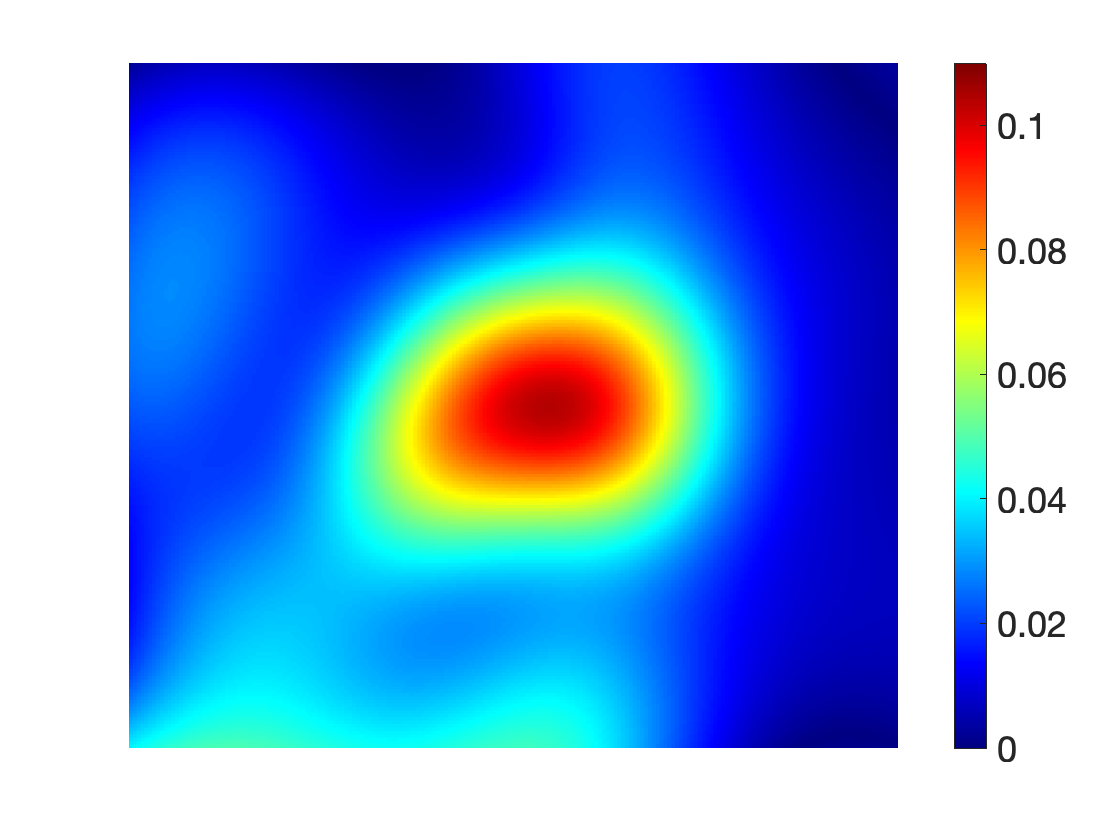}\\
(a) $q^\dag$  & (b) $\hat q$ & (c) $|\hat q-q^\dag|$
\end{tabular}
\caption{The reconstructions for Example \ref{exam:neu2} with exact data $($top$)$ and  noisy data $(\delta=10\%$, bottom$)$.}
\label{fig:neu2}
\end{figure}

The third example is about recovering a 3D conductivity coefficient.
\begin{example}\label{exam:neu2}
$\Omega=(0,1)^3$, $q^\dag=1 + 0.3e^{-20(x_1-0.5)^2-20(x_2-0.5)^2-20(x_3-0.5)^2}$, and $u^\dag=\sum_{i=1}^3(x_i+\frac13x_i^3)$.
\end{example}

Fig. \ref{fig:neu2} shows the reconstruction on a 2D cross section at $x_3=0.5$, for exact data and noisy data ($\delta=10\%$). The location at which the peak conductivity occurs and the overall Gaussian feature of $q^\dag$ are well recovered in both cases, except for a very mild flattening of the Gaussian bump when $10\%$ data noise is present. The relative $L^2(\Omega)$-error $e(\hat q)$ is 1.49e-2 and 3.42e-2,
for exact and noisy data, respectively. The reconstruction error in the noisy case is higher than that for exact data, but still quite acceptable. This again shows the high robustness of the DNN approach with respect to data noise.

\begin{figure}[htbp]
\centering
\setlength{\tabcolsep}{0em}
\begin{tabular}{ccc}
\includegraphics[width=0.32\textwidth]{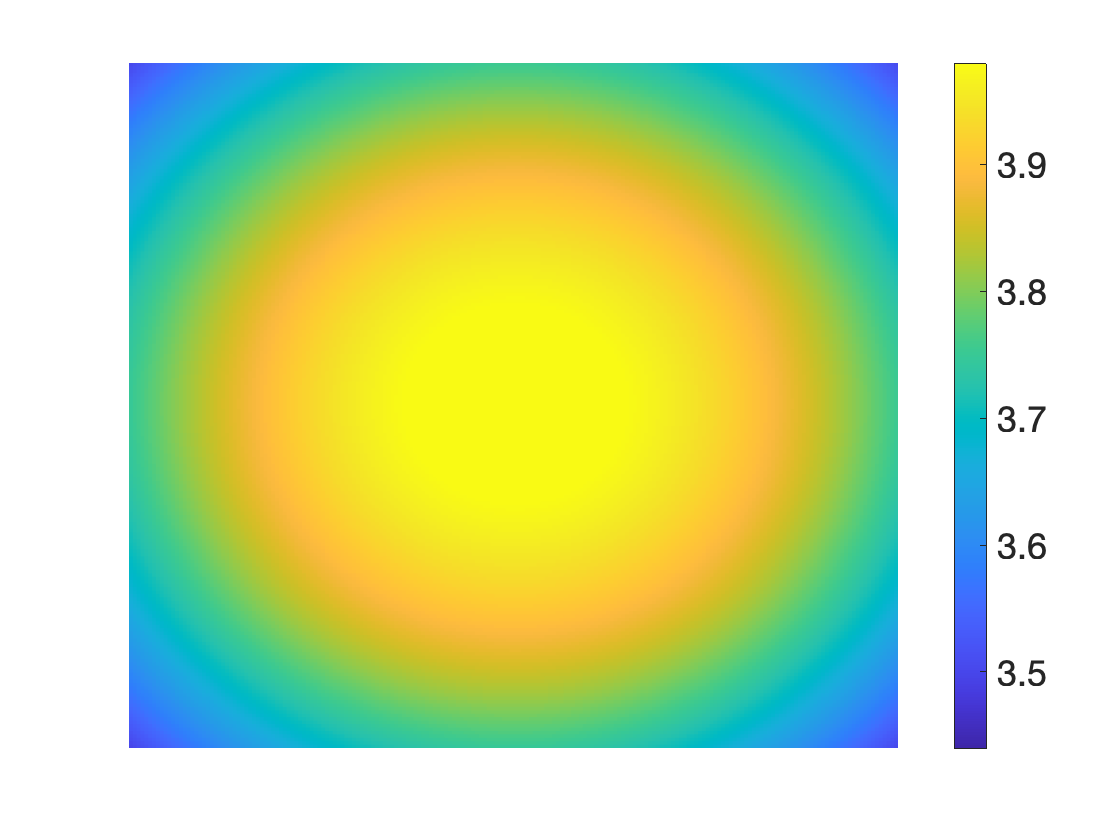} &
\includegraphics[width=0.32\textwidth]{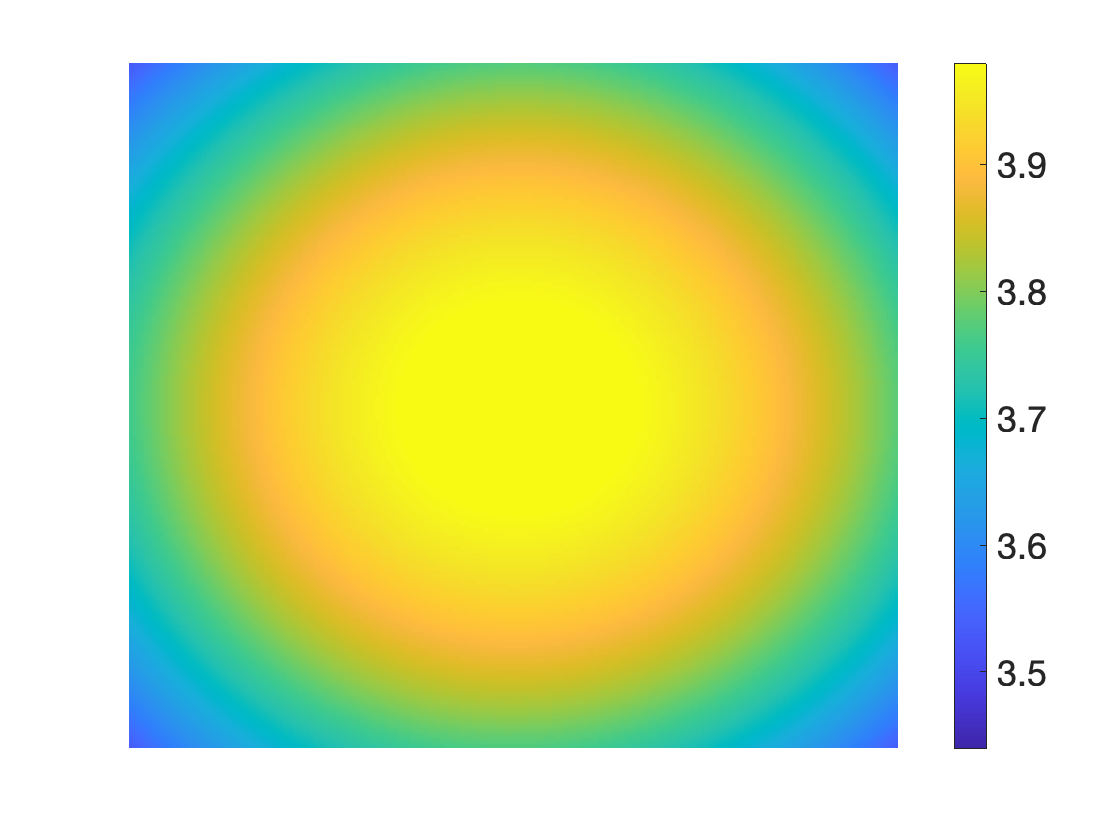} &
\includegraphics[width=0.32\textwidth]{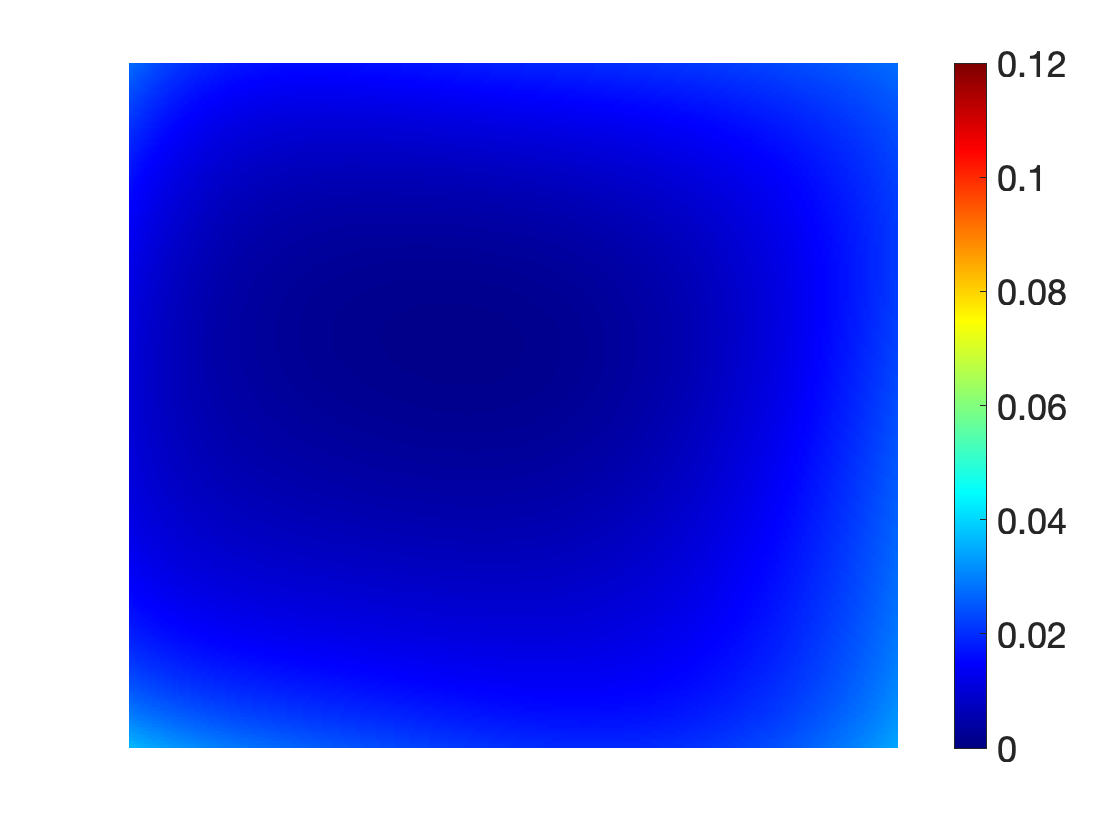}\\
\includegraphics[width=0.32\textwidth]{neudim5ex.png} &
\includegraphics[width=0.32\textwidth]{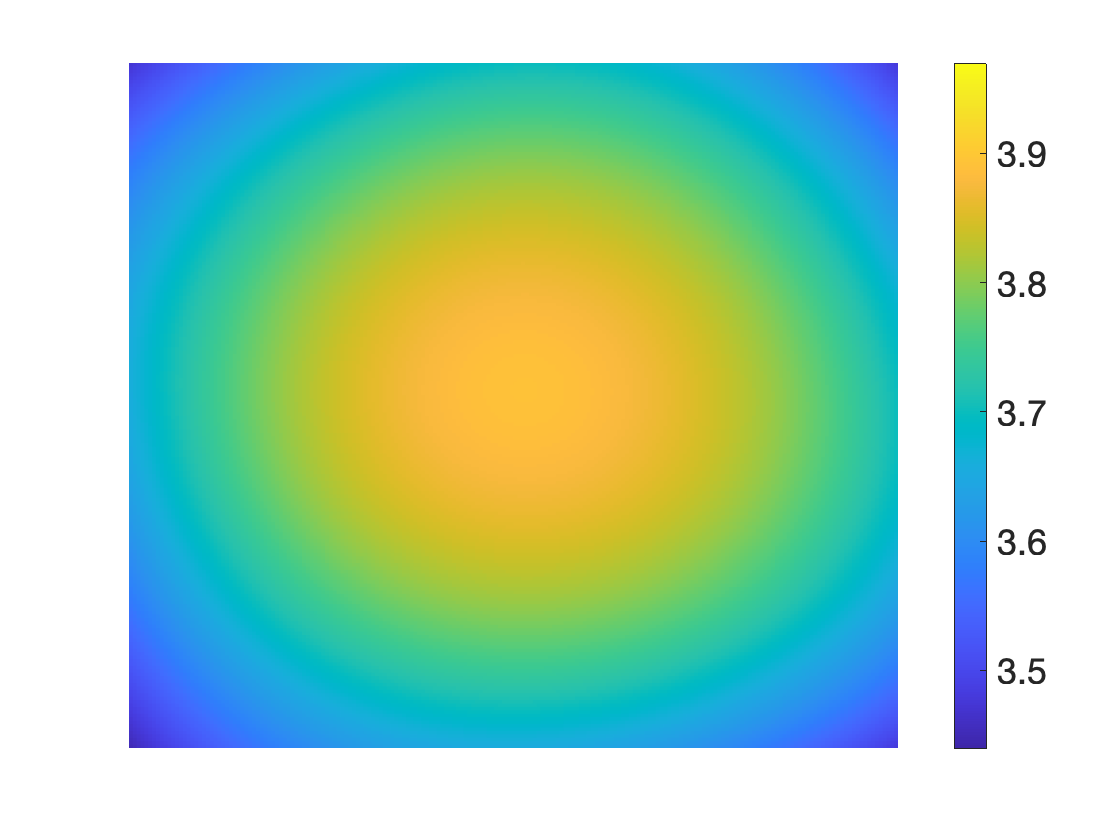} &
\includegraphics[width=0.32\textwidth]{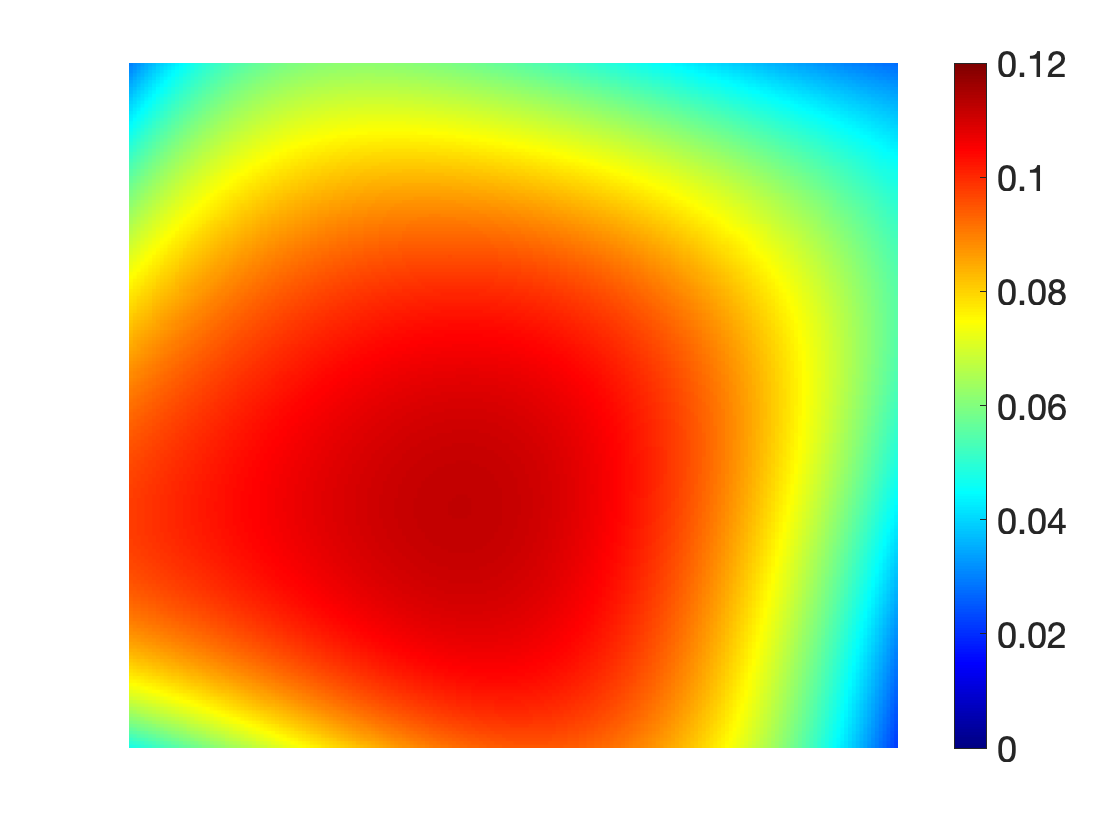}\\
(a) $q^\dag$  & (b) $\hat q$ & (c) $|\hat q-q^\dag|$
\end{tabular}
\caption{The reconstructions for Example \ref{exam:neudim5} with exact data $($top$)$ and noisy data $(\delta=10\%$, bottom$)$.}
\label{fig:neudim5}
\end{figure}

The fourth example is about recovering a conductivity in five dimension.
\begin{example}\label{exam:neudim5}
The domain $\Omega=(0,1)^5$, $q^\dag=1-(x_1-0.5)^2-(x_2-0.5)^2+\cos(\pi(x_3+1.5))+\cos(\pi(x_4+1.5))+\cos(\pi(x_5+1.5))$, and $u^\dag=\sum_{i=1}^5(x_i+\frac13x_i^3)$.
\end{example}

Fig. \ref{fig:neudim5} shows the reconstruction on a 2D cross section at
$x_3=x_4=x_5=0.5$. The relative $L^2(\Omega)$-error $e(\hat q)$ for exact and noisy data (10\%) is 3.44e-3
and 2.27e-2, respectively. The DNN approximations  capture the overall feature of the exact conductivity $q^\dag$
for both exact and noisy data, showing again remarkable robustness of the method. Except for the slight under-estimation in peak conductivity values, the result does not deteriorate much for up to $10\%$ noise. Moreover, it also shows the advantage of using DNNs: it can solve inverse problems for high-dimensional PDEs, which is not easily
tractable for more traditional approaches, e.g.,  FEM. Surprisingly, the approach can produce accurate results without using too many sampling points in the domain $\Omega$, despite the high-dimensionality of the problem.

Next, we present a more ill-posed problem of recovering the conductivity from partial interior data.
\begin{example}\label{exam:neupartial2d}
The domain $\Omega=(0,1)^2$, the measurement $\nabla z^\delta$ on the region $\omega=\Omega\setminus(0.2,0.8)^2$, $q^\dag=2+0.5\sin(2\pi x_1)\sin(2\pi x_2)$, and $u^\dag=x_1+x_2+\frac{1}{3}(x_1^3+x_2^3)$.
\end{example}

\begin{figure}[htbp]
\centering
\setlength{\tabcolsep}{0em}
\begin{tabular}{ccc}
\includegraphics[width=0.32\textwidth]{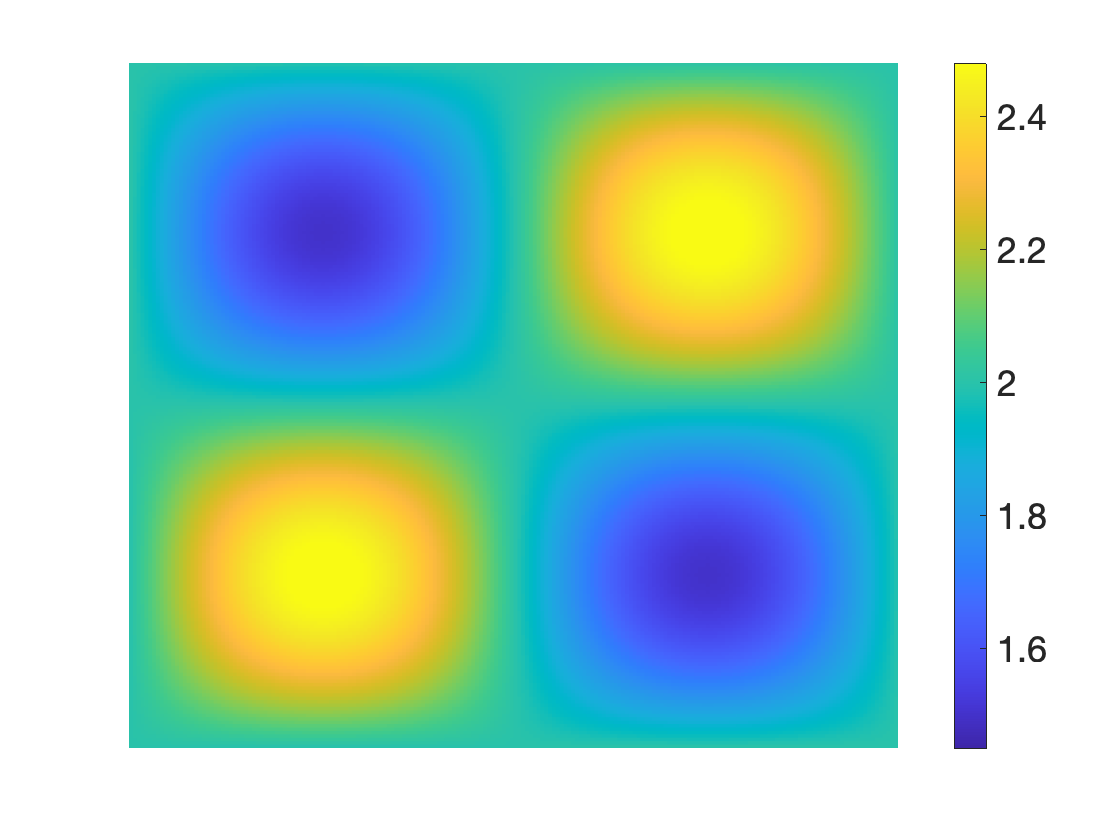} &
\includegraphics[width=0.32\textwidth]{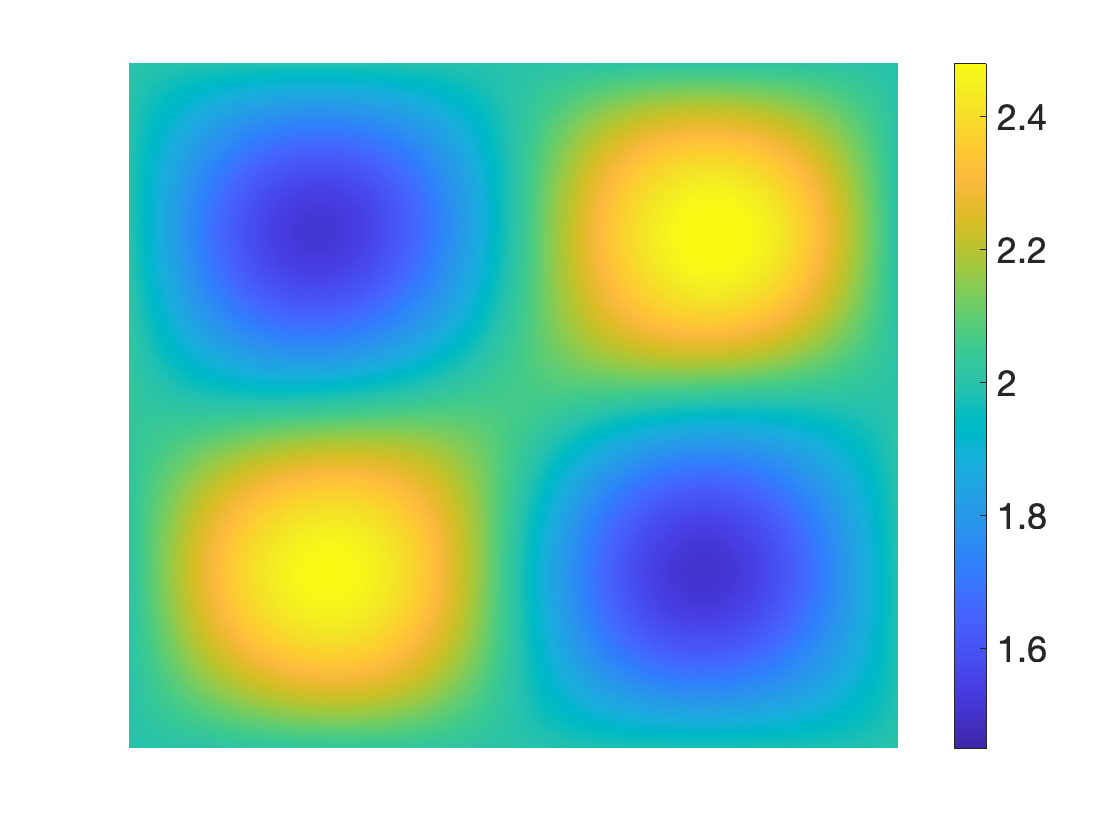} &
\includegraphics[width=0.32\textwidth]{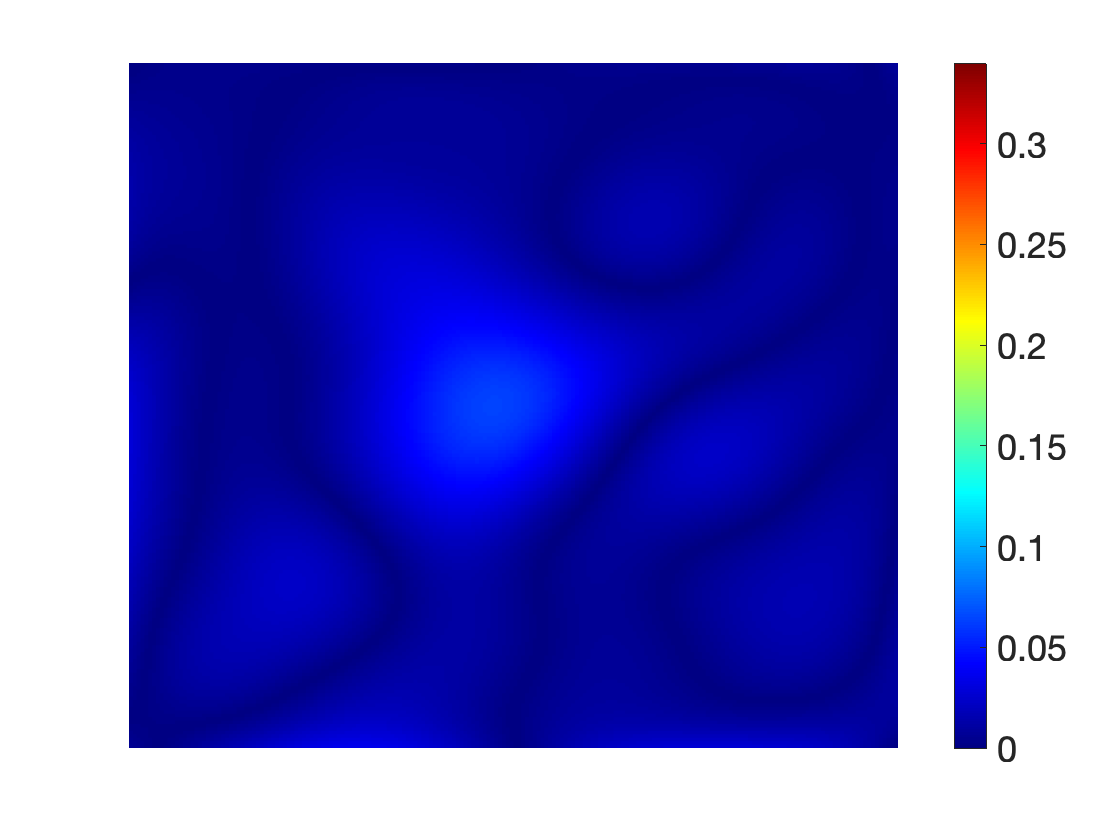}\\
\includegraphics[width=0.32\textwidth]{neupartial2dex.png} &
\includegraphics[width=0.32\textwidth]{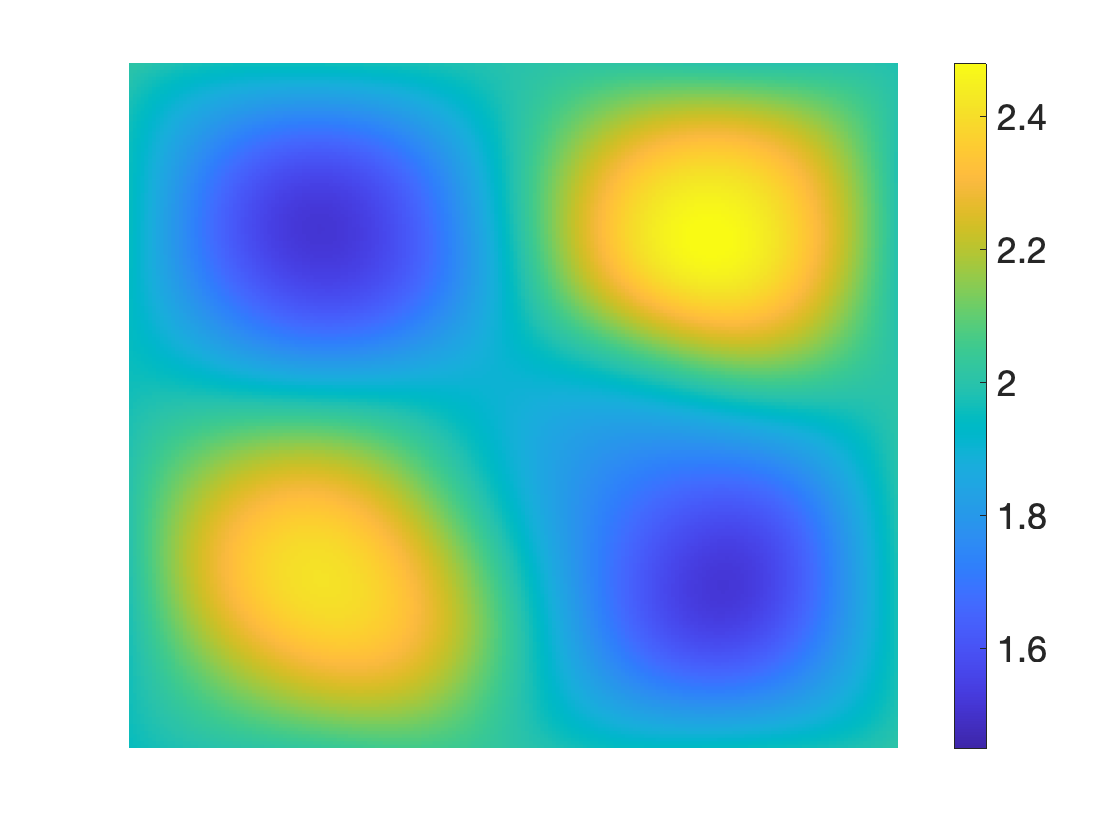} &
\includegraphics[width=0.32\textwidth]{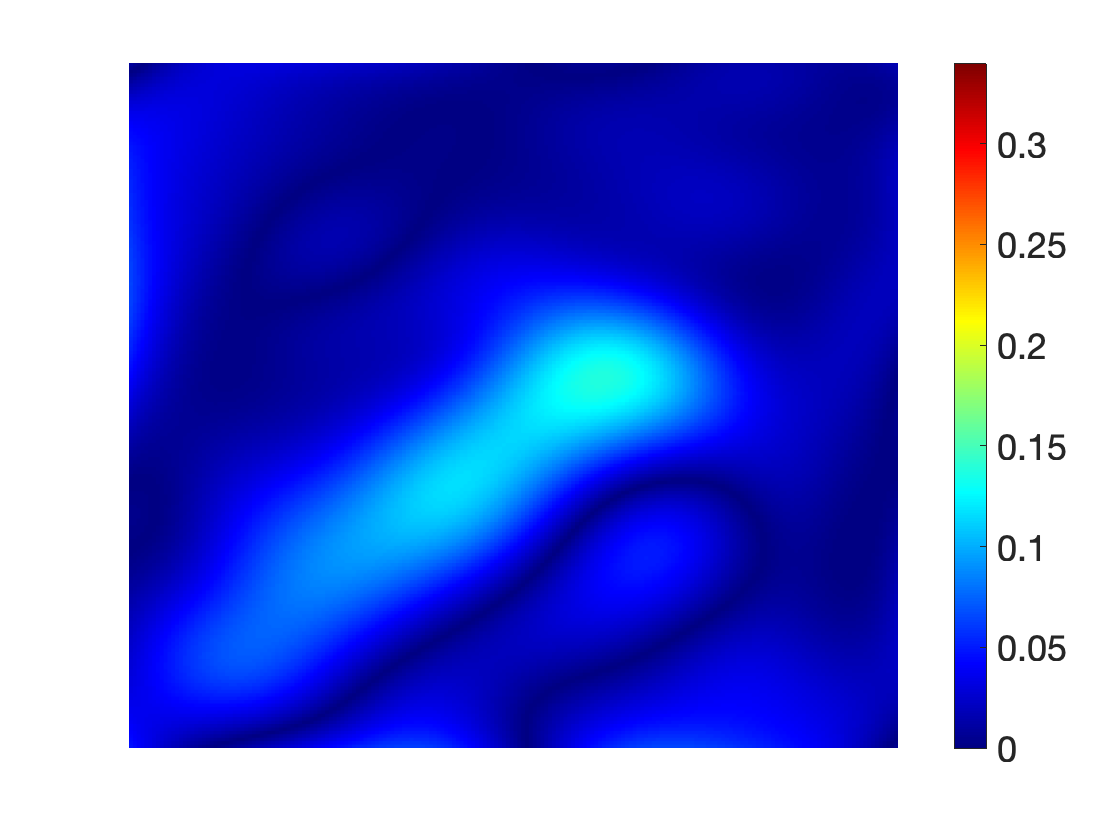}\\
\includegraphics[width=0.32\textwidth]{neupartial2dex.png} &
\includegraphics[width=0.32\textwidth]{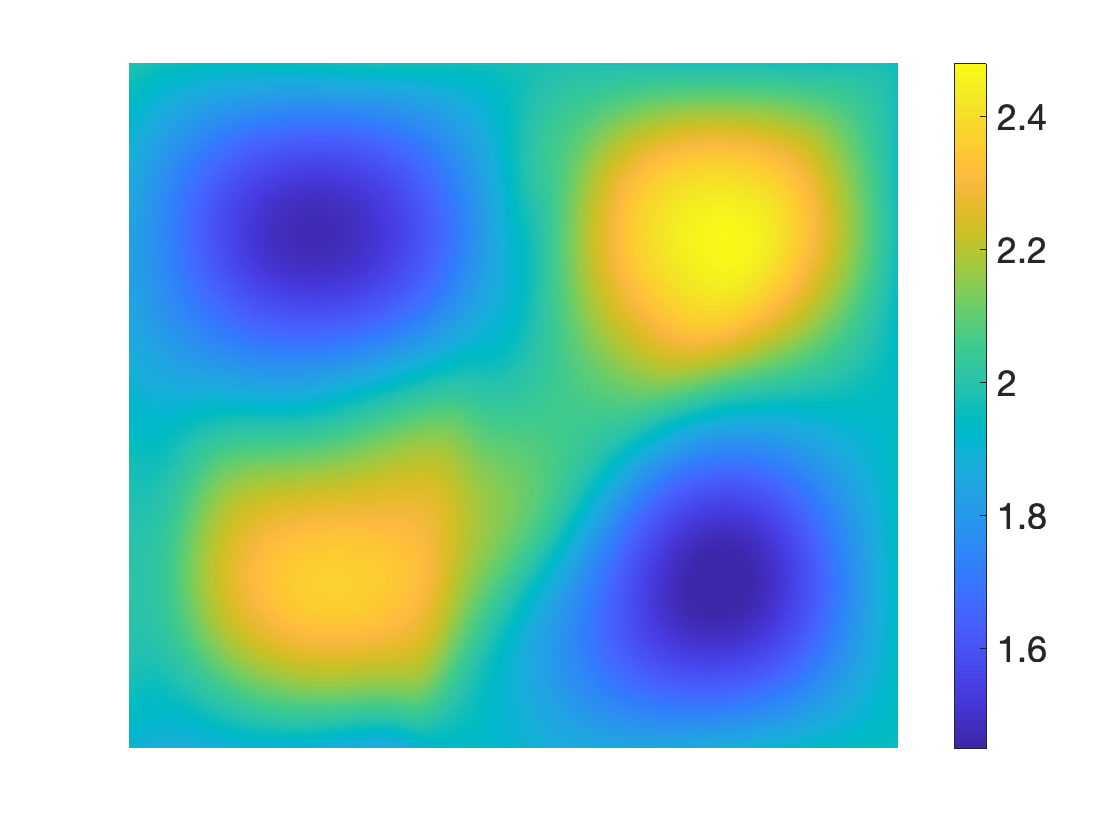} &
\includegraphics[width=0.32\textwidth]{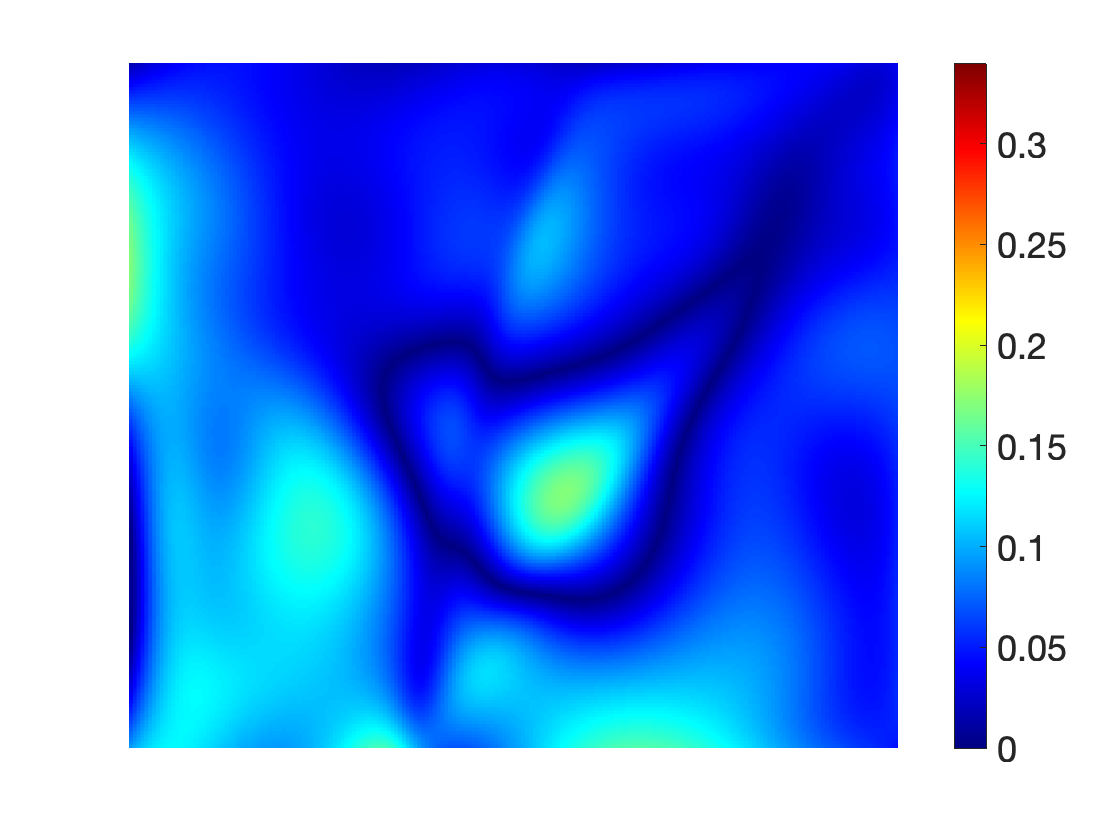}\\
(a) $q^\dag$  & (b) $\hat q$ & (c) $|\hat q-q^\dag|$
\end{tabular}
\caption{The reconstructions for Example \ref{exam:neupartial2d} using DNN method with exact data $($top$)$ and noisy data $(\delta=5\%,10\%$, middle, bottom$)$.}
\label{fig:neupartial2d}
\end{figure}

\begin{figure}[htbp]
\centering
\setlength{\tabcolsep}{0em}
\begin{tabular}{ccc}
\includegraphics[width=0.32\textwidth]{neupartial2dex.png} &
\includegraphics[width=0.32\textwidth]{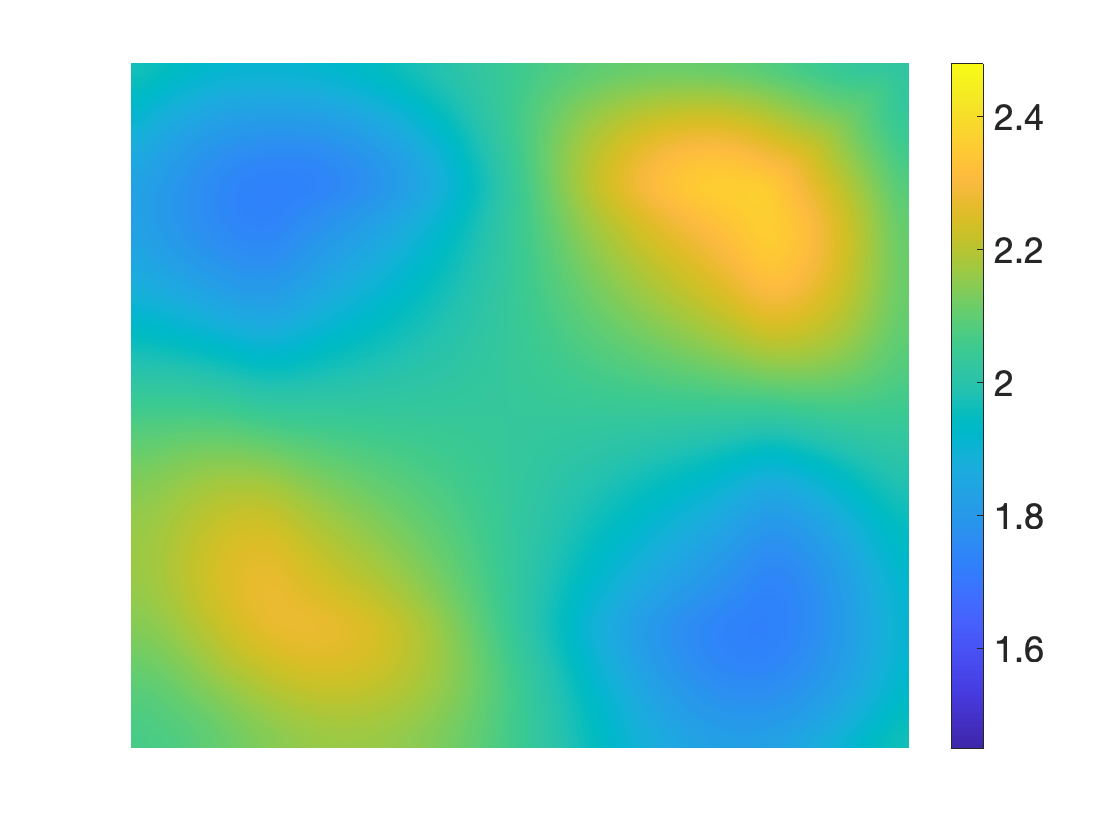} &
\includegraphics[width=0.32\textwidth]{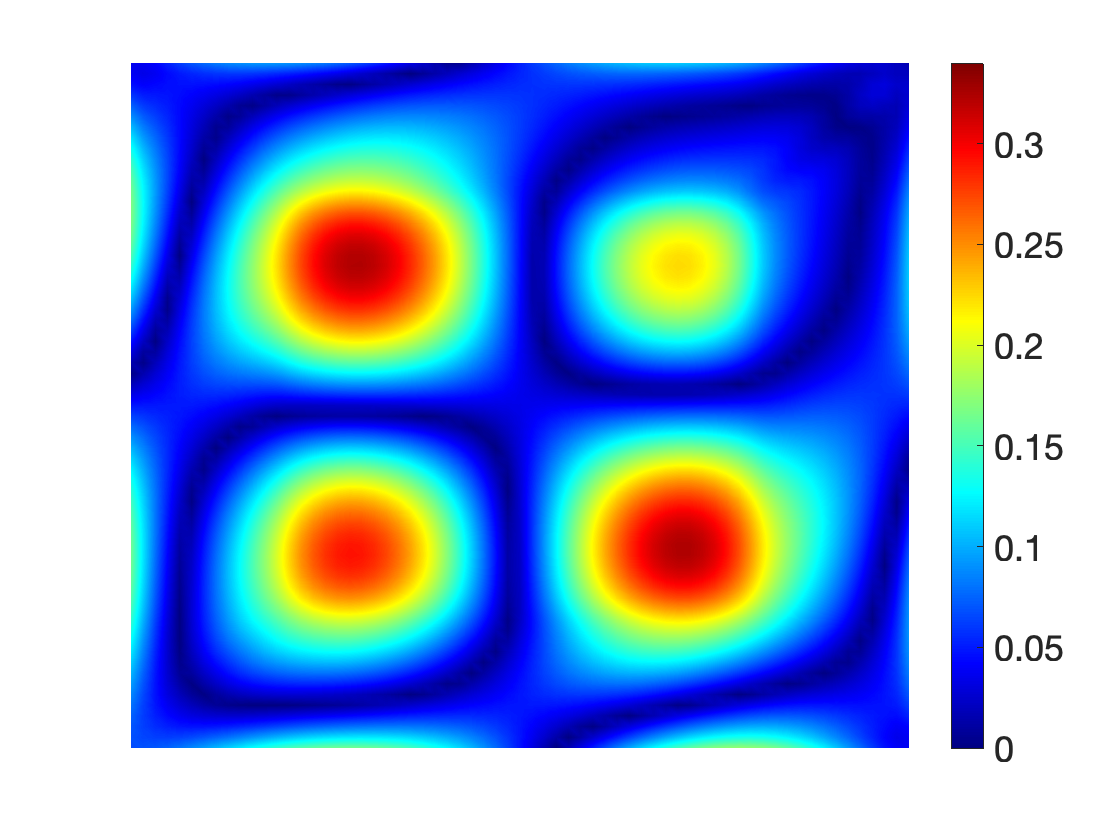}\\
\includegraphics[width=0.32\textwidth]{neupartial2dex.png} &
\includegraphics[width=0.32\textwidth]{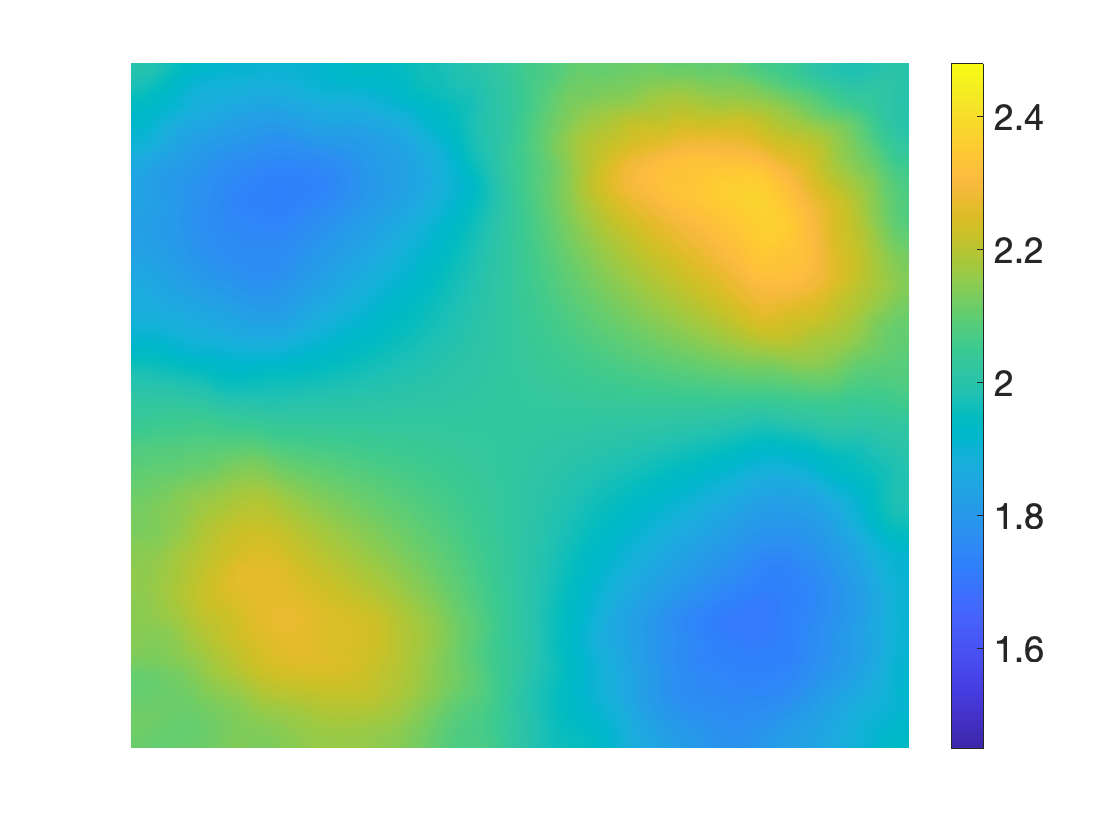} &
\includegraphics[width=0.32\textwidth]{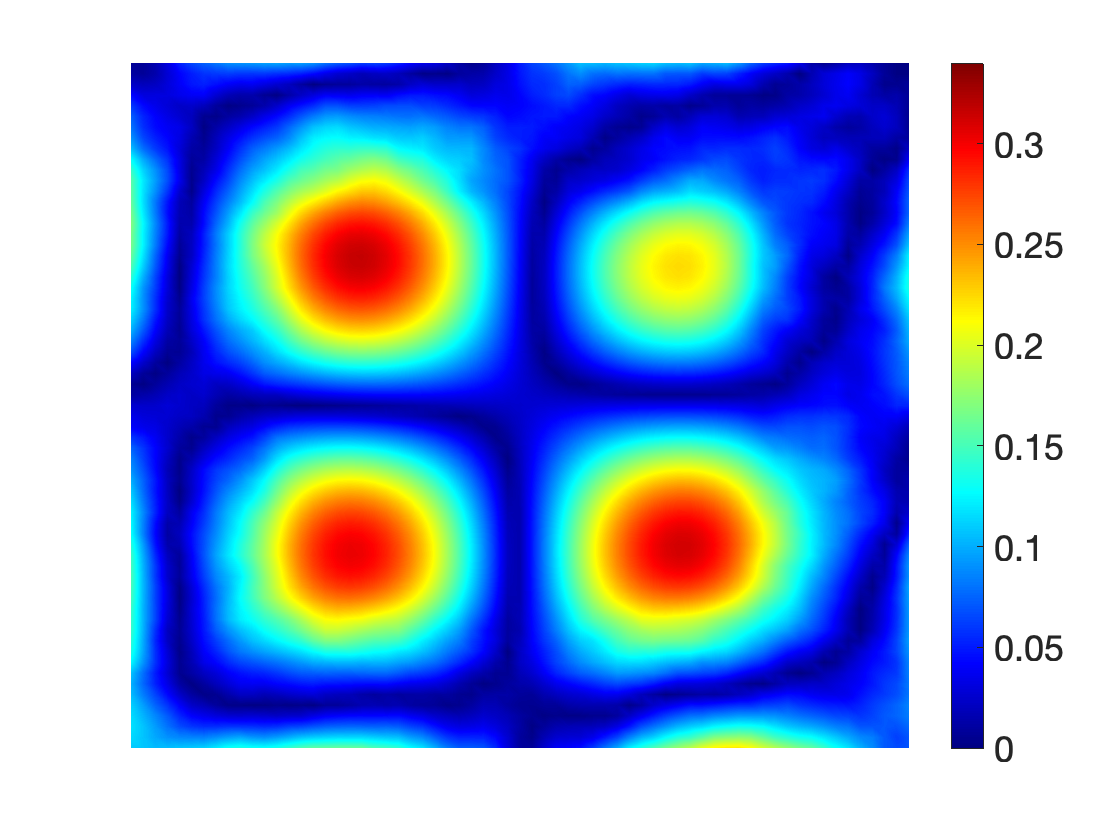}\\
\includegraphics[width=0.32\textwidth]{neupartial2dex.png} &
\includegraphics[width=0.32\textwidth]{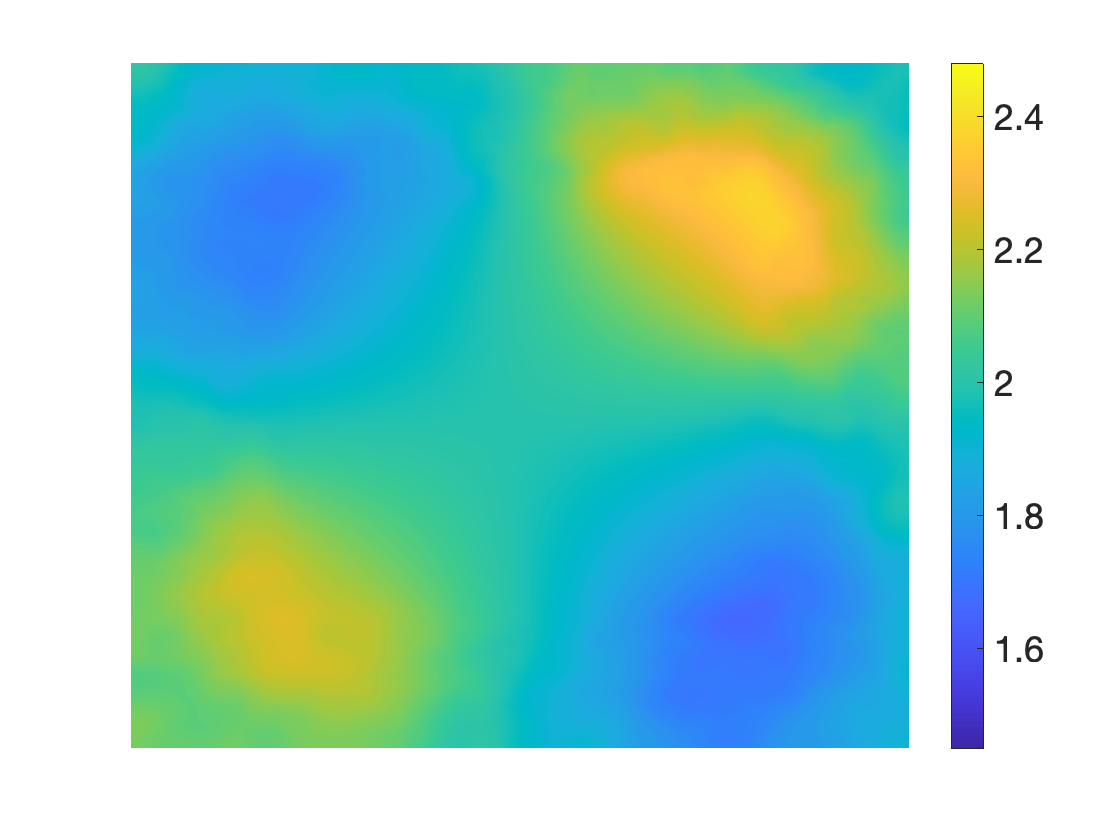} &
\includegraphics[width=0.32\textwidth]{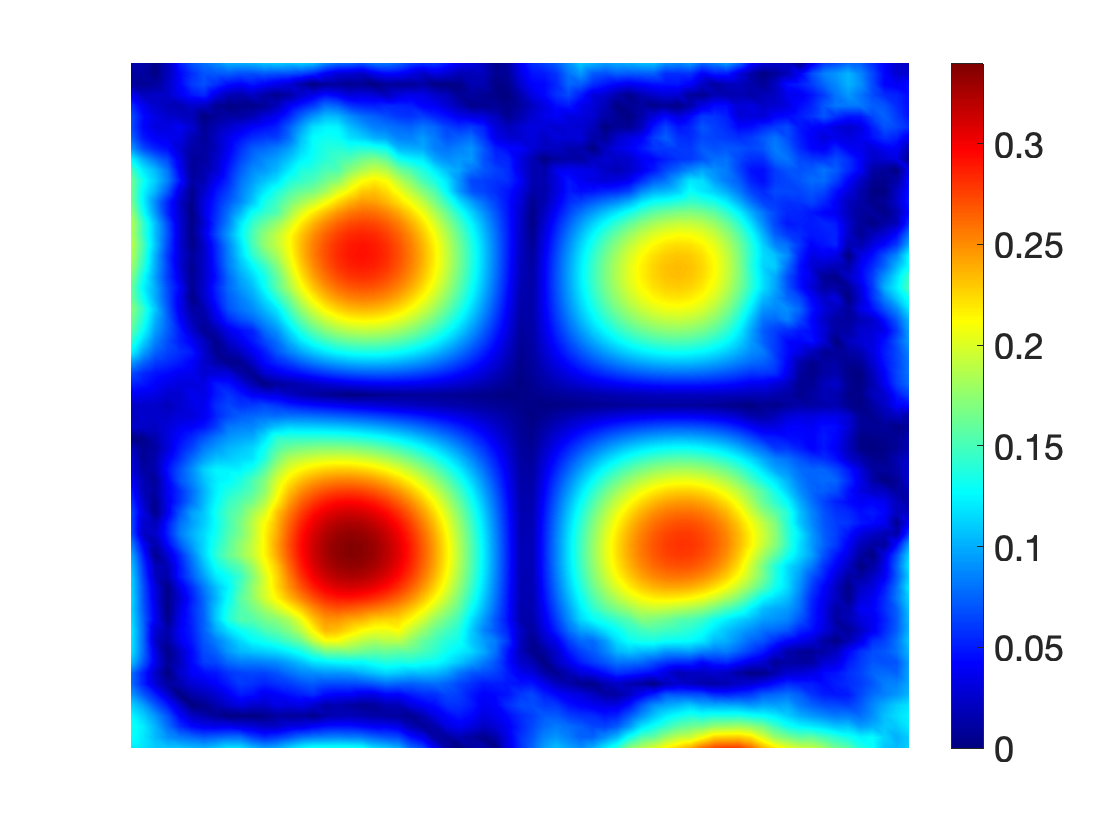}\\
(a) $q^\dag$  & (b) $\hat q$ & (c) $|\hat q-q^\dag|$
\end{tabular}
\caption{The reconstructions for Example \ref{exam:neupartial2d} using FEM with exact data $($top$)$ and noisy data $(\delta=5\%,10\%$, middle, bottom$)$.}
\label{fig:neupartial2dfem}
\end{figure}

Since the data $\nabla z^\delta$ is only given over a subdomain $\omega$ - a very narrow band near the domain boundary $\partial\Omega$ with band width $0.2$, we modify the loss in \eqref{eqn:loss-Neum} as
\begin{equation}\label{eqn:obj-neupartial}
\begin{aligned}
J_{\bsgamma}(\theta,\kappa)=&\|\sigma_\kappa- \P01(q_\theta)\nabla z^\delta\|_{L^2(\omega)}^2+\gamma_\sigma\|\nabla\cdot\sigma_\kappa+f\|_{L^2(\Omega)}^2 \\ &+\gamma_b\|\sigma_\kappa\cdot \n-g\|^2_{L^2(\partial\Omega)}+\gamma_q
    \| \nabla q_\theta\|_{L^2(\Omega)}^2.
\end{aligned}
\end{equation}
The classic FEM can also be used to construct a discrete loss. The reconstructions by both approaches are included to shed light into their stability for partial data. Fig. \ref{fig:neupartial2d} shows that the result by the DNN is accurate on the whole domain $\Omega$, including the central region on which no observational data $\nabla z^\delta$ is given, and it does not change much for noisy data. In contrast, the reconstruction by the FEM exhibits poorer stability, in that the conductivity features in the subdomain $\Omega\setminus\omega$ cannot be recovered at all (regardless of the noise levels), cf. Fig. \ref{fig:neupartial2dfem}. This example shows not only the robustness of the DNN approach, but also its potential for EIT type problems. We leave a comprehensive evaluation on EIT to future work.

The last example is about recovering a 3D conductivity from partial internal data.
\begin{example}\label{exam:neupartial3d}
The domain $\Omega=(0,1)^3$, the measurement $\nabla z^\delta$ on the region $\omega=\Omega\setminus (0.2,0.8)^3$, $q^\dag=2+0.25\cos(\pi (x_1+1.5))+0.5\cos(\pi (x_2+1.5))+0.5x_3^2$, and $u^\dag=x_1+x_2+x_3+\frac{1}{3}(x_1^3+x_2^3+x_3^3)$.
\end{example}

\begin{figure}[htbp]
\centering
\setlength{\tabcolsep}{0em}
\begin{tabular}{ccc}
\includegraphics[width=0.32\textwidth]{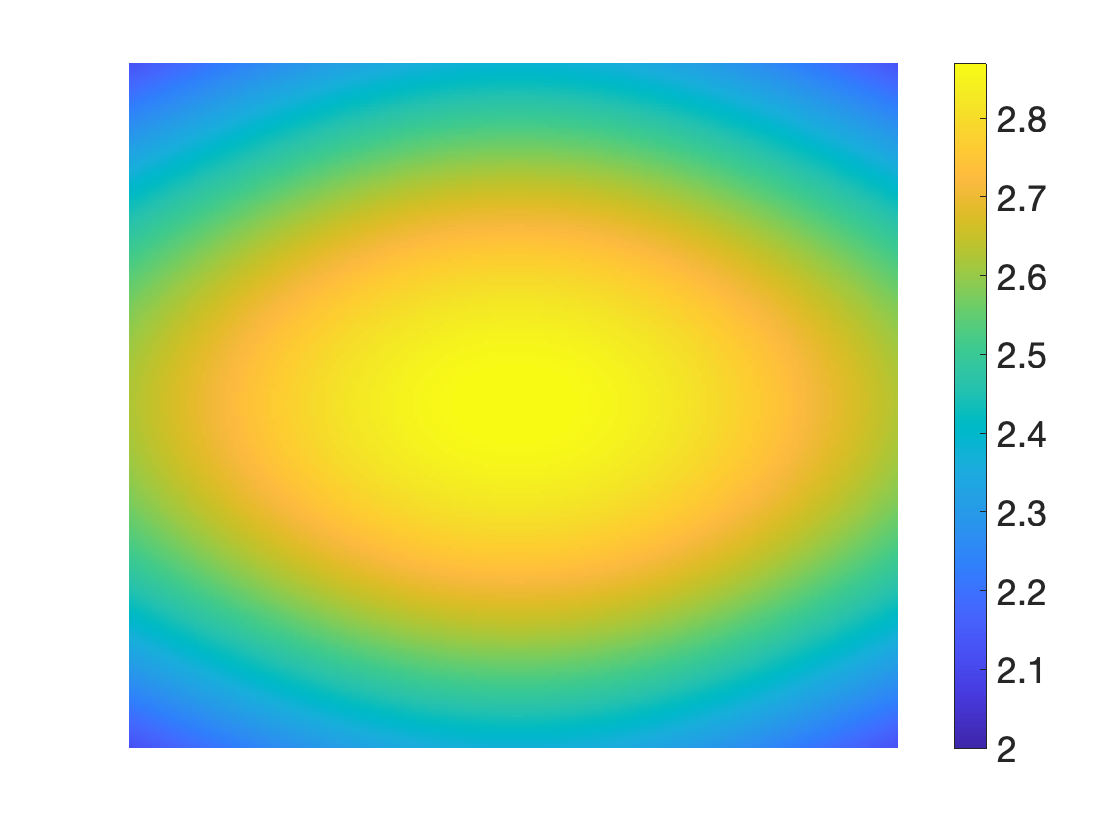} &
\includegraphics[width=0.32\textwidth]{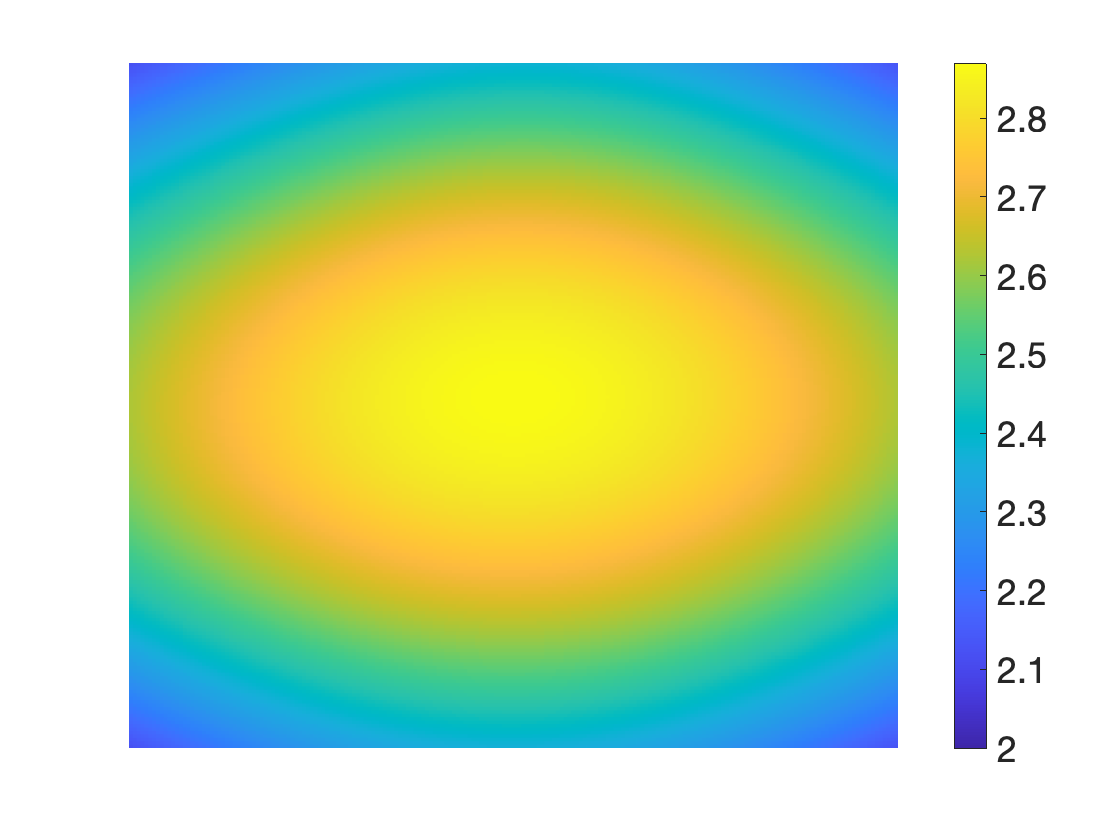} &
\includegraphics[width=0.32\textwidth]{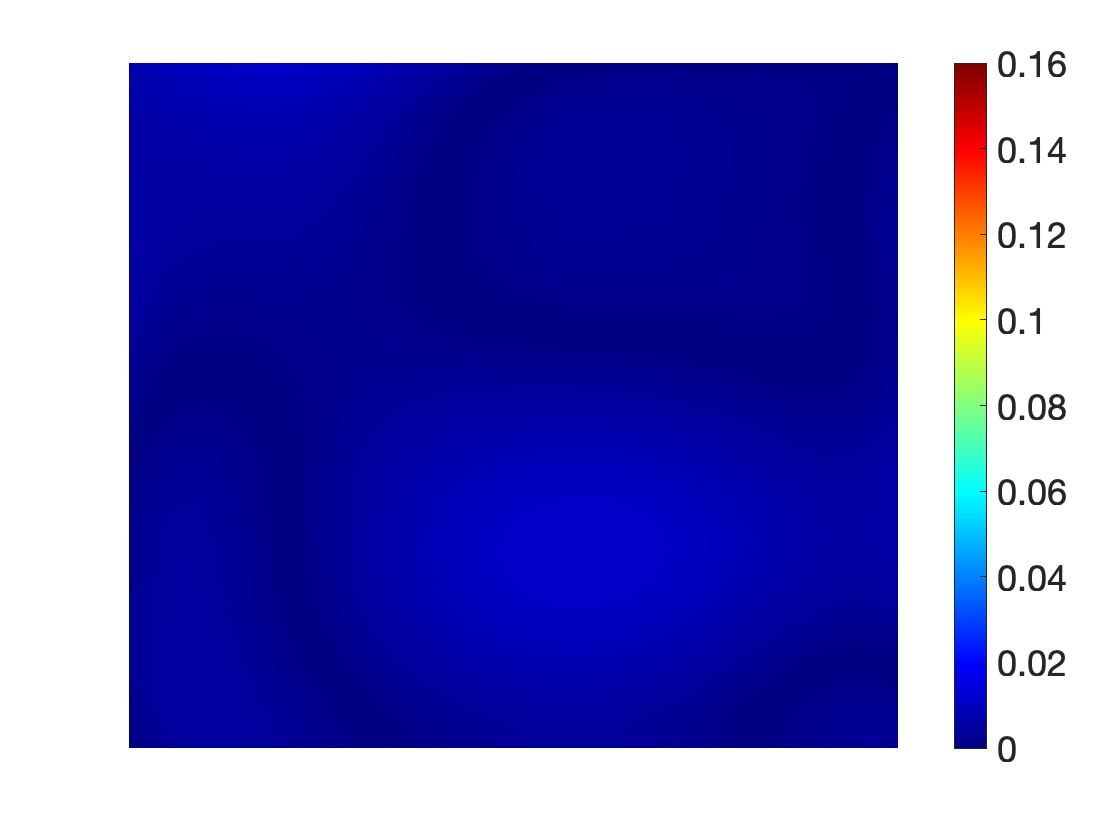}\\
\includegraphics[width=0.32\textwidth]{neupartial3dex.png} &
\includegraphics[width=0.32\textwidth]{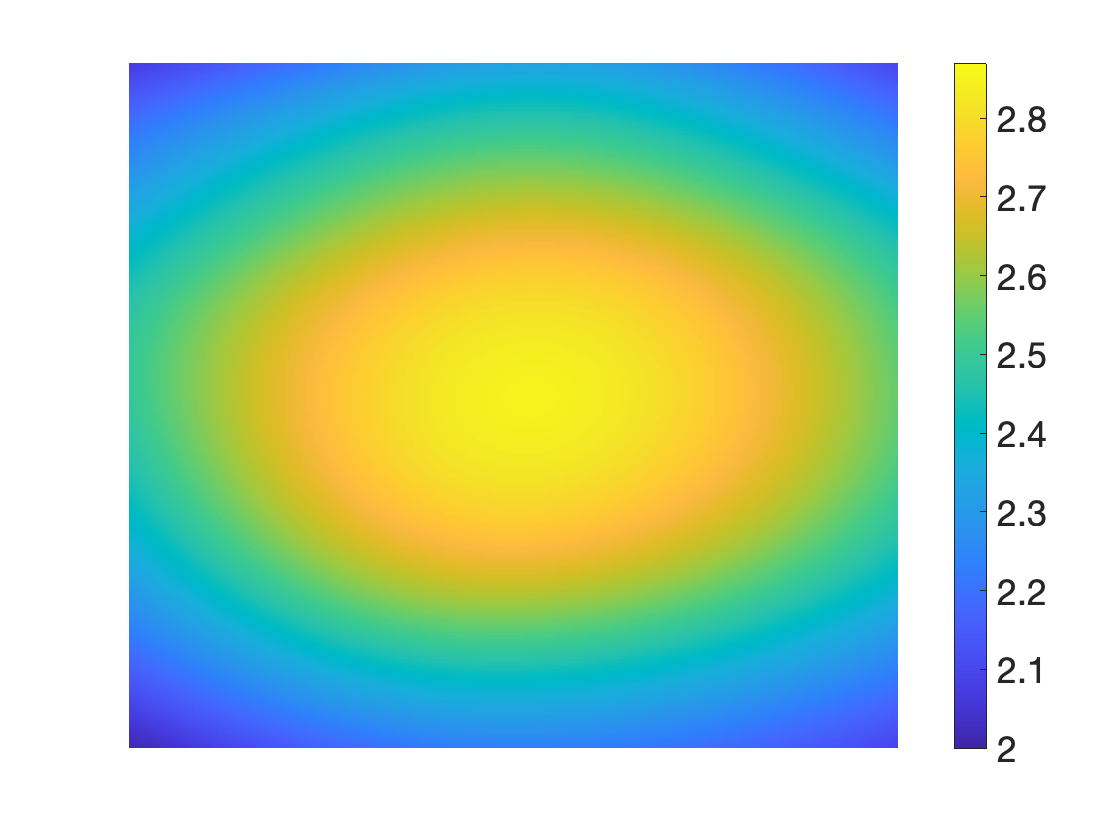} &
\includegraphics[width=0.32\textwidth]{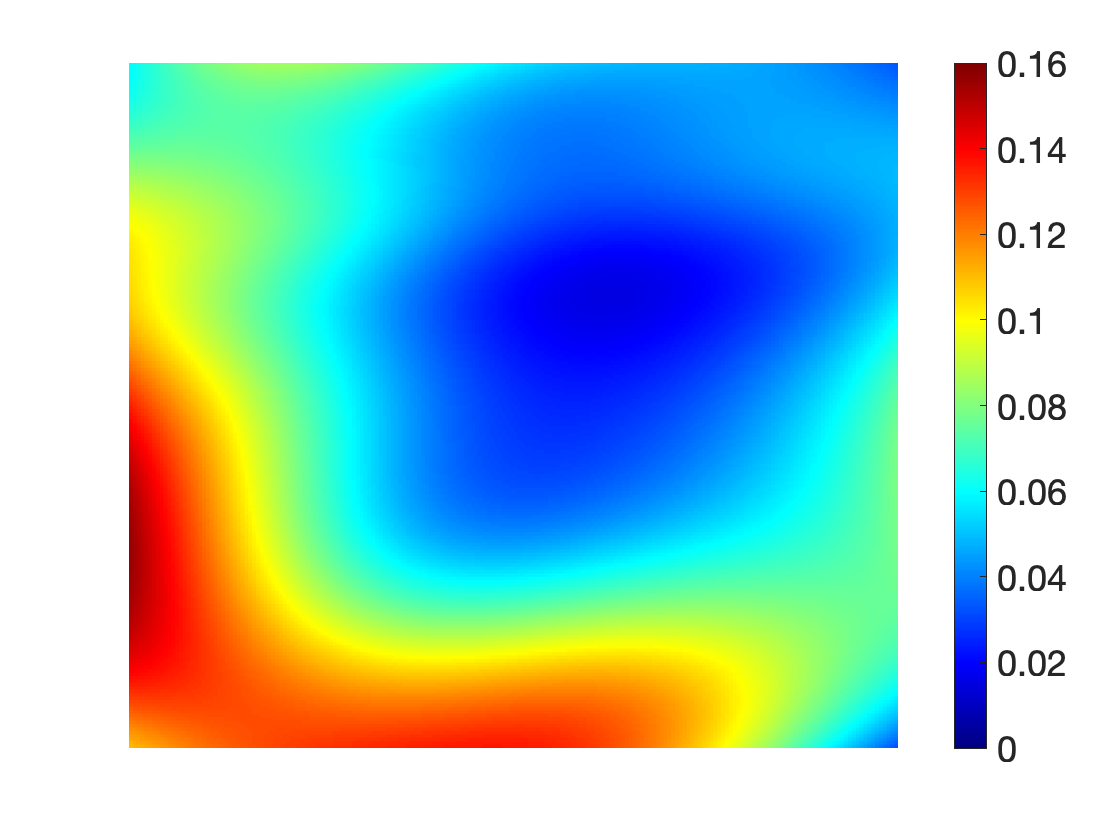}\\
(a) $q^\dag$  & (b) $\hat q$ & (c) $|\hat q-q^\dag|$
\end{tabular}
\caption{The reconstructions for Example \ref{exam:neupartial3d} using DNN method with exact data $($top$)$ and noisy data $(\delta=10\%$, bottom$)$.}
\label{fig:neupartial3d}
\end{figure}

Fig. \ref{fig:neupartial3d} shows the reconstruction on a 2D cross section at $x_3=0.5$. The
relative $L^2(\Omega)$-error $e(\hat q)$ for exact and noisy data is 1.82e-3 and 2.87e-2,
respectively. Note that the reconstruction is accurate over the whole domain $\Omega$,
including the central region $\Omega\setminus \omega$ where there is no observational data,
and the result does not deteriorate much for up to $10\%$ noise.

\subsection{The Dirichlet problem}\label{sec:diri}
In this part, we present a set of experiments for the Dirichlet case.
The first example is about recovering a smooth 2D conductivity.

\begin{figure}[htbp]
\centering
\setlength{\tabcolsep}{0em}
\begin{tabular}{ccc}
\includegraphics[width=0.32\textwidth]{diri1ex.png} &
\includegraphics[width=0.32\textwidth]{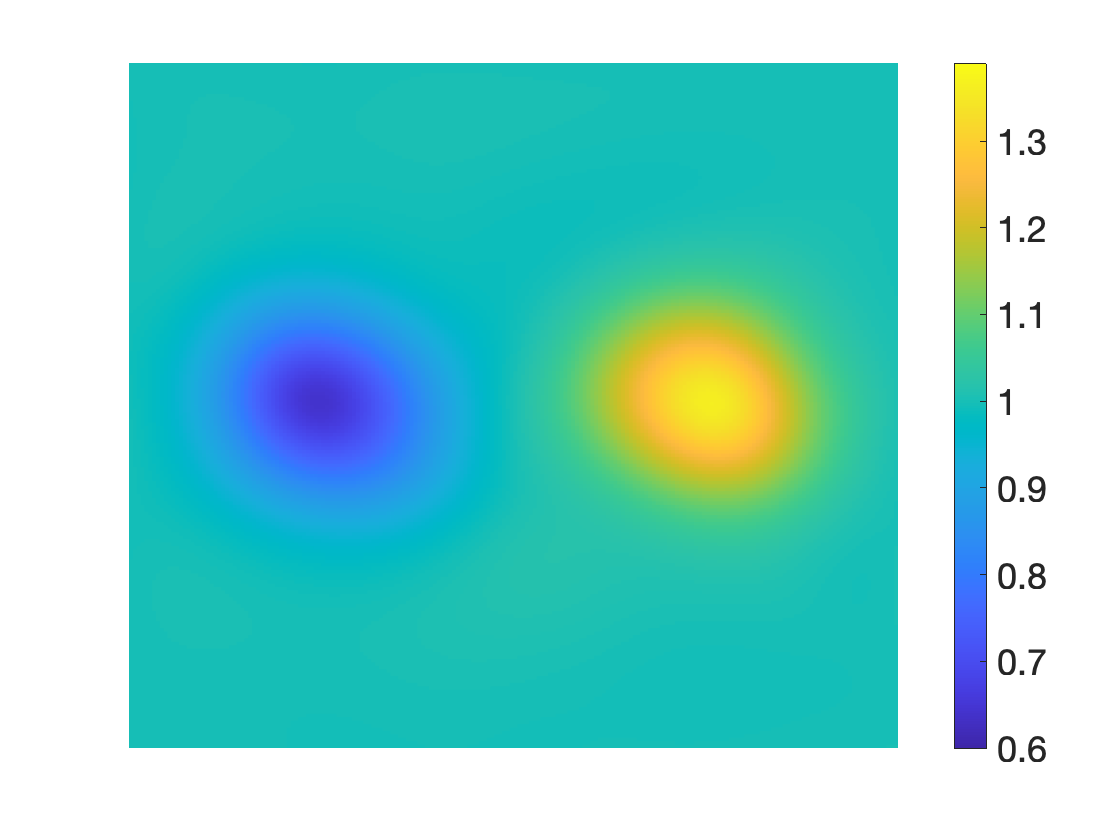} &
\includegraphics[width=0.32\textwidth]{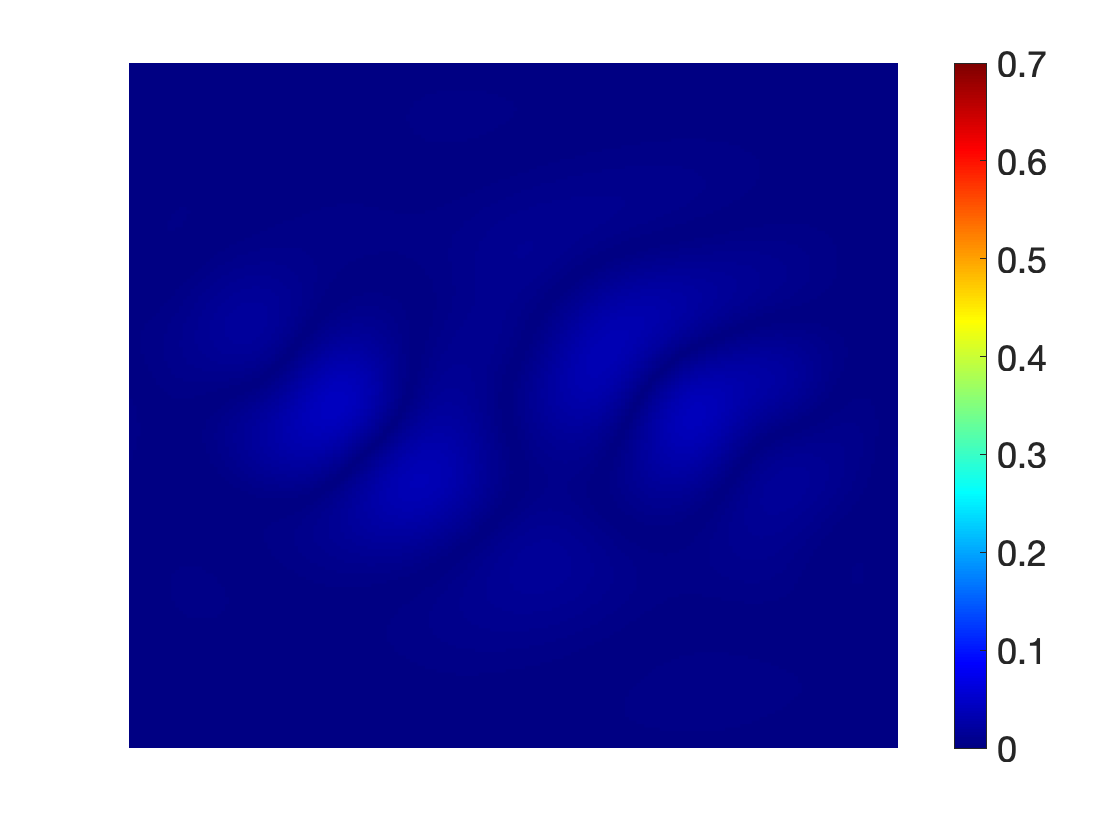}\\
\includegraphics[width=0.32\textwidth]{diri1ex.png} &
\includegraphics[width=0.32\textwidth]{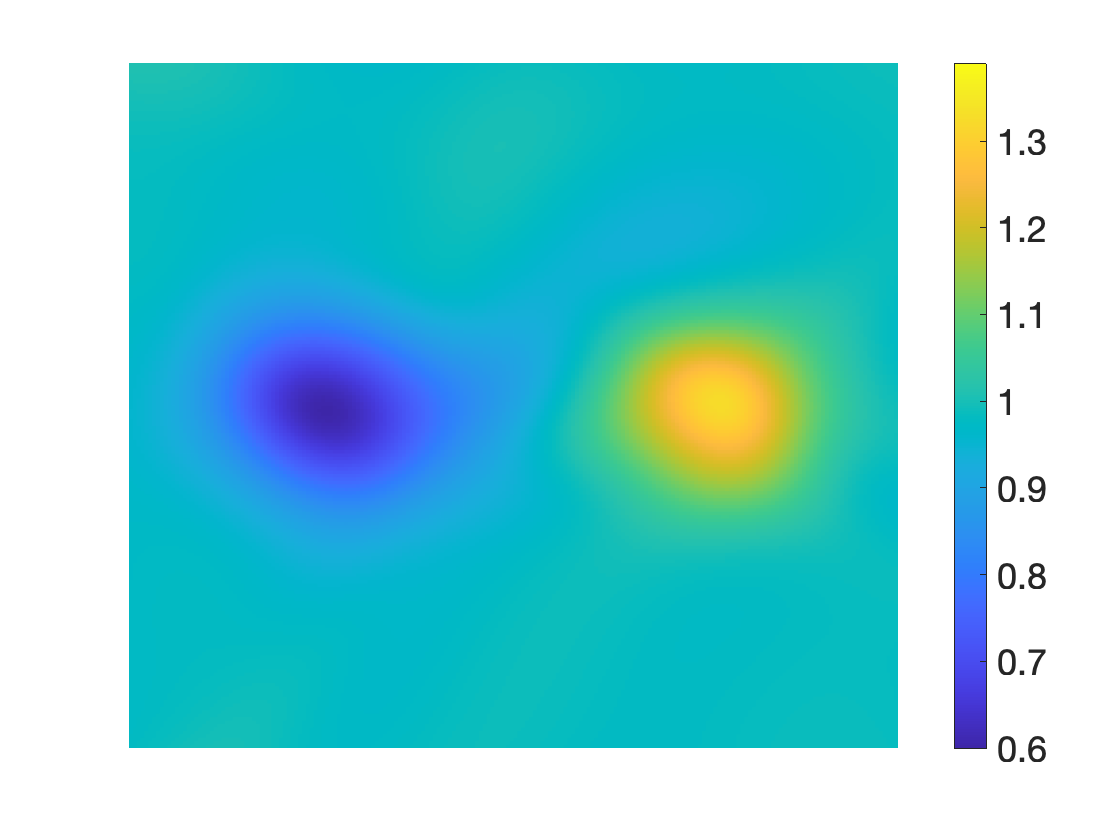} &
\includegraphics[width=0.32\textwidth]{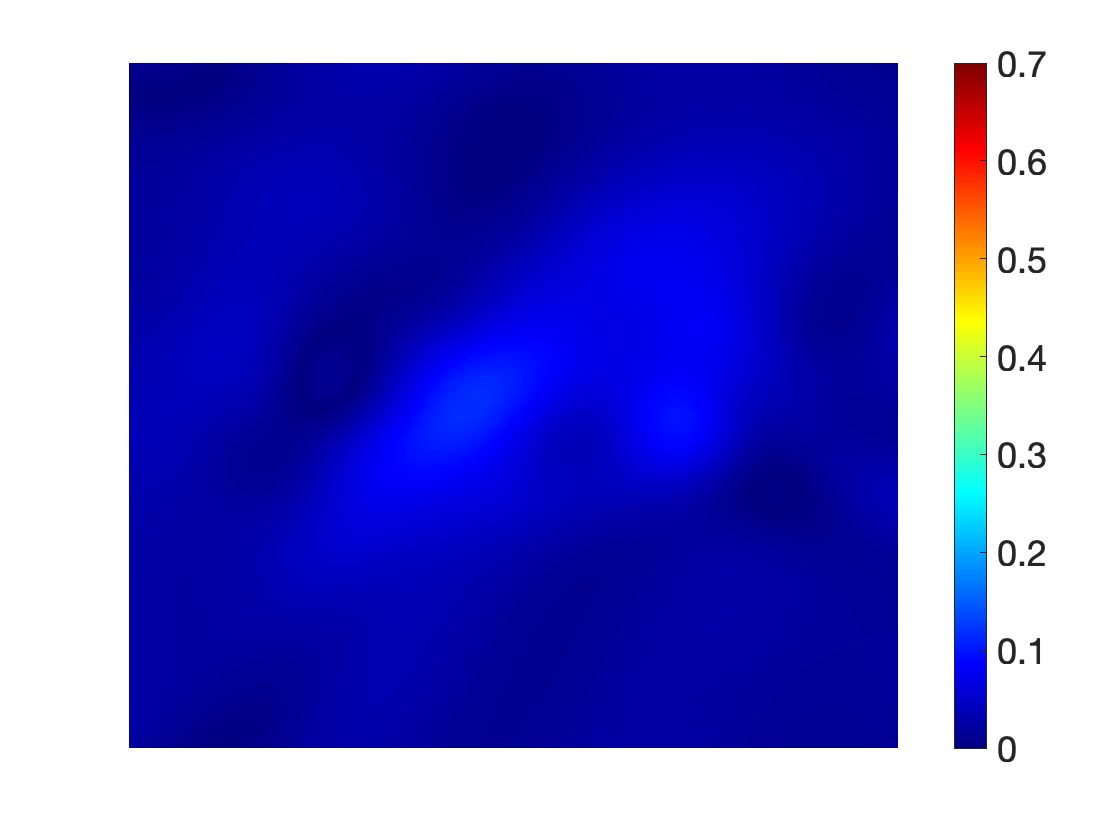}\\
\includegraphics[width=0.32\textwidth]{diri1ex.png} &
\includegraphics[width=0.32\textwidth]{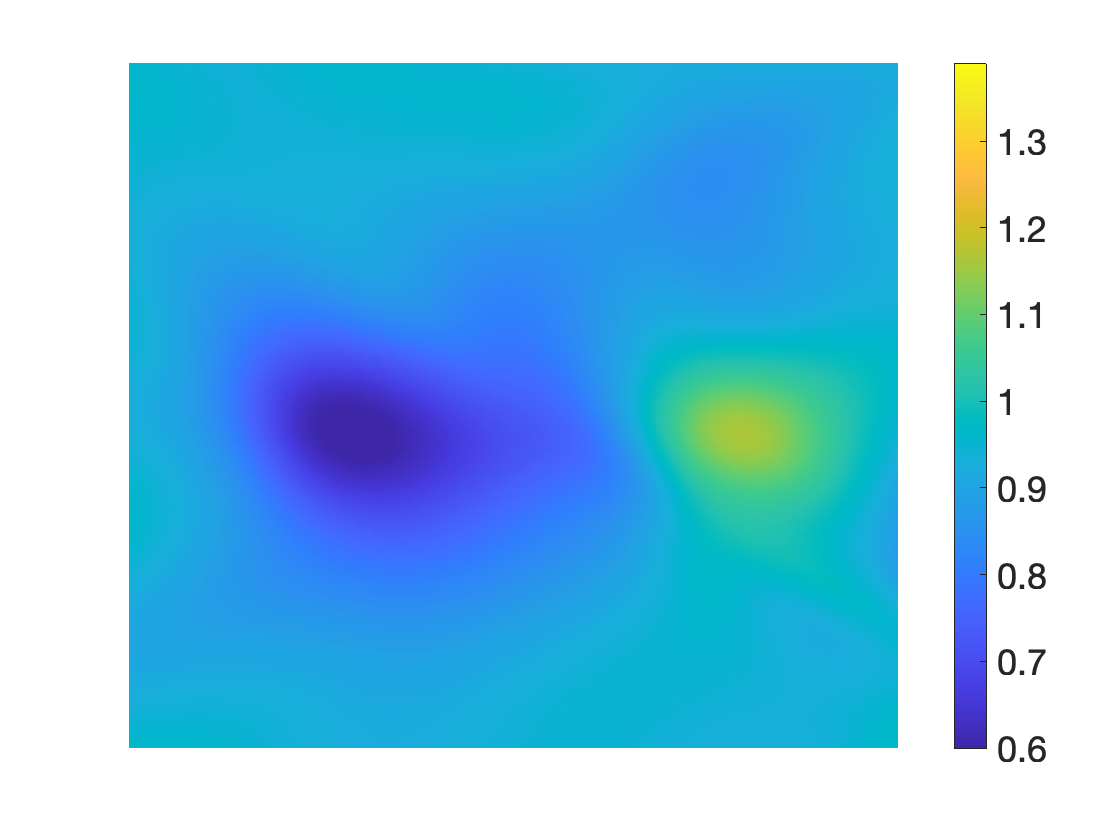} &
\includegraphics[width=0.32\textwidth]{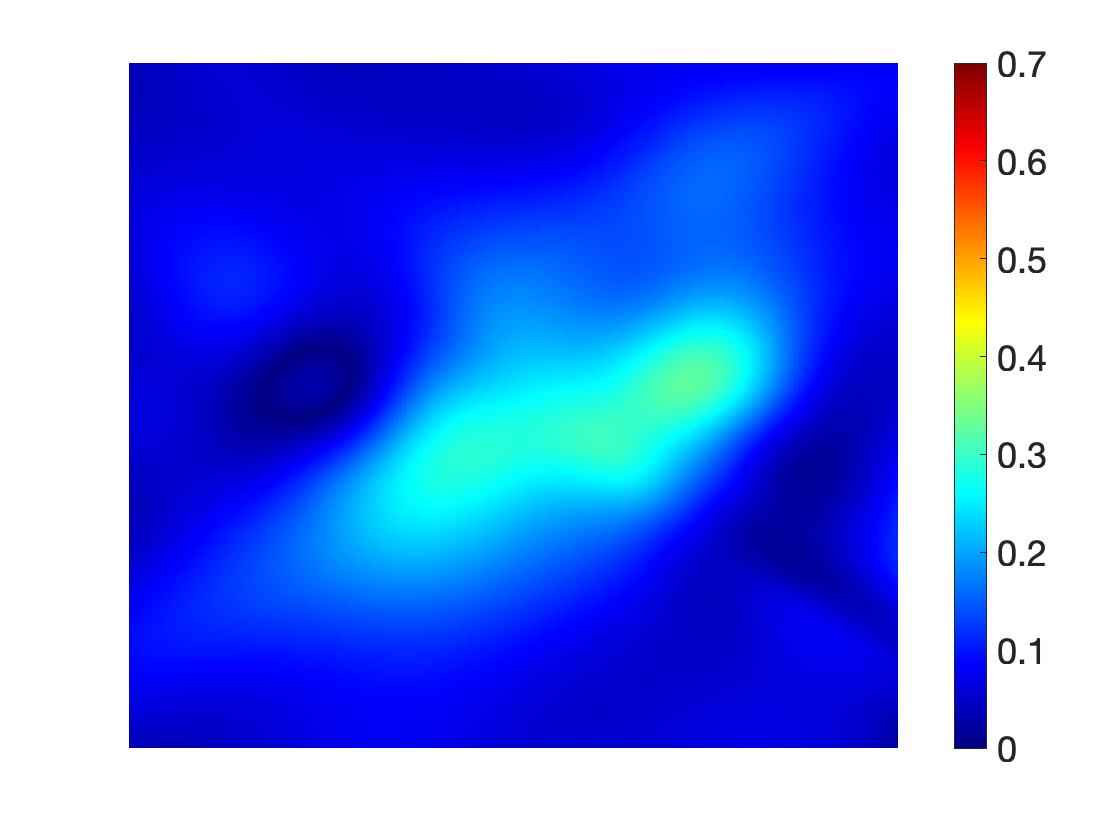}\\
(a) $q^\dag$  & (b) $\hat q$ & (c) $|\hat q-q^\dag|$
\end{tabular}
\caption{The reconstructions for Example \ref{exam:diri1} using DNN method with exact data $($top$)$ and noisy data $(\delta=10\%,\ 20\%$, middle, bottom$)$.}
\label{fig:diri1}
\end{figure}

\begin{figure}[htbp]
\centering
\setlength{\tabcolsep}{0em}
\begin{tabular}{ccc}
\includegraphics[width=0.32\textwidth]{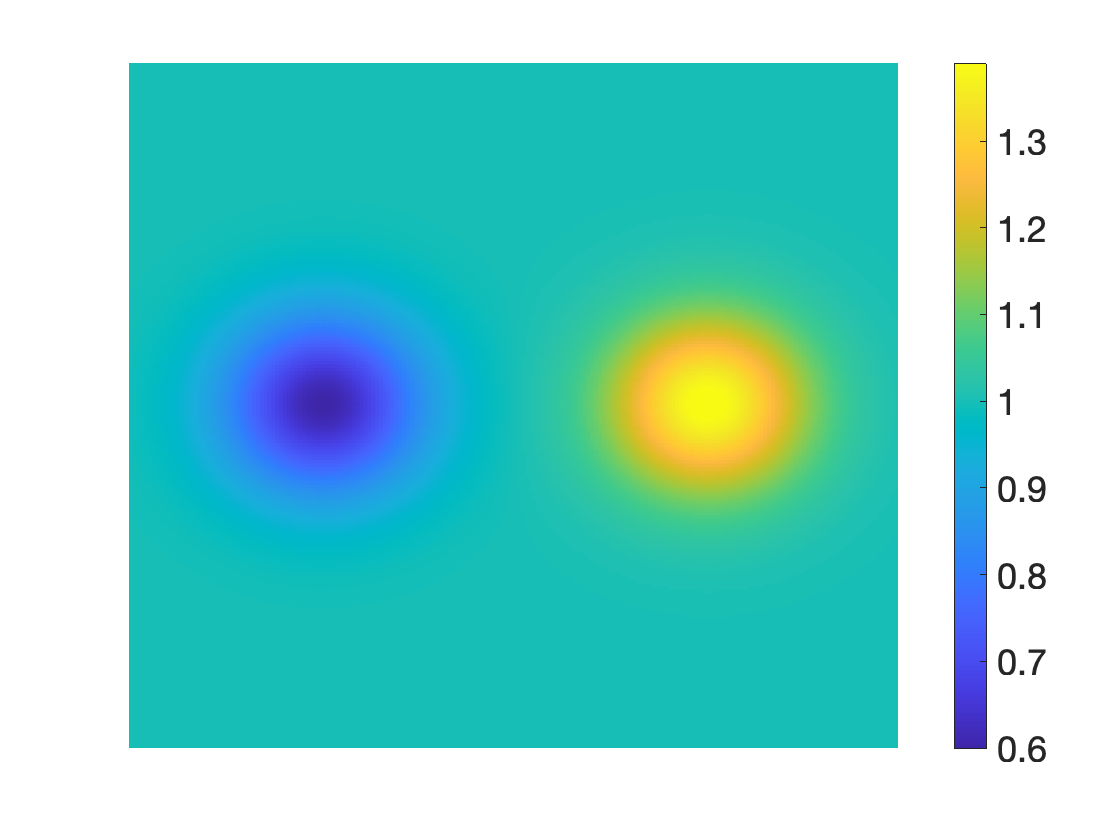} &
\includegraphics[width=0.32\textwidth]{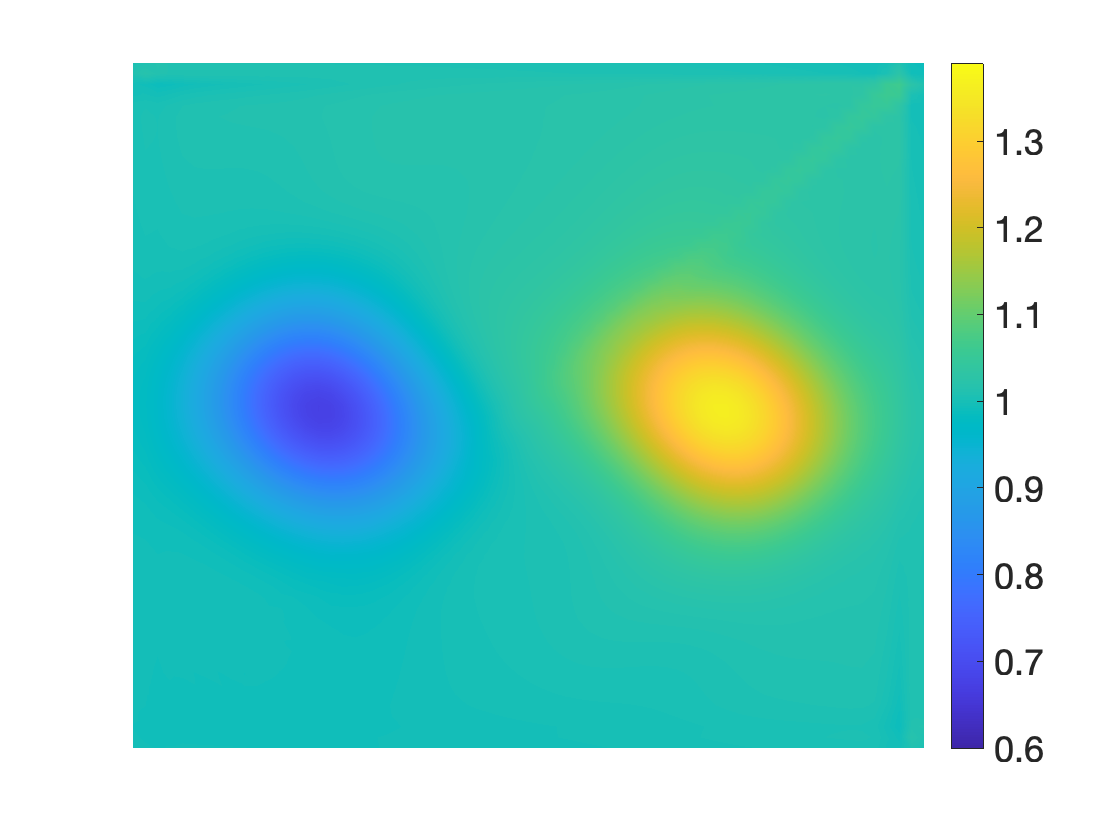} &
\includegraphics[width=0.32\textwidth]{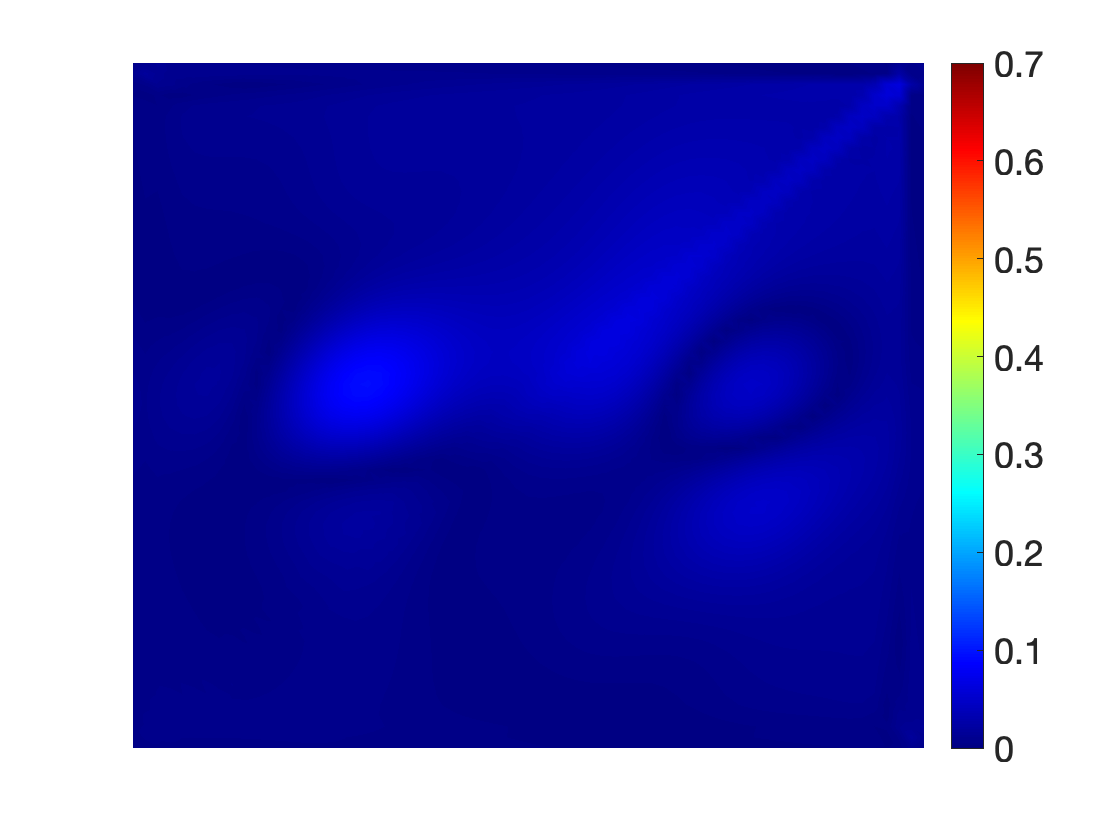}\\
\includegraphics[width=0.32\textwidth]{pics/diri1ex.png} &
\includegraphics[width=0.32\textwidth]{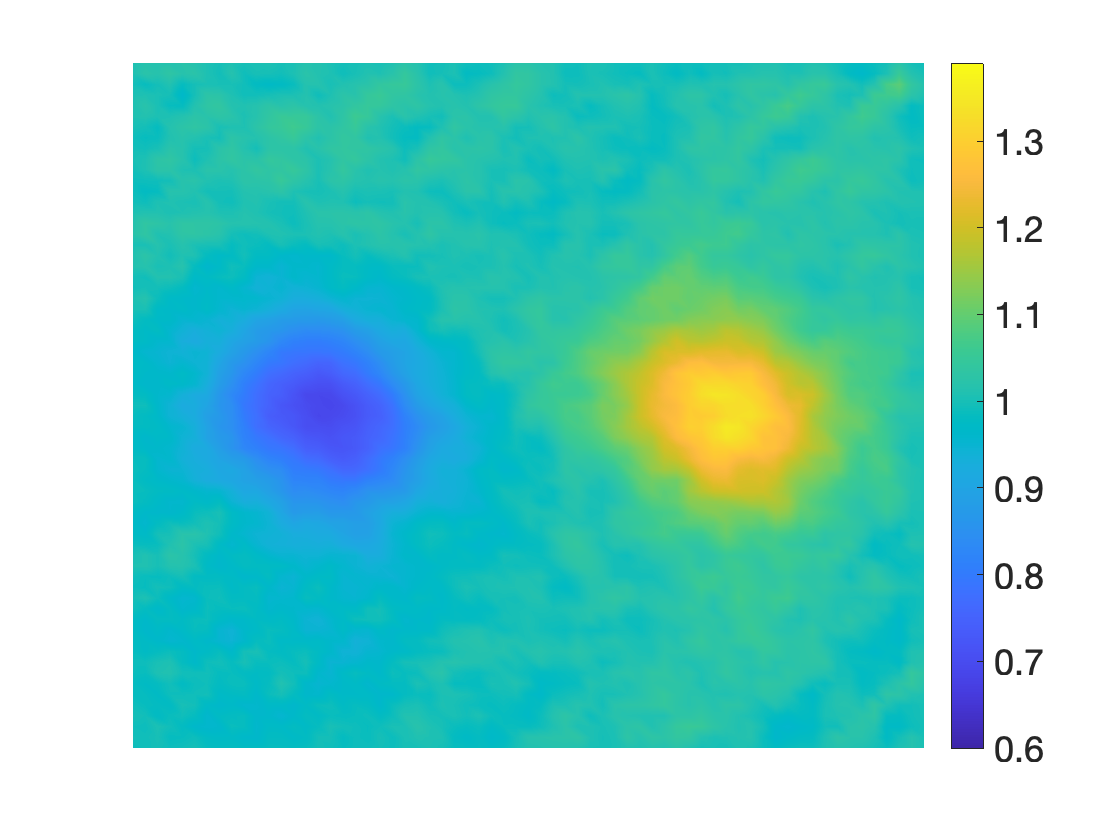} &
\includegraphics[width=0.32\textwidth]{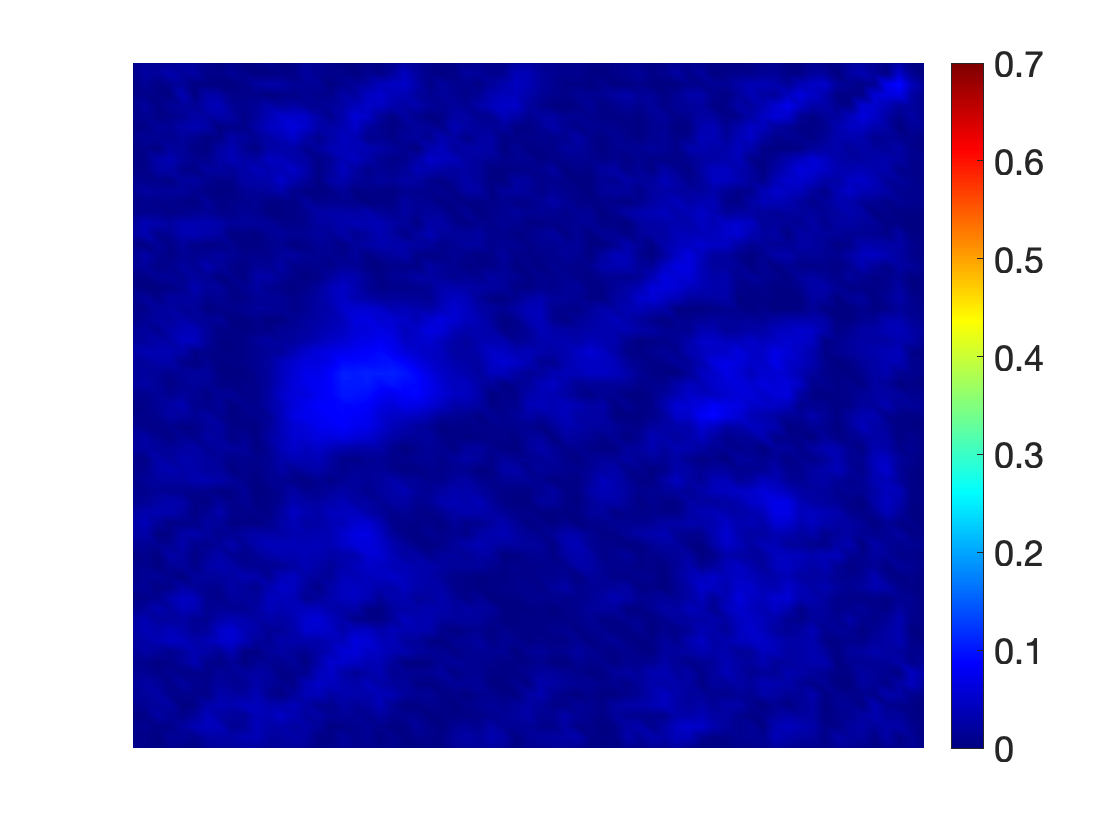}\\
\includegraphics[width=0.32\textwidth]{pics/diri1ex.png} &
\includegraphics[width=0.32\textwidth]{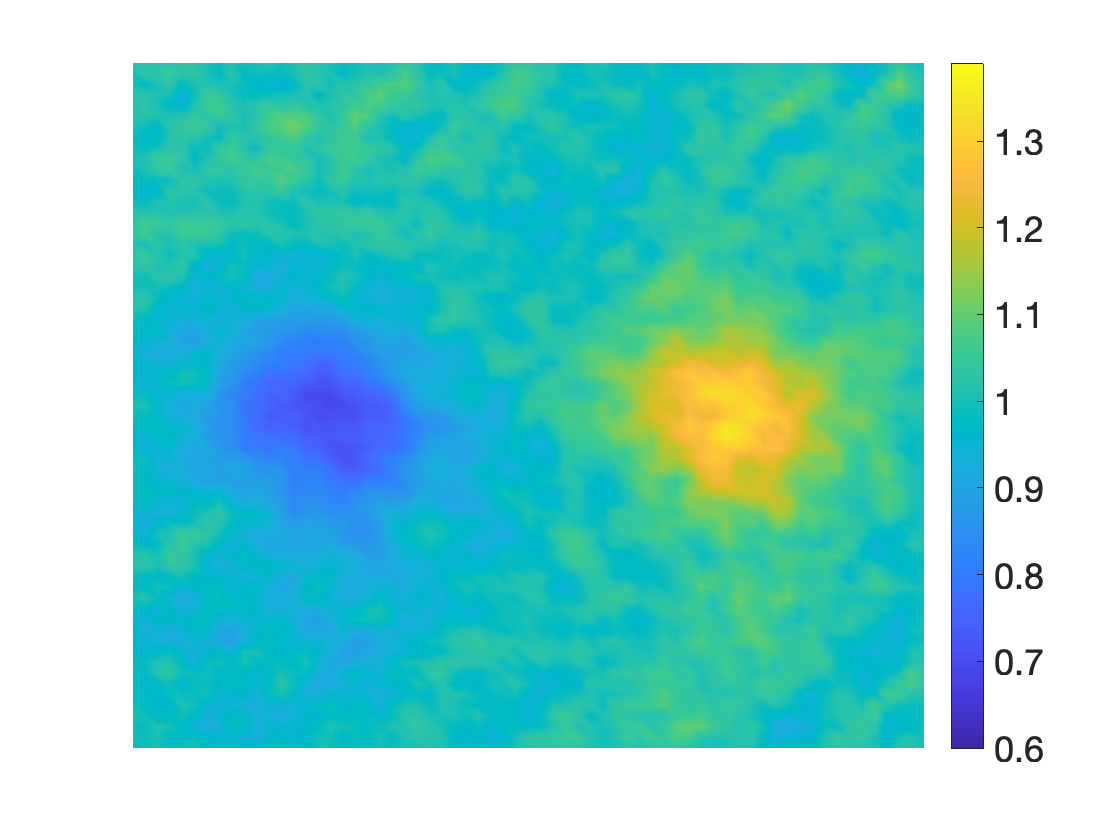} &
\includegraphics[width=0.32\textwidth]{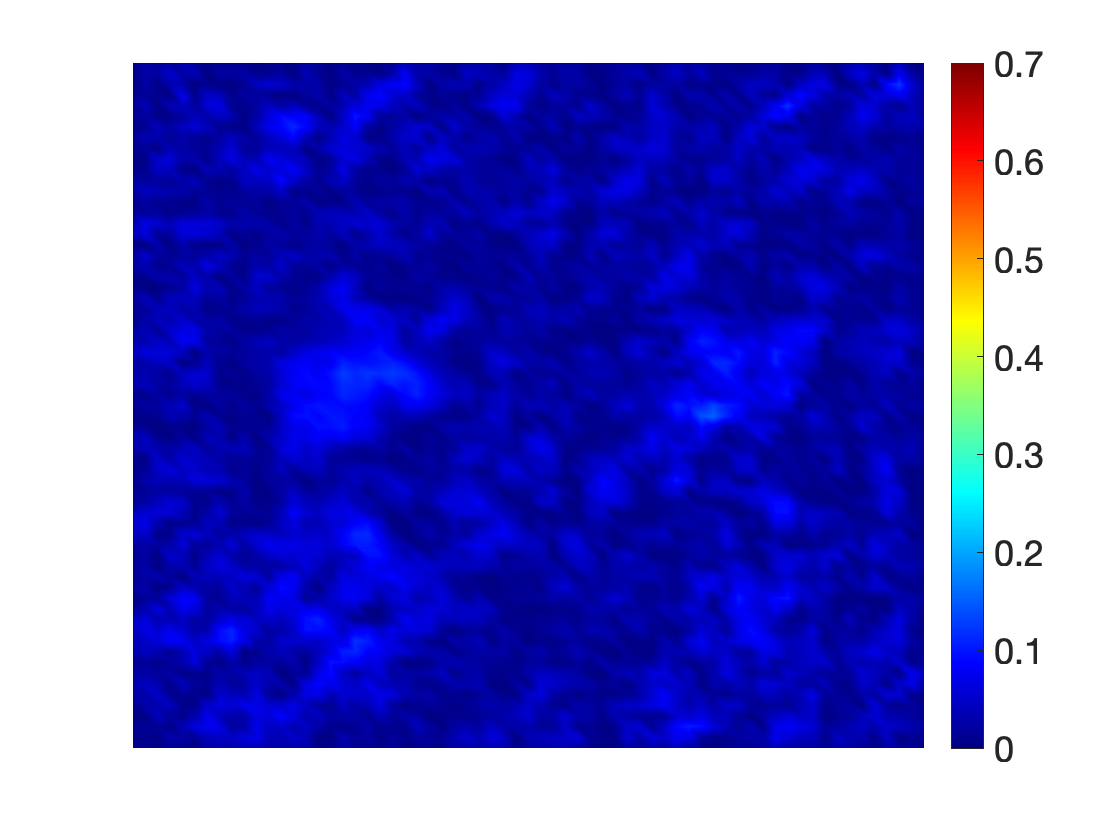}\\
(a) $q^\dag$  & (b) $\hat q$ & (c) $|\hat q-q^\dag|$
\end{tabular}
\caption{The reconstructions for Example \ref{exam:diri1} using FEM  with exact data $($top$)$ and noisy data $(\delta=10\%,\ 20\%$, middle, bottom$)$. }
\label{fig:diri1fem}
\end{figure}

\begin{example}\label{exam:diri1}
The domain $\Omega=(-1,1)^2$, exact conductivity $q^\dag = 1 +s_1(x_1, x_2) + s_2(x_1, x_2)$,
with $s_1=0.4e^{-15(x_1-0.5)^2-15x_2^2}$ and
$s_2 = -0.4e^{-15(x_1+0.5)^2-15x_2^2}$, and $u^\dag=x_1+x_2+\frac{1}{3}(x_1^3+x_2^3)$.
\end{example}

Figs. \ref{fig:diri1}-\ref{fig:diri1fem} show the reconstructions for the exact and noisy data by the DNN approach and classic FEM, respectively. The relative $L^2(\Omega)$-errors of the DNN approximations are 9.72e-3 and 3.78e-2 for exact and noisy data, which are comparable to that by the FEM approximation (2.32e-2 and 2.66e-2).
Like before, we observe excellent robustness of the DNN approach with respect to the presence of noise - when $\delta=10\%$, there is almost no degradation in the reconstruction
quality; when $\delta=20\%$, the bump features are generally well resolved but the magnitudes of the bumps are under-estimated, whereas the FEM reconstruction contains large oscillations in the background. Tables \ref{diri1table1}-\ref{diri1table3} show the impact of the regularization parameter $\gamma_q$, DNN architecture and sampling points on the reconstruction quality. The table shows that these parameter do not affect much the reconstruction error as long as they fall within a suitable range.

\begin{table}[htp!]
  \centering
  \caption{The variation of the relative $L^2(\Omega)$ error $e(\hat q)$ with respect to various algorithmic parameters. }
\begin{threeparttable}
\subfigure[$e$ v.s. $\gamma_q$ and $\delta$ using DNN method\label{diri1table1}]{\begin{tabular}{c|ccc}
\toprule
  $\gamma_q\backslash\delta$&   0\% &  1\%& 10\%\\
\midrule
     1.00e-2 &  2.58e-2 & 2.65e-2 & 4.71e-2 \\
     1.00e-3 &  1.19e-2 & 1.01e-2 & 4.48e-2\\
     1.00e-4 &  1.13e-2 & 8.22e-3 & 4.51e-2\\
     1.00e-5 &  9.72e-3 & 8.13e-3 & 3.78e-2\\
\bottomrule
\end{tabular}} \\
\subfigure[$e$ v.s. $W$ and  $L$\label{diri1table2}]{
\begin{tabular}{c|ccc}
\toprule
${W}\backslash L$&   5 &  10& 20\\
\midrule
     4 & 6.72e-2 & 8.79e-2 & 3.78e-1 \\
     12& 1.94e-2 & 7.87e-3 & 1.68e-2\\
     26& 7.76e-3 & 8.20e-3 & 6.59e-3\\
     40& 5.86e-3 & 6.11e-3 & 5.90e-3\\
\bottomrule
\end{tabular}}\quad
\subfigure[$e$ v.s. $n_r$ and $n_b$\label{diri1table3}]{
\begin{tabular}{c|cccc}
\toprule
$n_b\backslash n_r$&   5000 &  10000& 20000&40000\\
\midrule
     500&  1.17e-2 &1.08e-2 &1.06e-2 &1.31e-2 \\
     1000& 8.66e-3 &9.09e-3 &1.16e-2 &1.16e-2\\
     4000& 1.01e-2 &1.18e-2 &1.03e-2 &9.72e-3\\

\bottomrule
\end{tabular}}
\end{threeparttable}
\end{table}

\begin{figure}[htbp]
\centering
\setlength{\tabcolsep}{0em}
\begin{tabular}{ccc}
\includegraphics[width=0.32\textwidth]{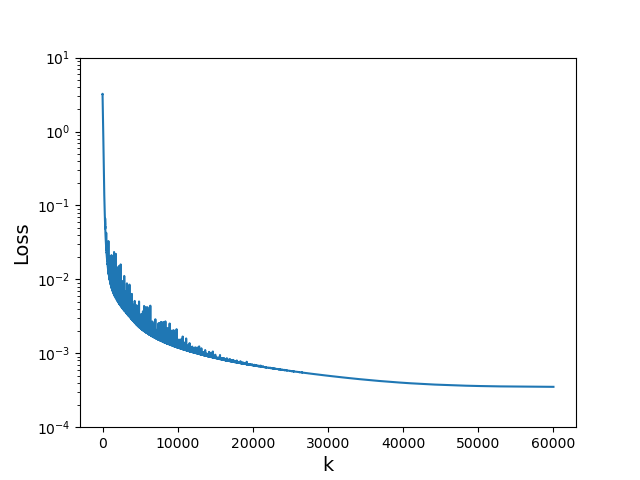} &
\includegraphics[width=0.32\textwidth]{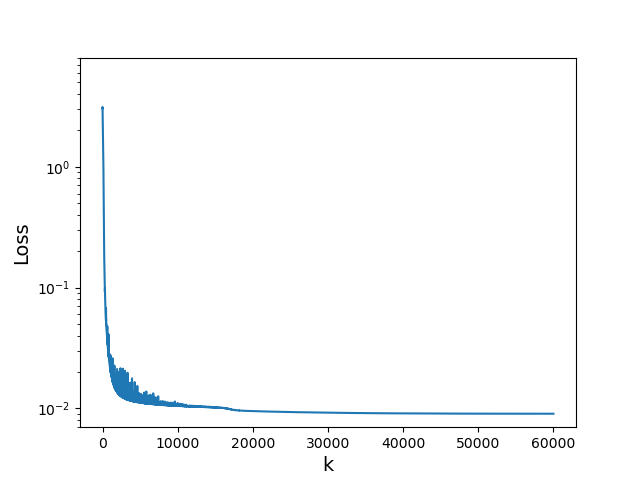} &
\includegraphics[width=0.32\textwidth]{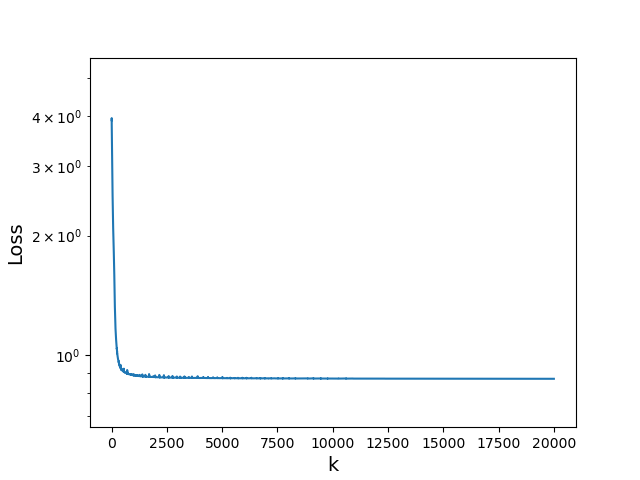}\\
\includegraphics[width=0.32\textwidth]{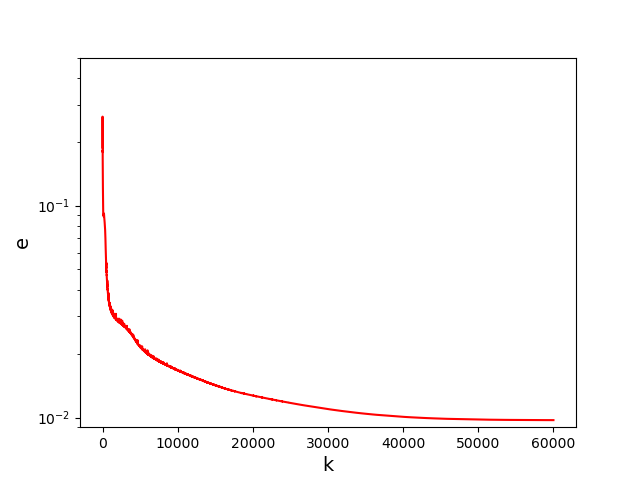} &
\includegraphics[width=0.32\textwidth]{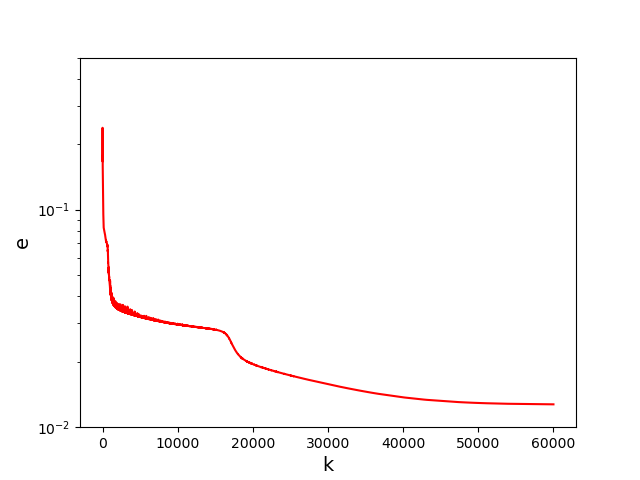} &
\includegraphics[width=0.32\textwidth]{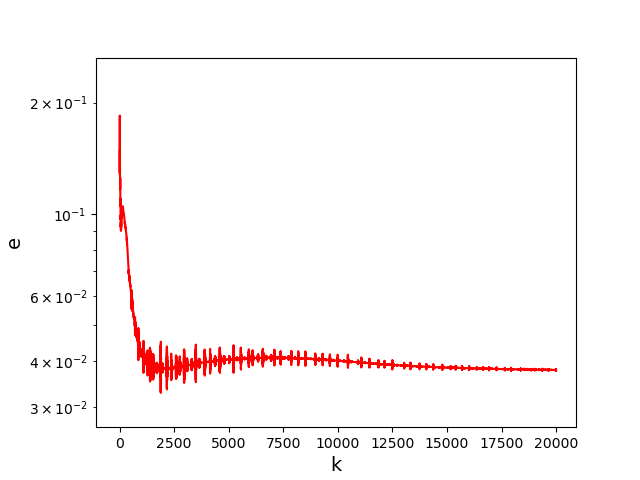}\\
(a) $\delta=0\%$  & (b) $\delta=1\%$ & (c) $\delta=10\%$
\end{tabular}
\caption{The variation of the loss $($top$)$ and relative $L^2(\Omega)$ error $e(\hat q)$ $($bottom$)$ during the training process for Example \ref{exam:diri1} at three different noise levels.}
\label{fig:diri1losse}
\end{figure}

\begin{figure}[htbp!]
\centering
\setlength{\tabcolsep}{0em}
\begin{tabular}{ccc}
\includegraphics[width=0.32\textwidth]{diri1ex.png} &
\includegraphics[width=0.32\textwidth]{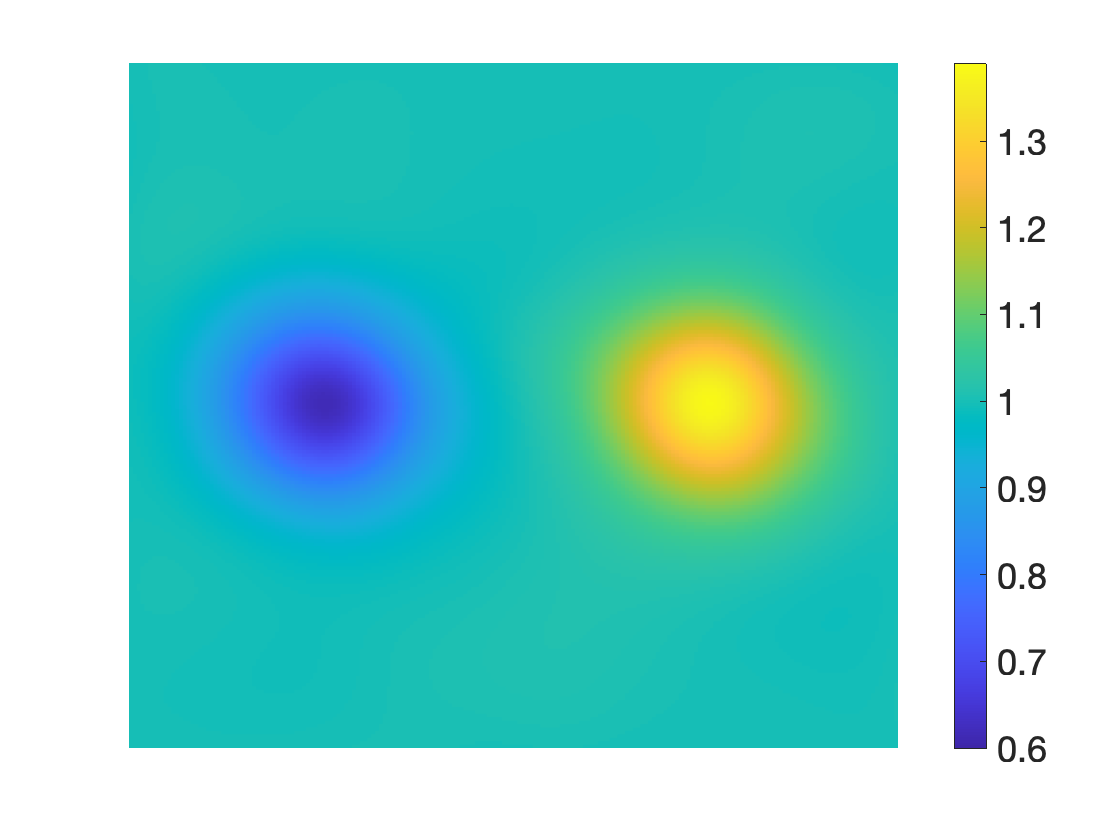} &
\includegraphics[width=0.32\textwidth]{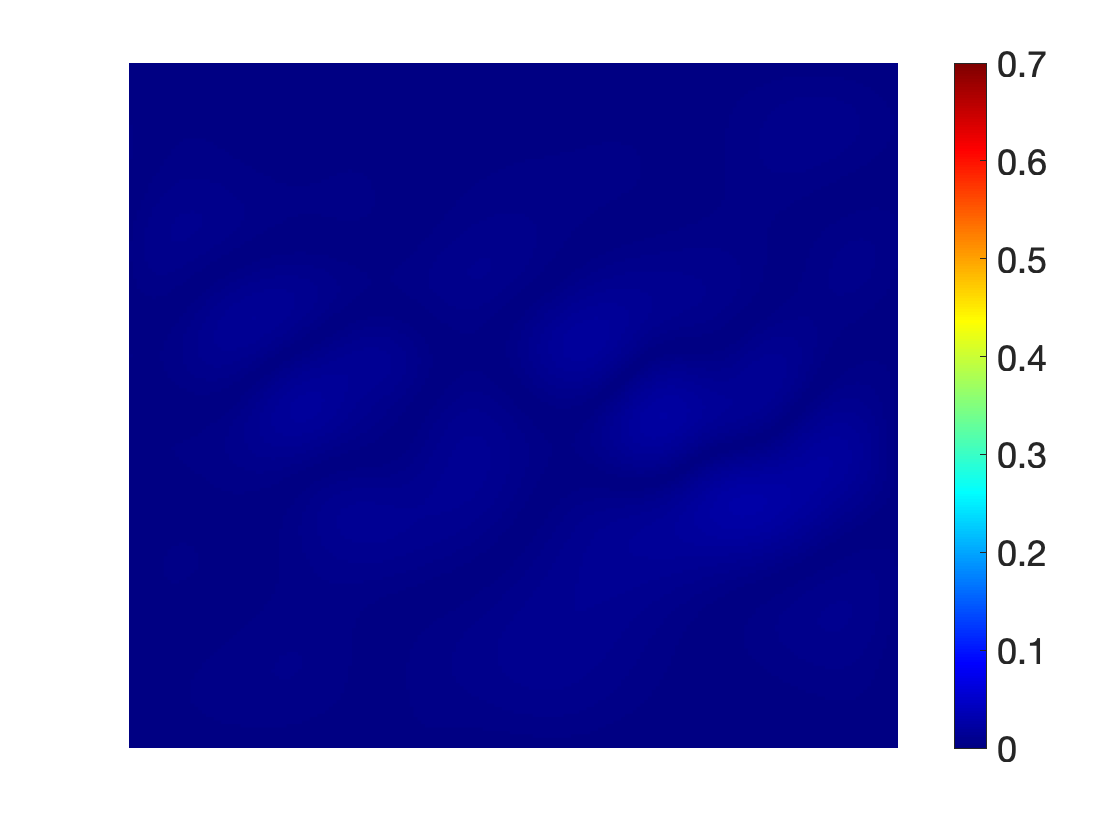}\\
\includegraphics[width=0.32\textwidth]{diri1ex.png} &
\includegraphics[width=0.32\textwidth]{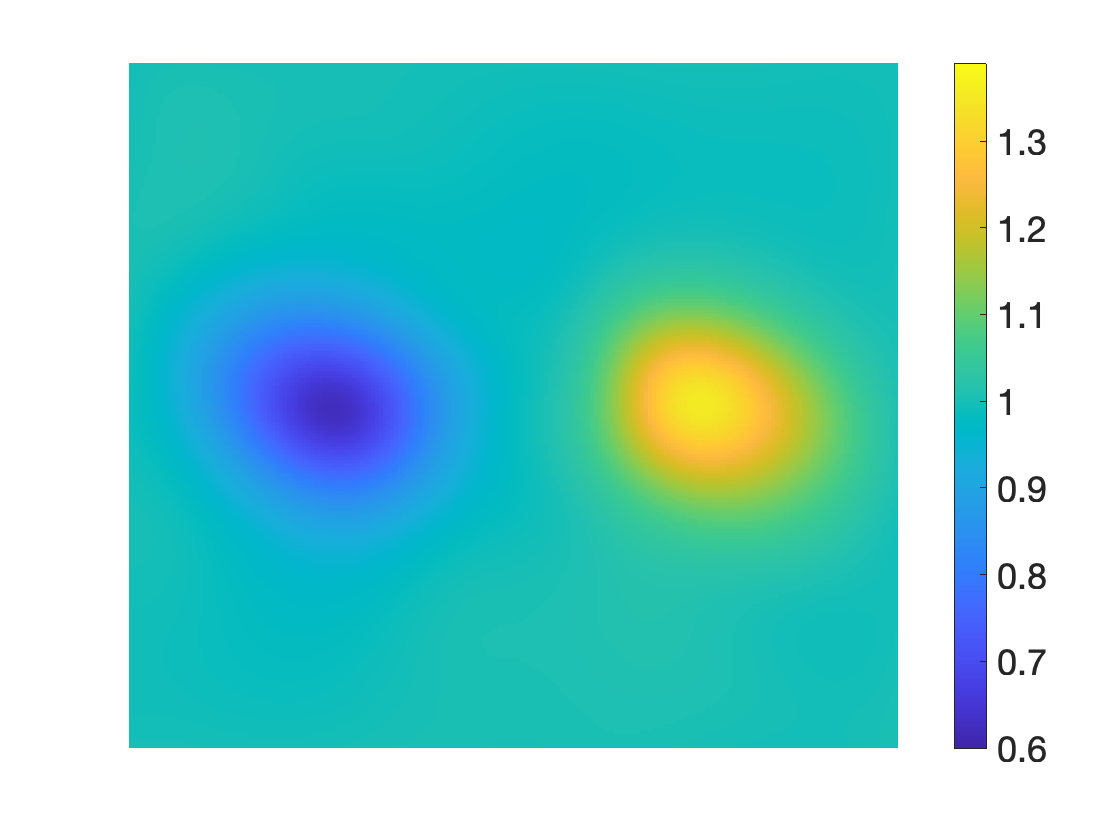} &
\includegraphics[width=0.32\textwidth]{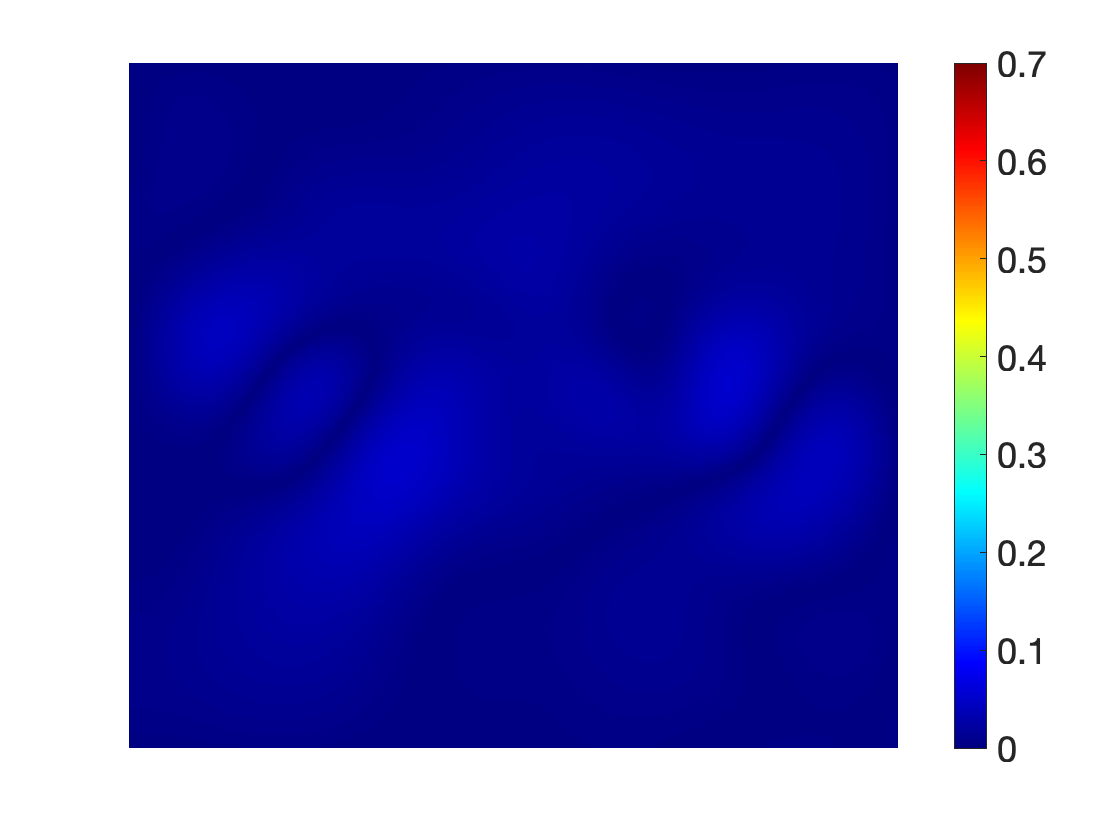}\\
\includegraphics[width=0.32\textwidth]{diri1ex.png} &
\includegraphics[width=0.32\textwidth]{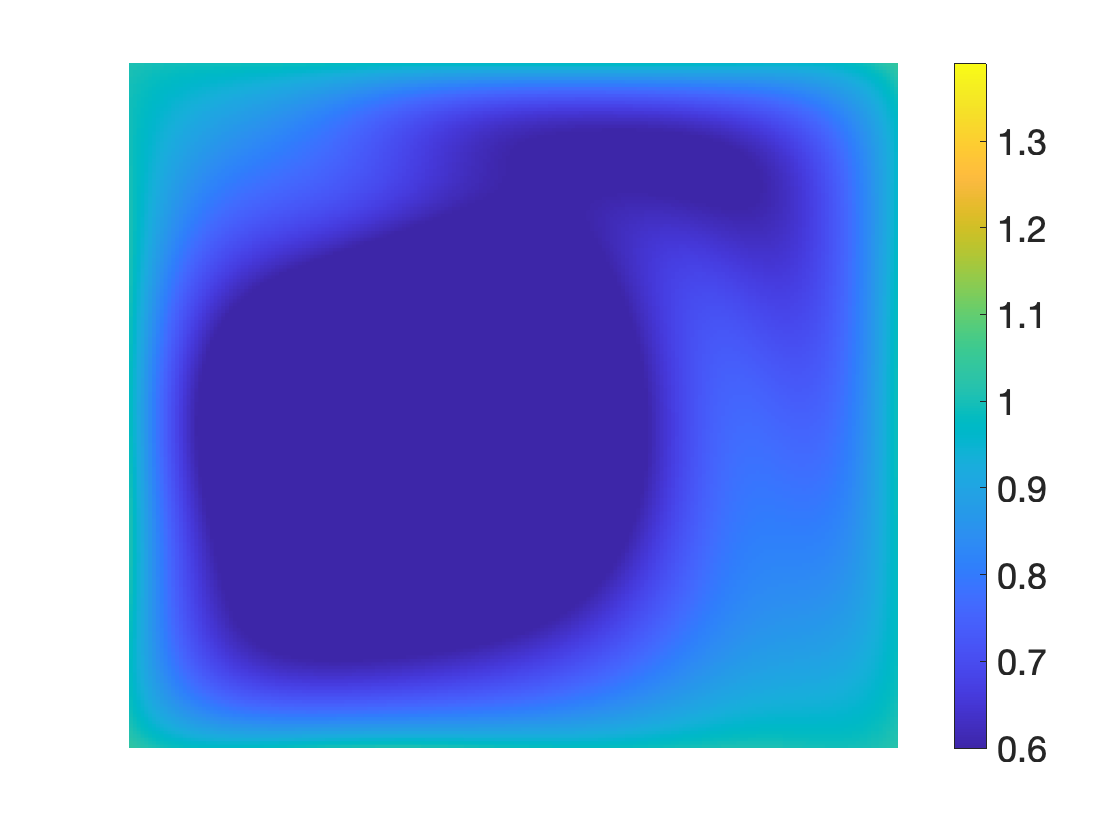} &
\includegraphics[width=0.32\textwidth]{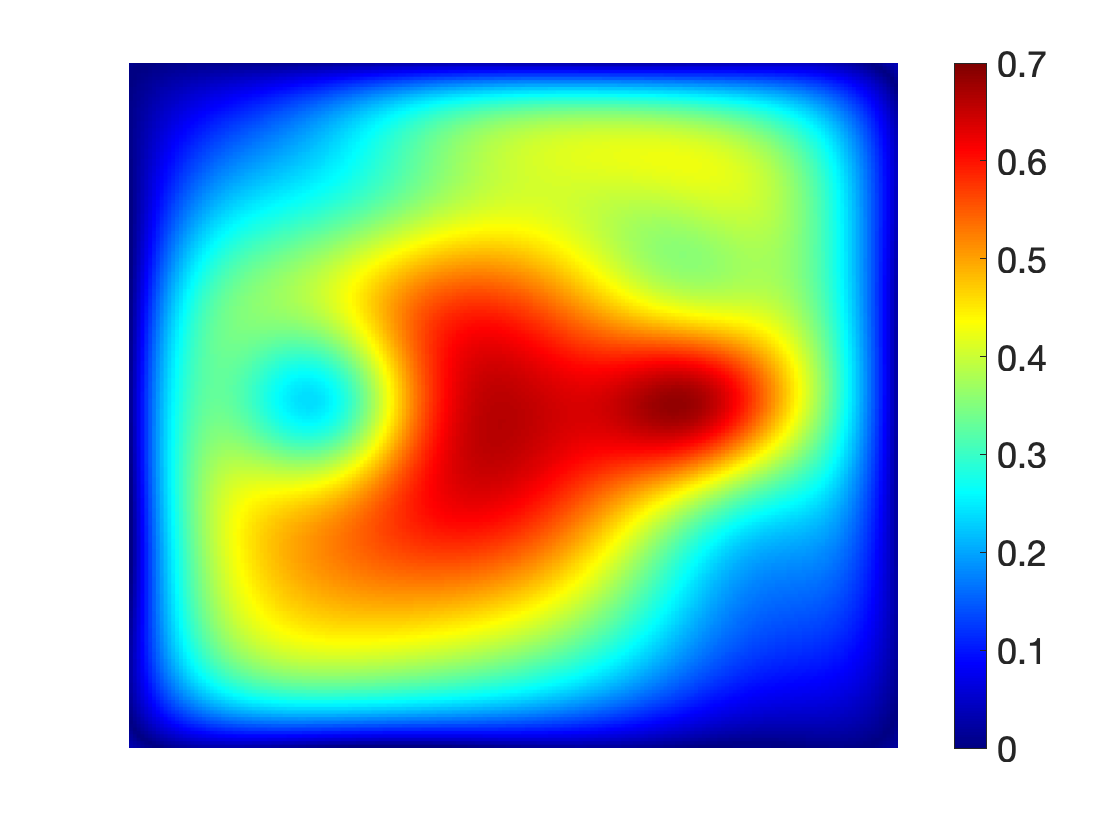}\\
(a) $q^\dag$  & (b) $\hat q$ & (c) $|\hat q-q^\dag|$
\end{tabular}
\caption{The reconstructions for Example \ref{exam:diri1} using loss \eqref{eqn:diriloss1} with exact data $($top$)$ and noisy data $($$\delta=1\%$, middle and $\delta=10\%$, bottom$)$.}
\label{fig:diri1sb}
\end{figure}

The convergence behavior of the optimizer in Fig. \ref{fig:diri1losse} shows that
the loss $J_{\bsgamma}$ and  error $e(\hat q)$ stagnate at a comparable level regardless of
the noise level. In light of Remark \ref{rmk:loss}, Fig. \ref{fig:diri1sb} shows the results of minimizing the loss \eqref{eqn:diriloss1} with exact data (top), $1\%$ noise (middle) and $10\%$ noise (bottom). The accuracy of reconstructions with exact data and $1\%$ noise is satisfactory, but the DNN fails to learn accurately using the loss \eqref{eqn:diriloss1} when the data $z^\delta$ becomes more noisy. In view of these experimental evidences, in practice one should prefer the loss \eqref{eqn:obj-Diri1} over the loss \eqref{eqn:diriloss1}.

The second example is about recovering a nearly piecewise constant conductivity.
\begin{example}\label{exam:diridisctn}
The domain $\Omega=(0,1)^2$, $q^\dag=1+0.3/(1+\exp(400((x_1-0.65)^2+2(x_2-0.7)^2-0.15^2)))$, and $u^\dag=\sin(\pi x_1)\sin(\pi x_2)$.
\end{example}

\begin{figure}[htbp]
\centering
\setlength{\tabcolsep}{0em}
\begin{tabular}{ccc}
\includegraphics[width=0.33\textwidth]{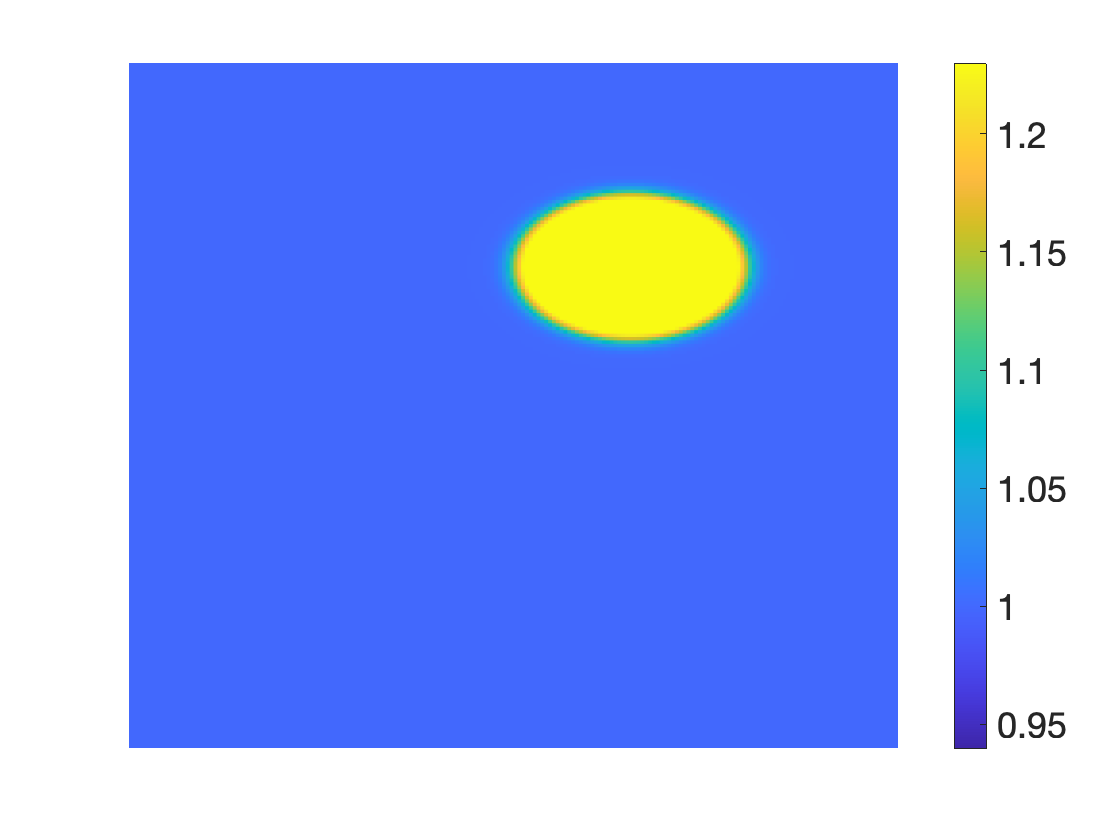} &
\includegraphics[width=0.33\textwidth]{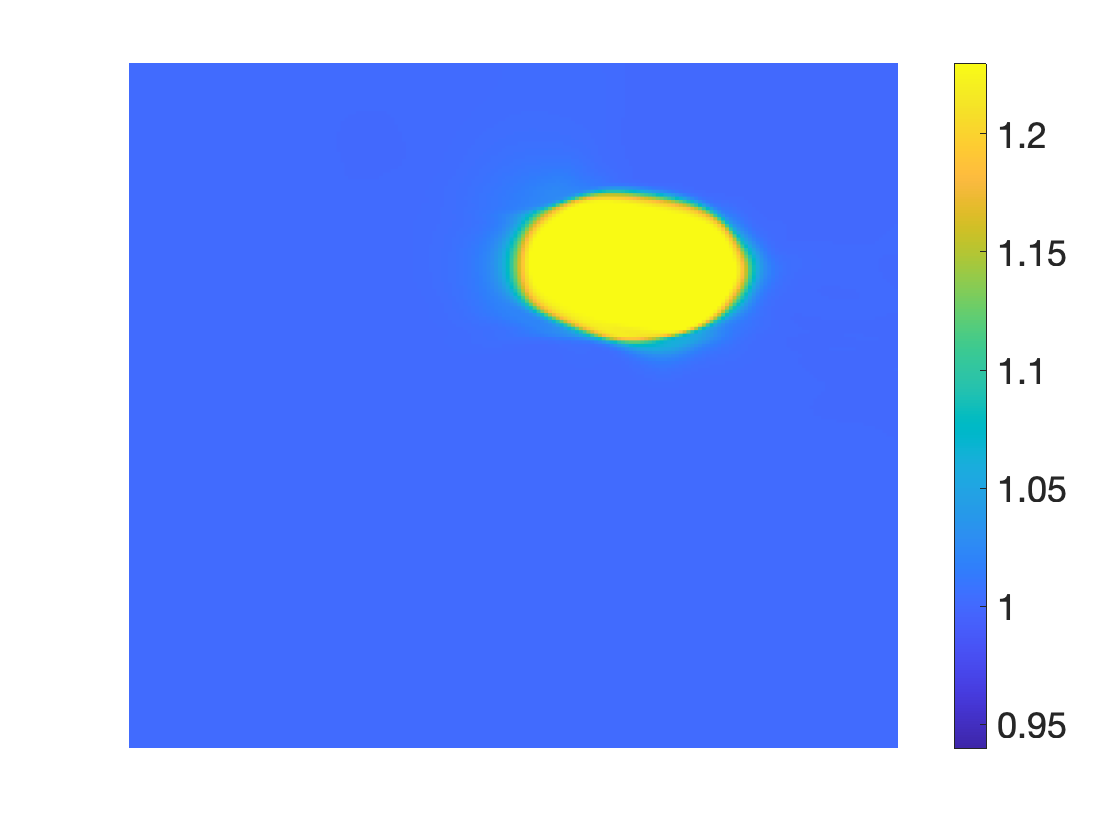} &
\includegraphics[width=0.33\textwidth]{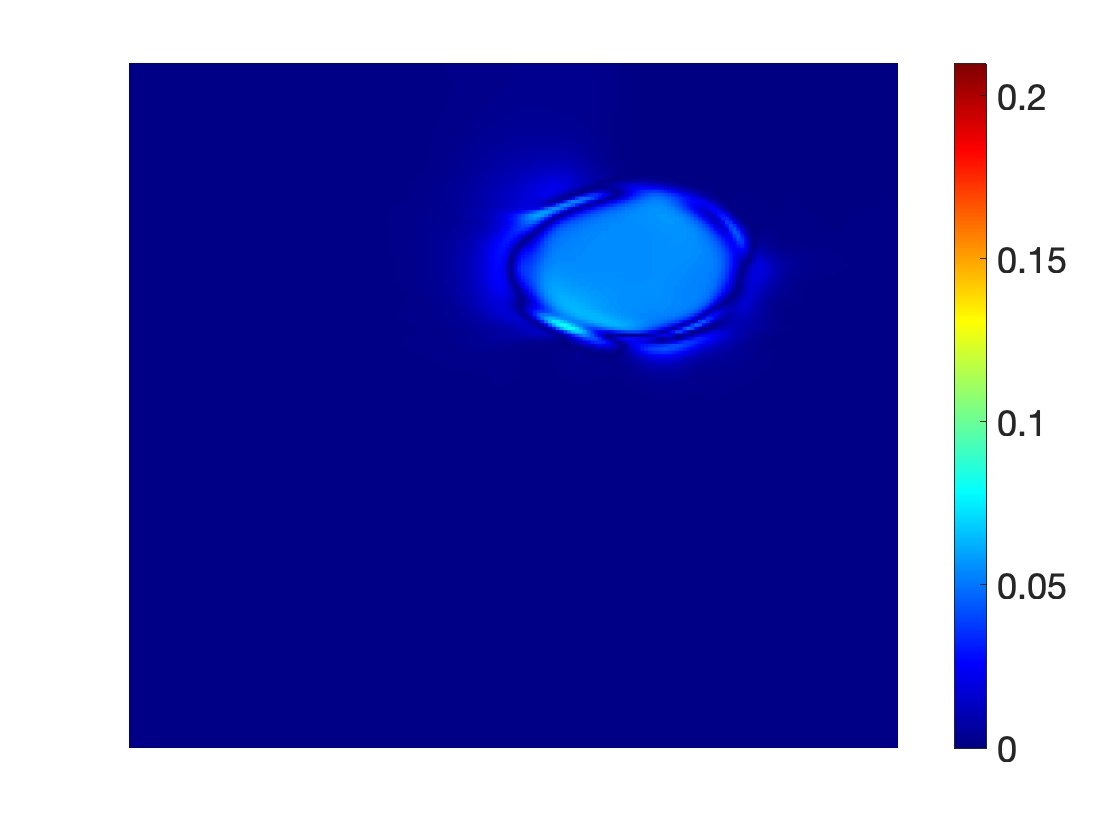}\\
\includegraphics[width=0.33\textwidth]{diridisctnex.png} &
\includegraphics[width=0.33\textwidth]{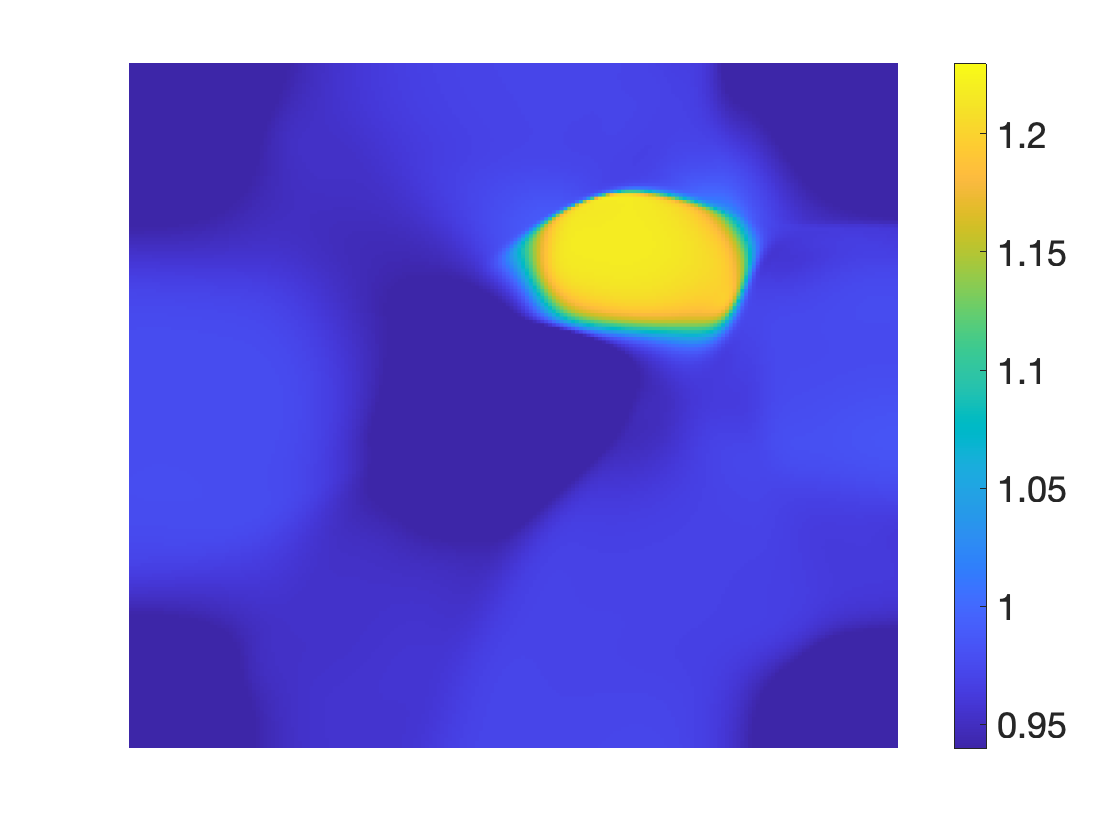} &
\includegraphics[width=0.33\textwidth]{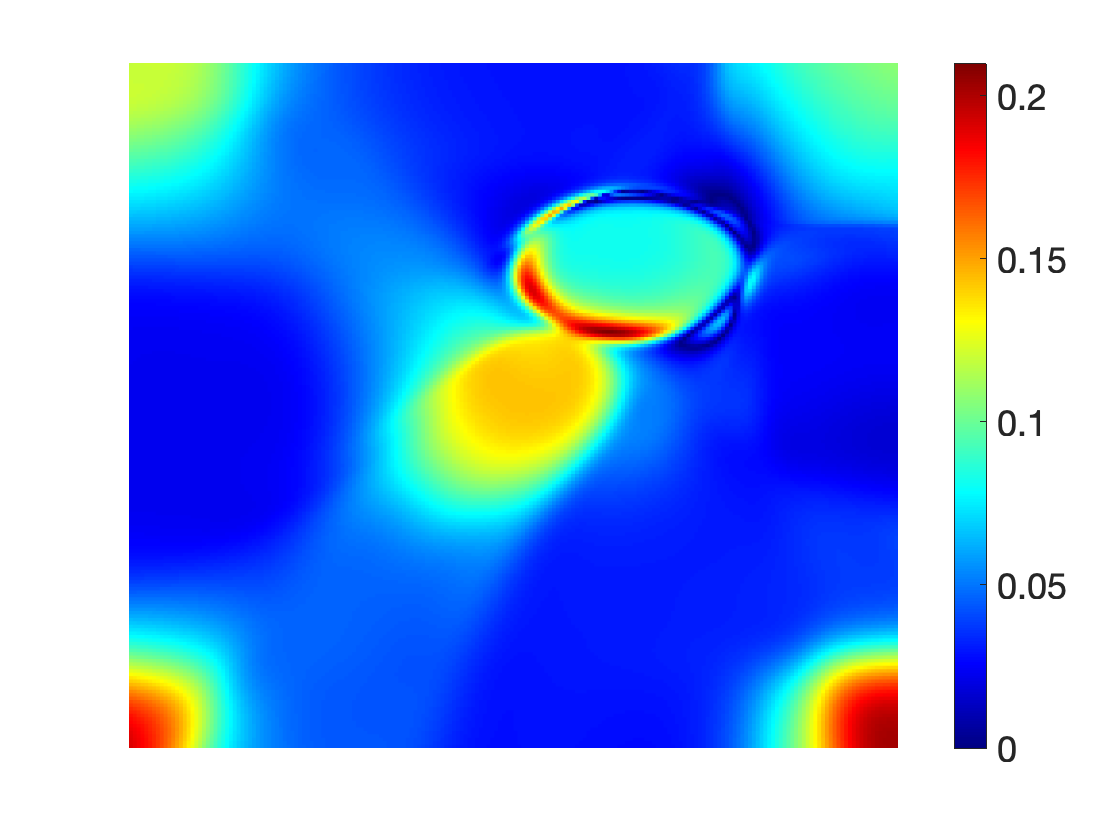}\\
(a) $q^\dag$  & (b) $\hat q$ & (c) $|\hat q-q^\dag|$
\end{tabular}
\caption{The reconstructions for Example \ref{exam:diridisctn} with exact data $($top$)$ and noisy data $(\delta=10\%$, bottom$)$.}
\label{fig:diridisctn}
\end{figure}

Note that in this example,
$\nabla u$ vanishes at the four corners of the domain $\Omega$ and the point $(\frac{1}{2},\frac{1}{2})$. Similar to Example \ref{exam:discon}, an additional total variation penalty,
i.e., $\gamma_{tv}| q|_{\rm TV}$ ($\gamma_{tv}=0.01$), is added to the loss $J_{\bsgamma}(\theta,\kappa)$ in order to promote piecewise constancy of the reconstruction. Fig. \ref{fig:diridisctn} shows the result using the loss
\eqref{eqn:obj-Diri1} with
exact (top) and noisy (bottom, $\delta=10\%$) data. The reconstruction is accurate for exact data, and in the presence of $10\%$ noise, its accuracy deteriorates only near the points where $\nabla u$ vanishes (see the reconstruction plot and the pointwise error plot). This is expected since substituting $\nabla u =0$ back into the formulation \eqref{eqn:obj-Diri1} leads to a loss nearly independent of $q_{\theta}$, which causes the inaccuracy in the reconstruction in the neighborhood of these points.

\begin{figure}[htbp]
\centering
\setlength{\tabcolsep}{0em}
\begin{tabular}{ccc}
\includegraphics[width=0.33\textwidth]{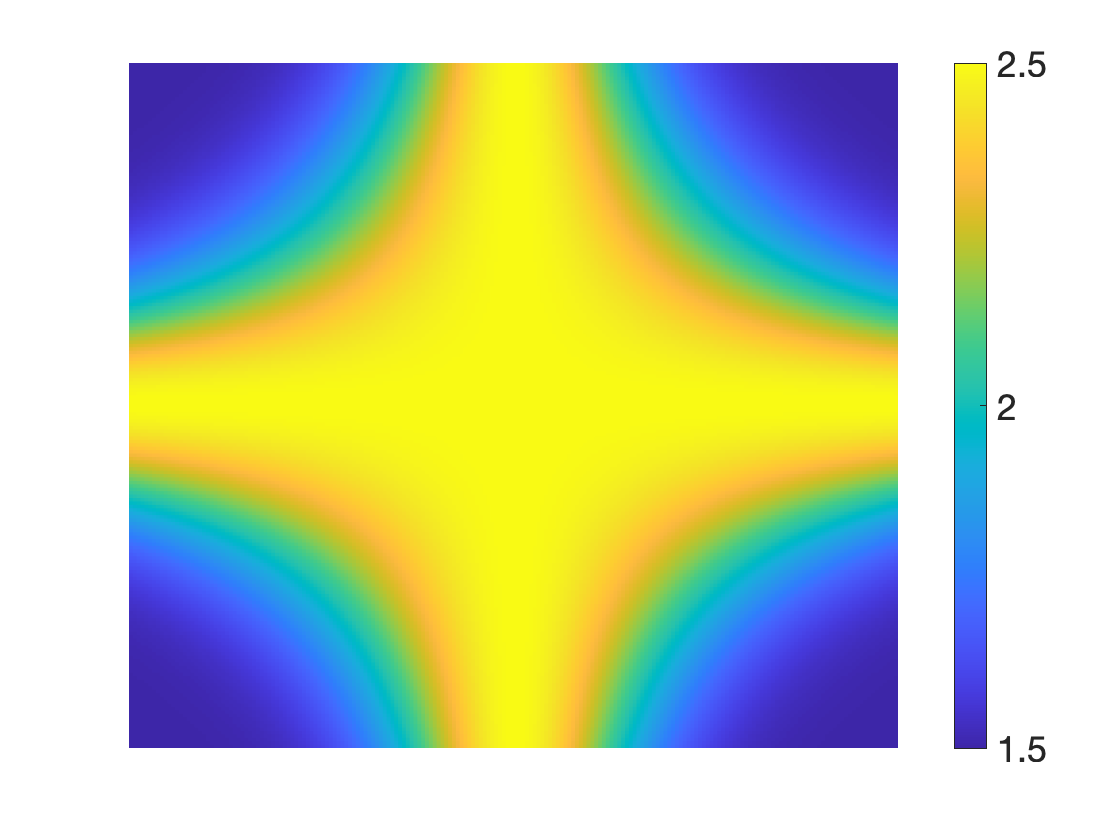} &
\includegraphics[width=0.33\textwidth]{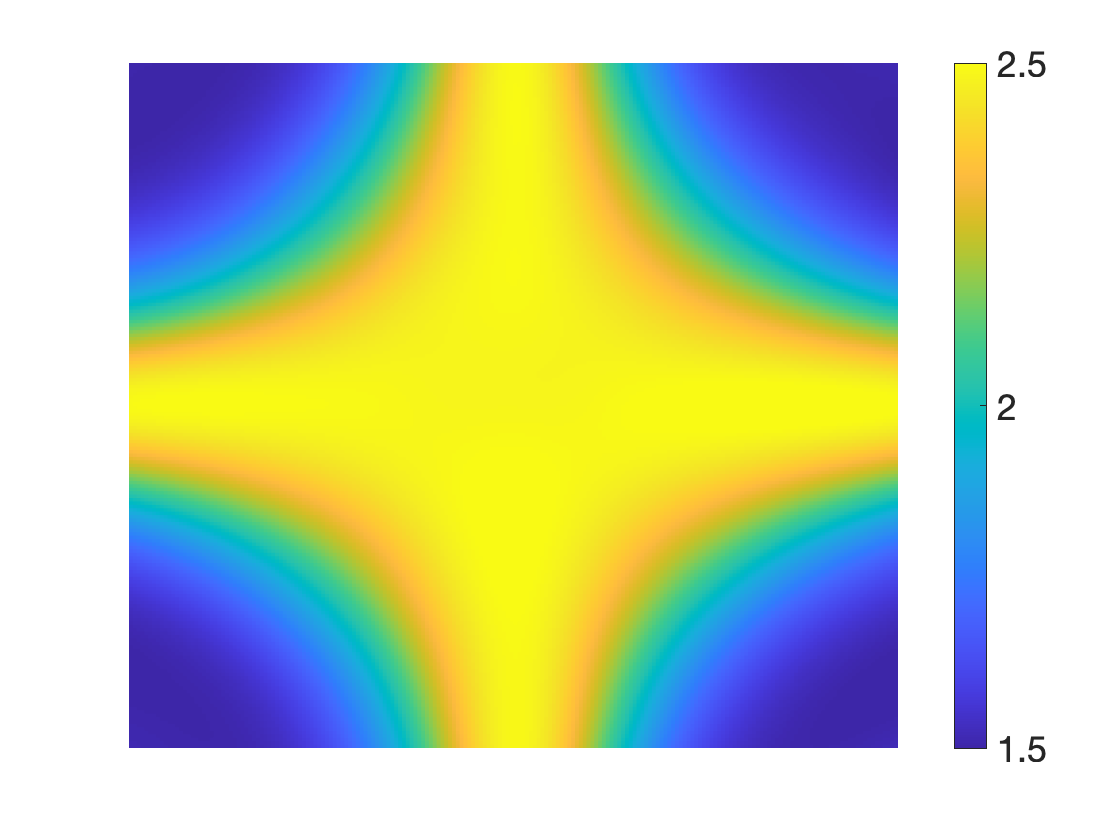} &
\includegraphics[width=0.33\textwidth]{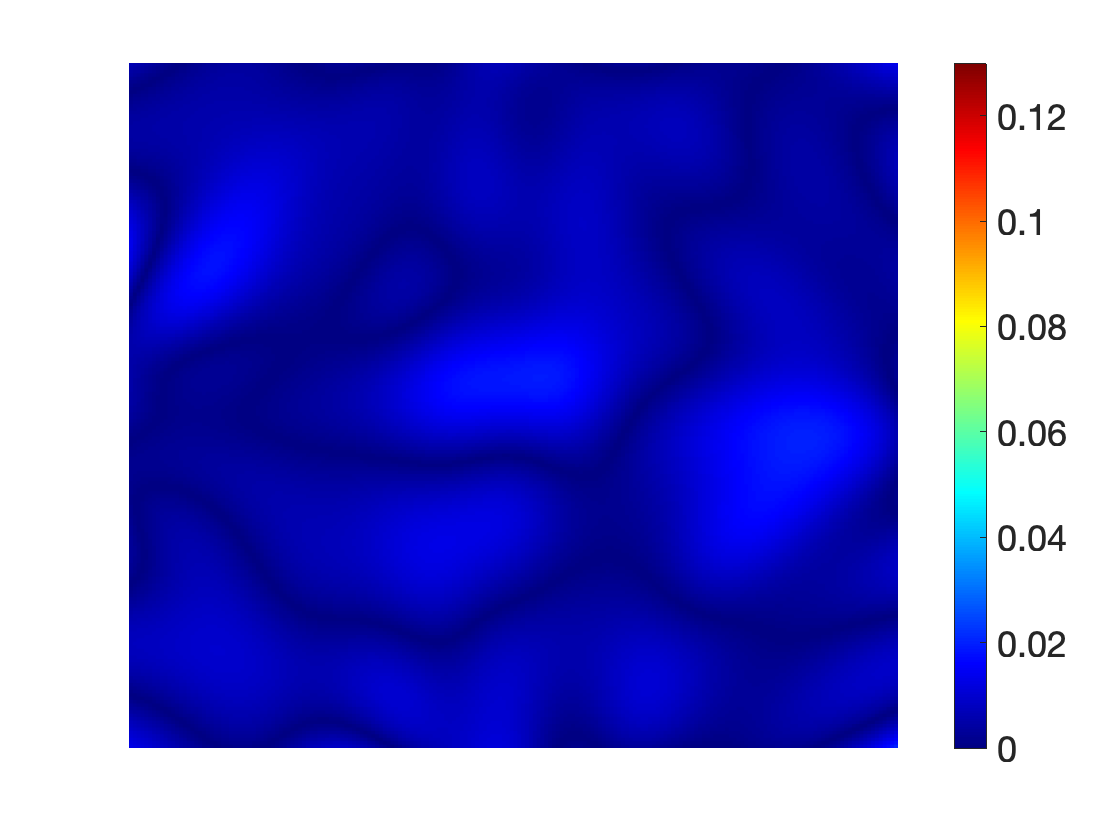}\\
\includegraphics[width=0.33\textwidth]{diri2ex.png} &
\includegraphics[width=0.33\textwidth]{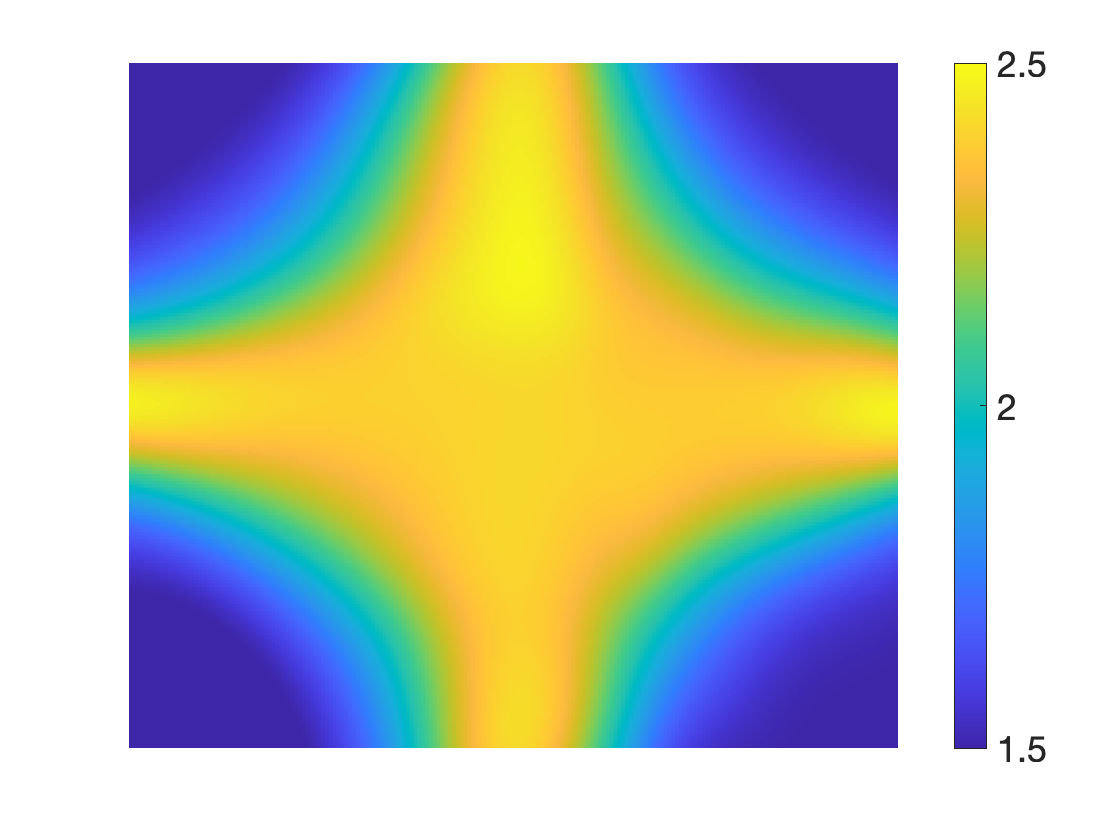} &
\includegraphics[width=0.33\textwidth]{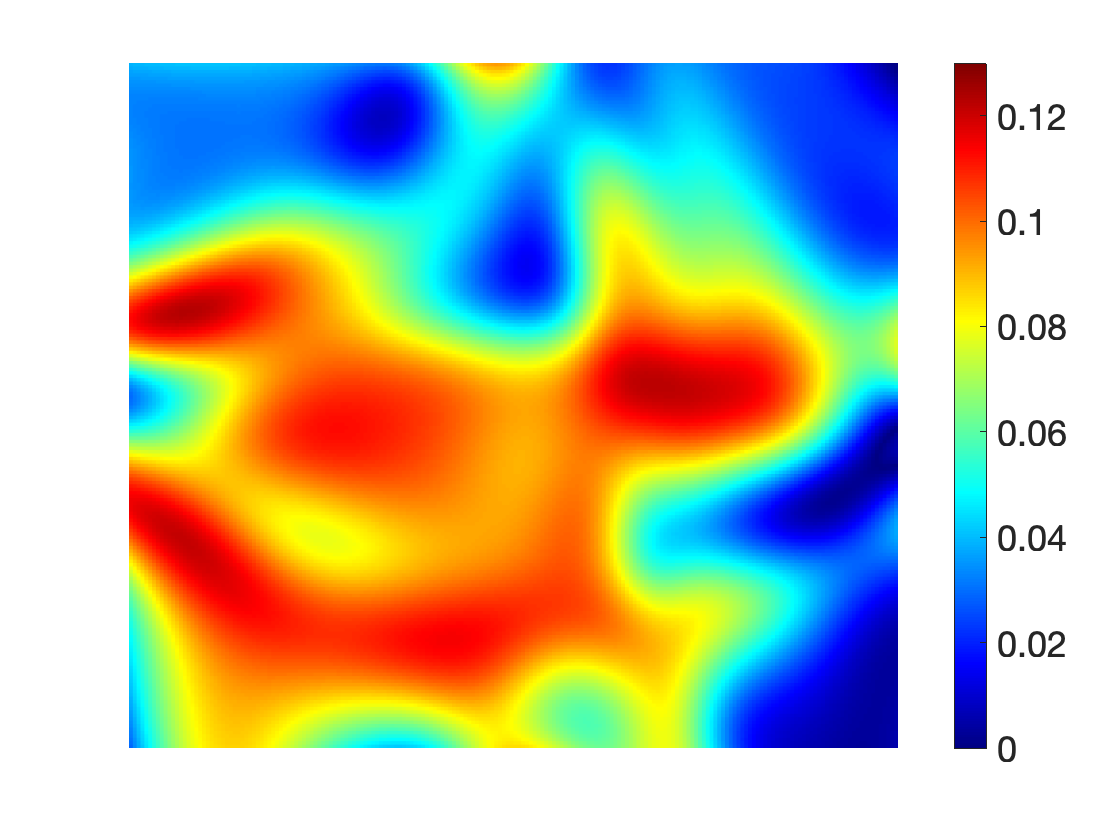}\\
(a) $q^\dag$  & (b) $\hat q$ & (c) $|\hat q-q^\dag|$
\end{tabular}
\caption{The reconstructions for Example \ref{exam:diri2} with exact data $($top$)$ and noisy data $(\delta=10\%$, bottom$)$.}
\label{fig:diri2}
\end{figure}

The third example is about recovering a 3D conductivity.
\begin{example}\label{exam:diri2}
The domain $\Omega=(0,1)^3$, $q^\dag=1+e^{-(12(x_1-0.5)(x_2-0.5))^2}+x_3$, and $u^\dag=\sum_{i=1}^3(x_i+\frac13x_i^3)$.
\end{example}

Fig. \ref{fig:diri2} show the reconstruction using the loss
\eqref{eqn:obj-Diri1} on a 2D cross section at $x_3=0.5$, with
exact (top) and noisy (bottom, $\delta=10\%$) data, with the relative $L^2(\Omega)$-error being 3.28e-3 and 3.48e-2, respectively, indicating the excellent robustness of the approach with respect to data noise: there
is only very mild deterioration in the reconstruction $\hat q$.

\begin{figure}[htbp]
\centering
\setlength{\tabcolsep}{0em}
\begin{tabular}{ccc}
\includegraphics[width=0.32\textwidth]{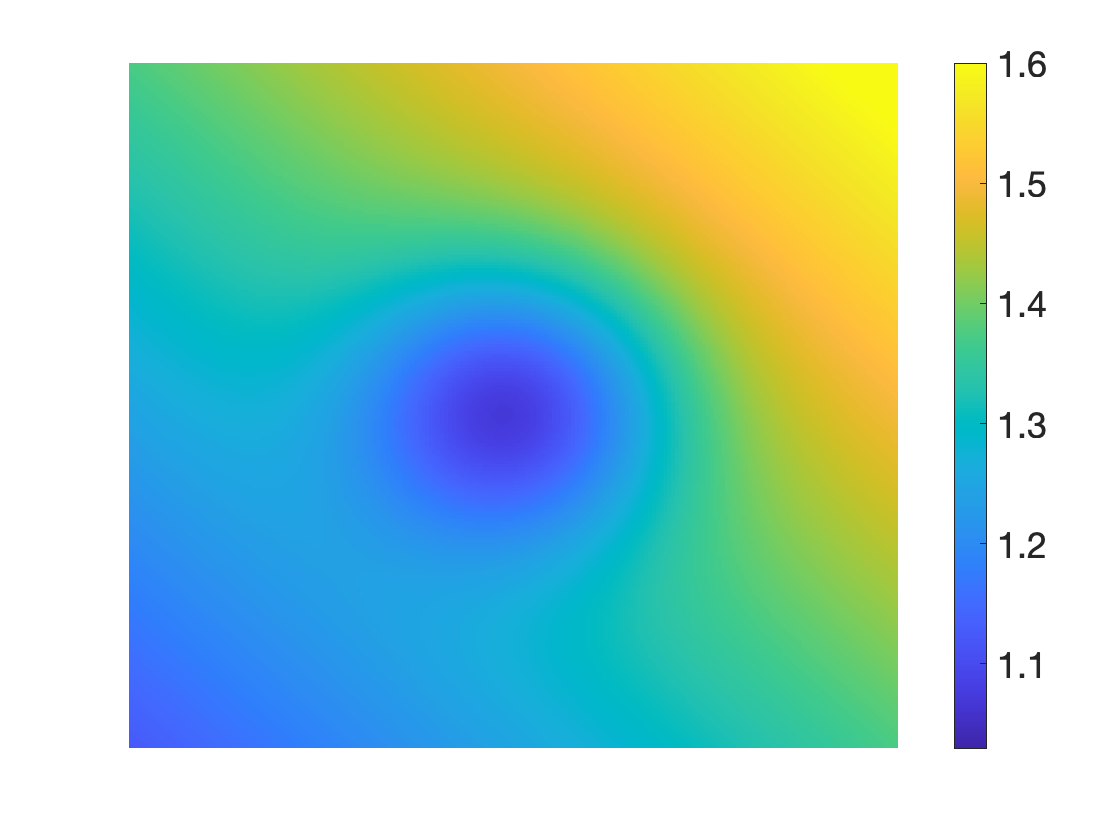} &
\includegraphics[width=0.32\textwidth]{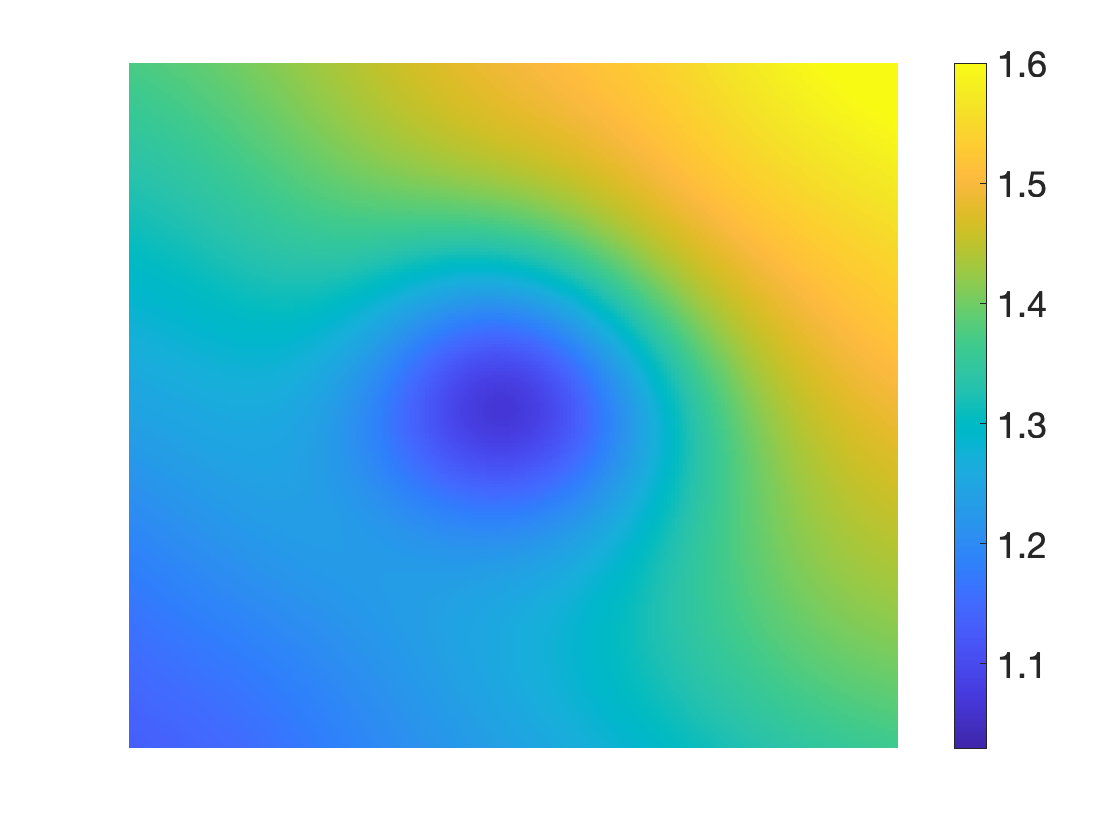} &
\includegraphics[width=0.32\textwidth]{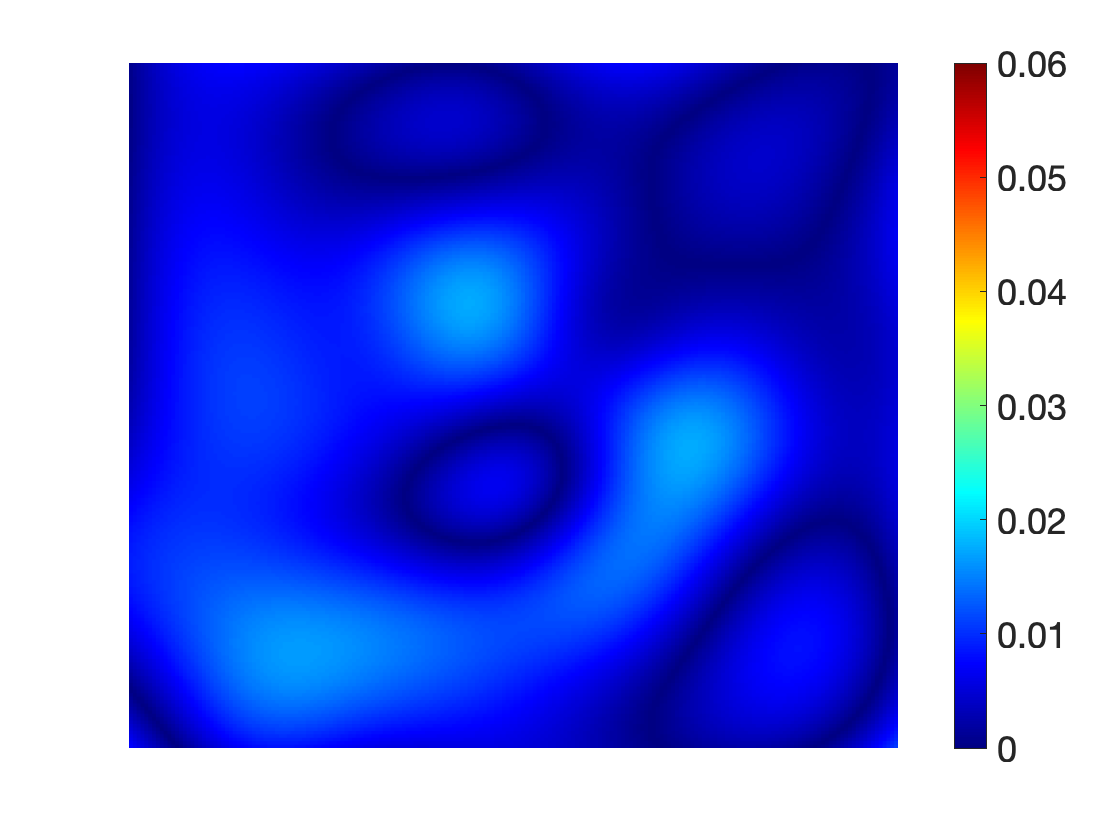}\\
\includegraphics[width=0.32\textwidth]{diridim5ex.png} &
\includegraphics[width=0.32\textwidth]{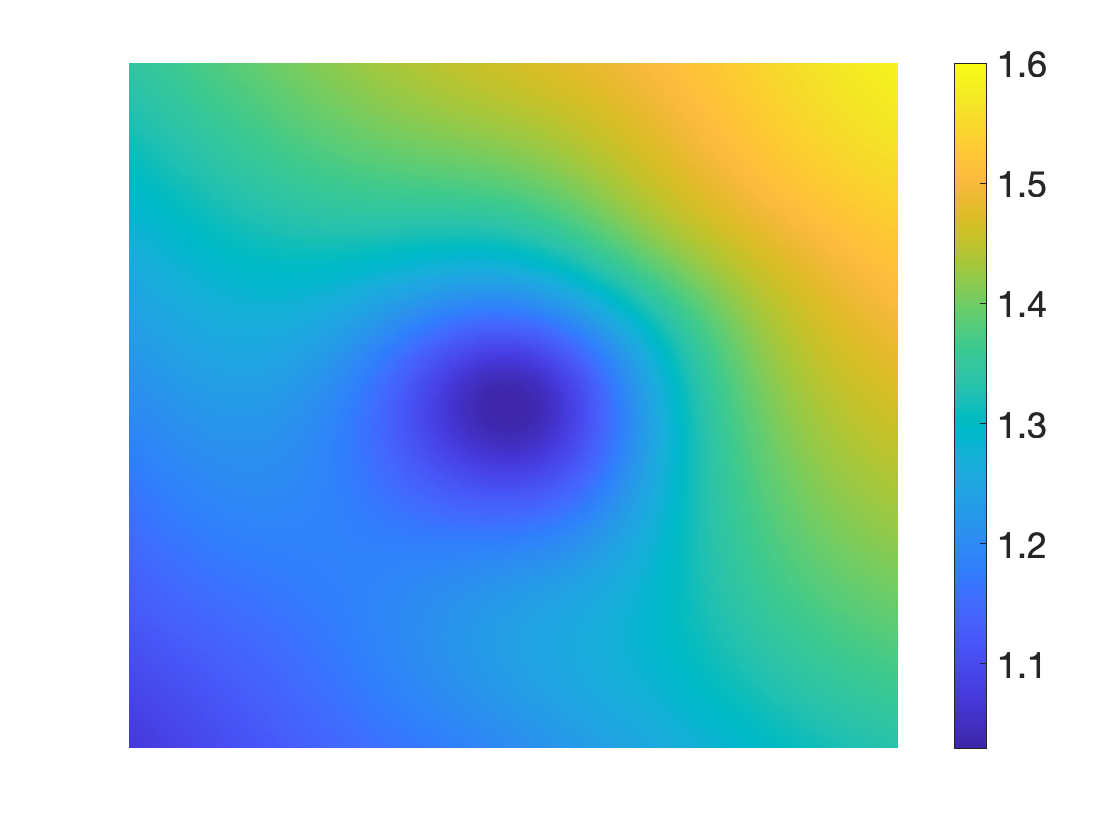} &
\includegraphics[width=0.32\textwidth]{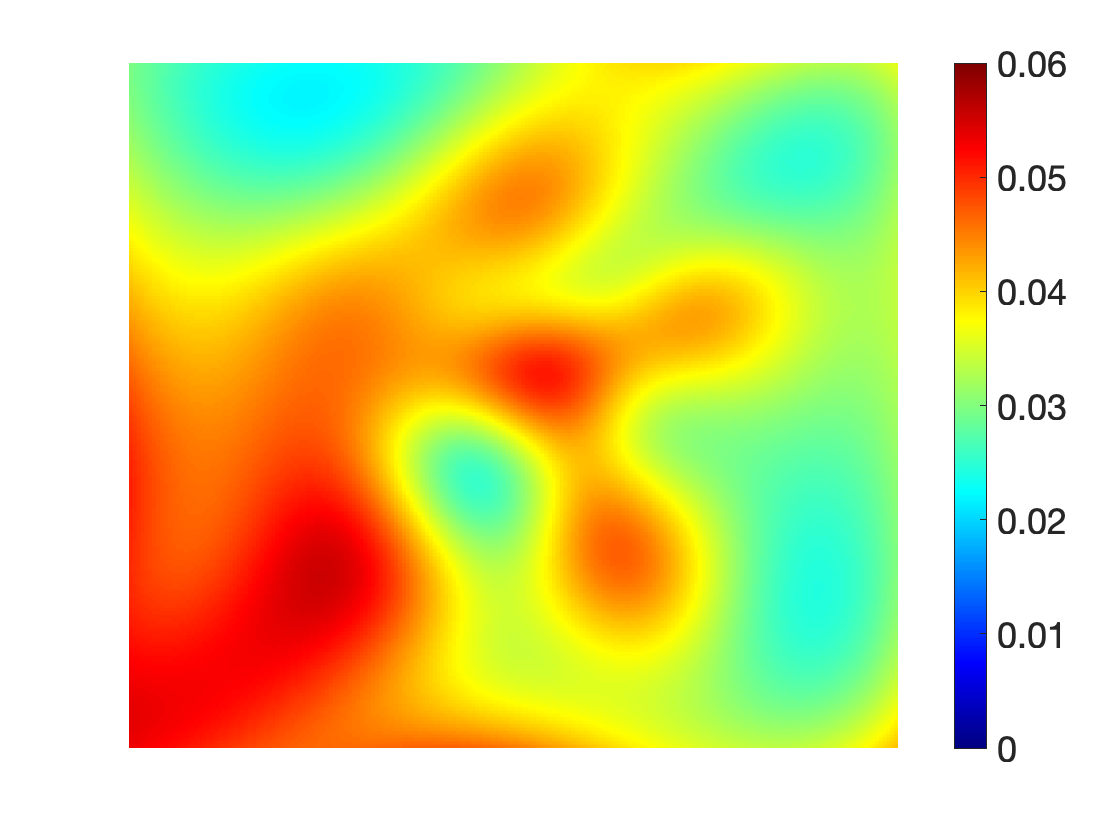}\\
(a) $q^\dag$  & (b) $\hat q$ & (c) $|\hat q-q^\dag|$
\end{tabular}
\caption{The reconstructions for Example \ref{exam:diridim5} with exact data $($top$)$ and noisy data $(\delta=10\%$, bottom$)$.}
\label{fig:diridim5}
\end{figure}

The fourth example is concerned with recovering a 5D conductivity function.
\begin{example}\label{exam:diridim5}
The domain $\Omega=(0,1)^5$, $q^\dag=1+0.5(x_1x_5+x_2x_4+x_3^2)-0.3e^{-25(x_1-0.5)^2-25(x_2-0.5)^2}$, and $u^\dag=\sum_{i=1}^5(x_i+\frac13x_i^3)$.
\end{example}

Fig. \ref{fig:diridim5} shows the reconstructions using the loss \eqref{eqn:obj-Diri1} on a 2D cross section with $x_3=x_4=x_5=0.5$, with exact (top) and noisy (bottom, $\delta=10\%$) data. The relative $L^2(\Omega)$-error is 5.78e-3 and 2.86e-2 for exact and noisy data, respectively. The features of the true conductivity $q^\dag$ have been successfully recovered and visually there is almost no difference between the results for exact and noisy data. This again shows the potential of the DNN approach for solving high-dimensional inverse problems.

Next we provide comparative experiments in the Dirichlet case.
\begin{example}\label{exam:diripartial2d}
The domain $\Omega=(0,1)^2$, the measurement $\nabla z^\delta$ on the region $\omega=\Omega\setminus (0.2,0.8)^2$, $q^\dag=2+0.5\sin(2\pi x_1)\sin(2\pi x_2)$, and $u^\dag=x_1+x_2+\frac{1}{3}(x_1^3+x_2^3)$.
\end{example}

\begin{figure}[htbp]
\centering
\setlength{\tabcolsep}{0em}
\begin{tabular}{ccc}
\includegraphics[width=0.32\textwidth]{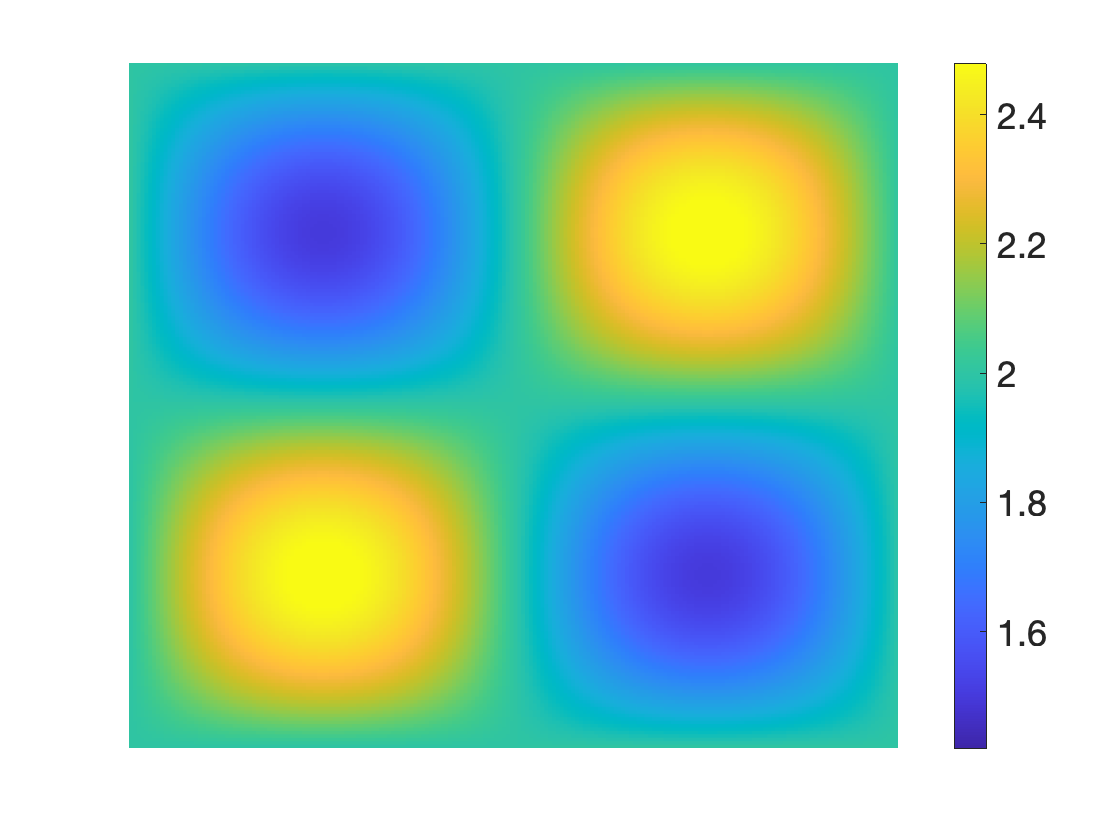} &
\includegraphics[width=0.32\textwidth]{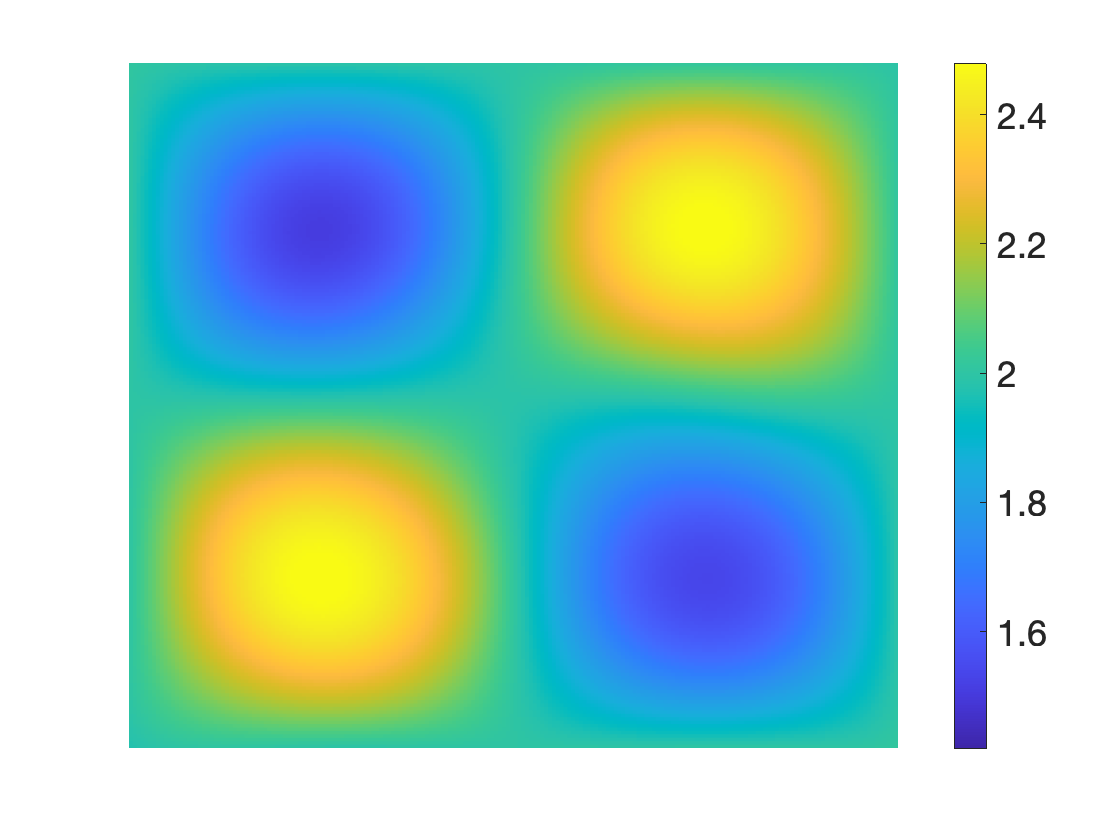} &
\includegraphics[width=0.32\textwidth]{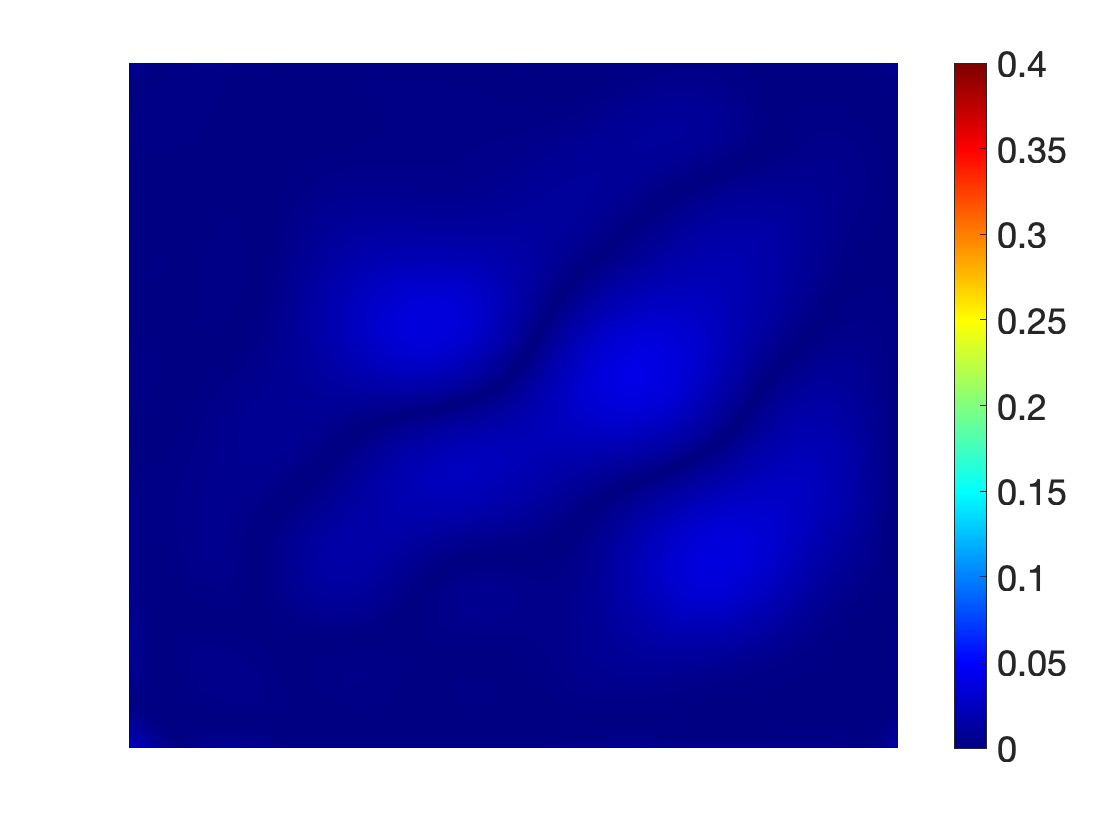}\\
\includegraphics[width=0.32\textwidth]{diripartial2dex.png} &
\includegraphics[width=0.32\textwidth]{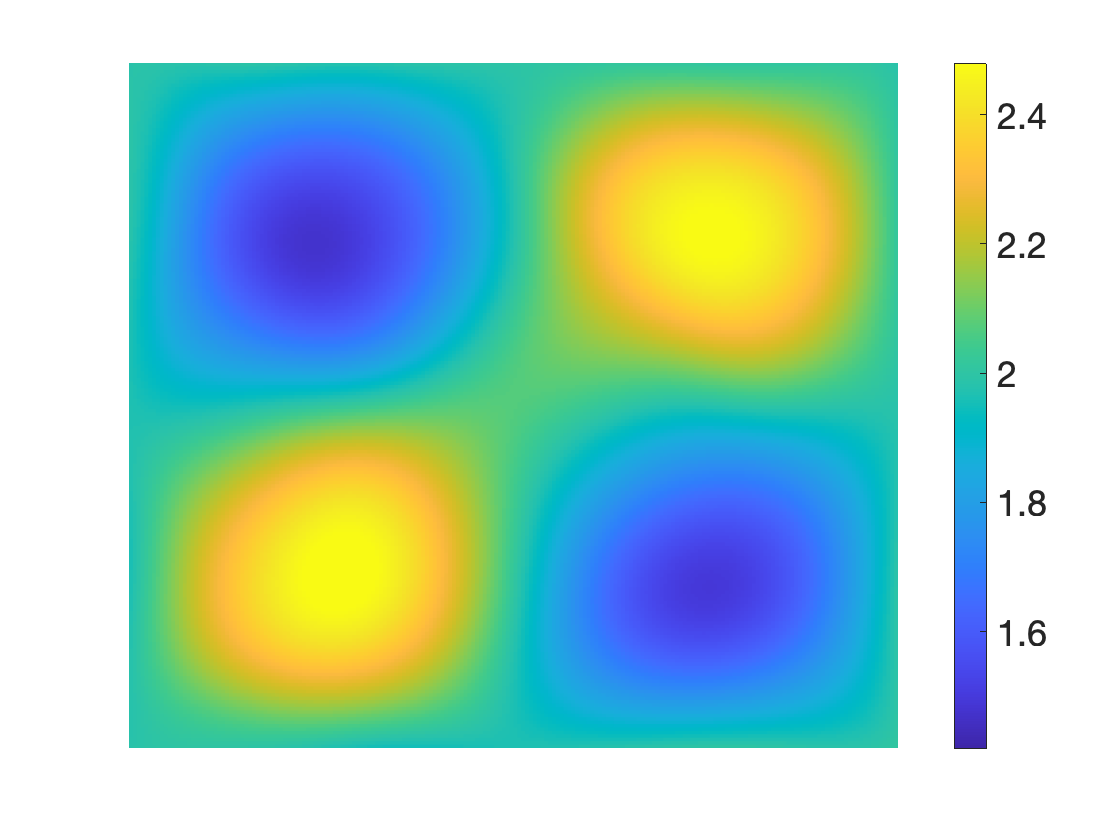} &
\includegraphics[width=0.32\textwidth]{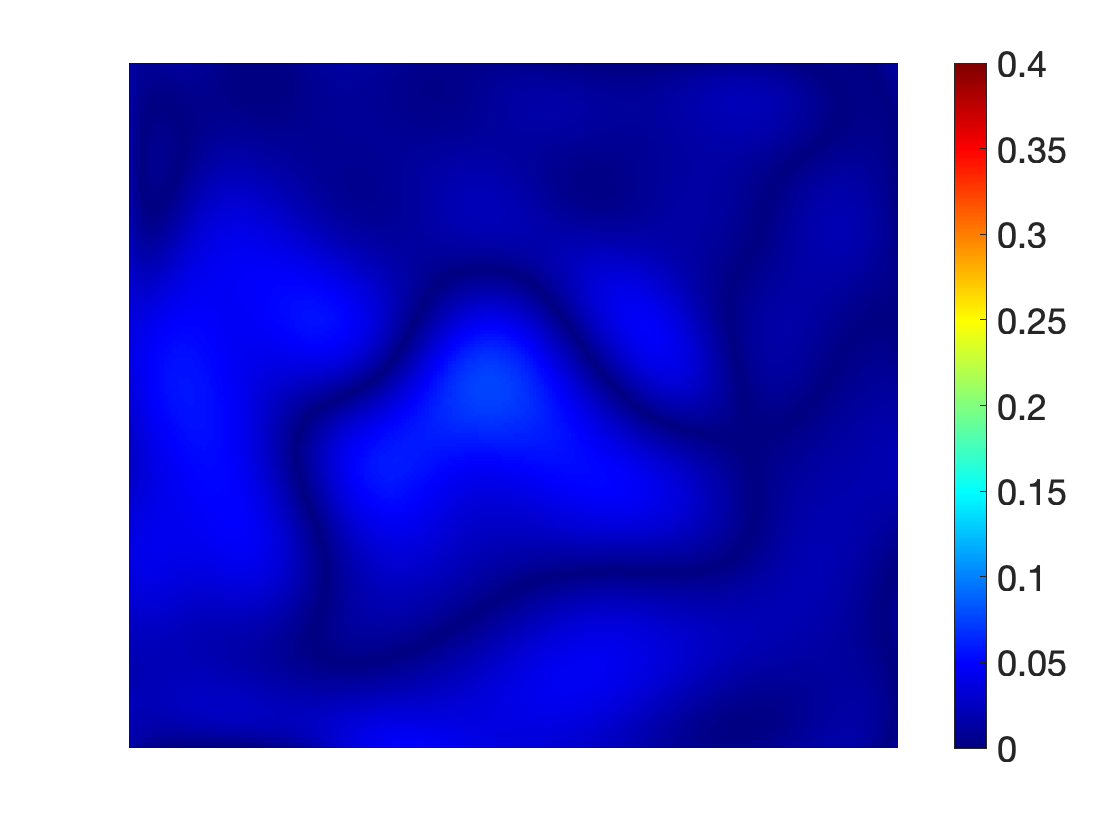}\\
\includegraphics[width=0.32\textwidth]{diripartial2dex.png} &
\includegraphics[width=0.32\textwidth]{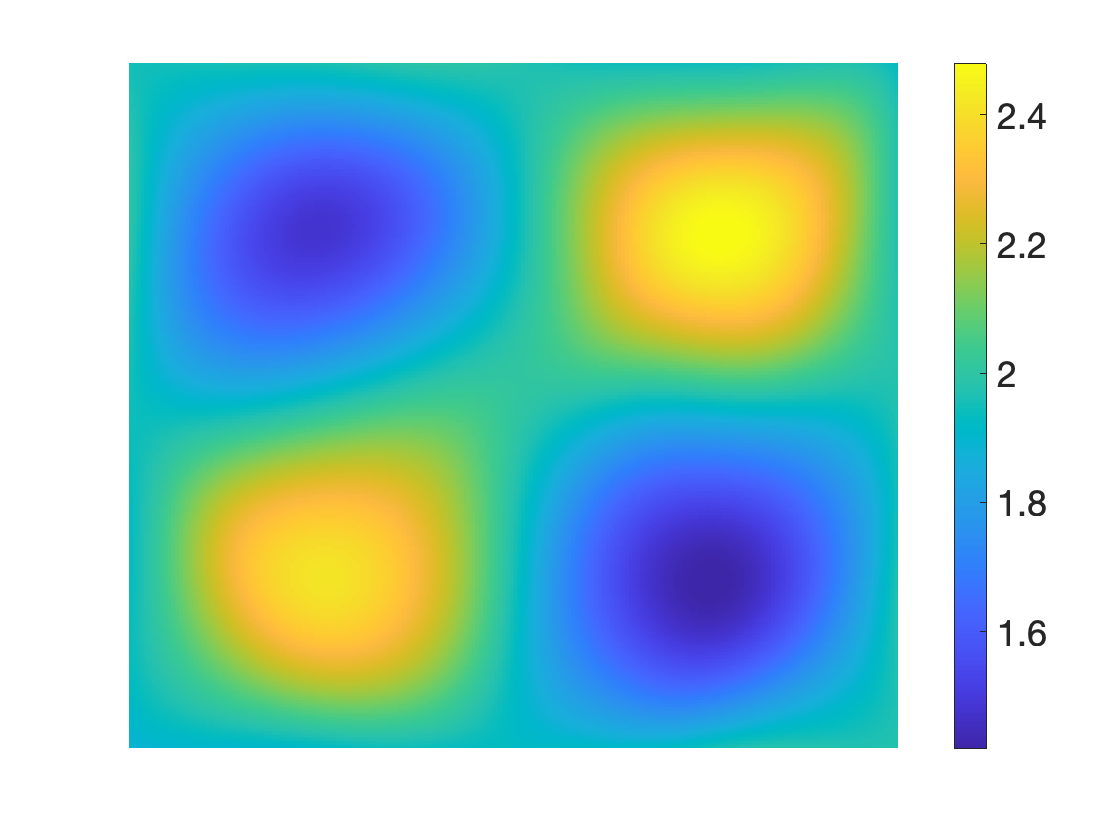} &
\includegraphics[width=0.32\textwidth]{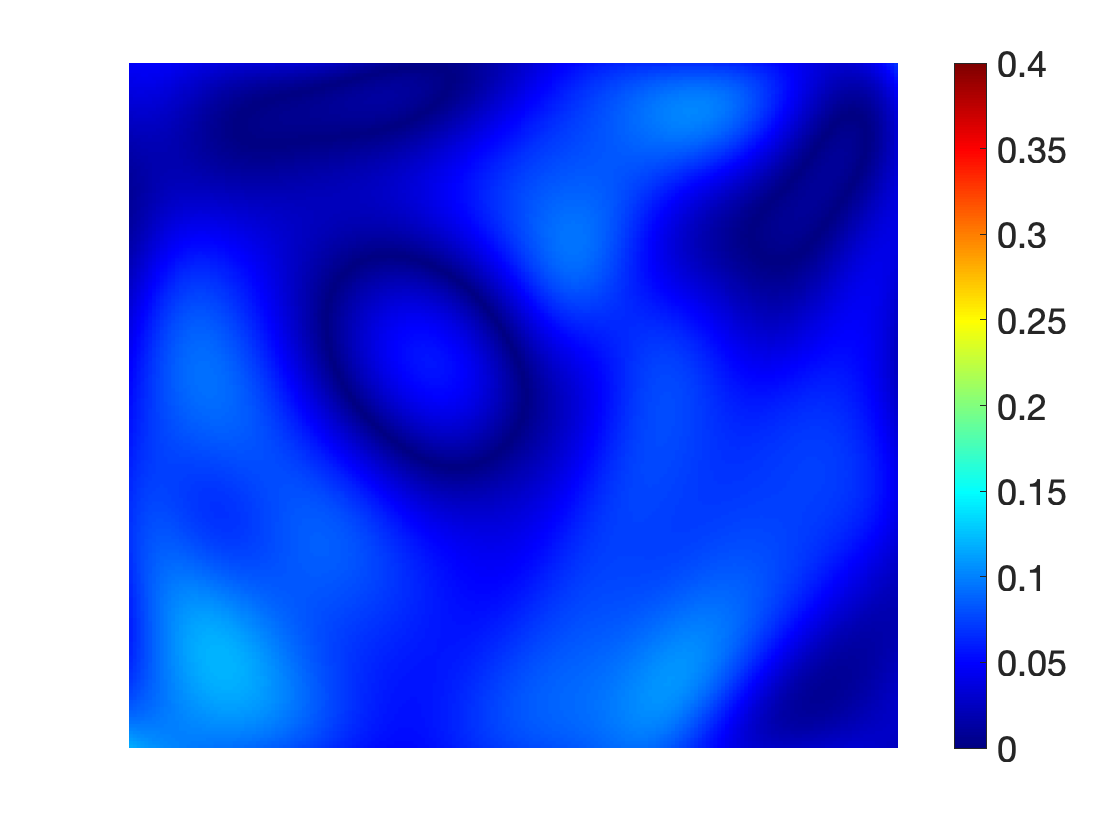}\\
(a) $q^\dag$  & (b) $\hat q$ & (c) $|\hat q-q^\dag|$
\end{tabular}
\caption{The reconstructions for Example \ref{exam:diripartial2d} using DNN method with exact data $($top$)$ and noisy data $(\delta=5\%,10\%$, middle, bottom$)$.}
\label{fig:diripartial2d}
\end{figure}

\begin{figure}[htbp]
\centering
\setlength{\tabcolsep}{0em}
\begin{tabular}{ccc}
\includegraphics[width=0.32\textwidth]{diripartial2dex.png} &
\includegraphics[width=0.32\textwidth]{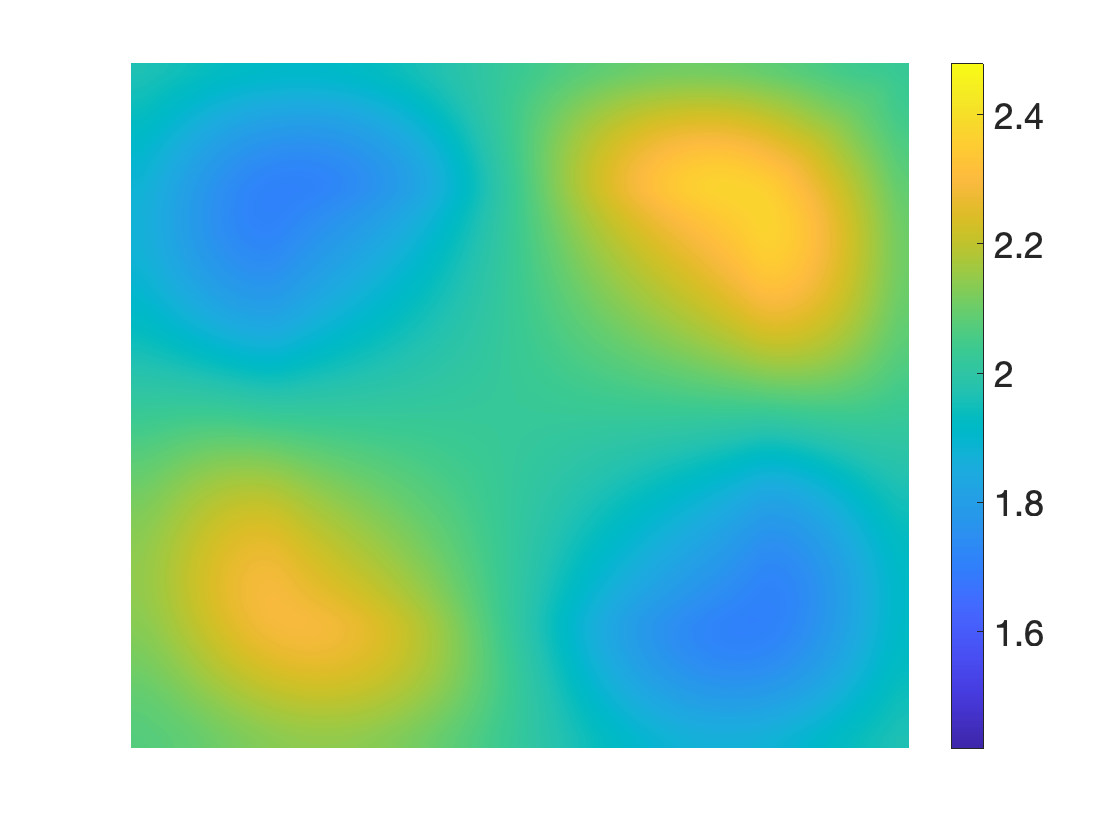} &
\includegraphics[width=0.32\textwidth]{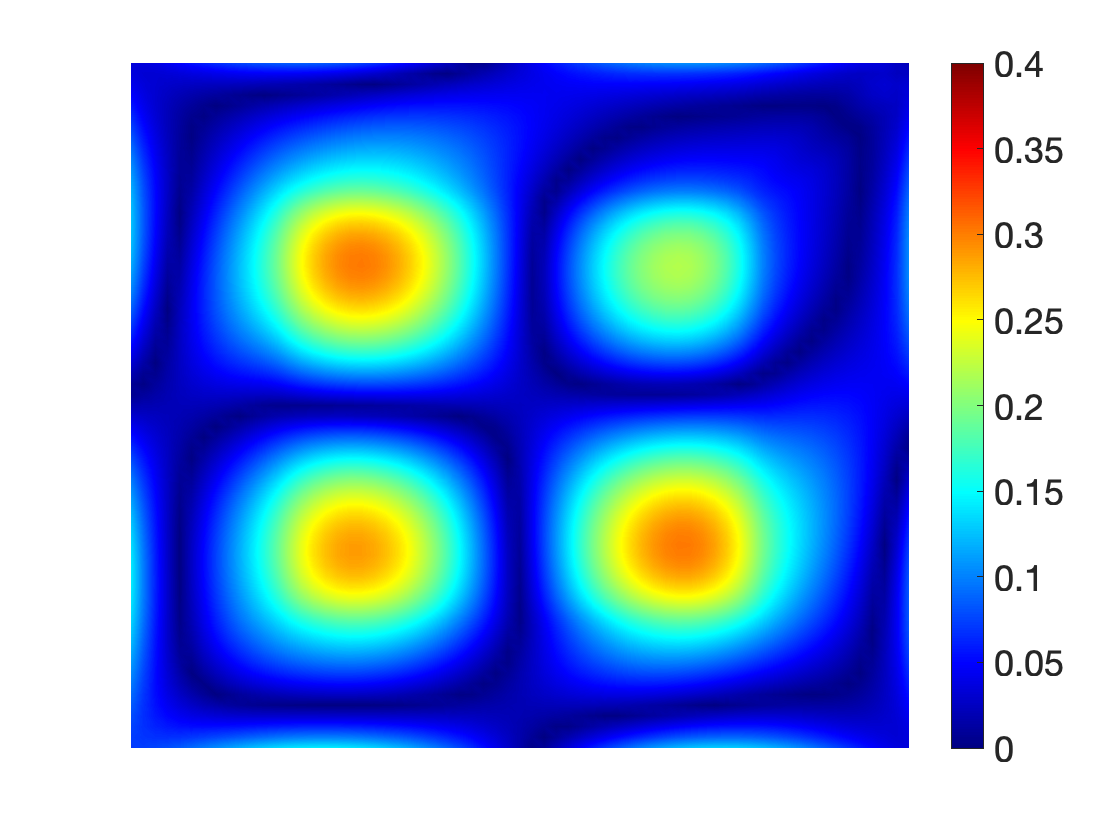}\\
\includegraphics[width=0.32\textwidth]{diripartial2dex.png} &
\includegraphics[width=0.32\textwidth]{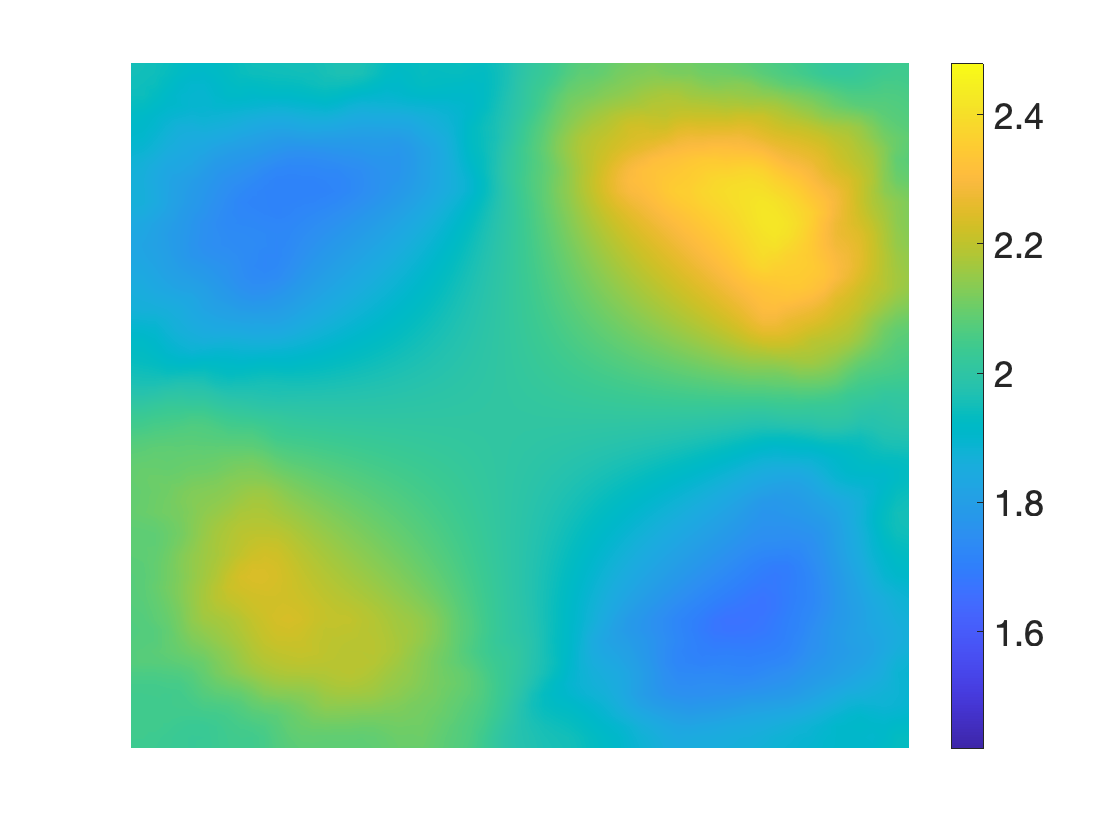} &
\includegraphics[width=0.32\textwidth]{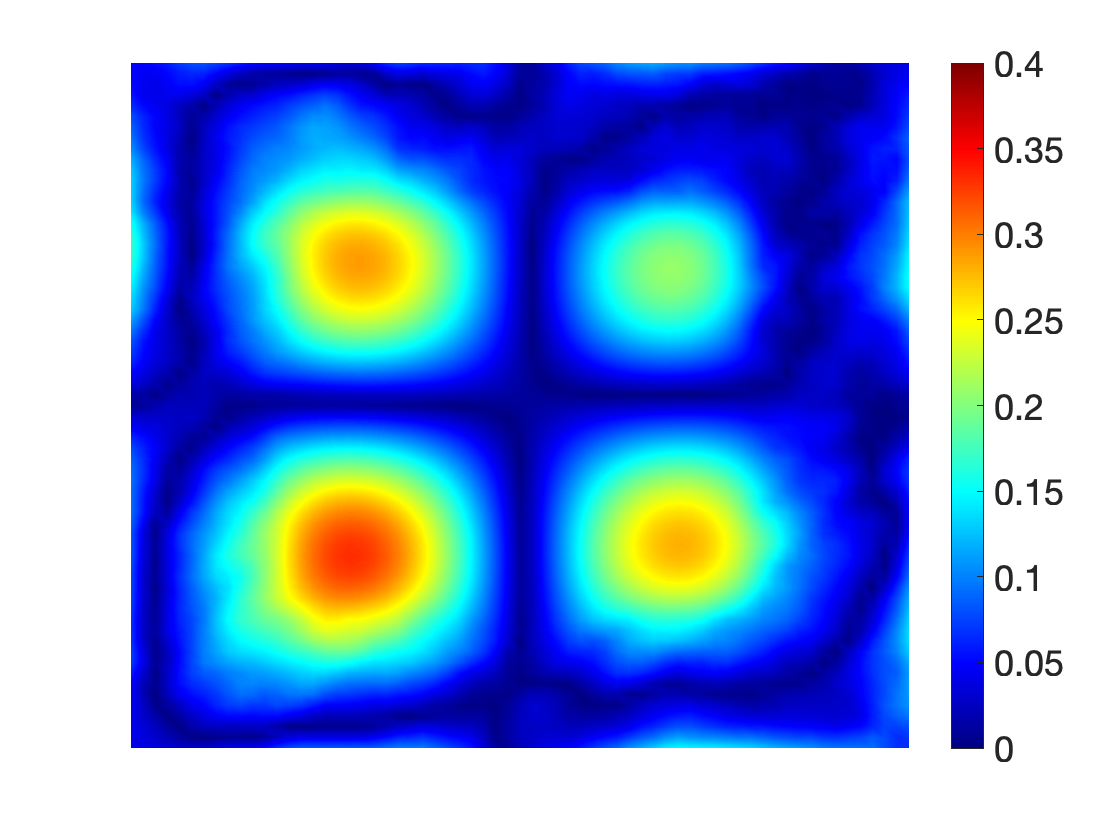}\\
\includegraphics[width=0.32\textwidth]{diripartial2dex.png} &
\includegraphics[width=0.32\textwidth]{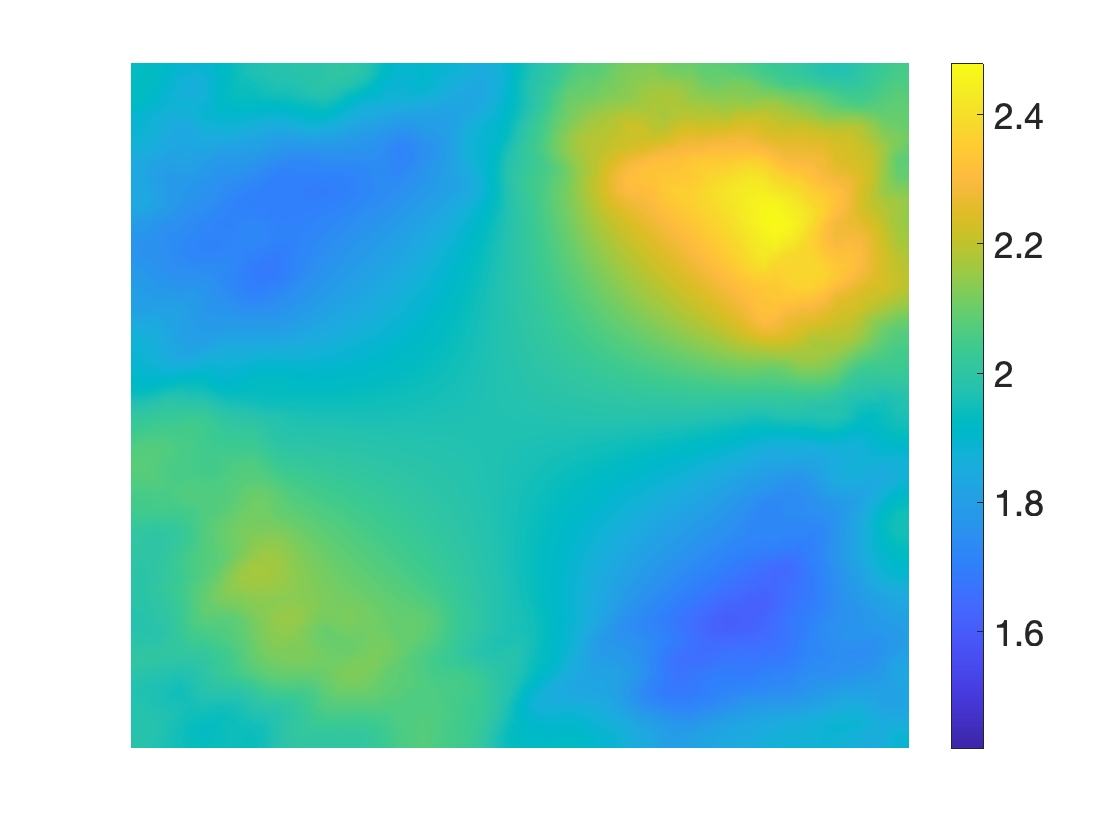} &
\includegraphics[width=0.32\textwidth]{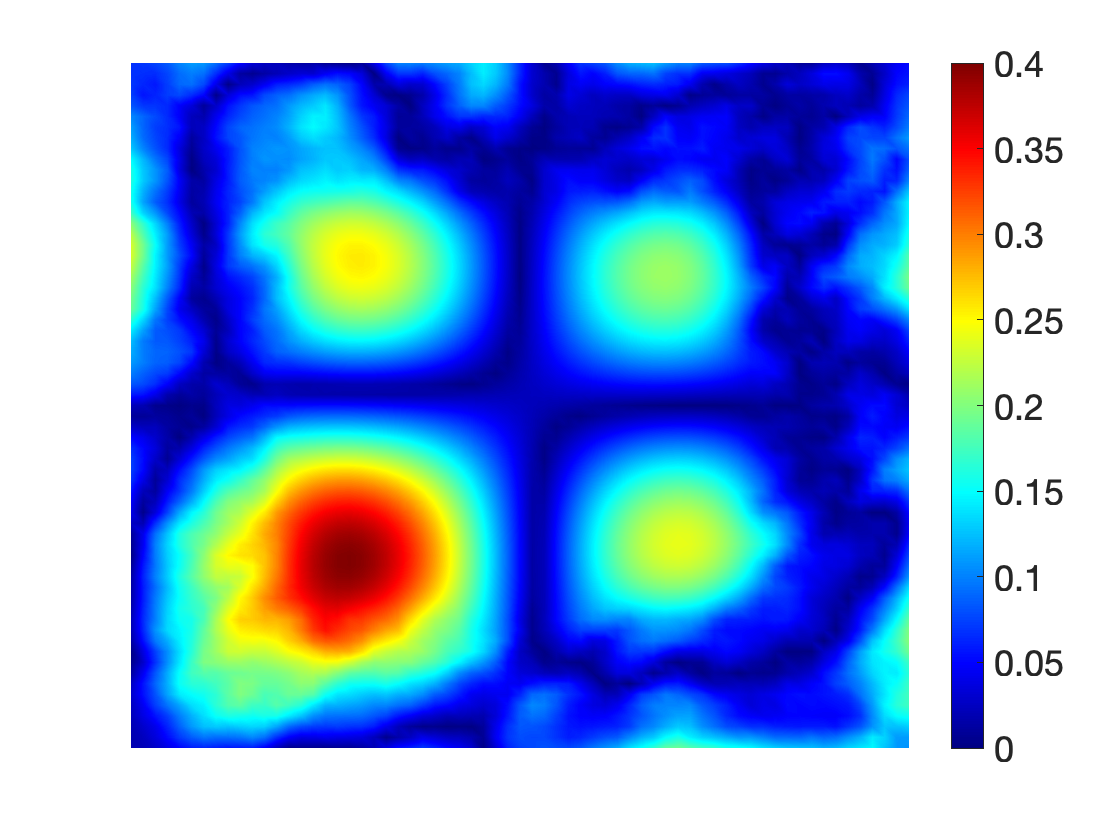}\\
(a) $q^\dag$  & (b) $\hat q$ & (c) $|\hat q-q^\dag|$
\end{tabular}
\caption{The reconstructions for Example \ref{exam:diripartial2d} using FEM with exact data $($top$)$ and noisy data $(\delta=5\%,10\%$, middle, bottom$)$.}
\label{fig:diripartial2dfem}
\end{figure}

The DNN approach minimizes the empirical version of the following loss:
\begin{equation*}
\begin{aligned}
J_{\bsgamma}(\theta,\kappa)
=&\|\sigma_\kappa-\P01(q_\theta)\nabla z^\delta\|_{L^2(\omega)}^2+\gamma_\sigma\|\nabla\cdot\sigma_\kappa+f\|_{L^2(\Omega)}^2\\
 &+\gamma_q\|\nabla q_\theta\|_{L^2(\Omega)}^2+\gamma_b\|\sigma_\kappa-q^\dag\nabla z^\delta\|^2_{L^2(\partial\Omega)}.
\end{aligned}
\end{equation*}
It yields accurate reconstruction for all noise levels, whereas the FEM approach fails to recover important features of the conductivity in the subdomain $\Omega\setminus\omega$, cf. Figs. \ref{fig:diripartial2d}-\ref{fig:diripartial2dfem}. This example again shows the superiority of the DNN approach over traditional FEM for partial internal data, and hence its potential for solving more ill-posed inverse problems, e.g., EIT.

The last example is about recovering a 3D conductivity from partial internal data.
\begin{example}\label{exam:diripartial3d}
The domain $\Omega=(0,1)^3$, the measurement $\nabla z^\delta$ on the subdomain $\Omega\setminus (0.2,0.8)^3$, $q^\dag=2+\sin(\pi x_1)\sin(\pi x_2)\sin(\pi x_3)$, and $u^\dag=\sum_{i=1}^3(x_i+\frac13x_i^3)$.
\end{example}
\begin{figure}[htbp]
\centering
\setlength{\tabcolsep}{0em}
\begin{tabular}{ccc}
\includegraphics[width=0.32\textwidth]{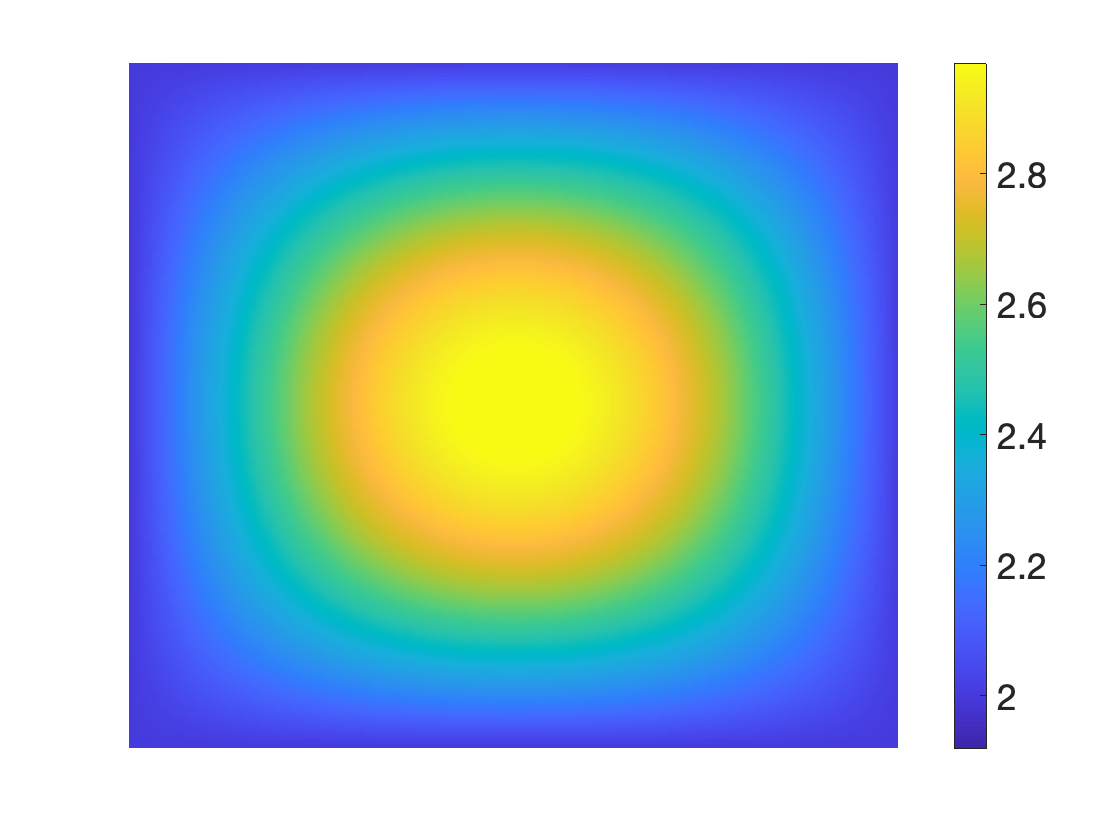} &
\includegraphics[width=0.32\textwidth]{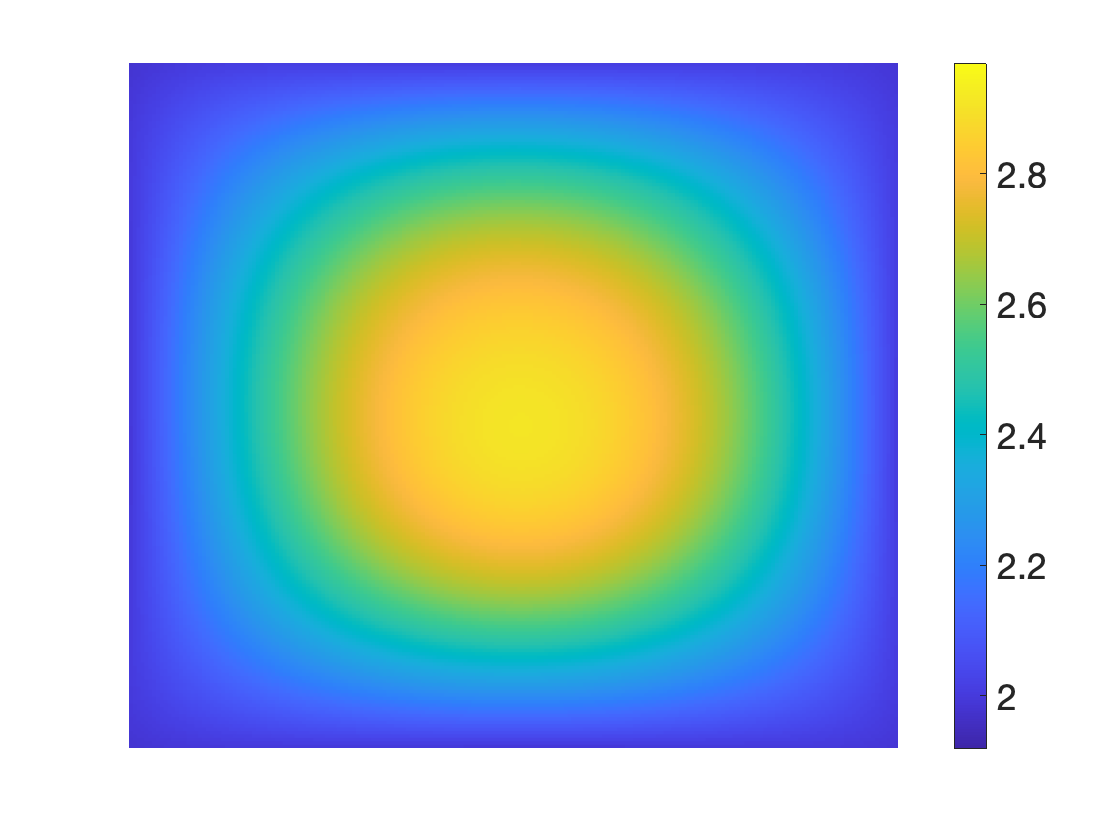} &
\includegraphics[width=0.32\textwidth]{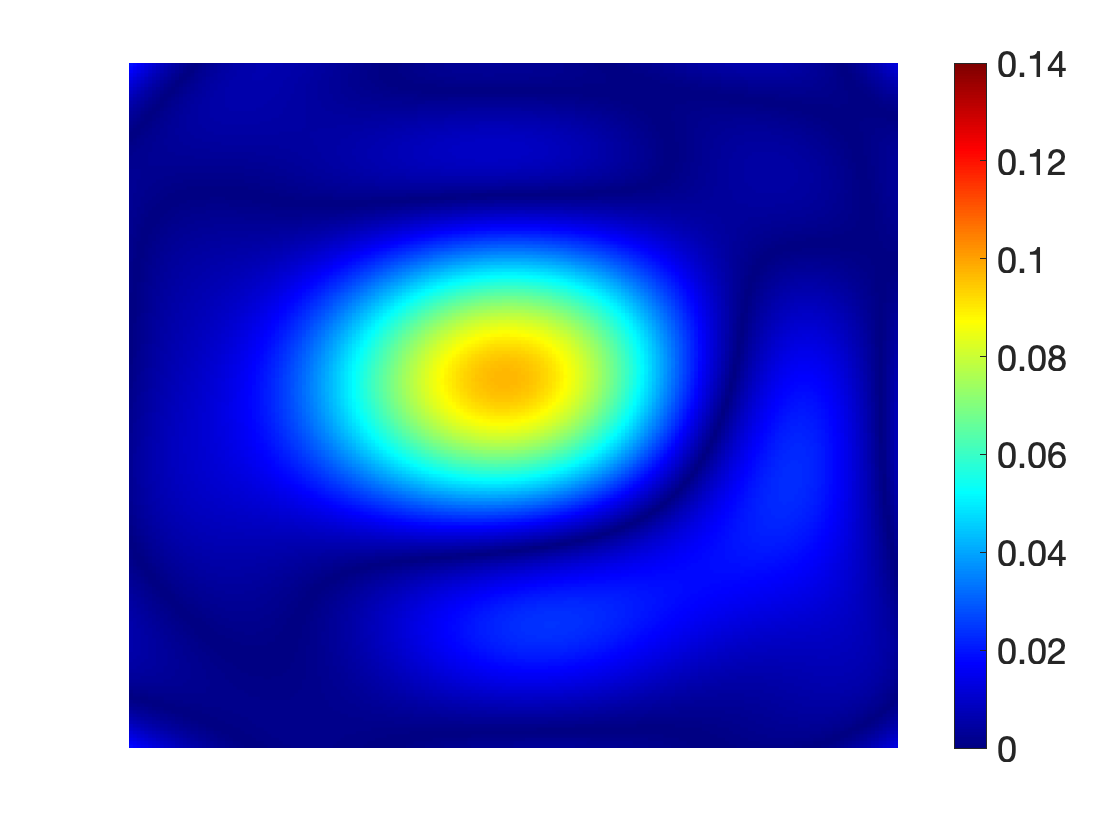}\\
\includegraphics[width=0.32\textwidth]{diripartial3dex.png} &
\includegraphics[width=0.32\textwidth]{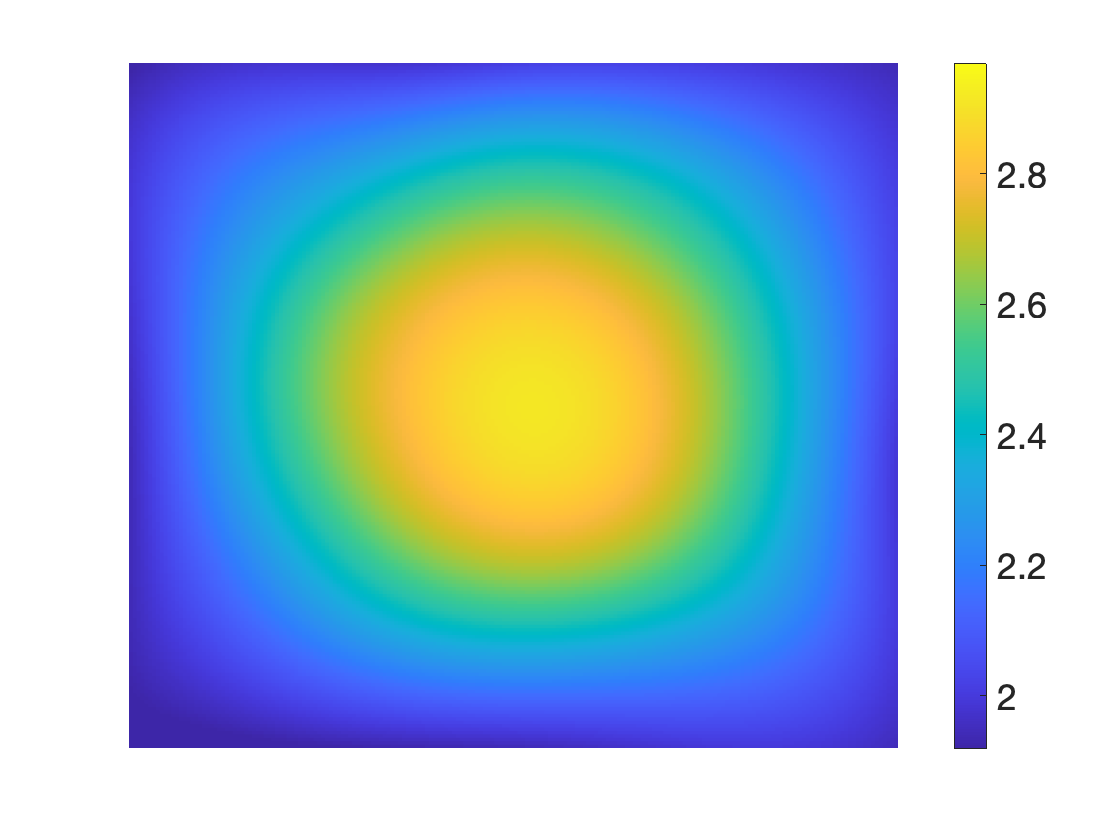} &
\includegraphics[width=0.32\textwidth]{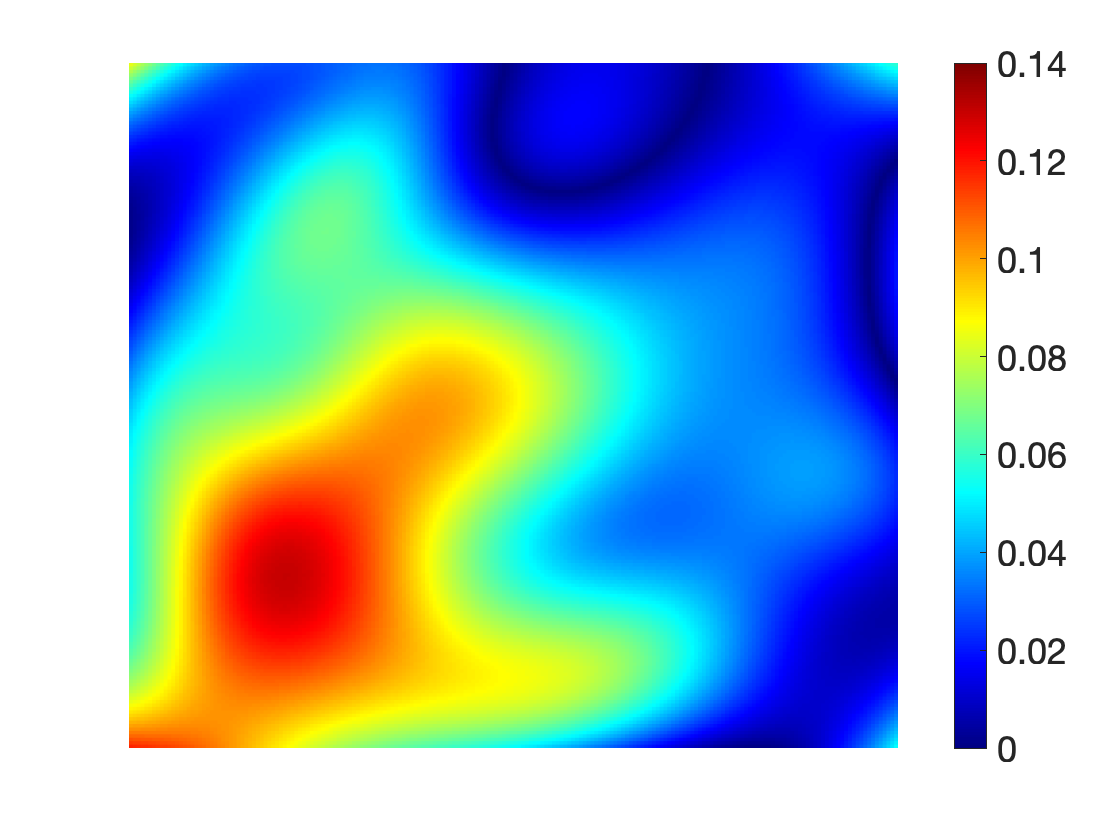}\\
(a) $q^\dag$  & (b) $\hat q$ & (c) $|\hat q-q^\dag|$
\end{tabular}
\caption{The reconstructions for Example \ref{exam:diripartial3d} using DNN method with exact data $($top$)$ and noisy data $(\delta=10\%$, bottom$)$.}
\label{fig:diripartial3d}
\end{figure}

Fig. \ref{fig:diripartial3d} shows the reconstruction on a 2D cross section at
$x_3=0.5$. The relative $L^2(\Omega)$-error $e(\hat q)$ for the exact and noisy data is 1.08e-2
and 2.52e-2, respectively. The reconstruction is accurate with both exact and $10\%$ noisy data over the domain $\Omega$.

\section{Conclusion}

In this work we have developed a novel approach for recovering the conductivity distribution in an elliptic PDE from one internal measurement. It combines the least-squares method and mixed formulation, and extends the work of Kohn and Lowe \cite{KohnLowe:1988} to the context of DNN approximations. By combining the stability argument with \textit{a priori} bounds on the loss, we have derived $L^2(\Omega)$ error estimates for the DNN approximations of the conductivity, thereby providing theoretical justifications of the approach. We have presented numerical experiments in high noise, high-dimensions and partial interior data, and the results show that the approach indeed holds big potentials for solving the concerned PDE inverse problem. Despite the compelling numerical evidence, there are still many interesting issues that deserve further investigations, which we briefly discuss below.

First, there are other types of measurements, e.g., power density or flux density in the (partial) interior domain, or boundary measurement (e.g., EIT or diffuse optical tomography), for which there are also diverse stability estimates. It is an interesting issue to extend the proposed approach to these important settings and to provide relevant theoretical underpinnings. Our preliminary experiments indicate that it is indeed promising for EIT type problems yet a more thorough evaluation is highly desirable. Relevant mathematical analysis requires new ideas.

Second, the formulation involves several hyper-parameters (e.g., $\gamma_\sigma$, $\gamma_b$ and $\gamma_q$) that have to be properly tuned. While the theoretical analysis sheds useful insights, the optimal tuple does require some tuning, and its efficient tuning represents one important task. One promising strategy is to adaptively learn via, e.g., hierarchical modeling (in the Bayesian framework) \cite{WangZabaras:2005} and bi-level type framework \cite{YuZhu:2020}. Further, the architectural choice of fully connected DNNs is important, e.g., depth, width and weight bound. The larger is the size, the more flexible is the DNN but the larger is the statistical error, leading to bias-variance type trade-off. In the lens of classical regularization, fixing the architectural parameters can also play a similar role to regularization by (nonlinear) projection / discretization, and its interplay with explicit regularization is an interesting theoretical question.

Third and last, from a practical aspect, the implementation of the proposed approach employs the off-shelf optimizer ADAM, which has been widely employed in practice. It often takes thousands of iterations to reach convergence, and the computational efficiency is not yet desirable. This is also widely observed for neural solvers for direct problems. Even worse, the optimizer tends to produce only a local minimum, resulting in pronounced (and often dominating) optimization errors. This observation greatly complicates the numerical verification of error estimates, which so far has not been achieved for neural solvers for PDEs and remains one outstanding challenge in the community \cite{SiegelHongJin:2023}. It is of great interest to analyze the optimization error and its impact on the reconstruction, especially its connection with the empirically observed robustness with respect to data noise.

\section*{Acknowledgements}
The authors are grateful to two anonymous referees for their constructive comments on the paper which have led to an improvement of the quality of the paper.

\appendix

\section{The proof of Theorem \ref{thm:stat-err}}
In this appendix, we prove Theorem \ref{thm:stat-err}, which gives high probability bounds on the statistical errors
$\Delta \mathcal{E}_{i}$, $i\in \{d,b,b',q,\sigma\}$. These bounds play a crucial role in proving Theorems \ref{thm:err-Neum-emp} and \ref{thm:err-Diri-emp}. By slightly abusing the notation, let $\mathcal{N}_\sigma\equiv \mathcal{N}(L_\sigma,W_\sigma,R_\sigma)$
and $\mathcal{N}_q=\mathcal{N}(L_q,W_q,R_q)$ be two DNN function
classes of given depth, width and parameter bound for approximating the flux $\sigma$ and conductivity $q$, respectively.
Then we define
\begin{align*}
  \mathcal{H}_d  &= \{h: \Omega\to \mathbb{R}|\,\, h(x) = \|\sigma_\kappa(x)-P_\mathcal{A}(q_\theta(x))\nabla z^\delta(x)\|_{\ell^2}^2, \sigma_\kappa\in \mathcal{N}_\sigma, q_\theta\in \mathcal{N}_q\},\\
  \mathcal{H}_\sigma & = \{h: \Omega\to \mathbb{R}|\,\, h(x) = |\nabla\cdot \sigma_\kappa(x)+f(x)|^2, \sigma_\kappa\in \mathcal{N}_\sigma\},\\
  \mathcal{H}_b & = \{h:\partial\Omega\to \mathbb{R}|\,\, h(x) = |\n \cdot \sigma_\kappa(x)-g(x)|^2, \sigma_\kappa\in \mathcal{N}_\sigma\},\\
  \mathcal{H}_{b'} & = \{h:\partial\Omega\to \mathbb{R}|\,\, h(x) = \|\sigma_\kappa(x)-q^\dag \nabla z^\delta\|_{\ell^2}^2, \sigma_\kappa\in \mathcal{N}_\sigma\},\\
  \mathcal{H}_q & = \{h:\Omega\to \mathbb{R}|\,\, h(x) =  \|\nabla q_\theta(x)\|_{\ell^2}^2, q_\theta\in \mathcal{N}_q\}.
\end{align*}

To this end, we employ Rademacher complexity \cite{BartlettMendelson:2002} and boundedness and
Lipschitz continuity of DNN functions and their derivatives with respect to the DNN parameters.
Rademacher complexity measures the complexity of a collection of functions
by the correlation between function values with Rademacher random variables.
\begin{definition}\label{def: Rademacher}
Let $\mathcal{F}$ be a real-valued function class defined on the domain $\Omega$ {\rm(}or the boundary $\partial\Omega${\rm)}, $\xi=\{\xi_j\}_{j=1}^n$ be i.i.d. samples from the distribution $\mathcal{U}(\Omega)$ {\rm(}or the distribution $\mathcal{U}(\partial\Omega)${\rm)}, and $\omega=\{\omega_j\}_{j=1}^n$ be i.i.d  Rademacher random variables with probability $P(\omega_j=1)=P(\omega_j=-1)=\frac12$. Then
	the Rademacher complexity $\mathfrak{R}_n(\mathcal{F})$ of the class $\mathcal{F}$ is defined by
	\begin{equation*}
		\mathfrak{R}_n(\mathcal{F})=\mathbb{E}_{\xi,\omega}\bigg{[}\sup_{v\in\mathcal{F}}\ n^{-1}\bigg{\lvert}\ \sum_{j=1}^{n}\omega_j v(\xi_j)\ \bigg{\rvert} \bigg{]}.
	\end{equation*}
\end{definition}

We use the following PAC-type generalization bound via Rademacher complexity
\cite[Theorem 3.1]{mohri2018foundations}. Note that \cite[Theorem 3.1]{mohri2018foundations}
only discusses the case of the function class ranging in $[0 ,1]$. For a general
bounded function, applying McDiarmid’s inequality and the original argument yields the following result.
\begin{lemma}\label{lem:PAC}
Let $X_1,\ldots,X_n$ be a set of i.i.d. random variables. Let $\mathcal{F}$
be a function class defined on $D$ such that $\sup_{v\in\mathcal{F}}\|v\|_{L^\infty(D)}\leq M_\mathcal{F}<\infty$. Then for any $\tau\in(0,1)$, with probability at least $1-\tau$,
\begin{equation*}
  \sup_{v\in \mathcal{F}}\bigg|n^{-1}\sum_{j=1}^n v(X_j)-\mathbb{E}[v(X)]\bigg| \leq 2\mathfrak{R}_n(\mathcal{F}) + 2M_\mathcal{F}\sqrt{\frac{\log\frac{1}{\tau}}{2n}}.
\end{equation*}
\end{lemma}

To apply Lemma \ref{lem:PAC}, we have to bound the Rademacher complexity of $\mathcal{H}_i$, $i\in\{d,\sigma,b,b',q\}$. This is achieved by combining Lipschitz continuity of DNN functions in $\mathcal{H}_i$ in the DNN parameter $\theta$ and $\kappa$, and Dudley's formula in Lemma \ref{lem:Dudley}. The next lemma gives useful boundedness and Lipschitz continuity of a $\tanh$-DNN function
class. Lemma \ref{lem:NN-Lip} also holds with
the $L^\infty(\partial\Omega)$ norm in place of the $L^\infty(\Omega)$ norm, since the argument depends only on the boundedness of $\rho=\tanh$ and
its derivative.
\begin{lemma}\label{lem:NN-Lip}
Let $\Theta$ be a parametrization with depth $L$ and width $W$, and $\theta=\{(A^{(\ell)},b^{(\ell)})_{\ell=1}^L\}, \tilde{\theta}=\{(\tilde{A}^{(\ell)},\tilde{b}^{(\ell)})_{\ell=1}^L\}\in\Theta$. Then for the DNN realizations
$v,\tilde{v}:\Omega\to\mathbb{R}$ of $\theta,\tilde{\theta}$ with $\|\theta\|_{\ell^\infty},
\|\tilde{\theta}\|_{\ell^\infty}\leq R$, the following estimates hold
\begin{enumerate}
\item[{\rm(i)}] $\|v\|_{L^\infty(\Omega)}\leq R(W+1), \quad \|\nabla v\|_{L^\infty(\Omega; \mathbb{R}^d)}\leq \sqrt{d}R^LW^{L-1}$;
\item[{\rm(ii)}] $\|v-\tilde{v}\|_{L^\infty(\Omega)}\leq 2LR^{L-1}W^{L}\|\theta-\tilde\theta\|_{\ell^\infty}$,\\ $\|\nabla (v-\tilde{v})\|_{L^\infty(\Omega; \mathbb{R}^d)}\leq  \sqrt{d}L^2R^{2L-2}W^{2L-2}\|\theta-\tilde\theta\|_{\ell^\infty}$.
\end{enumerate}
\end{lemma}
\begin{proof}
See \cite[p. 19 and p. 22 of Lemma 3.4]{jin2022imaging} and \cite[Remark 3.3]{jin2022imaging} for the proof.
\end{proof}

Next we appeal to reduction to
parameterisation, via Lipschitz continuity of the DNN functions.
\begin{lemma}\label{lem:fcn-Lip}
Let $c_z=\|\nabla z^\delta\|_{L^\infty(\Omega)}$ and $c_z'=\|\nabla z^\delta\|_{L^\infty(\partial\Omega)}$. For the function
classes $\mathcal{H}_i$, $i\in \{d,\sigma,b,b',q\}$, the functions are
uniformly bounded:
\begin{align*}
  \|h\|_{L^\infty(\Omega)} &\leq \left\{\begin{aligned}
    M_d&=2(dR_\sigma^{2}(W_\sigma+1)^{2}+c_1^2c_z^2), &&\quad h\in \mathcal{H}_d,\\
    M_\sigma &= 2(d^2R_\sigma^{2L_\sigma}W_\sigma^{2L_\sigma-2}+\|f\|^2_{L^\infty(\Omega)}), &&\quad h\in \mathcal{H}_\sigma,\\
    M_q &= dR_q^{2L_q}W_q^{2L_q-2}, && \quad h\in \mathcal{H}_q,
  \end{aligned}\right. \\
  \|h\|_{L^\infty(\partial\Omega)}& \leq \left\{\begin{aligned}
    M_b &= 2(dR_\sigma^{2}(W_\sigma+1)^2+\|g\|^2_{L^\infty(\partial\Omega)}), && \quad h\in \mathcal{H}_b,\\
    M_{b'} &= 2(dR_\sigma^{2}(W_\sigma+1)^{2}+c_1^2(c_z')^2),  && \quad h\in \mathcal{H}_{b'}.
  \end{aligned}\right.
\end{align*}
Moreover, the following Lipschitz continuity estimates in the DNN parameters hold
\begin{align*}
  \|h-\tilde h\|_{L^\infty(\Omega)} & \le  \Lambda_d (\|\theta -\tilde\theta\|_{\ell^\infty} + \|\kappa-\tilde\kappa\|_{\ell^\infty}),\quad \forall h,\tilde h\in \mathcal{H}_d,\\
  \|h-\tilde h\|_{L^\infty(\Omega)} & \leq \Lambda_\sigma \|\kappa-\tilde \kappa\|_{\ell^\infty},\quad \forall h,\tilde h\in \mathcal{H}_\sigma,\\
  \|h-\tilde h\|_{L^\infty(\partial\Omega)} & \leq \Lambda_b \|\kappa-\tilde \kappa\|_{\ell^\infty},\quad \forall h,\tilde h\in \mathcal{H}_b,\\
  \|h-\tilde h\|_{L^\infty(\partial\Omega)} & \leq \Lambda_{b'} \|\kappa-\tilde \kappa\|_{\ell^\infty},\quad \forall h,\tilde h\in \mathcal{H}_{b'},\\
  \|h_{\theta}-h_{\tilde \theta}\|_{L^\infty(\Omega)} & \leq \Lambda_q \|\theta-\tilde \theta\|_{\ell^\infty},\quad \forall h,\tilde h\in \mathcal{H}_q,
\end{align*}
with the constants $\Lambda_i$, $i\in \{d,\sigma,b,b',q\}$, given by
\begin{align*}
  \Lambda_d & =2\big(\sqrt{d}R_\sigma(W_\sigma+1)+c_1c_z\big)\max(2\sqrt{d}L_\sigma R_\sigma^{L_\sigma-1}W_\sigma^{L_\sigma},c_zL_qR_q^{L_q-1}W_q^{L_q}),\\
  \Lambda_\sigma & =  2(dR_\sigma^{L_\sigma}W^{L_\sigma-1}_{\sigma} +\|f\|_{L^\infty(\Omega)})dL_\sigma^2R_\sigma^{2L_\sigma-2}W_\sigma^{2L_\sigma-2},\\
  \Lambda_b & = 4(\sqrt{d}R_\sigma(W_{\sigma}+1)+\|g\|_{L^\infty(\partial\Omega)})\sqrt{d}L_\sigma R_\sigma^{L_\sigma-1}W_\sigma^{L_\sigma},\\
  \Lambda_{b'} & = 4\big(\sqrt{d}R_\sigma(W_\sigma+1)+c_1c_z'\big)\sqrt{d}L_\sigma R_\sigma^{L_\sigma-1}W_\sigma^{L_\sigma},\\
  \Lambda_q & = 2dL_q^2R_q^{3L_q-2}W_q^{3L_q-3}.
\end{align*}
\end{lemma}
\begin{proof}
The assertions follow directly from Lemma \ref{lem:NN-Lip}. For $h_{\theta,\kappa}\in \mathcal{H}_d$, we have
\begin{align*}
  |h_{\theta,\kappa}(x)| &\leq 2(\|\sigma_\kappa(x)\|_{\ell^2}^2 +\|\P01(q_\theta(x))\nabla z^\delta(x)\|_{\ell^2}^2)
 \leq 2(dR_\sigma^{2}(W_\sigma+1)^{2}+c_1^2c_z^2).
\end{align*}
For any $h_{\theta,\kappa}, h_{\tilde\theta,\tilde\kappa}\in \mathcal{H}_d$, by completing the squares and the Cauchy-Schwarz inequality, meanwhile noting the stability of $\P01$ in \eqref{eqn:P01-stab} we have
\begin{align*}
    & h_{\theta,\kappa}(x) - h_{\tilde\theta,\tilde\kappa}(x) = \|\sigma_\kappa(x)-\P01(q_\theta(x))\nabla z^\delta(x)\|_{\ell^2}^2 -
  \|\sigma_{\tilde\kappa}(x)-\P01(q_{\tilde\theta}(x))\nabla z^\delta(x)\|_{\ell^2}^2\\
    =& \big(\sigma_\kappa(x)-\P01(q_\theta(x))\nabla z^\delta(x)+ \sigma_{\tilde\kappa}(x)-\P01(q_{\tilde \theta}(x))\nabla z^\delta(x),\\
    &\quad \sigma_\kappa(x)- \sigma_{\tilde\kappa}(x)+(-\P01(q_\theta(x))+\P01(q_{\tilde\theta}(x)))\nabla z^\delta(x)\big)\\
  \leq & 2\big(\sup_{\sigma_\kappa\in \mathcal{N}_\kappa}\|\sigma_\kappa(x)\|_{L^\infty(\Omega;\mathbb{R}^d)} + c_1c_z\big) \\ &\times \big(\|\sigma_\kappa(x)-\sigma_{\tilde\kappa}(x)\|_{L^\infty(\Omega;\mathbb{R}^d)}+c_z\|q_\theta(x)-q_{\tilde\theta}(x)\|_{L^\infty(\Omega)}\big).
\end{align*}
Then by Lemma \ref{lem:NN-Lip}, with $A=2\big(\sqrt{d}R_\sigma(W_\sigma+1)+c_1c_z\big)$, we deduce
\begin{align*}
      &|h_{\theta,\kappa}(x) - h_{\tilde\theta,\tilde\kappa}(x)|  \leq  A (2\sqrt{d}L_\sigma R_\sigma^{L_\sigma-1}W_\sigma^{L_\sigma}\|\kappa-\tilde\kappa\|_{\ell^\infty}+c_zL_qR_q^{L_q-1}W_q^{L_q}\|\theta-\tilde\theta\|_{\ell^\infty})\\
       \le& A\max
      (2\sqrt{d}L_\sigma R_\sigma^{L_\sigma-1}W_\sigma^{L_\sigma},c_zL_qR_q^{L_q-1}W_q^{L_q})
      (\|\kappa-\tilde\kappa\|_{\ell^\infty}+\|\theta-\tilde\theta\|_{\ell^\infty}).
\end{align*}
The remaining estimates follow similarly. This completes the proof of the lemma.
\end{proof}

Next we bound the Rademacher complexities $\mathcal{R}_n(\mathcal{H}_i)$, $i\in\{d,
\sigma,b,b',q\}$, using the concept of the covering number.
\begin{definition}
Let $\mathcal{G}$ be a
real-valued function class equipped with the metric $\rho$. A collection
of points $\{x_i\}_{i=1}^n \subset \mathcal{G}$ is called an
$\epsilon$-cover of $\mathcal{G}$ if for any $x\in \mathcal{G}$, there
exists at least one $i \in \{1,\dots,n\}$ such that $\rho(x, x_i) \leq \epsilon$. The
$\epsilon$-covering number $\mathcal{C}(\mathcal{G}, \rho, \epsilon)$ is the minimum
cardinality among all $\epsilon$-cover of $\mathcal{G}$.
\end{definition}

Then we can state Dudley's theorem \cite[Theorem 9]{lu2021priori} and \cite[Theorem 1.19]{wolf2018mathematical}.
\begin{lemma}\label{lem:Dudley}
Let $M_\mathcal{F}:=\sup_{f\in\mathcal{F}} \|f\|_{L^{\infty}(\Omega)}$, and $\mathcal{C}(\mathcal{F},\|\cdot\|_{L^{\infty}(\Omega)},\epsilon)$ be the covering number of the set $\mathcal{F}$. Then the Rademacher complexity $\mathfrak{R}_n(\mathcal{F})$ is bounded by
\begin{equation*}
\mathfrak{R}_n(\mathcal{F})\leq\inf_{0<s< M_\mathcal{F}}\bigg(4s\ +\ 12n^{-\frac12}\int^{M_\mathcal{F}}_{s}\big(\log\mathcal{C}(\mathcal{F},\|\cdot\|_{L^{\infty}(\Omega)},\epsilon)\big)^{\frac12}\ {\rm d}\epsilon\bigg).
\end{equation*}
\end{lemma}

Last, we give the proof of Theorem \ref{thm:stat-err}.
\begin{proof}
For any $n \in \mathbb{N }$, $R \in [1, \infty)$, $\epsilon \in  (0,1)$,
and $ B_R := \{x\in\mathbb{R}^n:\ \|x\|_{\ell^\infty}\leq R\}$, then $\log \mathcal{C}(B_R,\|\cdot\|_{\ell^\infty},\epsilon)\leq n\log (4R\epsilon^{-1}) $
\cite[Proposition 5]{CuckerSmale:2002}.
By the Lipschitz continuity of DNN functions with respect to the DNN parameters, the covering number of the corresponding function class can be bounded by that of the parametrization. Thus by Lemma \ref{lem:NN-Lip},
\begin{align*}
 \log\mathcal{C}(\mathcal{H}_d,\|\cdot\|_{L^{\infty}(\Omega)},\epsilon)&\leq \log \mathcal{C}(\Theta_\sigma\otimes \Theta_q,\|\cdot\|_{\ell^\infty},\Lambda_d^{-1}\epsilon)\\
&\leq (N_{\kappa}+N_\theta)\log(4\max(R_q,R_\sigma)\Lambda_d\epsilon^{-1}),
\end{align*}
with $\Theta_\sigma$ and $\Theta_q$ denoting the neural network parameter sets for $\sigma$ and $q$, respectively, and the constant $\Lambda_d:=cR_\sigma W_\sigma \max(L_\sigma R_\sigma^{L_\sigma-1}W_\sigma^{L_\sigma},L_qR_q^{L_q-1}W_q^{L_q})$,
cf. Lemma \ref{lem:fcn-Lip}. By Lemma \ref{lem:fcn-Lip}, we also have $M_d
=cR_\sigma^{2}W_\sigma^2$. Then letting $s=n^{-\frac12}$
in Lemma \ref{lem:Dudley} gives
\begin{align*}
    &\quad \mathfrak{R}_n(\mathcal{H}_d)\leq4n^{-\frac12}+12n^{-\frac12}\int^{M_d}_{n^{-\frac12}}{\big((N_{\kappa}+N_\theta) \mbox{log}(4\max(R_q,R_\sigma)\Lambda_d\epsilon^{-1})\big)}^{\frac12}\ {\rm d}\epsilon\\
			&\leq4n^{-\frac12}+12n^{-\frac12}M_d{\big((N_{\kappa}+N_\theta) \mbox{log}(4\max(R_q,R_\sigma)\Lambda_dn^{\frac12})\big)}^\frac12 \\&
			\leq4n^{-\frac12}+cn^{-\frac12} R_\sigma^2 W_\sigma^2(N_\kappa+N_\theta)^{\frac12}\big(\log \max(R_q,R_\sigma)+\log \Lambda_d+\log n+ \log 2\big)^\frac12.
\end{align*}
Since $1\leq R_q,R_\sigma$, $1\leq W_q\leq N_\theta$, $1\leq W_\sigma\leq N_\kappa$ and $2\leq L_q, L_{\sigma}\leq c\log(d+2)$ (due to Lemma \ref{lem:tanh-approx}), we can bound the term $\log\Lambda_d$ by
\begin{equation*}
	\log\Lambda_d\leq c(\log R_\sigma +\log N_{\kappa} + \log R_q+\log N_\theta+\tilde c),
\end{equation*}
with the constants $c$ and $\tilde c$ depending on $c_1$, $c_z$, $d$, $L_q$ and $L_\sigma$ at most polynomially.
Hence, we have
\begin{equation*}
   \mathfrak{R}_n(\mathcal{H}_d)\leq c_{d} n^{-\frac12}R_\sigma^{2}N_\kappa^{2}(N_\kappa+N_\theta)^{\frac12}(\log^\frac12 R_\sigma +\log^\frac12 N_{\kappa} + \log^\frac12 R_q+\log^\frac12 N_\theta+\log^\frac12 n),
\end{equation*}
where $c_{d}>0$ depends on $d$, $c_1$ and $c_z$ at most polynomially.
Repeating the preceding argument leads to
\begin{align*}
\mathfrak{R}_n(\mathcal{H}_\sigma)&\leq c_\sigma n^{-\frac12}R_\sigma^{2L_\sigma}N_{\kappa}^{2L_\sigma-\frac32}\big(\log^\frac12R_\sigma+\log^\frac12N_{\kappa}+\log^\frac12n\big),\\
\mathfrak{R}_n(\mathcal{H}_b)&\leq c_bn^{-\frac12}R_\sigma^2 N_{\kappa}^{\frac52}\big(\log^\frac12 R_\sigma+\log^\frac12 N_{\kappa}+\log^\frac12 n\big),\\
\mathfrak{R}_n(\mathcal{H}_{b'})&\leq c_{b'}n^{-\frac12}R_\sigma^2 N_{\kappa}^{\frac52}\big(\log^\frac12 R_\sigma+\log^\frac12 N_{\kappa}+\log^\frac12 n\big),\\
\mathfrak{R}_n(\mathcal{H}_q)&\leq c_qn^{-\frac12}R_q^{2L_q}N_\theta^{2L_q-\frac32}\big(\log^\frac12R_q+\log^\frac12N_{\theta}+\log^\frac12n\big),
\end{align*}
where the involved constants depend on $d$ at most polynomially.
Last, the desired estimates follow from the PAC-type generalization  bound in Lemma \ref{lem:PAC}.
\end{proof}

\bibliographystyle{siam}
\bibliography{reference}

\end{document}